\documentclass[11pt,onecoulme]{article}
\usepackage{amsfonts}
\usepackage{mathrsfs}
\usepackage{amssymb}
\usepackage{color, amsmath,amssymb, amsfonts, amstext,amsthm, latexsym}

\usepackage{amssymb, epsfig, amssymb, latexsym}
\usepackage{amsmath}
\usepackage{graphicx}
\usepackage{longtable}

\numberwithin{equation}{section}

\allowdisplaybreaks

\newtheorem{definition}{Definition}[section]
\newtheorem{theorem}[definition]{Theorem}
\newtheorem{lemma}[definition]{Lemma}
\newtheorem{remark}[definition]{Remark}

\newtheorem{hyp}[definition]{Hypothesis}

\title{Viscosity Solutions to Second Order Path-Dependent  Hamilton-Jacobi-Bellman Equations in Hilbert Spaces \thanks{This work was partially supported by  the National Natural Science Foundation of China  (Grant No. 11401474), Shaanxi Natural Science Foundation
               (Grant No. 2017JM1016) and the Fundamental Research Funds for the Central Universities (Grant No. 2452019075).}}

\author{  Jianjun Zhou   \\
       College of Science,
             Northwest A\&F University,\\ Yangling 712100, Shaanxi, P. R.
             China\\
      \emph{E-mail:zhoujianjun@nwsuaf.edu.cn} }
          \date{}
\begin{document}

\maketitle

\pagestyle{plain}

\begin{abstract}
In this article, a notion of viscosity solutions is introduced for second order path-dependent Hamilton-Jacobi-Bellman (PHJB) equations associated with optimal control problems for path-dependent stochastic evolution equations in Hilbert spaces. We identify the value functional of optimal  control problems as unique viscosity solution to the associated PHJB equations. We also show that our notion of viscosity solutions is consistent with the corresponding notion of classical solutions, and satisfies a stability property.
\medskip

\end{abstract}


\section{Introduction}
\par
                         Let us consider a   controlled path-dependent stochastic evolution
                 equation (PSEE):
\begin{eqnarray}\label{state1}
\begin{cases}
            dX^{\gamma_t,u}(s)= AX^{\gamma_t,u}(s)ds+
           F\left(X_s^{\gamma_t,u},u(s)\right)ds+G\left(X_s^{\gamma_t,u},u(s)\right)dW(s),  \ \ s\in [t,T],\\
~~~~~X_t^{\gamma_t,u}=\gamma_t\in {\Lambda}_t,
\end{cases}
\end{eqnarray}
 for an unknown process $X^{\gamma_t,u}$ in a Hilbert space $H$. Here
%In the equation above,
 $T>0$ is an arbitrarily fixed finite time horizon;
 denote by $X^{\gamma_t,u}(s)$ the value of $X^{\gamma_t,u}$  at
 time $s$, and $X^{\gamma_t,u}_s$ the whole history path of $X^{\gamma_t,u}$ from time 0 to $s$; $\{W(t),t\geq0\}$ is a cylindrical Wiener process on  $(\Omega,{\cal {F}}, P)$
  with values in a
Hilbert space $\Xi$; $u(\cdot)=(u(s))_{s\in [t,T]}$ is  ${\cal{F}}^t$-progressively measurable  and take values in some metric space $(U,d)$ (we will say that $u(\cdot)\in {\cal{U}}[t,T]$) and ${\mathcal {F}}^{t}$ denotes the natural filtration of $W(s)-W(t)$, $s\in [t,T]$, augmented with the family ${\cal{N}}$ of $P$-null of $\cal{F}$.
 $A$ is the generator of
         a strongly continuous semigroup    of bounded linear operators $\{e^{tA}, t\geq0\}$ in Hilbert space
        $H$; $\Lambda_t$ denotes the set of all
 continuous $H$-valued functions defined over $[0,t]$ and ${\Lambda}=\bigcup_{t\in[0,T]}{\Lambda}_{t}$; $\gamma_t$ is an element of $\Lambda_t$ and denote by $\gamma_t(s)$ the value of $\gamma_t$  at
 time $s$.
We define a   norm on ${\Lambda}_t$  and a metric on ${\Lambda}$ as follows: for any $0\leq t\leq s\leq T$ and $\gamma_t,\eta_s\in {\Lambda}$,
\begin{eqnarray*}
   ||\gamma_t||_0:=\sup_{0\leq s\leq t}|\gamma_t(s)|,\ \ d_{\infty}(\gamma_t,\eta_s)%=d_{\infty}(\bar{\gamma}_{\bar{t}},\gamma_t)
   :=|t-s|
               +\sup_{0\leq \sigma\leq T}\left|e^{((\sigma-t)\vee0)A}\gamma_{t}(\sigma\wedge t)-e^{((\sigma-s)\vee0)A}\eta_{s}(\sigma\wedge s)\right|.
\end{eqnarray*}
 We assume the coefficients $F:\Lambda\times U\rightarrow H$ and  $G:\Lambda\times U\rightarrow L_2(\Xi,H)$   satisfy  Lipschitz condition under $||\cdot||_0$
                         with respect to  the path function.
\par
             One tries  to  maximize a cost functional of the form:
\begin{eqnarray}\label{cost1}
                     J(\gamma_t,u(\cdot)):=Y^{\gamma_t,u}(t),\ \ \ (t,\gamma_t)\in [0,T]\times {\Lambda},
\end{eqnarray}
 over  ${{\mathcal
                  {U}}}[t,T]$,
   where the process $Y^{\gamma_t,u}$ is defined by backward stochastic differential equation (BSDE):
    \begin{eqnarray}\label{fbsde1}
Y^{\gamma_t,u}(s)&=&\phi(X_T^{\gamma_t,u})+\int^{T}_{s}q(X_\sigma^{\gamma_t,u},Y^{\gamma_t,u}(\sigma),Z^{\gamma_t,u}(\sigma),u(\sigma))d\sigma\nonumber\\
                 &&-\int^{T}_{s}Z^{\gamma_t,u}(\sigma)dW(\sigma),\ \ \ a.s., \ \ \mbox{all}\ \ s\in [t,T].
\end{eqnarray}
Here  $q$ and $\phi$ are given real functionals on ${\Lambda}\times \mathbb{R}\times \Xi\times U$ and ${\Lambda}_T$, respectively,  and satisfy  Lipschitz condition under $||\cdot||_0$
                         with respect to  the path function.
             We define the value functional of the  optimal
                  control problem as follows:
\begin{eqnarray}\label{value1}
V(\gamma_t):=\mathop{\sup}\limits_{u(\cdot)\in{\cal{U}}[t,T]} Y^{\gamma_t,u}(t),\ \  (t,\gamma_t)\in [0,T]\times {\Lambda}.
\end{eqnarray}
  The goal of this article is to characterize this value functional  $V$. We    consider the following second order path-dependent Hamilton-Jacobi-Bellman (PHJB) equation:
  \begin{eqnarray}\label{hjb1}
\begin{cases}
\partial_tV(\gamma_t)+(A^*\partial_xV(\gamma_t),\gamma_t(t))_H+{\mathbf{H}}(\gamma_t,V(\gamma_t),\partial_xV(\gamma_t),\partial_{xx}V(\gamma_t))= 0,\ \ \  (t,\gamma_t)\in
                               [0,T)\times {\Lambda},\\
 V(\gamma_T)=\phi(\gamma_T), \ \ \ \gamma_T\in {\Lambda}_T;
 \end{cases}
\end{eqnarray}
      where
\begin{eqnarray*}
                                {\mathbf{H}}(\gamma_t,r,p,l)&=&\sup_{u\in{
                                         {U}}}[
                        (p,F(\gamma_t,u))_{H}+\frac{1}{2}\mbox{Tr}[ l G(\gamma_t,u)G^*(\gamma_t,u)]\\
                        &&\ \ \ \ \   +q(\gamma_t,r,pG(\gamma_t,u),u)], \ \ \ (t,\gamma_t,r,p,l)\in [0,T]\times{\Lambda}\times \mathbb{R}\times H\times {\cal{S}}(H).
\end{eqnarray*}
Here we let  $A^*$ the adjoint operator of $A$, $G^*$ the adjoint operator of $G$,   ${\cal{S}}(H)$ the space of bounded, self-adjoint operators on $H$ equipped with the operator norm,  $(\cdot,\cdot)_{H}$ the scalar product of $H$ and
                         $\partial_t,\partial_x$ and $\partial_{xx}$  the so-called pathwise  (or functional or Dupire; see \cite{dupire1, cotn0, cotn1}) derivatives, where $ \partial_t$  is known as horizontal derivative, while $\partial_x$ and $\partial_{xx}$  are first order and  second order vertical derivatives, respectively. 
\par
                         In this paper we will
                         develop a concept  of  viscosity solutions
                         to
                         PHJB equations on the space of  $H$-valued continuous paths %(see Definition \ref{definition4.1} for details)
                          and show that the value functional
                         $V$  defined in  (\ref{value1}) is   unique viscosity solution to the PHJB equation (\ref{hjb1}). In order to focus on our main well-posedness objective,
                         here we address
the Lipschitz case under $||\cdot||_0$.
\par
A definition of viscosity solutions for second order Hamilton-Jacobi-Bellman (HJB) equations in Hilbert spaces has been introduced  in Lions  \cite{lio1, lio3} for the case without unbounded term,
  and in \'{S}wi\c{e}ch \cite{swi} for the case with unbounded linear term. For the later case, it was based on the notion of the so-called $B$-continuous viscosity solutions which was introduced for first-order equations by
 Crandall and  Lions in \cite{cra6, cra7}. In  earlier paper Lions \cite{lio2},
 a specific second-order HJB equation for an optimal control of a Zakai
equation was studied by  also using  some ideas of the $B$-continuous viscosity solutions. Based on the maximum principle of Lions  \cite{lio3}, the first comparison theorem for $B$-continuous viscosity sub/supersolutions was proved in \'{S}wi\c{e}ch \cite{swi}.  We refer to   the monograph of  Fabbri,  Gozzi, and \'{S}wi\c{e}ch \cite{fab1} for a detailed account
for the theory of viscosity solutions.   One of the structural assumptions is
that the state space is a separable Hilbert space,  excluding for instance the metric space $\Lambda$ (notice
that in this paper we do not directly generalize those results to $\Lambda$,
as we adopt pathwise, rather than Fr\'{e}chet, derivatives on $\Lambda$).
\par
Fully nonlinear   first order path-dependent  Hamilton-Jacobi equations in finite dimensional spaces were deeply studied by Lukoyanov  \cite{luk}. By  formulating  the maximum/minimum condition on the subset of absolutely
continuous paths, the existence and uniqueness theorems of viscosity solutions
were obtained  when Hamilton function $\mathbf{H}$ is $d_p$-locally Lipschitz continuous in the path function. In our paper \cite{zhou}, we extended the results in \cite{luk} to Lipschitz continuous case under $||\cdot||_0$.
For minimax  solutions,  Bayraktar and Keller \cite{bay}  proposed the  notion for a
class of fully nonlinear PHJB equations with nonlinear,
monotone, and coercive operators on Hilbert spaces and proved
the existence, uniqueness and stability.
\par
Concerning the second order  case, In \cite{peng2.5}, Peng made the first attempt to extend Crandall-Lions framework  to path-dependent case in  finite dimensional spaces. By the left frozen maximization principle,  a comparison principle was proved under a technical  condition (16) in \cite{peng2.5} which is not easy to be satisfied in  practice. Tang and Zhang \cite{tang1} improved   Peng \cite{peng2.5} and  identified  the value functional of  optimal  control problems as
             a viscosity solution to the path-dependent Bellman equations by restricting
semi-jets on a space of $\alpha$-H\"{o}lder continuous paths.  However, the uniqueness result was not given.  Our paper \cite{zhou3} studied  a class of infinite-horizon optimal control problems for stochastic differential equations with delays
 and identified  the value functional  %of the optimal control problems
as a unique viscosity solution to the associated second order HJB equation. The Crandall-Lions definition was not further investigated, as the fact that the supremum norm $||\cdot||_0$ is not G${\hat{a}}$teaux differentiable
was perceived as an almost
insurmountable obstacle.
\par
 At the same time, when the space $H$ is finit-dimensiaonal, a new concept of viscosity solutions
was introduced by Ekren, Keller, Touzi
 and Zhang \cite{ekren1} in the semilinear context, and further extended to
 fully nonlinear parabolic
equations by Ekren, Touzi, and Zhang \cite{ekren3, ekren4},  elliptic
equations by Ren \cite{ren}, obstacle problems by Ekren \cite{ekren0}, and degenerate second-order
equations by Ren, Touzi, and Zhang \cite{ren1} and Ren and Rosestolato \cite{ren2}.
Cosso, Federico, Gozzi, Rosestolato and Touzi \cite{cosso} studied a class of semilinear second order PHJB equations with a linear unbounded operator on Hilbert
spaces. This new notion adopted is different from the
Crandall-Lions definition as the tangency condition is not pointwise but in the sense of
expectation with respect to an appropriate class of probability measures. This modification simplifies a lot the proof of uniqueness and  is successfully applied to the semilinear case, as
this does not require anymore the passage through the maximum principle of Lions  \cite{lio3}. For instance, the comparison theorem in  Cosso, Federico, Gozzi, Rosestolato and Touzi \cite{cosso}, even in the Markovian case, extends the result in Crandall-Lions framework. However, in the  nonlinear case, this new notion faces some obstacles, which limits its applications. The comparison theorems in Ekren,
Touzi \& Zhang  %\cite{ekren2},
\cite{ekren3} and \cite{ekren4}, Ekren \cite{ekren0} and Ren \cite{ren}, in particular,  require that the Hamilton function $\mathbf{H}$ is uniformly nondegenerate. In Ren,  Touzi and  Zhang \cite{ren1} and Ren and Rosestolato \cite{ren2}, the degenerate case was studied, but in order to apply these results one has to require that   the coefficients $F$, $G$, $q$ and $\phi$ are  $d_p$-uniformly continuous with respect to  the path function. Our context does not fall into the latter framework as in our case the coefficients are
required to have continuity properties under supremum norm $||\cdot||_0$.
\par
At present, some authors turn to study viscosity solutions of PHJB equations in  Crandall-Lions framework by defining a smooth functional and applying the Borwein-Preiss  variational principle.
         Cosso and  Russo \cite{cosso1} studied
 the existence and comparison theorem of viscosity solutions for the path-dependent heat equation. Our paper \cite{zhou5} and \cite{zhou6}  introduced   notions of viscosity solutions  for second order PHJB
          equations  in finite dimensional spaces and  first order PHJB
          equations in Hilbert spaces, respectively,  and identified the value functional  as
             unique viscosity solution to the associated PHJB equations. None of the results  above are directly applicable to our situation, as in our case, it is needed to adapt the maximum principle to the infinite-dimensional context.
             \par
              As mentioned earlier, In the Markovian case, the core of  Crandall-Lions viscosity solutions is $B$-continuity. Specifically, in order to overcome the difficulties caused by unbounded operator $A$, it is necessary to assume that the coefficients satisfy $B$-continuity.
              Then ones can  study the $B$-continuity of value function and  obtain the comparison theorem. The main objective of this paper is to extend  Crandall-Lions viscosity solutions to  our infinite-dimensional path-dependent context  without $B$-continuity assumption %under natural assumptions
              on the coefficients.
              \par

          Our core results are as follows.   We define a functional $\Upsilon^3:\Lambda\rightarrow \mathbb{R}$ by
          $$
          \Upsilon^3(\gamma_t)=S(\gamma_t)+3|\gamma_t(t)|^6, \  ~~ \gamma_t\in \Lambda,
          $$
          and
 \begin{eqnarray*}
S(\gamma_t)=\begin{cases}
            \frac{(||\gamma_{t}||_{0}^6-|\gamma_{t}(t)|^6)^3}{||\gamma_{t}||^{12}_{0}}, \
         ~~ ||\gamma_{t}||_{0}\neq0; \\
0, \ ~~~~~~~~~~~~~~~~~~~ ||\gamma_{t}||_{0}=0,
\end{cases}
\  ~~ \gamma_t\in \Lambda.
\end{eqnarray*}
%where $\gamma_t(t)$ denotes the value of
%$\gamma_t$ at time $t$.
   We show that it is equivalent to $||\cdot||_0^6$ and study its  regularity in the horizontal/vertical sense mentioned above. % regular enough  to be applied to the functional formula. % (see Theorem \ref{theoremito} and Lemma \ref{theoremS}).
     % and study
%its  properties (see Lemma \ref{theoremS000}).
This key functional  is the starting point for the proof of stability and uniqueness results. The uniqueness property is derived from the comparison theorem. We list some  points about  the proof of  the comparison theorem.
  \par
  (a)
            For every fixed $(t,\gamma_t)\in [0,T)\times\Lambda_t$,     $f(\eta_s):=\Upsilon^3(\eta_s-\gamma_{t,s,A})$, $(s,\eta_s)\in[t,T]\times \Lambda$ can be used as a test function in our definition of viscosity solutions as we show that $f$ satisfies  a functional It\^o inequality. This is important as  functional $\Upsilon^3$ is equivalent to $||\cdot||_0^6$. Then we can define an auxiliary function $\Psi$ which includes the  functional $f$.  By this, we only need to study the continuity  under $||\cdot||_0$ rather than  $B$-continuity of value functional. In particular, the comparison theorem is  established 
            when the coefficients satisfy   Lipschitz assumption under $||\cdot||_0$.
           % with  Lipschitz assumption under $||\cdot||_0$ on the coefficients. %under natural assumptions
              %on the coefficients case under $||\cdot||_0$  without $B$-continuity assumption %under natural assumptions
%              on the coefficients. under natural assumptions on the coefficients. % ( not under $B$-continuity assumptions on the coefficients like the Markovian case) comparison,
            We emphasize that,
with respect to the standard viscosity solution theory in infinite dimensional spaces, the
$B$-continuity assumption on the coefficients is   bypassed in our framework.
\par
                 (b)
 We  use $\Upsilon^3$ to  define a  smooth  gauge-type function $\overline{\Upsilon}^3:\Lambda\times \Lambda\rightarrow \mathbb{R}$ by
                 $$
                 \overline{\Upsilon}^3(\gamma_t,\eta_s)=\overline{\Upsilon}^3(\eta_s,\gamma_t):=\Upsilon^3(\eta_s-\gamma_{t,s,A})+|s-t|^2, \ \  0\leq t\leq s\leq T, \ \gamma_t,\eta_s\in \hat{\Lambda},
                 $$
where $\gamma_{t,s,A}\in \Lambda_s$, and $\gamma_{t,s,A}(\sigma)=\gamma_t(\sigma), \ \sigma\in [0,t]$ and  $\gamma_{t,s,A}(\sigma)=e^{(\sigma-t)A}\gamma_t(t), \ \sigma\in (t,s]$.  Then we can apply
                 a modification of Borwein-Preiss variational principle (see Theorem 2.5.2 in  Borwein \& Zhu  \cite{bor1}) to get a maximum of a perturbation of the auxiliary function $\Psi$.
            \par
            (c)   Unfortunately, the second order vertical derivative $\partial_{xx}S$ is not equal to $\mathbf{0}$. % then if  we apply Theorem 8.3 of \cite{cran2} as in Zhou \cite{zhou3}, we can not get the corresponding results.     % which makes it difficult to apply Theorem 8.3 of \cite{cran2}.
  To apply the maximum principle (see Theorem 8.3 of \cite{cran2}), a stronger convergence property of auxiliary functional is needed. %the more convergence of auxiliary functional is needed.
 By doing more detailed calculations, we obtain the expected convergence property of auxiliary functional  and prove the comparison theorem.
\par
An important consequence of (a) is that our notion of viscosity solutions is meaningful  even in the
Markovian  case. More precisely, in the
Markovian  case, the functional $f$ defined in (a) reduces to    $f(s,y):=|y-e^{(s-t)A}x|_H^2$, $(s,y)\in[t,T]\times H$,  for every fixed $(t,x)\in [0,T)\times H$.  The $B$-continuity assumption on the coefficients can be   bypassed with this functional.
\par
Regarding existence, we  show that the value functional  $V$  defined in  (\ref{value1}) is  a viscosity solution to the PHJB
                         equation given in  (\ref{hjb1}) by functional It\^o  formula, functional It\^o inequality  and dynamic programming principle.
                         Such a formula was firstly provided  in
                          Dupire  \cite{dupire1}
   (see Cont and  Fourni$\acute{e}$ \cite{cotn0}, \cite{cotn1} for a more general and systematic research). In this paper we extend the functional It\^o
  formula to  infinite dimensional spaces. We also provide a functional It\^o inequality for $f$ defined in (a) from functional It\^o
  formula.
  \par
  Finally, following Dupire  \cite{dupire1}, we study PHJB equation (\ref{hjb1}) in the metric space $(\Lambda,d_\infty)$ in the present paper, while some literatures (for example,  \cite{cosso}, \cite{cosso1}) study PHJB equations in a complete pseudometric space. The reason we do this is that it is convenient to define $\gamma_{t,\cdot, A}$ for every $(t,\gamma_t)\in [0,T)\times \Lambda$ in our framework, which is useful in defining metric $d_\infty$ and test functional $f$.
\par
                 The outline of this article is  as follows. In the following
              section,  we introduce the notations used throughout the paper and review the
              background for BSDEs. We  give a  functional It\^o formula which  will be used to prove the existence of
  viscosity solutions to PHJB equation (\ref{hjb1}). We also present a modification of
          Borwein-Preiss variational principle and  a smooth functional $S$ which are the key to prove the stability and uniqueness results of viscosity solutions.
   Section 3 is devoted to studying  PSEE (\ref{state1}).  A functional It\^o inequality for $f$ defined in (a) is also provided by  functional It\^o
  formula.
   In Section 4, we introduce  preliminary results on path-dependent stochastic  optimal control problems. We  give  the dynamic programming principle,  which will be used in the following sections.
  In Section 5, we define classical and viscosity solutions to
             PHJB equations (\ref{hjb1}) and  prove  that the value functional $V$ defined by (\ref{value1}) is a viscosity solution.   We also show
             the consistency with the notion of classical solutions and the stability result.  Finally,  the uniqueness of viscosity solutions to  (\ref{hjb1}) is proven in section 6.

\section{Preliminary work}  \label{RDS}
%%%%%%%%%%%%%%%%%%%%%%%%%%%%%%%%%%%%%%%%%%%%%%%%%%%%%%%%%%%%%%%%%
\par
2.1. ${Notations \ and \ Spaces}$.
 We list some notations that are used in this paper. Let $\Xi$, $K$ and $H$
         denote real separable Hilbert spaces, with scalar products
         $(\cdot,\cdot)_\Xi$, $(\cdot,\cdot)_K$ and $(\cdot,\cdot)_H$, respectively.  We use the symbol $|\cdot|$ to denote the norm in
         various spaces, with a subscript if necessary.
         $L(\Xi,H)$ denotes the space of all
         bounded linear operators from $\Xi$ into $H$; the subspace of
         Hilbert-Schmidt operators, with the Hilbert-Schmidt norm, is
         denoted by $L_2(\Xi,H)$.  Let
         $L(H)$ denotes the space of all
         bounded linear operators from $H$ into itself. Denote by ${\cal{S}}(H)$ the Banach space of bounded and self-adjoint
operators in the Hilbert space $H$ endowed with the operator norm.
            The operator $A$ is the generator of a strongly continuous
                 semigroup $\{e^{tA}, t\geq0\}$ of bounded linear operators in the
                 Hilbert space $H$.
                 The domain of the
         operator $A$ is denoted by ${\mathcal {D}}(A)$. $A^*$ denotes  the adjoint operator of $A$ with domain ${\mathcal {D}}(A^*)$.  %denotes the domain of the
%         operator $A^*$.
 %The domain of a linear (unbounded)
%         operator $A$ is denoted ${\mathcal {D}}(A)$.
         Let % $n$ be a positive integer and
$T>0$ be a fixed number.  For each  $t\in[0,T]$,
          define
         $\hat{\Lambda}_t:=D([0,t];H)$ as  the set of c$\grave{a}$dl$\grave{a}$g  $H$-valued
         functions on $[0,t]$.
       We denote $\hat{\Lambda}^t=\bigcup_{s\in[t,T]}\hat{\Lambda}_{s}$  and let  $\hat{\Lambda}$ denote $\hat{\Lambda}^0$.
       \par
%A very important remark on the notation $\hat{\Lambda}$:
As in Dupire \cite{dupire1}, we will denote elements of $\hat{\Lambda}$ by lower
case letters and often the final time of its domain will be subscripted, e.g. $\gamma\in \hat{\Lambda}_t\subset \hat{\Lambda}$ will be
denoted by $\gamma_t$. Note that, for any $\gamma\in  \hat{\Lambda}$, there exists only one $t$ such that $\gamma\in  \hat{\Lambda}_t$. For any $0\leq s\leq t$, the value of
$\gamma_t$ at time $s$ will be denoted by $\gamma_t(s)$. Moreover, if
a path $\gamma_t$ is fixed, the path $\gamma_t|_{[0,s]}$, for $0\leq s\leq t$, will denote the restriction of the path  $\gamma_t$ to the interval
$[0,s]$.
\par
        Following Dupire \cite{dupire1},  for $ x\in H,\gamma_t\in \hat{\Lambda}_t$, $0\leq t\leq \bar{t}\leq T$, we define $\gamma^x_{t}\in\hat{\Lambda}_t$ and $\gamma_{t,\bar{t}}, \gamma_{t,\bar{t},A}\in \hat{\Lambda}_{\bar{t}}$ as
\begin{eqnarray*}
  \gamma^x_{t}(s)&=&\gamma_t(s),\ \ s\in [0,t); \ \ \ \gamma^x_{t}(t)=\gamma_t(t)+x;\\
  \gamma_{t,\bar{t}}(s)&=&\gamma_t(s),\ \ s\in [0,t];\ \ \ \gamma_{t,\bar{t}}(s)=\gamma_t(t), \ \ s\in (t,\bar{t}];\\
  \gamma_{t,\bar{t},A}(s)&=&\gamma_t(s),\ \ s\in [0,t];\ \ \ \gamma_{t,\bar{t},A}(s)=e^{(s-t)A}\gamma_t(t), \ \ s\in (t,\bar{t}].
\end{eqnarray*}
 We define a norm on $\hat{\Lambda}_t$  and a metric on $\hat{\Lambda}$ as follows: for any $0\leq t\leq s\leq T$ and $\gamma_t,\eta_s\in \hat{\Lambda}$,
\begin{eqnarray}\label{2.1}
   ||\gamma_t||_0:=\sup_{0\leq s\leq t}|\gamma_t(s)|,\ \ \  \ d_{\infty}(\gamma_t,\eta_s):=|t-s|+||\gamma_{t,T,A}-\eta_{s,T,A}||_0.
               %+\sup_{0\leq s\leq T}|e^{((s-t)\vee0)A}\gamma_{t}(s\wedge t)-e^{((s-\bar{t})\vee0)A}\bar{\gamma}_{\bar{t}}(s\wedge \bar{t})|.
\end{eqnarray}
%Here and in the sequel,
% for notational simplicity,
%we use $||\gamma_{t}-\bar{\gamma}_{\bar{t}}||_0$ to denote $||\gamma_{t,\bar{t}}-\bar{\gamma}_{\bar{t}}||_0$.
Then $(\hat{\Lambda}_t, ||\cdot||_0)$ is a Banach space, and $(\hat{\Lambda}^t, d_{\infty})$ is a complete metric space by Lemma 5.1 in \cite{zhou6}.
\par Now we define the pathwise  derivatives of Dupire \cite{dupire1}.
 \begin{definition}\label{definitionc0} (Pathwise derivatives)
       Let $t\in [0,T)$ and  $f:\hat{\Lambda}^t\rightarrow \mathbb{R}$.
\begin{description}
        \item{(i)}  Given $(s,\gamma_s)\in [t,T)\times \hat{\Lambda}^t$, the horizontal derivative of $f$ at $\gamma_s$ (if the corresponding limit exists and is finite) is defined as
        \begin{eqnarray}\label{2.3}
               \partial_tf(\gamma_s):=\lim_{h\rightarrow0,h>0}\frac{1}{h}\left[f(\gamma_{s,s+h})-f(\gamma_s)\right].
\end{eqnarray}
For the final time $T$, the horizontal derivative of $f$ at $\gamma_T\in\hat{\Lambda}^t$ (if the corresponding limit exists and is finite) is defined as
$$
\partial_tf(\gamma_T):=\lim_{s<T, s\uparrow T}\partial_tf(\gamma_T|_{[0,s]}).
$$
If the above limit exists and is finite for every $(s,\gamma_s)\in [t,T]\times \Lambda^t$, the functional $\partial_tf:\hat{\Lambda}^t\rightarrow \mathbb{R}$ is called the horizontal derivative of $f$ with domain $\hat{\Lambda}^t$.
       \par
      \item{(ii)}  Given $(s,\gamma_s)\in [t,T]\times \hat{\Lambda}^t$,
      if there exists a $B\in H$ such that
$$
\lim_{|h|\rightarrow0}\frac{\left|f(\gamma_t^{h})-f(\gamma_t)-(B,h)_H\right|}{|h|}=0,
$$
we say $\partial_{x}f(\gamma_t):=B$ as the first order vertical derivative  of $f$ at $\gamma_s$. If $\partial_xf$ exists  for every $(s,\gamma_s)\in [t,T]\times \Lambda^t$, the map $\partial_xf:\hat{\Lambda}^t\rightarrow H$ is called the first order vertical  derivative of $f$ with domain $\hat{\Lambda}^t$.
\par
If   there exists a $B_1\in {\cal{S}}(H)$ such that
$$
\lim_{|h|\rightarrow0}\frac{\left|\partial_{x}f(\gamma_t^{h})-\partial_{x}f(\gamma_t)-B_1h\right|}{|h|}=0,
$$
we say $\partial_{xx}f(\gamma_t):=B_1$ as the second order vertical derivative  of $f$ at $\gamma_s$.
If $\partial_{xx}f$ exists  for every $(s,\gamma_s)\in [t,T]\times \Lambda^t$, the map $\partial_{xx}f:\hat{\Lambda}^t\rightarrow H$ is called the second order vertical  derivative of $f$ with domain $\hat{\Lambda}^t$.
%(if the corresponding limit exists and is finite) is defined as
%      \begin{eqnarray}\label{2.20jia}
% \partial_{x}f(\gamma_s):=(\partial_{x_1}f(\gamma_s),\partial_{x_2}f(\gamma_s),\ldots, \partial_{x_d}f(\gamma_s))
%\end{eqnarray}
%where
%        \begin{eqnarray}\label{2.2}
% \partial_{x_i}f(\gamma_s):=\lim_{h\rightarrow0}\frac{1}{h}\bigg{[}f(\gamma_s^{he_i})-f(\gamma_s)\bigg{]},\ \
% %\partial_{x_ix_j}u:=\partial_{x_i}(\partial_{x_j}u),\ \
% i=1,2,\ldots,d,
%\end{eqnarray}
%with $e_1,e_2,\ldots,e_d$ is the standard  orthonormal basis of $R^d$.
%If the above limit exists and is finite for every $(s,\gamma_s)\in [t,T]\times \Lambda^t$, the map $\partial_xf:=(\partial_{x_1}f,\partial_{x_2}f,\ldots, \partial_{x_d}f):\hat{\Lambda}^t\rightarrow R^d$ is called the vertical  derivative of $f$ with domain $\hat{\Lambda}^t$.
%\par
%We take the convention that $\gamma_s$ is column vector, but $\partial_{x}f$ denotes row vector.
\end{description}
\end{definition}
%Following Dupire \cite{dupire1}, we define spatial derivatives of $f:\hat{\Lambda}\rightarrow R$ in the standard sense: if there exists a $B\in H$ such that
%$$
%\lim_{|h|\rightarrow0}\frac{|f(\gamma_t^{h})-f(\gamma_t)-(B,h)_H|}{|h|}=0,
%$$
%we say $\partial_{x}f(\gamma_t)=B$; if   there exists a $B_1\in L(H)$ such that
%$$
%\lim_{|h|\rightarrow0}\frac{|\partial_{x}f(\gamma_t^{h})-\partial_{x}f(\gamma_t)-B_1h|}{|h|}=0,
%$$
%we say $\partial_{xx}f(\gamma_t)=B_1$;
%and the right time-derivative of $f$, if exists and is finite, as:
%\begin{eqnarray}\label{2.3}
%               \partial_tf(\gamma_t):=\lim_{l\rightarrow0,l>0}\frac{1}{l}\bigg{[}f(\gamma_{t,t+l})-f(\gamma_t)\bigg{]}, \ t<T.
%\end{eqnarray}
%For the final time $T$, we define
%$$
%\partial_tu(\gamma_T):=\lim_{t<T,t\uparrow T}\partial_tu(\gamma_t).
%$$
% We take the convention that $\gamma_t$ is column vector, but $\partial_{x}u$ denotes row vector and $\partial_{xx}u$ denotes $d\times d$-matrix.
\begin{definition}\label{definitionc}
       Let $t\in[0,T)$ and $f:\hat{\Lambda}^t\rightarrow \mathbb{R}$ be given.
\begin{description}
        \item{(i)}
                 We say $f\in C^0(\hat{\Lambda}^t)$ if $f$ is continuous in $\gamma_s$  on $\hat{\Lambda}^t$ under $d_{\infty}$.
                  \item{(ii)}
                 We say $f\in C_p^0(\hat{\Lambda}^t)\subset C^0(\hat{\Lambda}^t)$ if $f$ grows  in a polynomial way.
\par
       \item{(iii)}  We say %$f\in C^{0,2}(\hat{\Lambda}^t)\subset C^0(\hat{\Lambda}^t)$ if   $\partial_{x}f$ and $\partial_{xx}f$ exist and are continuous in $\gamma_s$  on $\hat{\Lambda}^t$ under $d_{\infty}$, and
        $f\in C^{1,2}(\hat{\Lambda}^t)\subset C^{0,2}(\hat{\Lambda}^t)$ if  $\partial_tf$ exists and is continuous in $\gamma_s$  on $\hat{\Lambda}^t$ under $d_{\infty}$.
       %\par
%       \item{(iv)}  We say $f\in C^{1,2}(\hat{\Lambda}^t)\subset C^0(\hat{\Lambda}^t)$ if  $\partial_tf$, $\partial_{x}f$ and $\partial_{xx}f$ exist and are continuous in $\gamma_s$  on $\hat{\Lambda}^t$ under $d_{\infty}$.
       \par
       \item{(iv)} We say %$f\in C^{0,2}_p(\hat{\Lambda}^t)\subset C^{0,2}(\hat{\Lambda}^t)$ if $f$ and all of its vertical derivatives grow  in a polynomial way, and
        $f\in C^{1,2}_p(\hat{\Lambda}^t)\subset C^{1,2}(\hat{\Lambda}^t)$ if $f$ and all of its derivatives grow  in a polynomial way.
       %\par
%       \item{(vi)} We say $f\in C^{1,2}_p(\hat{\Lambda}^t)\subset C^{1,2}(\hat{\Lambda}^t)$ if $f$ and all of its derivatives grow  in a polynomial way.
\end{description}
\end{definition}
\par
Let $\Lambda_t:= C([0,t],H)$ be the set of all continuous $H$-valued functions defined over $[0,t]$. We denote ${\Lambda}^t=\bigcup_{s\in[t,T]}{\Lambda}_{s}$  and let  ${\Lambda}$ denote ${\Lambda}^0$.
 Clearly, $\Lambda:=\bigcup_{t\in[0,T]}{\Lambda}_{t}\subset\hat{\Lambda}$, and each $\gamma\in \Lambda$ can also be viewed as an element of $\hat{\Lambda}$. $(\Lambda_t, ||\cdot||_0)$ is a
 Banach space, and $(\Lambda^t,d_{\infty})$ is a complete metric space by Lemma 5.1 in \cite{zhou6}.
 $f:\Lambda^t\rightarrow R$ and $\hat{f}:\hat{\Lambda}^t\rightarrow R$ are called consistent
  on $\Lambda^t$ if $f$ is the restriction of $\hat{f}$ on $\Lambda^t$.
\begin{definition}\label{definitionc2}
       Let  $t\in [0,T)$ and  $f:{\Lambda}^t\rightarrow \mathbb{R}$  be given.
\begin{description}
\item{(i)}
                 We say $f\in C^0({\Lambda}^t)$ if $f$ is continuous in $\gamma_s$  on $\Lambda^t$ under $d_{\infty}$. %For simplicity, we denote $C_0^0({\Lambda})$ by $C^0({\Lambda})$.
 \item{(ii)}
                 We say $f\in C_p^0({\Lambda}^t)\subset C^0({\Lambda}^t)$ if $f$ grows  in a polynomial way.
       \item{(iii)} We say  %$f\in C^{0,2}_p({\Lambda}^t)$ if
       % there exists $\hat{f}\in C^{0,2}_p(\hat{{\Lambda}}^t)$ which is consistent with $f$ on $\Lambda^t$, and
       $f\in C^{1,2}_p({\Lambda}^t)$ if
        there exists $\hat{f}\in C^{1,2}_p(\hat{{\Lambda}}^t)$ which is consistent with $f$ on $\Lambda^t$
        % \item{(iv)} We say $f\in C^{1,2}_p({\Lambda}^t)$ if
%        there exists $\hat{f}\in C^{1,2}_p(\hat{{\Lambda}}^t)$ which is consistent with $f$ on $\Lambda^t$.
       % \par
%       \item{(iv)} We say $f\in C^{1,2}_p({\Lambda}^t)\subset C^{1,2}({\Lambda}^t)$ if $f$ and all of its derivatives grow  in a polynomial way.
\end{description}
\end{definition}
  %\par
%           Let $\Omega:=\{\omega\in C([0,T],R^n):\omega(0)={\mathbf{0}}\}$, the set of continuous functions with initial value ${\mathbf{0}}$,
%           $W$ the canonical process, $P$ the Wiener measure, ${\cal {F}}$ the  Borel $\sigma$-field over $\Omega$, completed with
%            respect to the Wiener measure $P$ on this space. Then  $(\Omega,{\cal {F}},P)$ is a complete   space.  Here and in the sequel,
% for notational simplicity,
%we use $\mathbf{0}$ to denote vectors or matrices with appropriate dimensions whose components are all equal to 0. By $\{{\cal {F}}_t\}_{0\leq t\leq T}$ we denote  the filtration generated by $\{W(t),0\leq t\leq T\}$, augmented
%       with the family $\mathcal {N}$ of $P$-null of ${\cal {F}}$.
% The filtration $\{{\cal {F}}_t\}_{0\leq t\leq T}$ satisfies the
%       usual conditions.
       \par
           %Let $(\Omega,{\cal {F}},P)$ be a complete space with a filtration
%           $\{{\mathcal {F}}_{t}\}_{t\geq0}$ which satisfies the usual condition, i.e., $\{{\mathcal {F}}_{t}\}_{t\geq0}$
%            is a right continuous increasing family of sub $\sigma$-algebra of $\mathcal
%            {F}$ and ${\mathcal {F}}_0$ contains all $P$-null sets
%            of $\mathcal{F}$.  %$\{W(t)\}_{t\geq0}$ be a standard Wiener process  defined on
%        $(\Omega,{\cal {F}},P)$, and with values in $\Xi$.
               By a cylindrical Wiener process defined on a  complete probability  space
        $(\Omega,{\cal {F}},P)$, and with values in a Hilbert space $\Xi$,
        we mean a family $\{W(t),t\geq0\}$ of linear mappings $\Xi\rightarrow
        L^2(\Omega)$ such that for every $\xi, \eta \in \Xi$,
        $\{W(t)\xi,t\geq0\}$ is a real Wiener process and
       ${E}(W(t)\xi\cdot W(t)\eta)=(\xi,\eta)_\Xi t$.  ${\mathcal {F}}_{t}$, $t\in [0,T]$, will denote the natural filtration of $W$, augmented with the family ${\cal{N}}$ of $P$-null of $\cal{F}$:
       $$
              {\mathcal{F}}_{t}=\sigma(W(s):s\in[0,t])\vee \mathcal
             {N}.
       $$
      The filtration
           $\{{\mathcal {F}}_{t}\}_{0\leq t\leq T}$  satisfies the usual condition.
      %  In the following,
%        $\{W(t),t\geq0\}$ is a cylindrical Wiener process adapted
%       to the filtration $\{{\mathcal {F}}_{t}\}_{t\geq0}$.
   For every $[t_1,t_2]\subset[0,T]$, we also use the
       notation:
       $$
              {\mathcal{F}}_{t_1}^{t_2}=\sigma(W(s)-W(t_1):s\in[t_1,t_2])\vee \mathcal
             {N}.
       $$
        We also write ${\cal{F}}^t$ for $\{{\cal{F}}^s_t, t\leq s\leq T\}$.
        \par
        By $\mathcal{P}$  we denote the predictable $\sigma$-algebra generated by predictable processes and by $\mathcal{B}(\Theta)$ we denote, the Borel $\sigma$-algebra of any topological space $\Theta$.

Next we define several classes of stochastic processes with in a Hilbert space $K$.

 $\bullet L^{2}_{\mathcal{P}}(\Omega\times [0,T];K)$ denotes   the space of equivalence classes of  processes $y\in L^2(\Omega\times [0,T];K)$,  admitting a predictable version.  $ L^{2}_{\mathcal{P}}(\Omega\times [0,T];K)$
  is endowed with the norm
$$|{y}|^{2}=\mathbb{E}\int_{0}^{T}|{y}(t)|^{2}dt.$$

$\bullet L^{p}_{\mathcal{P}}(\Omega;L^{2}([0,T];K))$, defined for  $p\in[1,\infty)$, denotes   the space of equivalence classes of  processes $\{y(t),t\in [0,T]\}$, with values in $K$ such that  the norm
$$|{y}|^{p}=\mathbb{E}\left(\int_{0}^{T}|{y}(t)|^{2}dt\right)^{p/2}$$
is finite, and $\{y(t),t\in [0,T]\}$ admits a predictable version.

$\bullet L^{p}_{\mathcal{P}}(\Omega;C([0,t];K))$, defined for  $p\in[1,\infty)$ and $t\in (0,T]$, denotes  the space of predictable processes
$\{y(s),s\in[0,t]\}$ with continuous paths in $K$, such that  the norm
$$|{y}|^{p}=\mathbb{E}\sup_{s\in[0,t]}|y(s)|^{p}$$
is finite. Elements of $L^{p}_{\mathcal{P}}(\Omega;C([0,t];K))$ are identified up to indistinguishability.

\vbox{}
2.2. \emph{Functional   It\^o formula}.
Assume that $\vartheta\in L^{p}_{\mathcal{P}}(\Omega\times [0,T];H)$, $\varpi\in L^{p}_{\mathcal{P}}(\Omega\times [0,T];L_2(\Xi,H))$ for some $p>2$, and $(\theta,\gamma_{\theta})\in [0,T)\times \Lambda$, %$X(t)$ is an ${\cal{F}}_t$-measurable $H$-valued random variable,
then the following process
\begin{eqnarray}\label{formular1}
                 X(s)=e^{(s-\theta)A}\gamma_{\theta}(\theta)+\int^{s}_{\theta}e^{(s-\sigma)A}\vartheta(\sigma)d\sigma+\int^{s}_{\theta}e^{(s-\sigma)A}\varpi(\sigma)dW(\sigma),\ s\in [\theta,T],
\end{eqnarray}
and  $X(s)=\gamma_\theta(s),\ s\in [0,\theta)$
is well defined and $\mathbb{E}\sup_{s\in[0,T]}|X(s)|^{p}<\infty$ (see Proposition  7.3 in \cite{da}).
\begin{lemma}\label{theoremito}
\ \
Suppose  $f\in C_p^{1,2}(\hat{\Lambda}^{{t}})$ and $A^*\partial_xf\in C_p^0(\hat{\Lambda}^{{t}})$
for some $t\in[\theta,T)$. Then, under the above conditions, $P$-a.s.,  for all $s\in [t,T]$:
\begin{eqnarray}\label{statesop0}
                 f(X_s)&=&f(X_{t})+\int_{t}^{s}[\partial_tf(X_\sigma)+(A^*\partial_xf(X_\sigma,X(\sigma))_H
                 +(\partial_xf(X_\sigma),\vartheta(\sigma))_H\nonumber\\
                 &&+\frac{1}{2}\mbox{Tr}(\partial_{xx}f(X_\sigma)\varpi(\sigma)\varpi^*(\sigma))]d\sigma+\int^{s}_{t}(\partial_xf(X_\sigma),\varpi(\sigma)dW(\sigma))_H.
\end{eqnarray}
 Here and in the following, for every $s\in [0,T]$, $X(s)$ denotes  the value of $X$  at
 time $s$, and $X_s$ the whole history path of $X$ from time 0 to $s$.
\end{lemma}
The proof is  similar to Theorem 4.1 in  Cont \& Fournie \cite{cotn1} (see also Dupire \cite{dupire1}).
%the case without path-dependent (see  Proposition 3.3 in  Buckdahn and Li \cite{buck1}).
For the convenience of readers, here we give its proof.
\par
{\bf  Proof}. \ \    We can assume that the process $X(s),\ s\in[0,T]$ is
               bounded. This can be shown by localization. Namely
               for arbitrary constant $C> ||\gamma_t||_0$ define a stopping time
               $\tau_C$:
$$
               \tau_C=\inf\{s\in[\theta,T]: |X(s)|\geq C\}
$$
                  with the convention that $T_C=T$
                  if this set is empty. If one defines
$$
                  \vartheta_C(s)={\mathbf{1}}_{[0,\tau_C]}\vartheta(s),\ \
                  \varpi_C(s)={\mathbf{1}}_{[0,\tau_C]}\varpi(s),\ \
                 s\in[0,T],
$$
   and $X^{C}(s)=e^{(s-\tau_C)^+A}X (s\wedge C)$, i.e.
\begin{eqnarray*}
\begin{cases}
         X^{C}(s)=X(s)\ \ \ \ \ \ \ \ \ \ \ \ \ \ \ \ \ \ \ \ \mbox{if}\ s\leq \tau_C, \\
          X^{C}(s)=e^{(s-\tau_C)A}X(\tau_C)\ \ \ \ \ \ \  \mbox{if}\ s>\tau_C.
\end{cases}
\end{eqnarray*}
                 % and $X^{C}(s)=X(s\wedge \tau_C)$,  % i.e.
%\begin{eqnarray*}
%\begin{cases}
%         X^{n,A}(s)=X(s)\ \ \ \ \ \ \ \ \ \ \ \ \ \ \ \ \ \ \ \mbox{if}\ s\leq \tau_n, \\
%          X^{n,A}(s)=e^{(s-\tau_n)A}X(\tau_n)\ \ \ \ \ \ \  \mbox{if}\ s>\tau_n.
%\end{cases}
%\end{eqnarray*}
                  %It follows from  that %Corollary 2.5 in [31] that

                  It follows from  that %Corollary 2.5 in [31] that
\begin{eqnarray*}
             X^{C}(s)=e^{(s-\theta)A}\gamma_\theta(\theta)+\int_{\theta}^{s}{e^{(s-\sigma)A}}\vartheta_C(\sigma)d\sigma
             +\int_{\theta}^{s}{e^{(s-\sigma)A}}\varpi_C(\sigma)dW(\sigma), \ \
             s\in [\theta,T],
\end{eqnarray*}
and $X^{C}(s)=\gamma_\theta(s),  \ s\in [0,\theta)$.
                  If the formula (\ref{statesop0}) is true for $\vartheta_C$, $\varpi_C$ and
                  $X^{C}$ for arbitrary $C>0$, then, %using again
%                  Corollary 2.5 in [31],
it is true in the general
                  case.
\par
  For any  $s\in[t,T]$, denote $X^n(\sigma)=X{\mathbf{1}}_{[0,t)}(\sigma)+\sum^{2^n-1}_{i=0}X(t_{i+1}){\mathbf{1}}_{[t_i,t_{i+1})}(\sigma)+X(s){\mathbf{1}}_{\{s\}}(\sigma)$, $\sigma\in [0,s]$.  %  which is a  piecewise constant approximation of $X$.
Here $t_i=t+\frac{i(s-t)}{2^n}$. For every $(\sigma,\gamma_\sigma)\in [0,T]\times\hat{\Lambda}$, define $\gamma_{\sigma-}$  by
$$
               \gamma_{\sigma-}(l)=\gamma_{\sigma}(l),\ \ l\in [0,\sigma),\ \ \mbox{and}\  \ \gamma_{\sigma-}(\sigma)=\lim_{l\uparrow \sigma}\gamma_{\sigma}(l).
$$ We start with the decomposition
\begin{eqnarray}\label{decom}
&&f({X^n_{t_{1}-}})-f({X^n_{t_{0}-}})=f({X^n_{t_{1}-}})-f({X^n_{t_0}}),\nonumber\\
                &&f({X^n_{t_{i+1}-}})-f({X^n_{t_{i}-}})=f({X^n_{t_{i+1}-}})-f({X^n_{t_{i}}})
                   +f({X^n_{t_{i}}})-f({X^n_{t_{i}-}}),\ i\geq1.
\end{eqnarray}
 Let $\psi(\sigma)=f({X^n_{t_{i},t_i+\sigma}})$, we have $f({X^n_{t_{i+1}-}})-f({X^n_{t_{i}}})=\psi(h)-\psi(0)$, where $h=\frac{s-t}{2^n}$.
 Let $\psi_{t^+}$ denote the right derivative of $\psi$, then
 $$
 \psi_{t^+}(l)=\lim_{\delta>0,\delta\rightarrow0}\frac{\psi_{t^+}(l+\delta)-\psi_{t^+}(l)}{\delta}=\lim_{\delta>0,\delta\rightarrow0}\frac{f({X^n_{t_{i},t_i+l+\delta}})-f({X^n_{t_{i},t_i+l}})}{\delta}
 =\partial_tf(X_{t_i,t_i+l}^n),\ \ l\in[0,h].
 $$
 By
 $$
 d_\infty({X^n_{t_{i},t_i+l_1}},{X^n_{t_{i},t_i+l_2}})=|l_1-l_2|+\sup_{0\leq\sigma\leq |l_1-l_2|}|X^n(t_{i+1})-e^{\sigma A}X^n(t_{i+1})|, \ \ \ l_1, l_2\in [0,h],
$$
 and  $f\in C^{1,2}_p(\hat{\Lambda}^{t})$, we have $\psi$ and $\psi_{t^+}$ is continuous on $[0,h]$, therefore,
 $$
                                f({X^n_{t_{i+1}}}_{-})-f({X^n_{t_{i}}})=\psi(h)-\psi(0)=\int^{h}_{0}\psi_{t^+}(l)dl=\int^{t_{i+1}}_{t_i}\partial_tf(X_{t_i,l}^n)dl, \ i\geq 0.
$$
%
%
%
% Since $f\in C^{1,2}_p(\hat{\Lambda}^{t})$, the right derivative of $\psi$ denoted by $\psi_{t^+}$ is continuous on $[0,h]$ and $\psi_{t^+}(\sigma)=\partial_tf(X_{t_i,t_i+\sigma}^n)$, therefore,
% $$
%                                f({X^n_{t_{i+1}-}})-f({X^n_{t_{i}}})=\psi(h)-\psi(0)=\int^{h}_{0}\psi_{t^+}(\sigma)d\sigma=\int^{t_{i+1}}_{t_i}\partial_tf(X_{t_i,\sigma}^n)d\sigma,
%                                 \ \ i\geq0.
%$$
The term $f({X^n_{t_{i}}})-f({X^n_{t_{i}-}})$ in (\ref{decom}) can be written $\pi(X(t_{i+1})-X(t_i))-\pi(0)$, where
$\pi(l)=f({X^n_{t_{i}-}}+l{\mathbf{1}}_{\{t_i\}})$. Since $f\in C^{1,2}_p(\hat{\Lambda}^{t})$, $\pi$ is a $C^2$ function and $\nabla_x\pi(l)=\partial_xf({X^n_{t_{i}-}}+l{\mathbf{1}}_{\{t_i\}})$, $\nabla_x^2\pi(l)=\partial_{xx}f({X^n_{t_{i}-}}+l{\mathbf{1}}_{\{t_i\}})$.
Applying the It\^o formula (see %Theorem 4.17 in \cite{da}
Proposition 1.165 in \cite{fab1}
) to $\pi$ between $0$ and $h$ and the continuous process  $(X(t_{i}+s)-X(t_i))_{s\geq0}$, yields:
\begin{eqnarray}\label{ito}
                              &&f({X^n_{t_{i}}})-f({X^n_{t_{i}-}})=\pi({X(t_{i+1})-X(t_i)})-\pi(0)\nonumber\\
                              &=&\int^{t_{i+1}}_{t_i}[(A^*\partial_xf({X^n_{t_{i}}}_{-}+(X(\sigma)-X(t_i)){\mathbf{1}}_{\{t_{i}\}}),X(\sigma))_H\nonumber\\
                              &&+
                              (\partial_xf({X^n_{t_{i}}}_{-}+(X(\sigma)-X(t_i)){\mathbf{1}}_{\{t_{i}\}}),\vartheta(\sigma))_H\nonumber\\
                              &&
                              +\frac{1}{2}\mbox{Tr}[\partial_{xx}f({X^n_{t_{i}-}}+(X(\sigma)-X(t_i)){\mathbf{1}}_{\{t_i\}})\varpi(\sigma)\varpi^*(\sigma)]]d\sigma\nonumber\\
                              &&
                              + \int^{t_{i+1}}_{t_i}(\partial_xf({X^n_{t_{i}-}}+(X(\sigma)-X(t_i)){\mathbf{1}}_{\{t_i\}}),\varpi(\sigma)dW(\sigma)), \ \  i\geq 1.
 \end{eqnarray}
 Summing over $i\geq 0$ and denoting $i(l)$ the index such that $l\in [t_{i(l)},t_{i(l)+1})$, we obtain
 \begin{eqnarray}\label{app}
                                && f(X^n_s)-f(X_{t})=f(X^n_s)-f(X^n_{t})\nonumber\\
                                &=&\int^{s}_{t}\partial_tf(X_{t_{i(\sigma)},\sigma}^n)d\sigma
                                +\int^{s\vee t_1}_{t_1}[(A^*\partial_xf({X^n_{t_{i(\sigma)}}}_{-}+(X(\sigma)-X(t_{i(\sigma)})){\mathbf{1}}_{\{t_{i(\sigma)}\}}),X(\sigma))_H\nonumber\\
                              &&~~+(\partial_xf({X^n_{t_{i(\sigma)}}}_{-}+(X(\sigma)-X(t_{i(\sigma)})){\mathbf{1}}_{\{t_{i(\sigma)}\}}),\vartheta(\sigma))_H\nonumber\\
                              &&~~
                              +\frac{1}{2}\mbox{Tr}[\partial_{xx}f({X^n_{t_{i(\sigma)}-}}+(X(\sigma)-X(t_{i(\sigma)})){\mathbf{1}}_{\{t_{i(\sigma)}\}})\varpi(\sigma)\varpi^*(\sigma)]]d\sigma\nonumber\\
                              &&
                              + \int^{s\vee t_1}_{t_{1}}(\partial_xf({X^n_{t_{i(\sigma)}-}}+(X(\sigma)-X(t_{i(\sigma)})){\mathbf{1}}_{\{t_{i(\sigma)}\}}),\varpi(\sigma)dW(\sigma))_H.
 \end{eqnarray}
  $f(X^n_s)$ %and $f(X^n_{t})$
  converges to $f(X_s)$ %and $f(X_{t})$
   almost surely.
   %, respectively.
    Since all approximations of $X$ appearing in the various integrals have a $||\cdot||_{0}$-distance from $X_s$ less than $||X^n_s-X_s||_0\rightarrow0$, $f\in C^{1,2}_p(\hat{\Lambda}^{t})$ and $A^*\partial_xf\in C_p^0(\hat{\Lambda}^{{t}})$
    imply  that the integrands appearing in the above integrals converge respectively to $\partial_tf(X_s), A^*\partial_xf(X_s),
    \partial_xf(X_s), \partial_{xx}f(X_s)$ as $n\rightarrow\infty$. By $X$ is bounded and $f\in C^{1,2}_p(\hat{\Lambda}^{t})$, the integrands in the various above integrals are bounded.  The dominated convergence and Burkholder-Davis-Gundy inequalities  for  the stochastic integrals % (see \cite{pro}, Chapter IV, Theorem 32)
then ensure that the Lebesgue integrals converge almost surely,
and the stochastic integral in probability, to the terms appearing in (\ref{statesop0}) as $n\rightarrow\infty$.\ \ $\Box$
\par
By the above Lemma, we have the following important results.
\begin{lemma}\label{0815lemma}
             Let $f\in C_p^{1,2}(\Lambda^t)$ and $\hat{f}\in C_p^{1,2}(\hat{\Lambda}^t)$ such that $\hat{f}$ is consistent with $f$ on $\Lambda^t$, then the following definition
             $$
             \partial_tf:=\partial_t\hat{f}, \ \ \ \partial_xf:=\partial_x\hat{f}, \ \ \ \partial_{xx}f:=\partial_{xx}\hat{f} \ \ \mbox{on} \ \Lambda^t
             $$
              is independent of the choice of $\hat{f}$. Namely, if there is another $\hat{f}'\in C_p^{1,2}(\hat{\Lambda}^t)$ such that $\hat{f}'$ is consistent with $f$ on $\Lambda^t$, then the derivatives of $\hat{f}'$
              coincide with those of $\hat{f}$ on $\Lambda^t$.
\end{lemma}
\par
{\bf  Proof}. \ \
 By the definition of the horizontal derivative, it is clear that $\partial_t\hat{f}(\gamma_s)=\partial_t\hat{f}'(\gamma_s)$ for every $(s,\gamma_s)\in [t,T]\times\Lambda^t$.
%
%
% For every $(t,\gamma_t)\in [0,T)\times\Lambda$, first let $A=\mathbf{0}$, $\vartheta=\mathbf{0}$, $\varpi=\mathbf{0}$ and $\theta=t$ in (\ref{formular1}), by the Lemma \ref{theoremito},
%$$
%                        \int^{s}_{t}\partial_t\hat{f}(\gamma_{t,\sigma})d\sigma=\int^{s}_{t}\partial_t\hat{f}'(\gamma_{t,\sigma})d\sigma, \ \ s\in [t,T].
%$$
%Then from the  regularity $\partial_t\hat{f},\partial_t\hat{f}'\in C^0(\hat{\Lambda}^t)$, it follows that $\partial_t\hat{f}(\gamma_t)=\partial_t\hat{f}'(\gamma_t)$.
  Next, let $A=\mathbf{0}$,  $\varpi=\mathbf{0}$,  $\theta=t$ and $\vartheta\equiv h\in H$ in (\ref{formular1}),  by the Lemma \ref{theoremito},
$$
                        \int^{s}_{t}(\partial_x\hat{f}(X_\sigma),h)_Hd\sigma=\int^{s}_{t}(\partial_x\hat{f}'(X_\sigma),h)_Hd\sigma, \ \ s\in [t,T].
$$
Here and
in the sequel, for notational simplicity, we use $\mathbf{0}$ to denote elements, operators, or paths  which are
identically equal to zero. %vectors, matrices, or
%paths with appropriate dimensions whose components are all equal to 0
 By the regularity $\partial_x\hat{f},\partial_x\hat{f}'\in C^0(\hat{\Lambda}^t;H)$ and the arbitrariness of $h\in H$, we  have  $\partial_x\hat{f}(\gamma_s)=\partial_x\hat{f}(\gamma_s)$ for every $(s,\gamma_s)\in [t,T]\times\Lambda^t$.
 Finally,  let $A=\mathbf{0}$,  $\vartheta=\mathbf{0}$,  $\theta=t$ and $\varpi\equiv a\in L_2(\Xi;H)$ in (\ref{formular1}),  by the Lemma \ref{theoremito},
 $$
                        \int^{s}_{t}\mbox{Tr}(\partial_{xx}\hat{f}(X_\sigma)aa^*)d\sigma =\int^{s}_{t}\mbox{Tr}(\partial_{xx}\hat{f}'(X_\sigma)aa^*)d\sigma, \ \ s\in [t,T].
$$
By the regularity $\partial_{xx}\hat{f},\partial_{xx}\hat{f}'\in C^0(\hat{\Lambda}^t;L(H))$ and the arbitrariness of $a\in L_2(\Xi;H)$, we also have  $\partial_{xx}\hat{f}(\gamma_s)=\partial_{xx}\hat{f}(\gamma_s)$ for every $(s,\gamma_s)\in [t,T]\times\Lambda^t$. \ \ $\Box$

\vbox{}
2.3. \emph{Backward stochastic differential
             equations}.
      We consider the backward stochastic differential
             equations (BSDEs) in a Hilbert space $K$:
\begin{eqnarray}\label{bsde}
                         Y(t)+\int^{T}_{t}Z(\sigma) dW_\sigma=\int^{T}_{t}f(\sigma,Y(\sigma),Z(\sigma))d\sigma
                         +\eta.\ \ \ \ \ \ 0\leq t\leq T,
\end{eqnarray}
                         for $t$ varying on the time interval $[0,T]$.
                        the mapping $f: \Omega \times [0,T]\times
                         K\times L_2(\Xi,K)\rightarrow K$ is assume to be measurable with
                         respect to ${\cal {P}}\times {\cal {B}} (K\times L_2(\Xi,K))$. %(by $\cal {P} $ we denote the predictable
                        %$\sigma$-algebra on $\Omega \times [0,T]$).
                        $\eta: \Omega \rightarrow K$ is
                       assumed to be ${\mathcal {F}}_T$-measurable.
                       $(f,\eta)$ are the parameters of (\ref{bsde}).
\par
                       We recall the well known results on
                       the existence and uniqueness of the BSDEs. %(see also El Karoui, Peng, Quenez \cite{el} and Peng \cite{peng11}).
\begin{lemma}\label{lemma2.5} (Proposition 4.3 in \cite{fuh0}) Assume that: (i) there
                                          exists $L>0$ such that
\begin{eqnarray*}
                           |f(\sigma,y_1,z_1)-f(\sigma,y_2,z_2)|\leq
                           L(|y_1-y_2|+|z_1-z_2|),
\end{eqnarray*}
                             $P$-a.s. for every $\sigma \in [0,T], y_1,y_2\in K, z_1,z_2\in
                             L_2(\Xi,K)$;
\par
       (ii) there exists $p\in[2,\infty)$ such that
\begin{eqnarray*}
                  \mathbb{E}\bigg{(}\int^{T}_{0}|f(\sigma,0,0)|^2d\sigma\bigg{)}^{\frac{p}{2}}<\infty,\
                  \ \ \ \ \mathbb{E}|\eta|^p<\infty.
\end{eqnarray*}
                    Then there exists a unique pair of processes $(Y,Z)$% $Y\in L^p_{\mathcal {P}}(\Omega;C([0,T];R^m))$, $Z\in
                % L^p_{\mathcal {P}}(\Omega;L^2([0,T];R^{m\times n}))$
                 such that  (\ref{bsde})
                    holds for $t\in[0,T]$ and
\begin{eqnarray}\label{lemma2.510}
                 \mathbb{E}\sup_{t\in[0,T]}|Y(t)|^p+\mathbb{E}\bigg{(}\int^{T}_{0}|Z(\sigma)|^2d\sigma\bigg{)}^{\frac{p}{2}}\leq
                 C_p\mathbb{E}\bigg{(}\int^{T}_{0}|f(\sigma,0,0)|^2d\sigma\bigg{)}^{\frac{p}{2}}+C_p\mathbb{E}|\eta|^p,
\end{eqnarray}
for some constant $C_p>0$ depending only on $p$, $L$, $T$.
Moreover, let two BSDEs of  parameters    $(\eta^1,f^1)$ and $(\eta^2,f^2)$ satisfy all the assumptions (i) and (ii).
                     % We denote by $(Y^1,Z^1)$ and  $(Y^2,Z^2)$)  their respective adapted solutions.
                      Then the difference of the solutions
                      $(Y^1,Z^1)$ and $(Y^2,Z^2)$ of BSDE (\ref{bsde}) with
                      the data $(\eta^1,f^1)$ and
                      $(\eta^2,f^2)$, respectively, satisfies the
                      following  estimate:
\begin{eqnarray}\label{lemma2.51}
                 &&\mathbb{E}\sup_{t\in[0,T]}|Y^1(t)-Y^2(t)|^p+\mathbb{E}\bigg{(}\int^{T}_{0}|Z^1(\sigma)-Z^2(\sigma)|^2d\sigma\bigg{)}^{\frac{p}{2}}\nonumber\\
                 &\leq&
                 C_p\mathbb{E}\bigg{(}\int^{T}_{0}|f^1(\sigma,Y^2(\sigma),Z^2(\sigma))-f^2(\sigma,Y^2(\sigma),Z^2(\sigma))|^2d\sigma\bigg{)}^{\frac{p}{2}}+C_p\mathbb{E}|\eta^1-\eta^2|^p.
\end{eqnarray}
               %   for some constant $C_p>0$ depending only on $p$, $L$, $T$.
                  \end{lemma}
                  We also have  the following comparison theorem on BSDEs in infinite dimensional spaces.
    \begin{lemma}\label{lemma2.70904}
                      (Theorem 2.7 in \cite{zhou1})   Let two BSDEs of  parameters    $(\eta^1,f^1)$ and $(\eta^2,f^2)$ satisfy all the assumptions of Lemma \ref{lemma2.5}.
                     Denote by $(Y^1,Z^1)$ and  $(Y^2,Z^2)$)  their respective adapted solutions.
                      %Denote by
%                      $(Y^1,Z^1)$ and $(Y^2,Z^2)$ of BSDE (\ref{bsde}) with
%                      the data $(\eta^1,f^1)$ and
%                      $(\eta^2,f^2)$,
%
%
%                       Assume $f^1$
%                     satisfies  the assumption in Lemma 2.5 and  $\eta^1,\eta^2\in L^2(\Omega,{\mathcal
%                      {F}}_T,P)$. We denote by $(Y^1,Z^1)$ the unique solution of Eq.(2.8) with
%                      parameters $(\eta^1,f^1)$, and denote
%                   by $(Y^2,Z^2)$ the solution of BSDE:
%$$
%                         Y^2_t+\int^{T}_{t}Z^2_\sigma dW_\sigma=\int^{T}_{t}g_\sigma d\sigma
%                         +\eta^2,\ \ \ \ \ \ t\in[0,T].
%$$
%                          where $g\in L^2_{\mathcal {P}}(\Omega\times[0,T],R)$ and $E[\sup_{t\in[0,T]}|Y_t|^2]<\infty$.
 If
 $$
                     \eta^1\geq \eta^2, P\mbox{-a.s., and}
                        \ \ f^1(t,Y^2_t,Z^2_t)\geq  f^2(t,Y^2_t,Z^2_t),\ \
                        dP\otimes dt\ \ \mbox{a.s.}
 $$
                     Then we have that $Y^1_t\geq Y^2_t$, a.s., for
                     all $t\in[0,T]$.
\par
                     Moreover, the comparison is strict: that is
$$
                       Y^1_0=Y^2_0\Leftrightarrow \eta^1=\eta^2,\ f^1(t,Y^2_t,Z^2_t)= f^2(t,Y^2_t,Z^2_t),\ \
                        dP\otimes dt\ \ \mbox{a.s.}
$$
              \end{lemma}

                  \vbox{}
2.4. \emph{Borwein-Preiss variational principle and functional $S$}.
In this subsection we introduce the Borwein-Preiss variational principle and functional $S$, which are the key to proving  the uniqueness and  stability  of viscosity solutions.
%We conclude this section with  the following three  lemmas which will be used to prove the uniqueness and  stability  of viscosity solutions.
\begin{definition}\label{gaupe}
              Let $t\in [0,T]$ be fixed.  We say that a continuous functional $\rho:\Lambda^t\times \Lambda^t\rightarrow [0,+\infty)$ is a {gauge-type function} provided that:
             \begin{description}
        \item{(i)} $\rho(\gamma_s,\gamma_s)=0$ for all $(s,\gamma_s)\in [t,T]\times \Lambda^t$,
        \item{(ii)} for any $\varepsilon>0$, there exists $\delta>0$ such that, for all $\gamma_s, \eta_l\in \Lambda^t$, we have $\rho(\gamma_s,\eta_l)\leq \delta$ implies that
        $d_\infty(\gamma_s,\eta_l)<\varepsilon$.
        \end{description}
\end{definition}
The following lemma is a modification of Borwein-Preiss variational principle (see Theorem 2.5.2 in  Borwein \& Zhu  \cite{bor1}). It will be used to
get a maximum of a perturbation of the auxiliary function in the proof of uniqueness. The proof  is completely similar to the finite dimensional case (see Lemma 2.12 in \cite{zhou5}). Here we omit it.
\begin{lemma}\label{theoremleft} %(Borwein-Preiss Variational Principle)
Let $t\in [0,T]$ be fixed and let $f:\Lambda^t\rightarrow \mathbb{R}$ be an upper semicontinuous functional  bounded from above. Suppose that $\rho$ is a gauge-type function
 and $\{\delta_i\}_{i\geq0}$ is a sequence of positive number, and suppose that $\varepsilon>0$ and $(t_0,\gamma^0_{t_0})\in [t,T]\times \Lambda^t$ satisfy
 $$
f\left(\gamma^0_{t_0}\right)\geq \sup_{(s,\gamma_s)\in [t,T]\times \Lambda^t}f(\gamma_s)-\varepsilon.
 $$
 Then there exist $(\hat{t},\hat{\gamma}_{\hat{t}})\in [t,T]\times \Lambda^t$ and a sequence $\{(t_i,\gamma^i_{t_i})\}_{i\geq1}\subset [t,T]\times \Lambda^t$ such that
  \begin{description}
        \item{(i)} $\rho(\gamma^0_{t_0},\hat{\gamma}_{\hat{t}})\leq \frac{\varepsilon}{\delta_0}$,  $\rho(\gamma^i_{t_i},\hat{\gamma}_{\hat{t}})\leq \frac{\varepsilon}{2^i\delta_0}$ and $t_i\uparrow \hat{t}$ as $i\rightarrow\infty$,
        \item{(ii)}  $f(\hat{\gamma}_{\hat{t}})-\sum_{i=0}^{\infty}\delta_i\rho(\gamma^i_{t_i},\hat{\gamma}_{\hat{t}})\geq f(\gamma^0_{t_0})$, and
        \item{(iii)}  $f(\gamma_s)-\sum_{i=0}^{\infty}\delta_i\rho(\gamma^i_{t_i},\gamma_s)
            <f(\hat{\gamma}_{\hat{t}})-\sum_{i=0}^{\infty}\delta_i\rho(\gamma^i_{t_i},\hat{\gamma}_{\hat{t}})$ for all $(s,\gamma_s)\in [\hat{t},T]\times \Lambda^{\hat{t}}\setminus \{(\hat{t},\hat{\gamma}_{\hat{t}})\}$.

        \end{description}
\end{lemma}
 \par
Define  $S:\hat{\Lambda}\rightarrow \mathbb{R}$ by, for every $(t,\gamma_t)\in [0,T]\times{\hat{\Lambda}}$,
 \begin{eqnarray*}
S(\gamma_t)=\begin{cases}
            \frac{(||\gamma_{t}||_{0}^6-|\gamma_{t}(t)|^6)^3}{||\gamma_{t}||^{12}_{0}}, \
         ~~ ||\gamma_{t}||_{0}\neq0; \\
0, \ ~~~~~~~~~~~~~~~~~~~ ||\gamma_{t}||_{0}=0.
\end{cases}
\end{eqnarray*}
\par
Define, for every $M\in \mathbb{R}$,
$$
         \Upsilon^M(\gamma_t):=S(\gamma_t)+M|\gamma_t(t)|^6, \ \  \gamma_t\in \hat{\Lambda};
$$
$$
         \Upsilon^M(\gamma_t,\eta_s)= \Upsilon^M(\eta_s,\gamma_t):=\Upsilon^M(\eta_s-\gamma_{t,s,A}), \ \  0\leq t\leq s\leq T, \ \gamma_t,\eta_s\in \hat{\Lambda};
$$
and
$$
\overline{\Upsilon}^M(\gamma_t,\eta_s)=\overline{\Upsilon}^M(\eta_s,\gamma_t):=\Upsilon^M(\eta_s,\gamma_t)+|s-t|^2, \ \  0\leq t\leq s\leq T, \ \gamma_t,\eta_s\in \hat{\Lambda}.
$$
%For every $M>0$, define $\Upsilon^M$ and $\overline{\Upsilon}^M$ by
%\begin{eqnarray*}
%         \Upsilon^M(\gamma_t,\eta_s)=S(\gamma_t-\eta_s)+M|\gamma_t(t)-\eta_s(s)|_{8}^8, \ \ \ (t,\gamma_t), (s,\eta_s)\in [0,T]\times\hat{\Lambda},
%         %\cases{\frac{(||\gamma_t-\eta_s||_{0,8}^8-|\gamma_t(t)-\eta_s(s)|_{8}^8)^3}{||\gamma_t-\eta_s|_{0,8}^{16}}
%%         +M|\gamma_t(t)-\eta_s(s)|_{8}^8, \ \mbox{if}\  ||\gamma_t-\eta_s||_{0,8}^8>0, \cr
%%          0 ~~~~~~~~~~~~~~~~~~~~~~~~~~~~~~~~~~~~~~~~~~~~~~~~~~~~~~~~~
%%           \ \mbox{if}\  ||\gamma_t-\eta_s||_{0,8}^8=0,}
%\end{eqnarray*}
%and
%\begin{eqnarray*}
%         \overline{\Upsilon}^M(\gamma_t,\eta_s)= \Upsilon^M(\gamma_t,\eta_s)+|s-t|^2 \ \ \ (t,\gamma_t), (s,\eta_s)\in [0,T]\times\hat{\Lambda}.
%\end{eqnarray*}
%For simplicity, we let $\Upsilon^M(\gamma_t)$ denote $\Upsilon^M(\gamma_t,\eta_s)$ when $\eta_s(l)\equiv{\mathbf{0}}$ for all $l\in [0,s]$.
%In the proof of uniqueness of viscosity solutions, in order to apply theorem 8.3 in [2], we need the following lemma.
The following    lemma  is used  to prove the  uniqueness and  stability  of viscosity solutions.
\begin{lemma}\label{theoremS}
%For every fixed $(\hat{t}, a_{\hat{t}}) \in [0,T)\times \hat{\Lambda}_{\hat{t}}$, define   $S_1^{a_{\hat{t}}}:{\hat{\Lambda}}^{\hat{t}}\rightarrow R$ by
% $$S_1^{a_{\hat{t}}}(\gamma_t):=S_1(\gamma_{t},a_{\hat{t}}),\ \ (t,\gamma_t)\in [\hat{t}, T]\times{\hat{\Lambda}}^{\hat{t}}.$$ Then
 $S(\cdot)\in C^{1,2}_p(\hat{\Lambda})$. Moreover, for every $M\geq1$,
\begin{eqnarray}\label{s0}
                \frac{8}{27}||\gamma_t||_{0}^6\leq   \Upsilon^M(\gamma_t)%+   |\gamma_t(t)-a_{\hat{t}}(\hat{t})|_8^8
                \leq (M+1)||\gamma_t||_{0}^6, \ \ (t,\gamma_t)\in [0, T]\times{\hat{\Lambda}}.
\end{eqnarray}
\end{lemma}
\par
   {\bf  Proof  }. \ \ % First, by the definition of $S_1^{a_{\hat{t}}}(\gamma_{t})$, it is clear that $S_1^{a_{\hat{t}}}(\gamma_{t})\in C^0(\hat{\Lambda}^{\hat{t}})$ and $\partial_tS_1^{a_{\hat{t}}}(\gamma_{t})=0$ for
%   $(t,\gamma_t)\in [\hat{t},T]\times\hat{\Lambda}^{\hat{t}}$.
   First, we prove $S\in C^0(\hat{\Lambda})$. %Let $\eta_s\rightarrow \gamma_t$ in $(\hat{\Lambda},d_\infty)$.
   For any $(t,\gamma_t), (s, \eta_s)\in [0,T]\times \hat{\Lambda}$, if $s\geq t$,
   $$
   |\gamma_t(t)-\eta_s(s)|\leq \left|\gamma_t(t)-e^{(s-t)A}\gamma_t(t)\right|+\left|e^{(s-t)A}\gamma_t(t)-\eta_s(s)\right|, %\rightarrow0, \ \ \mbox{as}\ d_\infty(\gamma_t,\eta_s)\rightarrow0,
   $$
   and
   \begin{eqnarray*}
   |||\gamma_t||_0-||\eta_s||_0|&\leq& ||\gamma_{t,s,A}||_0-||\gamma_t||_0+||\gamma_{t,s,A}-\eta_s||_0\\
   &\leq& \sup_{t\leq l\leq s}\left|(e^{(l-t)A}-I)\gamma_t(t)\right|+||\gamma_{t,s,A}-\eta_s||_0;
   %\rightarrow0 \ \ \mbox{as}\ d_\infty(\gamma_t,\eta_s)\rightarrow0.
\end{eqnarray*}
   %Thus we have $S(\eta_s)\rightarrow S(\gamma_t)$ as $\ d_\infty(\gamma_t,\eta_s)\rightarrow0$.
    %Otherwise, we may assume $s<t$, then
    if $s<t$, let $\gamma_t(t-):=\lim_{l\uparrow t}\gamma_t(l)$, we have
   \begin{eqnarray*}
   |\gamma_t(t)-\eta_s(s)|&\leq& \left|\gamma_t(t)-e^{(t-s)A}\eta_s(s)\right|+\left|e^{(t-s)A}(\eta_s(s)-\gamma_t(s))\right|+\left|e^{(t-s)A}(\gamma_t(s)-\gamma_t(t-))\right|\\
   &&+\left|(e^{(t-s)A}-I)\gamma_t(t-)\right|+|\gamma_t(t-)-\gamma_t(s)|+|\gamma_t(s)-\eta_s(s)|, %\rightarrow0, \ \ \mbox{as}\ d_\infty(\gamma_t,\eta_s)\rightarrow0,
 \end{eqnarray*}
 and
 \begin{eqnarray*}
   &&|||\gamma_t||_0-||\eta_s||_0|\leq ||\eta_{s,t,A}||_0-||\eta_s||_0+||\eta_{s,t,A}-\gamma_t||_0\\
   &\leq& \sup_{s\leq l\leq t}\left|(e^{(l-s)A}-I)\gamma_t(t)\right|+\sup_{s\leq l\leq t}\left|(e^{(l-s)A}-I)\right||\gamma_t(t)-\eta_s(s)|+||\eta_{s,t,A}-\gamma_t||_0.
   %\rightarrow0 \ \ \mbox{as}\ d_\infty(\gamma_t,\eta_s)\rightarrow0.
\end{eqnarray*}
    Then we have $S(\eta_s)\rightarrow S(\gamma_t)$ as $\eta_s\rightarrow\gamma_t$ under $d_\infty$. Thus $S\in C^0(\hat{\Lambda})$.
    Second, by the definition of $S(\cdot)$, it is clear that $\partial_tS(\gamma_{t})=0$ for all $(t,\gamma_t)\in [0,T]\times \hat{\Lambda}$.
   Next, we consider $ \partial_{x}S$. Clearly,
   \begin{eqnarray}\label{s10902jia}
   \partial_xS(\gamma_{0})=0,\  \ \gamma_0\in \Lambda_0.
\end{eqnarray}
For every $(t,\gamma_t)\in (0,T]\times \hat{\Lambda}$, let $||\gamma_t||_{0^-}=\sup_{0\leq s<t}|\gamma_t(s)|$.
Then,
 if $|\gamma_t(t)|<||\gamma_t||_{0^-}$,
\begin{eqnarray*}
  &&\lim_{|h|\rightarrow0}\frac{\bigg{|}S(\gamma_t^{h})-S(\gamma_{t})+\frac{18(||\gamma_t||^6_{0}-|\gamma_t(t)|^6)^2|\gamma_t(t)|^4(\gamma_t(t),h)_H}
   {||\gamma_t||^8_{0}}\bigg{|}}{|h|}\\
  &=&\lim_{|h|\rightarrow0}\frac{\left|\left(||\gamma_t||^6_{0}-|\gamma_t(t)+{h}|^6\right)^3
   -{\left(||\gamma_t||^6_{0}-|\gamma_t(t)|^6\right)^3}+18\left(||\gamma_t||^6_{0}-|\gamma_t(t)|^6\right)^2|\gamma_t(t)|^4(\gamma_t(t),h)_H \right|}{|h|\times||\gamma_t||^{12}_{0}}\nonumber\\
  % &=&\lim_{|h|\rightarrow0}\frac{{|-(2||\gamma_t||^2_{0}-|\gamma_t(t)+{h}|^2-|\gamma_t(t)|^2)(2(\gamma_t(t),h)_H+|h|^2)+4(||\gamma_t||^2_{0}-|\gamma_t(t)|^2)(\gamma_t(t),h)_H|}
%  }{|h|\times||\gamma_t||^2_{0}}\nonumber\\
   &=&0.
   \end{eqnarray*}
   Thus,
\begin{eqnarray}\label{s1}
   \partial_xS(\gamma_{t})=-\frac{18\left(||\gamma_t||^6_{0}-|\gamma_t(t)|^6\right)^2
    |\gamma_t(t)|^4\gamma_t(t)}{||\gamma_t||^{12}_{0}}.
\end{eqnarray}
 If  $|\gamma_t(t)|>||\gamma_t||_{0^-}$,
\begin{eqnarray}\label{s2}
   \partial_{x}S(\gamma_{t})=0;
   \end{eqnarray}
if  $|\gamma_t(t)|=||\gamma_t||_{0^-}\neq0$,
since
\begin{eqnarray}\label{jiaxis}
||\gamma^{h}_t||_{0}^6-|\gamma_t(t)+{h}|^6=
\begin{cases}
0,\ \ \ \ \ \ \ \ \ \  \ \ \ \ \ \ \ \ \ \ \  \ \ \ \ \  \ \ \ \
|\gamma_t(t)+h|\geq |\gamma_t(t)|,\\
  %|\gamma_t(t)-a_{\hat{t}}(\hat{t})|_8^4-|\gamma_t(t)+{he_i}-a_{\hat{t}}(\hat{t})|_8^4,
  |\gamma_t(t)|^6-|\gamma_t(t)+{h}|^6,\ \ \ \ \ |\gamma_t(t)+h|<|\gamma_t(t)|,\end{cases}
\end{eqnarray}
we have%, %if $|\gamma_t(t)-a_{\hat{t}}(\hat{t})|=||\gamma_t-a_{\hat{t}}||_{0^-}\neq0$,
\begin{eqnarray}\label{s3}
 0\leq\lim_{|h|\rightarrow0}\frac{\left|S(\gamma_t^{h})-S(\gamma_t)\right|}{|h|}
  \leq\lim_{|h|\rightarrow0}\frac{{
  ||\gamma_t(t)|^6-|\gamma_t+{h}|^6|^3}}{|h|||\gamma^{h}_t||_{0}^{12}} =0; \ \ \
   \end{eqnarray}
    if $|\gamma_t(t)|=||\gamma_t||_{0^-}=0$,
\begin{eqnarray}\label{ss4}
 \partial_{x}S(\gamma_{t})=0.
% \lim_{h\rightarrow0}\frac{S_1(\gamma_t^{he_i},a_{\hat{t}})-S_1(\gamma_t,a_{\hat{t}})}{h}
%  %\leq  \lim_{h\rightarrow0}\bigg{|}\frac{{(|\gamma_t(t)|^2-|\gamma_t(t)+{he_i}|^2)^3}}{h(1+|\gamma_t|_C^4)} \bigg{|}
%  =\lim_{h\rightarrow0}\frac{0}{h^5} =0;
   \end{eqnarray}
From (\ref{s10902jia}), (\ref{s1}), (\ref{s2}), (\ref{s3}) and (\ref{ss4}) we obtain that, for all $(t,\gamma_t)\in [0, T]\times{\hat{\Lambda}}$,
\begin{eqnarray*}
    \partial_{x}S(\gamma_t)=\begin{cases}-\frac{18\left(||\gamma_t||^6_{0}-|\gamma_t(t)|^6\right)^2|\gamma_t(t)|^4\gamma_t(t)}
   {||\gamma_t||^{12}_{0}}, \ \ \ \ \ \ \  ||\gamma_t||_{0}\neq0,\\
    0, ~~~~~~~~~~~~~~~~~~~~~~~~~~~~~~~~~~~~~~~~ ~||\gamma_t||_{0}=0.
    \end{cases}
\end{eqnarray*}
It is clear that $\partial_{x}S\in C^0(\hat{\Lambda})$. %Thus, we have show that $S_1^{a_{\hat{t}}}\in C^{1}({\hat{\Lambda}}^{\hat{t}})$.
\par
We now consider $\partial_{xx}S$.   Clearly,
   \begin{eqnarray}\label{s10902jia1}
   \partial_{xx}S(\gamma_{0})=0,\  \ \gamma_0\in \Lambda_0.
\end{eqnarray}
For every $(t,\gamma_t)\in (0,T]\times \hat{\Lambda}$, since
\begin{eqnarray*}
&&\left(||\gamma_t||^6_{0}-|\gamma_t(t)+h|^6\right)^2|\gamma_t(t)+h|^4(\gamma_t(t)+h)-\left(||\gamma_t||^6_{0}-|\gamma_t(t)|^6\right)^2|\gamma_t(t)|^4\gamma_t(t)\\
&=&\left(||\gamma_t||^6_{0}-|\gamma_t(t)+h|^6\right)^2|\gamma_t(t)+h|^4h+[4\left(||\gamma_t||^6_{0}-|\gamma_t(t)+h|^6\right)^2|\gamma_t(t)|^2(\gamma_t(t),h)\\
&&-12\left(||\gamma_t||^6_{0}-|\gamma_t(t)+h|^6\right)|\gamma_t(t)+h|^8(\gamma_t(t),h)+o(h)]\gamma_t(t),
\end{eqnarray*}
we have
if $|\gamma_t(t)|<||\gamma_t||_{0^-}$,
\begin{eqnarray}\label{s5}
  && \partial_{xx}S(\gamma_t)\nonumber\\
   &=&\frac{216\left(||\gamma_t||^6_{0}-|\gamma_t(t)|^6\right)|\gamma_t(t)|^{8}(\gamma_t(t),\cdot)\gamma_t(t)-72\left(||\gamma_t||^6_{0}-|\gamma_t(t)|^6\right)^2|\gamma_t(t)|^2(\gamma_t(t),\cdot)\gamma_t(t)}
   {{||\gamma_t||_{0}^{12}}}\nonumber\\
  % &&-\frac{24(||\gamma_t-a_{\hat{t}}||_{0,8}^6-|\gamma_t(t)-a_{\hat{t}}(\hat{t})|^6)^2(\gamma^i_t(t)-a^i_{\hat{t}}(\hat{t}))(\gamma_t^j(t)-a^j_{\hat{t}}(\hat{t}))}{{||\gamma_t-a_{\hat{t}}||_{0,8}^6}}\nonumber\\
   &&-18\frac{\left(||\gamma_t||^6_{0}-|\gamma_t(t)|^6\right)^2|\gamma_t(t)|^4I}{{||\gamma_t||_{0}^{12}}};
   \end{eqnarray}
   if $|\gamma_t(t)|>||\gamma_t||_{0^-}$,
\begin{eqnarray}\label{s4}
   \partial_{xx}S(\gamma_t)=0;
   \end{eqnarray}
 if $|\gamma_t(t)|=||\gamma_t||_{0^-}\neq0$, by (\ref{jiaxis}),
we have
\begin{eqnarray}\label{s6666}
 0\leq\lim_{|h|\rightarrow0}\bigg{|}\frac{\partial_{x}(\gamma_t^{h})-\partial_{x}S(\gamma_t)}{|h|}\bigg{|}\leq\lim_{|h|\rightarrow0}18\bigg{|}\frac{{\left(|\gamma_t(t)|^6-|\gamma_t(t)+h|^6\right)^2|\gamma_t(t)+h|^6(\gamma_t(t)+h)}}
  {|h|||\gamma^{h}_t||_{0}^{12}} \bigg{|}
  =0;
   \end{eqnarray}
    if $|\gamma_t(t)|=||\gamma_t||_{0^-}=0$,
\begin{eqnarray}\label{ss42}
\partial_{xx}S(\gamma_t)=0.
% 0\leq\lim_{h\rightarrow0}\bigg{|}\frac{S_1(\gamma_t^{he_i},a_{\hat{t}})-S_1(\gamma_t,a_{\hat{t}})}{h}\bigg{|}
%  %\leq  \lim_{h\rightarrow0}\bigg{|}\frac{{(|\gamma_t(t)|^2-|\gamma_t(t)+{he_i}|^2)^3}}{h(1+|\gamma_t|_C^4)} \bigg{|}
%  \leq\lim_{h\rightarrow0}\bigg{|}\frac{h^5}{h^3} \bigg{|}=0;
   \end{eqnarray}
%
%
% From $\partial_{x_i}S_1(\gamma_t)=0$ when $|\gamma_t(t)|\geq||\gamma_t||_{0^-}$, it follows that
   Combining (\ref{s5}), (\ref{s4}), (\ref{s6666}) and (\ref{ss42}) we obtain, for all $(t,\gamma_t)\in [\hat{t}, T]\times{\hat{\Lambda}}^{\hat{t}}$,
   \begin{eqnarray*}
    \partial_{xx}S(\gamma_t)=\begin{cases}\frac{216\left(||\gamma_t||^6_{0}-|\gamma_t(t)|^6\right)|\gamma_t(t)|^{8}(\gamma_t(t),\cdot)\gamma_t(t)-72\left(||\gamma_t||^6_{0}-|\gamma_t(t)|^6\right)^2
    |\gamma_t(t)|^2(\gamma_t(t),\cdot)\gamma_t(t)}
   {{||\gamma_t||_{0}^{12}}}\\
-18\frac{\left(||\gamma_t||^6_{0}-|\gamma_t(t)|^6\right)^2|\gamma_t(t)|^4I}{{||\gamma_t||_{0}^{12}}},~~~~~~~~~~~~~~~~~~~~~~~~~~~~~~~~~~~~~~~~~~~~~~~~~~ ||\gamma_t||_{0}\neq0,\\
    0, ~~~~~~~~~~~~~~~~~~~~~~~~~~~~~~~~~~~~~~~~~~~~~~~~~~~~~~~~~~~~~~~~~~~~~~~~~~~~~~~~~||\gamma_t||_{0}=0.
    \end{cases}
\end{eqnarray*}
It is clear that $\partial_{xx}S\in C^0(\hat{\Lambda})$. By  simple calculation, we can see that  $S$ and all of its derivatives grow  in a polynomial way.
 Thus, we have show that $S\in C^{1,2}_{p}(\hat{\Lambda})$.
\par
Now we prove (\ref{s0}). It is clear that, for every $M\geq1$,
$$
                      \Upsilon^M(\gamma_t)%+   M|\gamma_t(t)-a_{\hat{t}}(\hat{t})|_8^4
                      \leq (M+1)||\gamma_t||_{0}^6, \ \ (t,\gamma_t)\in [0, T]\times{\hat{\Lambda}}.
$$
On the other hand, for all $(t,\gamma_t)\in [0, T]\times{\hat{\Lambda}}$ and  every $M\geq1$,
$$
                       \Upsilon^M(\gamma_t)%+  M |\gamma_t(t)-a_{\hat{t}}(\hat{t})|_8^4
                       \geq \frac{M}{3}||\gamma_t||_{0}^6, \ \  \mbox{if }\ \ ||\gamma_t||_{0}^6-|\gamma_t(t)|^6\leq\frac{2}{3}||\gamma_t||_{0}^6,
$$
 and
 $$
                       \Upsilon^M(\gamma_t)%+  M |\gamma_t(t)-a_{\hat{t}}(\hat{t})|_8^4
                       \geq \frac{8}{27}||\gamma_t||_{0}^6,
                       \ \  \mbox{if }\ \ ||\gamma_t||_{0}^6-|\gamma_t(t)|^6>\frac{2}{3}||\gamma_t||_{0}^6.
$$
Thus, we  have (\ref{s0}) holds true.
 The proof is now complete. \ \ $\Box$
\par
In the proof of uniqueness of viscosity solutions, in order to apply theorem 8.3 in [2], we also need the following lemma.
%The proof of the following Lemma
Its proof is completely similar to the finite dimensional case (see Lemma 2.13 in \cite{zhou5}). Here we omit it.
\begin{lemma}\label{theoremS00044} For $M\geq3$, we have that
\begin{eqnarray}\label{jias5}
2^5\Upsilon^M(\gamma_t)+2^5\Upsilon^M(\eta_t)\geq\Upsilon^M(\gamma_t+\eta_t), \ \ (t,\gamma_t, \eta_t)\in [0,T]\times\hat{ \Lambda}\times \hat{\Lambda}.
\end{eqnarray}
\end{lemma}

\section{ Path-dependent stochastic evolution equations.}
\par
In this section, we consider the controlled state
             equation (\ref{state1}).
We introduce the admissible control. Let $t,s$ be two deterministic times, $0\leq t\leq s\leq T$.
\begin{definition}
                An admissible control process $u(\cdot)=\{u(r),  r\in [t,s]\}$ on $[t,s]$  is an ${\cal{F}}^t$-progressing measurable process taking values in some  metric space $(U,d)$. The set of all admissible controls on $[t,s]$ is denoted by ${\cal{U}}[t,s]$. We identify two processes $u(\cdot)$ and $\tilde{u}(\cdot)$ in ${\cal{U}}[t,s]$
                and write $u(\cdot)\equiv\tilde{u}(\cdot)$ on $[t,s]$, if $P(u(\cdot)=\tilde{u}(\cdot) \ a.e. \ \mbox{in}\ [t,s])=1$.
\end{definition}
First, we describe some continuous properties of the solutions of  state equation
(\ref{state1}).  We assume the following. % Let us assume that functionals
 %$ F:{\Lambda}\times U\rightarrow H$ and $ G:{\Lambda}\times U\rightarrow L_2(\Xi,H)$ satisfy the
%following assumption.
      % For any $t\in [0,T]$, denote by ${\cal{H}}^2(t,T)$ the space of all ${{\cal{F}}}^t$-adapted, $R^d$-valued  processes $(Y(s))_{t\leq s\leq T}$ such that
%       $||Y||^2=E[\int^{T}_{t}|Y(s)|^2ds]<\infty$ and by ${\cal{S}}^2(t,T)$ the space of all ${{\cal{F}}}^t$-adapted, $R$-valued continuous processes
%       $(Y(s))_{t\leq s\leq T}$ such that
%       $||Y||^2=E[\sup_{t\leq s\leq T}|Y(s)|^2]<\infty$.
\begin{hyp}\label{hypstate}
\begin{description}
     \item{(i)}
    The operator $A$ is the generator of
         a $C_0$  semigroup  $\{e^{tA}, t\geq0\}$ of bounded linear operator in Hilbert space
        $H$.
         \item{(i')}
                 The operator $A$ is the generator of a $C_0$ contraction
                 semigroup $\{e^{tA}, t\geq0\}$ of bounded linear operators in the
                 Hilbert space $H$.
                 \par
        \item{(ii)} $F:{\Lambda}\times U\rightarrow H$ and $G:{\Lambda}\times U\rightarrow L_2(\Xi,H)$ are continuous, and
        %For every fixed $\gamma_t\in\Lambda$, $F(\gamma_t,\cdot)$ and
%        $G(\gamma_t,\cdot)$ are continuous in $u$.
%\par
%       \item{(iii)}
                 there exists a constant $L>0$ such that, for all $(t,\gamma_t,u)$,  $ (s,\eta_t,u) \in [0,T]\times {\Lambda}\times U$,
                 %$(t,\gamma_t,\zeta_T,u)$,  $ (s, \eta_s,\zeta'_Tu) \in [0,T]\times {\Lambda}\times{\Lambda_T}\times U$,
      \begin{eqnarray}\label{assume1111}
                &&|F(\gamma_t,u)|^2\vee|G(\gamma_t,u)|^2_{L_2(\Xi,H)}\leq
                 L^2(1+||\gamma_t||^2_0),\nonumber\\
                 &&|F(\gamma_t,u)-F(\eta_t,u)|\vee|G(\gamma_t,u)-G(\eta_t,u)|_{L_2(\xi,H)}\leq
                 L||\gamma_t-\eta_t||_0.
             \end{eqnarray}
             % \begin{eqnarray}\label{asumme2222}
%               |q(\gamma_t,u)-q(\eta_s,u)|\leq Ld_\infty(\gamma_t,\eta_s), \ \ |q(\gamma_t,u)|\leq L(1+||\gamma_t||_0);
%\end{eqnarray}
% \begin{eqnarray}\label{asumme3333}
%                 |\phi(\zeta_T)-\phi(\zeta'_T)|\leq L||\zeta-\zeta'_T||_0, \ \ \ |\phi(\zeta_T)|\leq L(1+||\zeta||_0).
%\end{eqnarray}
\end{description}
\end{hyp}
\par
              We say that $X$ is a mild solution of equation $(\ref{state1})$ with initial data $\xi_t\in L_{\cal{P}}^p(\Omega;C([0,t];H))$ if it is a continuous, $\{{\cal{F}}_t\}_{t\geq0}$-predictable process with values in $H$, and it satisfies:  $P$-a.s.,
\begin{eqnarray*}
            X(s)=e^{(s-t)A}\xi_t(t)+\int_{t}^{s}{e^{(s-\sigma)A}}F(X_\sigma,u(\sigma))d\sigma+\int_{t}^{s}{e^{(s-\sigma)A}}G(X_\sigma,u(\sigma))dW(\sigma),  \ s\in [t,T],
\end{eqnarray*}
where $X(s)=\xi_t(s),  \ s\in[0,t)$. To emphasize dependence on initial data and control, we denote the solution by $X^{\xi_t,u}(\cdot)$. Note that, if $u(\cdot)\in {\mathcal{U}}[t,T]$ and $\xi_t=\gamma_t \in \Lambda_t$, then $X^{\gamma_t,u}(s)$ is
${\cal{F}}^T_t$-measurable, hence independent of ${\cal{F}}_t$.
\par
The following lemma is standard; see, for example, Theorem 3.6 in \cite{ro}. We include the proof for completeness and because it will be useful in the following.
 % but we do not find it in the existing literature. For the convenience of readers, we give its proof.
\begin{lemma}\label{lemmaexist}
\ \ Assume that Hypothesis \ref{hypstate} (i) and (ii)  hold. Then for every $p>2$, $u(\cdot)\in {\cal{U}}[0,T]$,
$\xi_t\in L^{p}_{\mathcal{P}}(\Omega;C([0,t];H))$, (\ref{state1}) admits a
unique mild solution $X^{\xi_t,u}$.  Moreover, if we let  $X^{\xi'_t,u}$  be the solutions of  (\ref{state1})
 corresponding $\xi'_t\in L^{p}_{\mathcal{P}}(\Omega;C([0,t];H))$ and $u(\cdot)\in {\cal{U}}[0,T]$. Then the following estimates hold:
\begin{eqnarray}\label{state1est}
               %\mathbb{E}||X_T^{\xi_t,u}-X_T^{\xi'_t,u}||^p_0\leq C_p\mathbb{E}||\xi_t-\xi'_t||^p_0,\ \ \ \
                \mathbb{E}||X_T^{\xi_t,u}||^p_0\leq C_p(1+\mathbb{E}||\xi_t||^p_0).
\end{eqnarray}
              The constant $C_p$ depending only on  $p$, $T$, $L$ and
              $M_1=:\sup_{s\in [0,T]}|e^{sA}|$.
\par
Finally,  let $A_\mu=\mu A(\mu I-A)^{-1}$ be the Yosida approximation of $A$ and
let $X^\mu$ be the solution of the following:
\begin{eqnarray}\label{07162}
X^\mu(s)=e^{(s-t)A_\mu}\xi_t(t)+\int^{s}_{t}e^{(s-\sigma)A_\mu}F(X^{\mu}_\sigma,u(\sigma))d\sigma+\int^{s}_{t}e^{(s-\sigma)A_\mu}G(X^{\mu}_\sigma,u(\sigma))dW(\sigma),  s\in [t,T]; \ \
\end{eqnarray}
and $  X^\mu(s)=\xi_t(s), \ s\in[0,t)$.
Then
\begin{eqnarray}\label{0717}
                       \lim_{\mu\rightarrow\infty}\mathbb{E}\sup_{s\in [t,T]}|X^{\xi_t,u}(s)-X^\mu(s)|^p=0.
\end{eqnarray}
\end{lemma}
{\bf  Proof  }. \ \
We define a mapping $\Phi$ from $L^p_{\cal{P}}(\Omega; C([0,T];H))$ to itself by the formula
\begin{eqnarray*}
                 \Phi(X)(s)=e^{(s-t)A}\xi_t(t)+\int^{s}_{t}e^{(s-\sigma)A}F(X_\sigma,u(\sigma))d\sigma+\int^{s}_{t}e^{(s-\sigma)A}G(X_\sigma,u(\sigma))dW(\sigma),\  s\in [t,T],
\end{eqnarray*}
$$
                 \Phi(X)(s)=\xi_t(s),\ \ s\in [0,t),
$$
and show that it is a contraction, under an equivalent norm $||X||=\left(\mathbb{E}\sup_{s\in[0,T]}e^{-\beta sp}|X(s)|^p\right)^{\frac{1}{p}}$, where $\beta>0$ will be chosen later.
\par
We will use the so called factorization method (see Theorem 5.2.5 in \cite{da}). Let us take %$p>2$ and
 $\alpha\in (0,1)$ such that
$$
                      \frac{1}{p}<\alpha<\frac{1}{2}\ \ \mbox{and let}\ \ c_{\alpha}^{-1}=\int^{s}_{\sigma}(s-l)^{\alpha-1}(l-\sigma)^{-\alpha}dl.
$$
Then by the Fubini theorem and stochastic Fubini theorem, for $s\in [t,T]$,
\begin{eqnarray*}
                 \Phi(X)(s)&=&e^{(s-t)A}\xi_t(t)+c_{\alpha}\int^{s}_{t}\int^{s}_{\sigma}(s-l)^{\alpha-1}(l-\sigma)^{-\alpha}e^{(s-l)A}e^{(l-\sigma)A}dlF(X_\sigma,u(\sigma))d\sigma\\
                 &&+c_{\alpha}\int^{s}_{t}\int^{s}_{\sigma}(s-l)^{\alpha-1}(l-\sigma)^{-\alpha}e^{(s-l)A}e^{(l-\sigma)A}dlG(X_\sigma,u(\sigma))dW(\sigma)\\
                % &=&e^{(s-t)A}\xi_t(t)+c_{\alpha}\int^{s}_{t}(s-l)^{\alpha-1}e^{(s-l)A}\int^{l}_{t}(l-\sigma)^{-\alpha}e^{(l-\sigma)A} F(X_\sigma,u(\sigma))d\sigma dl\\
%                 &&+c_{\alpha}\int^{s}_{t}(s-l)^{\alpha-1}e^{(s-l)A}\int^{l}_{t}(l-\sigma)^{-\alpha}e^{(l-\sigma)A} G(X_\sigma,u(\sigma))dW(\sigma) dl\\
                 &=&e^{(s-t)A}\xi_t(t)+c_{\alpha}\int^{s}_{t}(s-l)^{\alpha-1}e^{(s-l)A}Y(l)dl,
\end{eqnarray*}
where
$$
Y(l)=\int^{l}_{t}(l-\sigma)^{-\alpha}e^{(l-\sigma)A} F(X_\sigma,u(\sigma))d\sigma+\int^{l}_{t}(l-\sigma)^{-\alpha}e^{(l-\sigma)A} G(X_\sigma,u(\sigma))dW(\sigma).
$$
By the H$\ddot{o}$lder inequality, setting %$M:=\sup_{s\in [0,T]}|e^{sA}|$,
$q=\frac{p}{p-1}$,
\begin{eqnarray*}
                 &&e^{-\beta s}\bigg{|}\int^{s}_{t}(s-l)^{\alpha-1}e^{(s-l)A}Y(l)dl\bigg{|}\\
                 &\leq& \bigg{(}\int^{s}_{t}e^{-q\beta (s-l)}(s-l)^{(\alpha-1)q}dl\bigg{)}^{\frac{1}{q}}
                 \bigg{(}\int^{s}_{t}e^{-p\beta l}\left|e^{(s-l)A}Y(l)\right|^pdl\bigg{)}^{\frac{1}{p}}\\
                 &\leq& M_1\bigg{(}\int^{T}_{0}e^{-q\beta l}l^{(\alpha-1)q}dl\bigg{)}^{\frac{1}{q}}
                 \bigg{(}\int^{T}_{t}e^{-p\beta l}|Y(l)|^pdl\bigg{)}^{\frac{1}{p}}.
\end{eqnarray*}
Then we get
\begin{eqnarray*}
                 ||\Phi(X)||\leq M_1(\mathbb{E}||\xi_t||^p_0)^{\frac{1}{p}}+M_1c_{\alpha}\bigg{(}\int^{T}_{0}e^{-q\beta l}l^{(\alpha-1)q}dl\bigg{)}^{\frac{1}{q}}
                 \bigg{(}\mathbb{E}\int^{T}_{t}e^{-p\beta l}|Y(l)|^pdl\bigg{)}^{\frac{1}{p}}.
\end{eqnarray*}
By the Burkholder-Davis-Gundy
                              inequalities and  Hypothesis \ref{hypstate} (i) and (ii), there exists a constant $c_p$ depending only on $p$  and it may vary from line to line such that
\begin{eqnarray*}
                &&\mathbb{E}|Y(l)|^p\\
                &\leq& c_p\mathbb{E}\bigg{(}\int^{l}_{t}(l-\sigma)^{-\alpha}\left|e^{(l-\sigma)A} F(X_\sigma,u(\sigma))\right|d\sigma\bigg{)}^p\\
                &&+c_p\mathbb{E}\bigg{(}\int^{l}_{t}(l-\sigma)^{-2\alpha}\left|e^{(l-\sigma)A} G(X_\sigma,u(\sigma))\right|^2_{L_2(\Xi,H)}d\sigma\bigg{)}^{\frac{p}{2}}\\
                &\leq&c_pM_1^pL^p\mathbb{E}\bigg{(}\int^{l}_{t}(l-\sigma)^{-\alpha}(1+||X_\sigma||_0)d\sigma\bigg{)}^p
                +c_pM_1^pL^p\mathbb{E}\bigg{(}\int^{l}_{t}(l-\sigma)^{-2\alpha}(1+||X_\sigma||_0)^2d\sigma\bigg{)}^{\frac{p}{2}}\\
                &\leq&c_pM_1^pL^p\mathbb{E}\sup_{\sigma\in [0,l]}[(1+||X_\sigma||_0)^pe^{-p\beta\sigma}]\bigg{[}\bigg{(}\int^{l}_{t}(l-\sigma)^{-\alpha}e^{\beta\sigma}d\sigma\bigg{)}^p
                +\bigg{(}\int^{l}_{t}(l-\sigma)^{-2\alpha}e^{2\beta\sigma}d\sigma\bigg{)}^{\frac{p}{2}}\bigg{]},
\end{eqnarray*}
which implies
\begin{eqnarray*}
                e^{-p\beta l}\mathbb{E}|Y(l)|^p
                %&\leq&c_pM_1^pL^p(1+||X||^p)\bigg{[}\bigg{(}\int^{l}_{t}(l-\sigma)^{-\alpha}e^{-\beta(l-\sigma)}d\sigma\bigg{)}^p
%                +\bigg{(}\int^{l}_{t}(l-\sigma)^{-2\alpha}e^{-2\beta(l-\sigma)}d\sigma\bigg{)}^{\frac{p}{2}}\bigg{]}\\
                \leq c_pM_1^pL^p(1+||X||^p)\bigg{[}\bigg{(}\int^{T}_{0}\sigma^{-\alpha}e^{-\beta\sigma}d\sigma\bigg{)}^p
                +\bigg{(}\int^{T}_{0}\sigma^{-2\alpha}e^{-2\beta\sigma}d\sigma\bigg{)}^{\frac{p}{2}}\bigg{]}.
\end{eqnarray*}
We conclude that
\begin{eqnarray*}
                 ||\Phi(X)||&\leq& M_1(\mathbb{E}||\xi_t||^p_0)^{\frac{1}{p}}+M_1^2Lc_{\alpha}(Tc_p(1+||X||^p))^{\frac{1}{p}}\bigg{(}\int^{T}_{0}e^{-q\beta l}l^{(\alpha-1)q}dl\bigg{)}^{\frac{1}{q}}\\
                 &&\times
                 \bigg{[}\bigg{(}\int^{T}_{0}\sigma^{-\alpha}e^{-\beta\sigma}d\sigma\bigg{)}^p
                +\bigg{(}\int^{T}_{0}\sigma^{-2\alpha}e^{-2\beta\sigma}d\sigma\bigg{)}^{\frac{p}{2}}\bigg{]}^{\frac{1}{p}}.
\end{eqnarray*}
 This show that $\Phi$ is a well defined mapping on $L^p_{\mathcal{P}}(\Omega; C([0,T];H))$. If $X,X'$ are processes belonging to this space, similar passages show that
 \begin{eqnarray*}
                 ||\Phi(X)-\Phi(X')||&\leq& M_1^2Lc_{\alpha}(Tc_p)^{\frac{1}{p}}||X-X'||\bigg{(}\int^{T}_{0}e^{-q\beta l}l^{(\alpha-1)q}dl\bigg{)}^{\frac{1}{q}}\\
                 &&\times
                 \bigg{[}\bigg{(}\int^{T}_{0}\sigma^{-\alpha}e^{-\beta\sigma}d\sigma\bigg{)}^p
                +\bigg{(}\int^{T}_{0}\sigma^{-2\alpha}e^{-2\beta\sigma}d\sigma\bigg{)}^{\frac{p}{2}}\bigg{]}^{\frac{1}{p}}.
\end{eqnarray*}
Therefore, for $\beta$ sufficiently large, the mapping is a contraction. In particular, we obtain $||X^{\xi_t,u}||\leq C_p(1+(\mathbb{E}||\xi_t||^p_0)^{\frac{1}{p}})$  %and $||X^{\xi_t,u}-X^{\xi'_t,u}||\leq C_p(\mathbb{E}||\xi_t-\xi'_t||^p_0)^{\frac{1}{p}}$,
 which prove the
estimate (\ref{state1est}).
\par
Finally, for some constant $c_p$ depending only on $p$ that may be vary from line to line,
\begin{eqnarray}\label{0716}
                            &&\mathbb{E}\sup_{r\in[0,s]}\left|X^{\xi_t,u}(r)-X^\mu(r)\right|^p\nonumber\\
                           &\leq&c_p\bigg{(}\mathbb{E}\int^{s}_{t}\sup_{\theta\in[0,T]}\left|(e^{\theta A}-e^{\theta A_\mu})F(X^{\xi_t,u}_\sigma,u(\sigma))\right|^pd\sigma
                           +L^pM^p_1\mathbb{E}\int^{s}_{t}\sup_{\theta\in[0,\sigma]}\left|X^{\xi_t,u}(\theta)-X^\mu(\theta)\right|^pd\sigma\nonumber\\
                            &&+\mathbb{E}\sup_{r\in[t,s]}\bigg{|}\int^{r}_{t}(e^{(r-\sigma)A}G(X^{\xi_t,u}_\sigma,u(\sigma))-e^{(r-\sigma)A_\mu}G(X^\mu_\sigma,u(\sigma)))dW(\sigma)\bigg{|}^p\bigg{)}\nonumber\\
                            &&+\mathbb{E}\sup_{r\in[0,s]}\left|(e^{rA}-e^{rA_\mu})\xi_t(t)\right|^p.
\end{eqnarray}
                             On the other hand, by the stochastic Fubini
                               theorem,
\begin{eqnarray*}
                                 &&\int^{r}_{t}\left(e^{(r-\sigma)A}G(X^{\xi_t,u}_\sigma,u(\sigma))-e^{(r-\sigma)A_\mu}G(X^\mu_\sigma,u(\sigma))\right)dW(\sigma)\\
                                % &=&c_\alpha\int^{r}_{0}\int^{r}_{\sigma}(r-\theta)^{\alpha-1}(\theta-\sigma)^{-\alpha}
%                                 (e^{(r-\theta)A}e^{(\theta-\sigma)A}G(X(\sigma))-e^{(r-\theta)A_\mu}e^{(\theta-\sigma)A_\mu}
%                                 G(X(\sigma)))d\theta dW(\sigma)\\
                                 &=&c_\alpha\int^{r}_{t}(r-\theta)^{\alpha-1}\left(e^{(r-\theta)A}Y(\theta)-e^{(r-\theta)A_\mu}Y^\mu(\theta)\right)d\theta,
\end{eqnarray*}
                              where
$$
                             Y(\theta)=\int^{\theta}_{t}(\theta-\sigma)^{-\alpha}e^{(\theta-\sigma)A}G(X^{\xi_t,u}_\sigma,u(\sigma))dW(\sigma),
                             $$
                             and
                             $$
                             Y^\mu(\theta)=\int^{\theta}_{t}(\theta-\sigma)^{-\alpha}e^{(\theta-\sigma)A_\mu}G(X^\mu_\sigma,u(\sigma))dW(\sigma).
$$
                              By the H\"{o}lder inequality and the Burkholder-Davis-Gundy
                             inequalities, we have
\begin{eqnarray}\label{07161}
                            &&\mathbb{E}\sup_{r\in[t,s]}\bigg{|}\int^{r}_{t}\left(e^{(r-\sigma)A}G(X^{\xi_t,u}_\sigma,u(\sigma))-e^{(r-\sigma)A_\mu}
                                              G(X^\mu_\sigma,u(\sigma))\right)dW(\sigma)\bigg{|}^p\nonumber\\
                          %  &=&c^p_\alpha
%                             E\sup_{r\in[0,s]}|\int^{r}_{0}(r-\theta)^{\alpha-1}(e^{(r-\theta)A}Y(\theta)-e^{(r-\theta)A_\mu}Y^\mu(\theta))d\theta|^p\\
                           %&\leq&c^p_\alpha\bigg{(}\int^{T}_{0}r^{q(\alpha-1)}dr\bigg{)}^{\frac{p}{q}}
%                             E
%                                               \int^{s}_{0}|e^{(r-\theta)A}Y(\theta)-e^{(r-\theta)A_\mu}Y^\mu(\theta)|^pd\theta\\
                             &\leq&c_pc^p_\alpha\bigg{(}\int^{T}_{0}r^{q(\alpha-1)}dr\bigg{)}^{\frac{p}{q}}\bigg{(}\mathbb{E}\int^{s}_{t}\sup_{r\in[0,T]}\left|\left(e^{rA}-e^{rA_\mu}\right)Y(\theta)\right|^pd\theta
                             +M^p_1\mathbb{E}\int^{s}_{t}|Y(\theta)-Y^\mu(\theta)|^pd\theta\bigg{)}\nonumber\\
                             &\leq&c_pc^p_\alpha\bigg{(}\int^{T}_{0}r^{q(\alpha-1)}dr\bigg{)}^{\frac{p}{q}}\bigg{(}\mathbb{E}\int^{s}_{t}\sup_{r\in[0,T]}\left|\left(e^{rA}-e^{rA_\mu}\right)Y(\theta)\right|^pd\theta
                             \nonumber\\
                             &&+
                                                       M_1^{2p}L^p\bigg{(}\int^{T}_{0}r^{-2\alpha}dr\bigg{)}^{\frac{p}{2}}\mathbb{E}\int^{s}_{t}\sup_{\sigma\in[0,\theta]}|X^{\xi_t,u}(\sigma)-X^\mu(\sigma)|^pd\theta\nonumber\\
                             &&+\mathbb{E}\int^{s}_{t}\bigg{(}\int^{\theta}_{t}(\theta-\sigma)^{-2\alpha}\sup_{r\in[0,T]}\left|\left(e^{rA}-e^{rA_\mu}\right)G(X^{\xi_t,u}_\sigma,u(\sigma))\right|^2
                             d\sigma\bigg{)}^{\frac{p}{2}}d\theta\bigg{)}.
\end{eqnarray}
                             From (\ref{0716}), (\ref{07161}) and  applying  Gronwall's inequality  and the Dominated
                             Convergence Theorem, it follows that
                                  $\mathbb{E}\sup_{0\leq s\leq
                                  T}|X^{\xi_t,u}(s)-X^\mu(s)|^p\rightarrow
                                  0$ as $\mu\rightarrow \infty$. \ \ $\Box$

  \par
  The next result contains the local boundedness and the continuity of the trajectory $X^{\xi_t,u}$.% and  value functional $V$. %In what follows, $C$ is an absolute
%constant, that can be different in different places.
  \begin{lemma}\label{lemmaexist111}
\ \ Assume that Hypothesis \ref{hypstate}  (i) and (ii)  hold. Then, for any $p>2$, $0\leq t\leq \bar{t}\leq T$, $\xi_t, \xi'_t\in L_{\cal{P}}^p(\Omega;C([0,t];H))$ and $u(\cdot)\in {\cal{U}}[0,T]$,
\begin{eqnarray}\label{2.6}
            \sup_{u(\cdot)\in {\cal{U}}[0,T]}\mathbb{E}\left|X^{\xi_t,u}(s)-e^{(s-t)A}\xi_t(t)\right|^p\leq
            C^1_p\left(1+\mathbb{E}||\xi_t||_0^p\right)|s-t|^{\frac{p}{2}}, \ \ \ s\in[t,T];
\end{eqnarray}
\begin{eqnarray}\label{0604}
\mathbb{E}\left|\left|X^{\xi'_t,u}_T-X^{\xi_{t,\bar{t},A},u}_T\right|\right|^p_0\leq C^1_p\left(1+\mathbb{E}||\xi'_t||^p_0\right)(\bar{t}-t)^{\frac{p}{2}}+C^1_p\mathbb{E}||\xi'_t-\xi_t||^p_0;
\end{eqnarray}
%\begin{eqnarray}\label{3.5}
% |V(\xi_t)|\leq C(1+||\xi_t||_0);
% \end{eqnarray}
%\begin{eqnarray}\label{0604001}
%|V(\xi_{t,\bar{t},A})-V(\xi'_t)|
%                        \leq C(1+||\xi'_t||_0)(\bar{t}-t)+C||\xi'_t-\xi_t||_0.
%\end{eqnarray}
 The constant $C^1_p$ depending only on  $p$, $T$, $L$ and $M_1$.
\end{lemma}
{\bf  Proof  }. \ \
 For any  $\xi_t\in\Lambda$, by (\ref{assume1111}) and (\ref{state1est}), we
                obtain the following result:
\begin{eqnarray*}
                   \mathbb{E}\left|X^{\xi_t,u}(s)-e^{(s-t)A}\xi_t(t)\right|^p
               \leq c_pL^pM^p_1[1+C_p\left(1+\mathbb{E}||\xi_t||^p_0\right)]|s-t|^{\frac{p}{2}}\left(1+|s-t|^{\frac{p}{2}}\right).
\end{eqnarray*}
Here and in the rest of this proof, $c_p$ denotes a positive constant, whose value depend only on $p$ and may vary from line to line.
    Taking the supremum in ${\cal{U}}[0,T]$, we obtain (\ref{2.6}).
  For any $0\leq t\leq \bar{t}\leq T$, $\xi_t, \xi'_t\in {\Lambda}$ and $u(\cdot)\in {\cal{U}}[0,T]$, by (\ref{assume1111}) and (\ref{state1est}), we have
  %for some absolute constant $c_p$ depending only on $p$, that can be different in different places,
%\begin{eqnarray*}
%             |X^{\xi_t,u}(s)-e^{A(s-t)}\xi_t(t))|\leq \int^{s}_{t}L(1+||X^{\xi_t,u}_r||_0)dr\leq C(1+||\xi_t||_0)(s-t).
%\end{eqnarray*}
\begin{eqnarray*}
             &&\mathbb{E}\sup_{\bar{t}\leq s\leq \sigma}\left|X^{\xi'_t,u}(s)-X^{\xi_{t,\bar{t},A},u}(s)\right|^p\\
             &\leq& c_pM^p_1\mathbb{E}|\xi'_t(t)-\xi_t(t)|^p+%\sup_{\bar{t}\leq s\leq \sigma}
             c_p\mathbb{E}\bigg{(}\int^{\bar{t}}_{t}\left|e^{(s-r)A}F(X^{\xi'_t,u}_r,u(r))\right|dr\bigg{)}^p\\
             &&+c_p\mathbb{E}\bigg{|}\int^{\bar{t}}_{t}e^{(s-r)A}G(X^{\xi'_t,u}_r,u(r))dW(r)\bigg{|}^p\\
             &&
                                                             +c_p\mathbb{E}\bigg{(}\sup_{\bar{t}\leq s\leq \sigma}\int^{s}_{\bar{t}}\left|e^{(s-r)A}\left(F(X^{\xi'_t,u}_r,u(r))-F(X^{\xi_{t,\bar{t},A},u}_r,u(r))\right)\right|dr\bigg{)}^p\\
             &&
                                                             +c_p\mathbb{E}\sup_{\bar{t}\leq s\leq \sigma}\bigg{|}\int^{s}_{\bar{t}}e^{(s-r)A}\left(G(X^{\xi'_t,u}_r,u(r))-G(X^{\xi_{t,\bar{t},A},u}_r,u(r))\right)dW(r)\bigg{|}^p\\
                                                            % &\leq& M|\xi'_t(t)-\xi_t(t)|+LM\int^{\bar{t}}_{t}(1+||X^{\xi'_t,u}_r||_0)dr
%                                                             +LM\int^{\sigma}_{\bar{t}}||X^{\xi'_t,u}_r-X^{\xi_{t,\bar{t},A},u}_r||_0dr\\
                                                              &\leq&c_pM^p_1\mathbb{E}|\xi'_t(t)-\xi_t(t)|^p+c_pL^pM^p_1[1+C_p(1+\mathbb{E}||\xi'_t||^p_0)](\bar{t}-t)^{\frac{p}{2}}\left(1+(\bar{t}-t)^{\frac{p}{2}}\right)\\
                                                              &&
                                                             +c_pL^pM^p_1\int^{\sigma}_{\bar{t}}\mathbb{E}\left|\left|X^{\xi'_t,u}_r-X^{\xi_{t,\bar{t},A},u}_r\right|\right|^p_0dr\\
                                                             &&+c_p\mathbb{E}\sup_{\bar{t}\leq s\leq \sigma}\bigg{|}\int^{s}_{\bar{t}}e^{(s-r)A}\left(G(X^{\xi'_t,u}_r,u(r)))-G(X^{\xi_{t,\bar{t},A},u}_r,u(r)\right)dW(r)\bigg{|}^p.
\end{eqnarray*}
Proceeding as in the proof of Lemma \ref{lemmaexist} we get  for $\frac{1}{p}<\alpha<\frac{1}{2}$ and for a suitable constant $c_p$
\begin{eqnarray*}
             &&\mathbb{E}\bigg{|}\sup_{\bar{t}\leq s\leq \sigma}\int^{s}_{\bar{t}}e^{(s-r)A}\left(G(X^{\xi'_t,u}_r,u(r)))-G(X^{\xi_{t,\bar{t},A},u}_r,u(r)\right)dW(r)\bigg{|}^p \\
             &\leq&c_pc^p_\alpha\bigg{(}\int^{T}_{0}r^{q(\alpha-1)}dr\bigg{)}^{\frac{p}{q}}
                                                       M_1^{2p}L^p\bigg{(}\int^{T}_{0}r^{2\alpha}dr\bigg{)}^{\frac{p}{2}}\int^{\sigma}_{\bar{t}}\mathbb{E}\left|\left|X^{\xi'_t,u}_r-X^{\xi_{t,\bar{t},A},u}_r\right|\right|^p_0dr.
\end{eqnarray*}
Thus,
%\begin{eqnarray*}
%             ||X^{\xi'_t,u}_\sigma-X^{\xi_{t,\bar{t},A},u}_\sigma||_0
%                                                              \leq M_1||\xi'_t-\xi_t||_0+LM_1(1+C_1(1+||\xi'_t||_0))(\bar{t}-t)
%                                                             +LM_1\int^{\sigma}_{\bar{t}}||X^{\xi'_t,u}_r-X^{\xi_{t,\bar{t},A},u}_r||_0dr.
%\end{eqnarray*}
\begin{eqnarray*}
             \mathbb{E}\left|\left|X^{\xi'_t,u}_\sigma-X^{\xi_{t,\bar{t},A},u}_\sigma\right|\right|^p_0
                                                              &\leq& c_pM_1^p\mathbb{E}||\xi'_t-\xi_t||^p_0+c_pL^pM^p_1(1+T^{\frac{p}{2}})[1+C_p(1+\mathbb{E}||\xi'_t||^p_0)](\bar{t}-t)^{\frac{p}{2}}\\
                                                             &&+c_pL^pM_1^p(1+c^p_\alpha M_1^p)\int^{\sigma}_{\bar{t}}\mathbb{E}\left|\left|X^{\xi'_t,u}_r-X^{\xi_{t,\bar{t},A},u}_r\right|\right|^p_0dr.
\end{eqnarray*}
Then, by Gronwall's inequality, we obtain (\ref{0604}). %Next, by (\ref{asumme2222}), (\ref{asumme3333}) and (\ref{0604}), we get
The lemma  is proved. \ \  $\Box$
\par
               By Lemmas \ref{theoremito}, % and \ref{theoremito2},
               we get the following
                lemma which is needed to prove the existence  of viscosity solutions.
\begin{lemma}\label{theoremito3}
\ \
Suppose $X^{\gamma_t,u}$ is a solution of (\ref{state1}) with initial data $\gamma_t\in \Lambda_t$, $f\in C_p^{1,2}({\Lambda}^{{\hat{t}}})$  and $A^*\partial_xf\in C_p^0({\Lambda}^{{\hat{t}}})$
for some $\hat{t}\in[t,T)$. Then, $P$-a.s., for any $s\in [\hat{t},T]$:
\begin{eqnarray}\label{statesop03}
                 &&f\left(X^{\gamma_t,u}_s\right)=f\left(X^{\gamma_t,u}_{\hat{t}}\right)+\int_{\hat{t}}^{s}\partial_tf\left(X^{\gamma_t,u}_\sigma\right)d\sigma
                 +\int^{s}_{\hat{t}}\left(A^*\partial_xf\left(X^{\gamma_t,u}_\sigma\right), X^{\gamma_t,u}\left(\sigma\right)\right)_H\nonumber\\
                 &&+\left(\partial_xf\left(X^{\gamma_t,u}_\sigma\right),F\left(X^{\gamma_t,u}_\sigma,u(\sigma)\right)\right)_H+\frac{1}{2}\mbox{Tr}\left(\partial_{xx}f\left(X^{\gamma_t,u}_\sigma\right)
                 G\left(X^{\gamma_t,u}_\sigma,u\left(\sigma\right)\right)G^*\left(X^{\gamma_t,u}_\sigma,u\left(\sigma\right)\right)\right)d\sigma\nonumber\\
                 &&+\int^{s}_{\hat{t}}\left(\partial_xf\left(X^{\gamma_t,u}_\sigma\right),G\left(X^{\gamma_t,u}_\sigma,u(\sigma)\right)dW(\sigma)\right)_H.
\end{eqnarray}
\end{lemma}

\par

                The following
                lemma is also needed to prove the existence  of viscosity solutions.
%\begin{lemma}\label{theoremS0000815}  Assume the Hypothesis \ref{hypstate} holds true and $g\in {\mathcal{G}}_t$ for some $t\in [0,T)$,  %and $g\in C^1(R)$ with $g'\geq0$,
%we have, $P$-a.s.,
%\begin{eqnarray}\label{jias510815jia}
%                  g\left(X^{\gamma_t,u}_s\right)
%                  &\leq& g\left(X^{\gamma_t,u}_t\right)+\int_{\hat{t}}^{s}\partial_t\varphi\left(X^{\gamma_t,u}_\sigma\right)d\sigma+\int^{s}_{{t}}\left(\partial_xg\left(X^{\gamma_t,u}_\sigma\right),
%                  F\left(X^{\gamma_t,u}_\sigma,u\left(\sigma\right)\right)\right)_H
%                  \nonumber\\
%                  &&+\frac{1}{2}\mbox{Tr}\left(\partial_{xx}g\left(X^{\gamma_t,u}_\sigma\right)G\left(X^{\gamma_t,u}_\sigma,u\left(\sigma\right)\right)
%                  G^*\left(X^{\gamma_t,u}_\sigma,u\left(\sigma\right)\right)\right)d\sigma\nonumber\\
%                  &&
%                  +\int^{s}_{t}\left(\partial_xg\left(X^{\gamma_t,u}_\sigma\right),G\left(X^{\gamma_t,u}_\sigma,u\left(\sigma\right)\right)dW\left(\sigma\right)\right)_H, \ \ s\in [t,T].
%\end{eqnarray}
%\end{lemma}\
%\par
\begin{lemma}\label{theoremito2}
\ \
Assume the Hypothesis \ref{hypstate} holds true. For every $t\in [0,T)$, $\eta_t\in \Lambda$ and $M\geq3$,  $P$-a.s.,  for all $s\in [t,T]$:
\begin{eqnarray}\label{jias510815jia}
                  &&\Upsilon^M(X^{\gamma_t,u}_s-\eta_{t,s,A})\nonumber\\
                  &\leq& \Upsilon^M(X^{\gamma_t,u}_t-\eta_t)+\int^{s}_{{t}}(\partial_x\Upsilon^M(X^{\gamma_t,u}_\sigma-\eta_{t,\sigma,A}),F(X^{\gamma_t,u}_\sigma,u(\sigma))_H\nonumber\\
                  &&+\frac{1}{2}\mbox{Tr}(\partial_{xx}\Upsilon^M(X^{\gamma_t,u}_\sigma-\eta_{t,\sigma,A})G(X^{\gamma_t,u}_\sigma,u(\sigma))
                  G^*(X^{\gamma_t,u}_\sigma,u(\sigma)))d\sigma\nonumber\\
                  &&
                  +\int^{s}_{t}(\partial_x\Upsilon^M(X^{\gamma_t,u}_\sigma-\eta_{t,\sigma,A}),G(X^{\gamma_t,u}_\sigma,u(\sigma))dW(\sigma))_H;
\end{eqnarray}
and 
\begin{eqnarray*}\label{jias510815jia11}
                  &&|X^{\gamma_t,u}(s)-e^{(s-t)A}\eta_{t}(t)|^6=|y(s)|^6\nonumber\\
                  &\leq& |y(t)|^6+6\int^{s}_{{t}}|y(s)|^4(y(s),F(X^{\gamma_t,u}_\sigma,u(\sigma))_H\nonumber\\
                  &&+\frac{1}{2}\mbox{Tr}((6|y(s)|^4I+24|y(s)|^2y(s)(y(s),\cdot))G(X^{\gamma_t,u}_\sigma,u(\sigma))
                  G^*(X^{\gamma_t,u}_\sigma,u(\sigma)))d\sigma\nonumber\\
                  &&
                  +6\int^{s}_{t}|y(s)|^4(y(s),G(X^{\gamma_t,u}_\sigma,u(\sigma))dW(\sigma))_H,
\end{eqnarray*}
where
 $y(s)=X^{\gamma_t,u}(s)-e^{A(s-t)}\eta_t(t),\  t\leq s\leq T$ and $y(s)=\gamma_t(s)-\eta_t(s),\  0\leq s<t$.
\end{lemma}
\par
   {\bf  Proof  }. \ \ %Let $A_\mu=\mu A(\mu I-A)^{-1}$ be the Yosida approximation of $A$ and
%let $X^\mu$ be the solution of the following:
%$$
%X^\mu(s)=e^{(s-t)A_\mu}\gamma_t(t)+\int^{s}_{t}e^{(s-\sigma)A_\mu}F(X^{\mu}_\sigma,u(\sigma))d\sigma+\int^{s}_{t}e^{(s-\sigma)A_\mu}G(X^{\mu}_\sigma,u(\sigma))dW(\sigma), \ s\in [t,T],
%$$
%and $ X^\mu(s)=\gamma_t(s), \ s\in[0,t)$.
 Let $X^\mu$ be the solution of equation (\ref{07162}) and
 define $y^\mu$ by $y^\mu(s)=X^\mu(s)-e^{(s-t)A_\mu}\eta_t(t),\ \ t\leq s\leq T$ and $y^\mu(s)=\gamma_t(s)-\eta_t(s),\ \ 0\leq s<t$, then $y^\mu$ be the solution of the following:
$$
                      y^\mu(s)=e^{(s-t)A_\mu}y^\mu_t(t)+\int^{s}_{t}e^{(s-\sigma)A_\mu}F(X^{\mu}_\sigma, u(\sigma))d\sigma+\int^{s}_{t}e^{(s-\sigma)A_\mu}G(X^{\mu}_\sigma,u(\sigma))dW(\sigma),\ s\in [t,T],
$$
and $y^\mu(s)=\gamma_t(s)-\eta_t(s), \ s\in[0,t)$.
  By Lemmas  \ref{theoremito} and \ref{theoremS}, we have, $P$-a.s.,  for all $s\in [t,T]$:  %for every $g\in C^1(R)$ with $g'\geq0$,
\begin{eqnarray*}
                  \Upsilon^M(y^\mu_s)&=&\Upsilon^M(y^\mu_{{t}})+\int^{s}_{{t}}(\partial_x\Upsilon^M(y^\mu_\sigma),A_\mu y^\mu(\sigma)+F(X^{\mu}_\sigma,u(\sigma))_H\\
                  &&+\frac{1}{2}\mbox{Tr}(\partial_{xx}\Upsilon^M(y^\mu_\sigma)G(X^{\mu}_\sigma,u(\sigma))G^*(X^{\mu}_\sigma,u(\sigma)))d\sigma\\
                  &&+\int^{s}_{t}(\partial_x\Upsilon^M(y^\mu_\sigma),G(X^{\mu}_\sigma,u(\sigma))dW(\sigma))_H.
\end{eqnarray*}
   Noting that $A$ is the infinitesimal generator of a $C_0$ contraction semigroup, we have, if $||y^\mu_\sigma||^2_{0}\neq0$,
    $$\bigg{(}6M|y^\mu(\sigma)|^4y^\mu(\sigma)-\frac{18\left(||y^\mu_\sigma||^6_{0}-|y^\mu(\sigma)|^6\right)^2|y^\mu(\sigma)|^4 y^\mu(\sigma)}{||y^\mu_\sigma||^{12}_{0}},A_\mu y^\mu(\sigma)\bigg{)}_H\leq 0,\ \ \mbox{for}\  M\geq3.$$
     Thus, $P$-a.s.,  for all $s\in [t,T]$:
   \begin{eqnarray*}
                  \Upsilon^M(y^\mu_s)&\leq& \Upsilon^M(y^\mu_{{t}})+\int^{s}_{{t}}(\partial_x\Upsilon^M(y^\mu_\sigma),F(X^{\mu}_\sigma,u(\sigma))_H \\
                  &&+\frac{1}{2}\mbox{Tr}(\partial_{xx}\Upsilon^M(y^\mu_\sigma)G(X^{\mu}_\sigma,u(\sigma))G^*(X^{\mu}_\sigma,u(\sigma)))d\sigma\\
                  &&
                  +\int^{s}_{t}(\partial_x\Upsilon^M(y^\mu_\sigma),G(X^{\mu}_\sigma,u(\sigma))dW(\sigma))_H.
\end{eqnarray*}
Letting $\mu\rightarrow\infty$, by (\ref{0717}), we obtain, $P$-a.s.,  for all $s\in [t,T]$:
\begin{eqnarray*}
                  \Upsilon^M(y_s)&\leq & \Upsilon^M(y_{t})+\int^{s}_{{t}}(\partial_x\Upsilon^M(y_\sigma),F(X^{\gamma_t,u}_\sigma,u(\sigma))_Hd\sigma\\
                  &&+\frac{1}{2}\mbox{Tr}(\partial_{xx}\Upsilon^M(y_\sigma)G(X^{\gamma_t,u}_\sigma,u(\sigma))
                  G^*(X^{\gamma_t,u}_\sigma,u(\sigma)))d\sigma\\
                  &&
                  +\int^{s}_{t}(\partial_x\Upsilon^M(y_\sigma),G(X^{\gamma_t,u}_\sigma,u(\sigma))dW(\sigma))_H,
\end{eqnarray*}
where
 $y(s)=X^{\gamma_t,u}(s)-e^{A(s-t)}\eta_t(t),\  t\leq s\leq T$ and $y(s)=\gamma_t(s)-\eta_t(s),\  0\leq s<t$.
  That is (\ref{jias510815jia}). By the similar (ever easier) process, we can show (\ref{jias510815jia11})  holds true.
The proof is now complete. \ \ $\Box$
\par
              According to Lemma  \ref{theoremito2}, the following
\begin{remark}\label{remarks}
 \begin{description}
  \rm{
 \item{(i)}      Since $||\cdot||_{0}^6$ is not belongs to $C^{1,2}_p({\Lambda})$, then, for every $a_{\hat{t}}\in \Lambda$, $||\gamma_t-a_{\hat{t},t,A}||_0^6$ cannot  appear  as an auxiliary functional in the proof of the
    uniqueness and stability of viscosity solutions. However, by the above lemma, we can replace $||\gamma_t-a_{\hat{t},t,A}||_{0}^6$
   with its equivalent functional  $\Upsilon^3(\gamma_t-a_{\hat{t},t,A})$. Therefore, we can get the uniqueness result of viscosity solutions when   coefficients satisfy  continuity  assumptions under $||\cdot||_0$.
 \item{(ii)}
       It follows from (\ref{s0}) that $\overline{\Upsilon}^3(\cdot,\cdot)$ is a gauge-type function.  We can apply it to Lemma \ref{theoremleft} to
get a maximum of a perturbation of the auxiliary functional in the proof of uniqueness.
% For notational simplicity,
%we use $S_1(\cdot)$ to denote $S_1(\cdot,\mathbf{0})$, which will be applied   to prove the uniqueness of viscosity solutions.
  }
\end{description}
\end{remark}

\section{ A DPP for optimal control problems.}
\par
In this section, we consider optimal control problem  (\ref{state1}) and  (\ref{cost1}).
  % Now, we describe some continuous properties of the solutions of  state equation
%(\ref{state1}) and cost equation (\ref{fbsde1}).
 First let us make the following assumptions:
 %we assume that functionals
% $
%        q: {\Lambda}\times R\times \Xi\times U\rightarrow R$ and $\phi: {\Lambda}_T\rightarrow R
%$ satisfy the
%following assumption.
      % For any $t\in [0,T]$, denote by ${\cal{H}}^2(t,T)$ the space of all ${{\cal{F}}}^t$-adapted, $R^d$-valued  processes $(Y(s))_{t\leq s\leq T}$ such that
%       $||Y||^2=E[\int^{T}_{t}|Y(s)|^2ds]<\infty$ and by ${\cal{S}}^2(t,T)$ the space of all ${{\cal{F}}}^t$-adapted, $R$-valued continuous processes
%       $(Y(s))_{t\leq s\leq T}$ such that
%       $||Y||^2=E[\sup_{t\leq s\leq T}|Y(s)|^2]<\infty$.
\begin{hyp}\label{hypcost}
$
        q: {\Lambda}\times \mathbb{R}\times \Xi\times U\rightarrow \mathbb{R}$ and $\phi: {\Lambda}_T\rightarrow \mathbb{R}$ are continuous and
%\begin{description}
%        \item{(i)}
%        For every fixed $(t,\gamma_t,y,z)\in[0,T]\times\Lambda\times R\times \Xi$,
%        $q(\gamma_t,y,z,\cdot)$ is continuous in $u$.
%\par
%       \item{(ii)}
                 there exists a  constant $L>0$
                 such that, for all $(t,\gamma_t,\eta_T,y,z,u)$,  $ (t, \gamma'_t,\eta'_T,y',z',u)
                 \in [0,T]\times\Lambda\times {\Lambda}_T\times \mathbb{R}\times \Xi\times U$,
                 % it holds that
      \begin{eqnarray*}
                  &&|q(\gamma_t,y,z,u)|\leq L(1+||\gamma_t||_0+|y|+|z|),
                 \\
               &&|q(\gamma_t,y,z,u)-q(\gamma'_{t},y',z',u)|\leq L(||\gamma_t-\gamma'_{t}||_0+|y-y'|+|z-z'|),
                 \\
                 &&|\phi(\eta_T)-\phi(\eta'_T)|\leq L||\eta_T-\eta'_{T}||_0.
\end{eqnarray*}
%Here, for simplicity, we denote %$\int^{t\vee t'}_{0}[\gamma_{t,t\vee t'}(\sigma)-\gamma'_{t',t\vee t'}(\sigma)]d\mu$ and
%$\int^{t\vee t'}_{0}|\gamma_{t,t\vee t'}(\sigma)-\gamma'_{t',t\vee t'}(\sigma)|d\mu$ by
%  %$\mu(\gamma_{t}-\gamma'_{t'})$ and
%  $\mu(|\gamma_{t}-\gamma'_{t'}|)$.
%  %, respectively.
%\end{description}
\end{hyp}
Combing Lemma \ref{lemma2.5} and Lemma \ref{lemmaexist111}, we obtain
\begin{lemma}\label{lemmaexist}
\ \ Assume that Hypothesis \ref{hypstate} (i), (ii)  and Hypothesis \ref{hypcost}  hold. Then for every %$u\in {\cal{U}}$,
$(t,\xi_t,u(\cdot))\in [0,T]\times L_{\cal{P}}^p(\Omega;C([0,t];H))\times {\cal{U}}[0,T]$ and $p>2$,  BSDE (\ref{fbsde1}) admits a unique
 pair of solutions $(Y^{\xi_t,u}, Z^{\xi_t,u})$.  Furthermore, let  $X^{\xi'_t,u}$ and $(Y^{\xi'_t,u}, Z^{\xi'_t,u})$ be the solutions of SDE (\ref{state1}) and BSDE (\ref{fbsde1})
 corresponding $(t,\xi'_t,u(\cdot))\in [0,T]\times L_{\cal{P}}^p(\Omega;C([0,t];H))\times {\cal{U}}[0,T]$. Then the following estimates hold:
\begin{eqnarray}\label{fbjia4}
                 \mathbb{E}\left[\sup_{t\leq s\leq T}\left|Y^{\xi_t,u}(s)-Y^{\xi'_t,u}(s)\right|^p\right]
                \leq C_p\mathbb{E}||\xi_t-\xi'_t||_0^p;
                \end{eqnarray}
                \begin{eqnarray}\label{fbjia5}
                \mathbb{E}\bigg{[}\sup_{t\leq s\leq T}\left|Y^{\xi_t,u}(s)\right|^p\bigg{]}
                +\bigg{(}\int^{T}_{t}\left|Z^{\xi_t,u}(s)\right|^2ds\bigg{)}^{\frac{p}{2}}\leq C_p(1+\mathbb{E}||\xi_t||_0^p).
                \end{eqnarray}
%\begin{eqnarray*}
%                E[\sup_{t\leq s\leq T}|X^{t,\xi,u}(s)-X^{t,x',u}(s)|^p]&\leq& C_p|x-x'|_C^p,\\
%                E[\sup_{t\leq s\leq T}|X^{t,\xi,u}(s)|^p]&\leq& C_p(1+|x|_C^p),\\
%               % E[|X_{r}^{t,\xi,u}-x|_C^p|{\cal{F}}_t]&\leq& C_p{(}1+|x|_C^p)(r-t)^{\frac{p}{2}},
%                E[\sup_{t\leq s\leq r}|X^{t,\xi,u}(s,)-x(0)|^p]&\leq& C_p{(}1+|x|_C^p)(r-t)^{\frac{p}{2}},
%\end{eqnarray*}
%and
%\begin{eqnarray*}
%               && E\bigg{[}\sup_{t\leq s\leq T}|Y^{t,\xi,u}(s)-Y^{t,x',u}(s)|^p\bigg{]}
%                \leq C_p|x-x'|_C^p,\\
%               && E\bigg{[}\sup_{t\leq s\leq T}|Y^{t,\xi,u}(s)|^p\bigg{]}
%                +\bigg{(}\int^{T}_{t}|Z^{t,\xi,u}(s)|^2ds\bigg{)}^{\frac{p}{2}}\leq C_p(1+|x|_C^p).
%               %&&E\bigg{[}\sup_{t\leq s\leq r}|Y^{t,\xi,u}(s)-Y^{t,\xi,u}(t)|^p\bigg{]}\leq C_p{(}1+|x|_C^p)(r-t)^{\frac{p}{2}}.
%\end{eqnarray*}
              The constant $C_p$ depending only on  $p$, $T$ and $L$.
\end{lemma}
{\bf  Proof}. \ \
                  Existence and uniqueness of the solution of the backward equation (\ref{fbsde1}) follows from  Lemma  \ref{lemma2.5}.
                  Using  inequalities (\ref{lemma2.51}) and (\ref{state1est}) we  get inequality  (\ref{fbjia4}).  Combing inequalities (\ref{lemma2.510}) and (\ref{state1est}), we obtain inequality (\ref{fbjia5}).
                   \ \ $\Box$
                   \par
The  $J(\xi_t,u(\cdot)):=J(\gamma_t,u(\cdot))|_{\gamma_t=\xi_t}$ and $Y^{\xi_t,u}(t)$, $(t,\xi_t)\in [0,T]\times L_{\cal{P}}^p(\Omega;C([0,t];H))$ and $p>2$, are related by the following
theorem.
 \begin{theorem}\label{theoremj=y}
\ \
 Under  Hypothesis \ref{hypstate} (i), (ii)  and Hypothesis \ref{hypcost},    for every $p>2$, $t\in [0,T]$, $u(\cdot)\in {\cal{U}}[t,T]$ and $\xi_t, \xi'_t\in L_{\cal{P}}^p(\Omega;C([0,t];H))$, we have
\begin{eqnarray}\label{j=y}
J(\xi_t,u(\cdot))=Y^{\xi_t,u}(t),
\end{eqnarray}
and %for a constant $C>0$
\begin{eqnarray}\label{0903jia}
|Y^{\xi_t,u}(t)-Y^{\xi'_t,u}(t)|\leq C_p^{\frac{1}{p}}||\xi-\xi'||_0.
\end{eqnarray}
\end{theorem}
{\bf  Proof}. \ \
                     Let $\{h^n_t\}$, $n\in N$,
                      be a dense subset of $\Lambda_t$, $B(h^n_t,\frac{1}{k})$
                      be the open sphere in $\Lambda_t$ with the radius
                      equal to $\frac{1}{k}$ and the center at the
                      point $h^n_t$.  Set $B_{n,k}:=B(h^n_t,\frac{1}{k})\setminus
                      \bigcup_{m<n}B(h^m_t,\frac{1}{k})$  and $A_{n,k}:=\{\omega\in \Omega|\xi_t(\omega)\in
                      B_{n,k}\}$. Then $\cup_{n=1}^{\infty}A_{n,k}=\Omega$ and the
                      sequence $f^k_t(\omega):= \Sigma^{\infty}_{n=1}h^n_t1_{
                      A_{n,k}}(\omega)$ is ${\mathcal
                      {F}}_{t}$-measurable and
                             converges to $\xi_t$ strongly and
                             uniformly.
                     \par
                     For every  $n$ and $u(\cdot)\in {\mathcal {U}}[t,T]$, we put $(X^{n}(s),Y^{n}(s),Z^{n}(s))=(X^{h^n_t,u}(s),Y^{h^n_t,u}(s),Z^{h^n_t,u}(s))$.
                     Then $X^{n}(s)$ is the solution of the PSEE
\begin{eqnarray}\label{stateint333}
             X^n(s)=e^{(s-t)A}h^n_t(t)+\int_{t}^{s}e^{(s-\sigma)A}F(X^n_\sigma,u(\sigma))d\sigma+\int_{t}^{s}e^{(s-\sigma)A}G(X^n_\sigma,u(\sigma))dW(\sigma), \ \
             s\in [t,T],\nonumber
\end{eqnarray}
       where $X^n_t=h^n_t$; and $(Y^{n}(s),Z^{n}(s))$ is the solution of the associated  BSDE
\begin{eqnarray}\label{fbsde333}
Y^{n}(s)=\phi(X_T^{n})+\int^{T}_{s}q(X_\sigma^{n},Y^{n}(\sigma),Z^{n}(\sigma),u(\sigma))d\sigma-\int^{T}_{s}Z^{n}(\sigma)dW(\sigma),\  \ s\in [t,T].\nonumber
\end{eqnarray}
The above two equations are multiplied by $1_{A_{n,k}}$  and summed up with respect to $n$. Thus, taking into account that $\sum_{n=1}^{\infty}\varphi(h^n_t)1_{A_{n,k}}=\varphi(\sum_{n=1}^{\infty}h^n_t1_{A_{n,k}})$, we obtain
\begin{eqnarray*}
             \sum_{n=1}^{\infty}1_{A_{n,k}}X^n(s)&=&\sum_{n=1}^{\infty}1_{A_{n,k}}e^{(s-t)A}h^n_t(t)+\int_{t}^{s}e^{(s-\sigma)A}F\left(\sum_{n=1}^{\infty}1_{A_{n,k}}X^n_\sigma, u(\sigma)\right)d\sigma\nonumber\\
             && ~~
             +\int_{t}^{s}e^{(s-\sigma)A}G\left(\sum_{n=1}^{\infty}1_{A_{n,k}}X^n_\sigma,u(\sigma)\right)dW(\sigma), \nonumber
\end{eqnarray*}
       and
\begin{eqnarray}\label{fbsde333}
&&\sum_{n=1}^{\infty}1_{A_{n,k}}Y^{n}(s)\nonumber\\
&=&\phi\left(\sum_{n=1}^{\infty}1_{A_{n,k}}X_T^{n}\right)+\int^{T}_{s}q\left(\sum_{n=1}^{\infty}1_{A_{n,k}}X_\sigma^{n},
\sum_{n=1}^{\infty}1_{A_{n,k}}Y^{n}(\sigma),\sum_{n=1}^{\infty}1_{A_{n,k}}Z^{n}(\sigma),u(\sigma)\right)d\sigma\nonumber\\
&&-\int^{T}_{s}\sum_{n=1}^{\infty}1_{A_{n,k}}Z^{n}(\sigma)dW(\sigma).\nonumber
\end{eqnarray}
Then the strong uniqueness property of the solution to the PSEE and the BSDE yields
\begin{eqnarray}\label{stateint333444}
             X^{f^k_t,u}(s)&=&\sum_{n=1}^{\infty}1_{A_{n,k}}X^n(s),\nonumber\\
             \left(Y^{f^k_t,u}(s),Z^{f^k_t,u}(s)\right)
             &=&\left(\sum_{n=1}^{\infty}1_{A_{n,k}}Y^{n}(s),\sum_{n=1}^{\infty}1_{A_{n,k}}Z^{n}(s)\right),\ \ s\in [t,T].
\end{eqnarray}
Finally, from $J(h^n_t,u(\cdot))=Y^n(t)$, $n\geq 1$, we deduce that
\begin{eqnarray}\label{stateint33344455}
             Y^{f^k_t,u}(t)
             =\sum_{n=1}^{\infty}1_{A_{n,k}}Y^{n}(t)=\sum_{n=1}^{\infty}1_{A_{n,k}}J(h^n_t,u(\cdot))=J\left(\sum_{n=1}^{\infty}1_{A_{n,k}}h^n_t,u(\cdot)\right)=J(f^k_t,u(\cdot)).
\end{eqnarray}
Consequently, from the estimate (\ref{fbjia4}),  we get
\begin{eqnarray*}
 \mathbb{E}|Y^{\xi_t,u}(t)-J(\xi_t,u(\cdot))|^p&\leq&  2^{p-1}\mathbb{E}|Y^{\xi_t,u}(t)-Y^{f^{k}_t,u}(t)|^p+2^{p-1} \mathbb{E}|J(f^k_t,u(\cdot))-J(\xi_t,u(\cdot))|^p\\
 &\leq& 2^p C_p\mathbb{E}||\xi_t-f^k_t||_0^p\rightarrow0\ \mbox{as}\ k\rightarrow \infty.
\end{eqnarray*}
Now let us prove $(\ref{0903jia})$. By (\ref{fbjia4}) and (\ref{j=y}),
\begin{eqnarray*}\label{0903jia}
|Y^{\xi_t,u}(t)-Y^{\xi'_t,u}(t)|&=&|J(\xi_t,u(\cdot))-J(\xi'_t,u(\cdot))|=|Y^{\gamma_t,u}(t)|_{\gamma_t=\xi_t}-Y^{\eta_t,u}(t)|_{\eta=\xi'_t}|\\
&=&|(Y^{\gamma_t,u}(t)-Y^{\eta_t,u}(t))|_{\gamma_t=\xi_t,\eta=\xi'_t}|\leq C_p^{\frac{1}{p}}||\xi_t-\xi'_t||_0.
\end{eqnarray*}
The proof is complete.\ \ $\Box$
\par
                  From the uniqueness of the solution of (\ref{fbsde1}), it
                  follows that
$$
                 Y^{\gamma_t,u}(t+\delta)=Y^{X^{\gamma_t,u}_{t+\delta},u}(t+\delta)=J\left(X^{\gamma_t,u}_{t+\delta},u(\cdot)\right),\ \ \mbox{a.s.}
$$
As we mention in Section 3,  for every $(t,\gamma_t,u(\cdot))\in [0,T]\times \Lambda\times {\cal{U}}[t,T]$, $X^{\gamma_t,u}(s)$ is ${\cal{F}}_t^s$-measurable for all $s\in[t,T]$, then $Y^{\gamma_t,u}(s)$ is ${\cal{F}}_t^s$-measurable for all $s\in[t,T]$.
%(see also Proposition 4.2 in El Karoui, Peng, Quenez \cite{el}).
%On the other hand, $Y^{t,x,x_0,u}(s)$ is ${\cal{F}}_s$-measurable for all $s\in[t,T]$.
In particular, $Y^{\gamma_t,u}(t)$ is deterministic. Then we have
%Formally,  under the assumptions  Hypothesis \ref{hypstate}, the  value functional $V(t,x)$
% defined by (\ref{value1})
% is  ${\cal{F}}_t$-measurable.
%However, we have %according to  Proposition 3.3 in  Buckdahn and Li \cite{buck1}, we can prove the following.
\begin{theorem}\label{valuedet} \ \  Suppose   Hypothesis \ref{hypstate} (i), (ii) and Hypothesis \ref{hypcost} hold true.
                      Then $V$ is a deterministic functional.
\end{theorem}
\par
The following  property of the value functional $V$ which we present is an immediate consequence of Lemma \ref{lemmaexist}.
\begin{lemma}\label{lemmavaluev}
\ \ Assume that  Hypothesis \ref{hypstate} (i), (ii)  and Hypothesis \ref{hypcost} hold, then,  %$V\in {{C}}^0(\Lambda)$ and
 %there exists a  constant $C>0$ such that, 
 for all $p>2$, $0\leq t\leq T$, $\gamma_t, \eta_t\in \Lambda_t$,
\begin{eqnarray}\label{valuelip}
|V(\gamma_t)-V(\eta_t)|\leq C_p^{\frac{1}{p}}||\gamma_t-\eta_t||_0; \ \ \ \ \
               |V(\gamma_t)|\leq C_p^{\frac{1}{p}}(1+||\gamma_t||_0).
\end{eqnarray}
%If the Hypotheses  \ref{hypstate} (ii) and   \ref{hypvalue} (ii) also hold, we have
%\begin{eqnarray}\label{valuelip22222}
%|V(\gamma_t)-V(\gamma'_{t'})|
%&\leq& C(|\gamma_t(t)-\gamma_{t'}'(t')|+|\gamma_t(0)-\gamma_{t'}'(0)|\nonumber\\
%               && +|t'-t|^{\frac{1}{2}}(1+||\gamma_t||_0
%                +||\gamma_{t'}'||_0)+|\gamma_{t,t'}-\gamma_{t'}'|_H).
%\end{eqnarray}
\end{lemma}
 \par
We also have that
\begin{lemma}\label{lemma3.6}
                        For all $t\in[0,T]$, $\xi_t\in L_{\cal{P}}^p(\Omega,{\mathcal {F}}_t,\Lambda_t)$ for any $p>2$, and $u(\cdot)\in{\mathcal
                        {U}}[t,T]$ we have
\begin{eqnarray}\label{3.14}
                                  V(\xi_t)\geq Y^{\xi_t,u}(t),\
                                  \ \mbox{a.s.},
\end{eqnarray}
                  and for any $\varepsilon>0$ there exists an
                  admissible control $u(\cdot)\in{\mathcal
                        {U}}[t,T]$ such that
\begin{eqnarray}\label{3.15}
                                  V(\xi_t)\leq Y^{\xi_t,u}(t)+\varepsilon,\
                                  \ \mbox{a.s.}
\end{eqnarray}
\end{lemma}
{\bf  Proof}. \ \
                    By Theorem \ref{theoremj=y} and  the definition of $V(\gamma_t)$ we have, for any $u(\cdot)\in{\mathcal
                        {U}}[t,T]$,
$$
                      V(\xi_t)=V(\eta_t)|_{\eta_t=\xi_t}=\mathop{\sup}\limits_{v(\cdot)\in{\mathcal
                      {U}}[t,T]}J(\eta_t,v(\cdot))|_{\eta_t=\xi_t}\geq
                      J(\eta_t,u(\cdot))|_{\eta_t=\xi_y}=Y^{\xi_t,u}(t).
$$
                      We now prove (\ref{3.15}). Let $\{h^n_t\}$, $n\in N$,
                      be a dense subset of $\Lambda_t$, $B(h^n_t,\frac{1}{k})$
                      be the open sphere in $\Lambda_t$ with the radius
                      equal to $\frac{1}{k}$ and the center at the
                      point $h^n_t$.  Set $B_{n,k}:=B(h^n_t,\frac{1}{k})\setminus
                      \bigcup_{m<n}B(h_m,\frac{1}{k})$. Then the
                      sequence $f^k_t(\omega):= \Sigma^{\infty}_{n=1}h^n_t1_{\{\xi_t\in
                      B_{n,k}\}}(\omega)$ is ${\mathcal
                      {F}}_{t}$-measurable and
                             converges to $\xi_t$ strongly and
                             uniformly. For
                             %$k>5\frac{C^{\frac{1}{2}}_2\vee C}{\varepsilon}$,
                             $k>3\frac{1+C_p^{\frac{1}{p}}}{\varepsilon}$,  we have
$$
                      ||f^k_t-\xi_t||_0\leq \frac{\varepsilon}{3}.
$$
                     Therefore, for every $u(\cdot)\in {\mathcal {U}}[t,T]$,\\
$$
                          \left|Y^{f^k_t,u}(t)-Y^{\xi_t,u}(t)\right|\leq
                          \frac{\varepsilon}{3},
$$
$$
                          \left|V(f^k_t)-V(\xi_t)\right|\leq
                          \frac{\varepsilon}{3}.
$$
 Moreover, for every $h^n_t\in \Lambda_t$,  by the definition of $V$
                                we can choose an admissible control
                                $v^n$ such that
$$
                          V(h^n_t)\leq
                          Y^{h^n_t, v^n}(t)+\frac{\varepsilon}{3}, \ \  P\mbox{-a.s.}.
$$
                    Denote %$A_n:=\{V(h^n_t)-Y^{h^n_t,v^{n,l}}(t)\leq\frac{\varepsilon}{3}\}$ and
                    $u(\cdot):=\sum^{\infty}_{n=1}v^n(\cdot)1_{\{\xi_t\in
                      B_{n,k}\}}$, then
\begin{eqnarray*}
                        Y^{\xi_t,u}(t)&\geq&
                        -|Y^{f^k_t,u}(t)-Y^{\xi_t,u}(t)|+Y^{f^k_t,u}(t)\geq-\frac{\varepsilon}{3}+\sum^{\infty}_{n=1}Y^{h^n_t,v^{n}}(t)1_{\{\xi_t\in
                      B_{n,k}\}}\\
                      %&\geq&-\frac{\varepsilon}{5}+E\sum^{\infty}_{n=1}(V(t, h^n_t)-\frac{\varepsilon}{5})1_{A_n\{\xi_t\in
%                      B_{n,k}\}}+E\sum^{\infty}_{n=1}Y^{t,h^n_t,v^n}(t)1_{\overline{A_n}\{\xi_t\in
%                      B_{n,k}\}}\\
                      &\geq&-\frac{2\varepsilon}{3}+\sum^{\infty}_{n=1}V(h^n_t)1_{\{\xi_t\in
                      B_{n,k}\}}=-\frac{2\varepsilon}{3}+V(f^k_t)\\
                      % &\geq&-\frac{\varepsilon}{3}+V(f^k_t)-\frac{\varepsilon}{3}+\sum^{\infty}_{n=1}(Y^{h^n_t,v^{n,l}}(t)-V(h^n_t))1_{\overline{A_n}\{\xi_t\in
%                      B_{n,k}\}}\\
                      &\geq&-\varepsilon+V(\xi_t),  \ \  P\mbox{-a.s.}.
\end{eqnarray*}
%Let $\overline{D_l}:=\cup_{n=1}^{\infty}\overline{A_n}$ and $D:=\cup_{l=1}^{\infty}D_l$, then
%$$
%1\geq P(D)=1-P(\overline{D})=1-P(\cap_{l=1}^{\infty}\overline{D_l})\geq1-P(\overline{D_l})\geq 1-\lim_{l\rightarrow\infty}\frac{\varepsilon}{l}=1,
%$$
% and
%\begin{eqnarray*}
%                        Y^{\xi_t,u}(t)\geq-\varepsilon+V(\xi_t), \ \ \mbox{on} \ \ {D}.
%\end{eqnarray*}
Thus (\ref{3.15}) holds.\ \ $\square$
\par
We now discuss a dynamic programming principle  (DPP) for the optimal control problem (\ref{state1}), (\ref{fbsde1}) and (\ref{value1}).
For this purpose, we define the family of backward semigroups associated with BSDE (\ref{fbsde1}), following the
idea of Peng \cite{peng11}.
\par
Given the initial path $(t,\gamma_t)\in [0,T)\times{\Lambda}$, a positive number $\delta\leq T-t$, an admissible control $u(\cdot)\in {\cal{U}}[t,t+\delta]$ and
a real-valued random variable $\zeta\in L^2(\Omega,{\cal{F}}_{t+\delta},P;\mathbb{R})$, we put
\begin{eqnarray}\label{gdpp}
                        G^{\gamma_t,u}_{s,t+\delta}[\zeta]:=\tilde{Y}^{\gamma_t,u}(s),\ \
                        \ \ \ \ s\in[t,t+\delta],
\end{eqnarray}
                        where $(\tilde{Y}^{\gamma_t,u}(s),\tilde{Z}^{\gamma_t,u}(s))_{t\leq s\leq
                        t+\delta}$ is the solution of the following
                        BSDE with the time horizon $t+\delta$:
\begin{eqnarray}\label{bsdegdpp}
\begin{cases}
d\tilde{Y}^{\gamma_t,u}(s) =-q(s,X^{\gamma_t,u}_s,\tilde{Y}^{\gamma_t,u}(s),\tilde{Z}^{\gamma_t,u}(s),u(s))ds+\tilde{Z}^{\gamma_t,u}(s)dW(s),\\
 ~\tilde{Y}^{\gamma_t,u}(t+\delta)=\zeta,
 \end{cases}
\end{eqnarray}
                  and $X^{\gamma_t,u}(\cdot)$ is the solution of SDE (\ref{state1}).
               % \par
%                  The proof of DPP for  stochastic delay optimal control on fixed time horizontal
%is similar to the Markovian case (see Buckdahn and Li \cite{buck1}). Here we omit it.
                  %
%
%
%                   Then, obviously, for the solution
%                   $({Y}^{\gamma_t,u}(s),$ ${Z}^{\gamma_t,u}(s))_{t\leq s\leq
%                        T}$ of BSDE (\ref{fbsde1}), the uniqueness of the BSDE yields
%\begin{eqnarray*}
%               J(\gamma_t,u)=Y^{\gamma_t,u}(t)=G^{\gamma_t,u}_{t,T}[\phi(X^{\gamma_t,u}_T)]=G^{\gamma_t,u}_{t,t+\delta}[Y^{\gamma_t,u}(t+\delta)]
%               =G^{\gamma_t,u}_{t,t+\delta}[J(X^{\gamma_t,u}_{t+\delta},u)].
%\end{eqnarray*}
\begin{theorem}\label{theoremddp} %(See Theorem 3.4 in \cite{tang1})
    Assume  Hypothesis \ref{hypstate} (i), (ii) and Hypothesis \ref{hypcost} hold true, the value functional
                              $V$ obeys the following DPP: for
                              any $\gamma_t\in {\Lambda_t}$ and $0\leq t<t+\delta\leq T$,
\begin{eqnarray}\label{ddpG}
                              V(\gamma_t)=\mathop{\sup}\limits_{u(\cdot)\in{\mathcal
                              {U}}[t,t+\delta]}G^{\gamma_t,u}_{t,t+\delta}[V(X^{\gamma_t,u}_{t+\delta})].
\end{eqnarray}
\end{theorem}
%\begin{theorem}\label{theorem3.8}
%                              Assume Hypotheses \ref{hypstate} and \ref{hypvalue} hold true, the value functional
%                              $V(\gamma_t)$ obeys the following DPP: for
%                              any $0\leq t<t+\delta\leq T$, $x\in
%                              {\Lambda_t}$.
%\begin{eqnarray}\label{3.18}
%                              V(\gamma_t)=\mbox{esssup}_{u\in{\mathcal
%                              {U}}[t,T]}G^{\gamma_t,u}_{t,t+\delta}[V(t+\delta,X^u_{t+\delta}(\gamma_t))].
%\end{eqnarray}
%\end{theorem}
\par
{\bf  Proof}. \ \
                     By the definition of $V(\gamma_t)$ we have
$$
                     V(\gamma_t)=\mathop{\sup}\limits_{u(\cdot)\in{\mathcal
                              {U}}[t,T]}G^{\gamma_t,u}_{t,T}[\phi(X^{\gamma_t,u}_T)]=\mathop{\sup}\limits_{u(\cdot)\in{\mathcal
                              {U}}[t,T]}G^{\gamma_t,u}_{t,t+\delta}\left[Y^{X^{\gamma_t,u}_{t+\delta},u}(t+\delta)\right].
$$
                     From (\ref{3.14}) and the comparison theorem (see Lemma \ref{lemma2.70904}) it
                     follows that
$$
                      V(\gamma_t)\leq \mathop{\sup}\limits_{u(\cdot)\in{\mathcal
                              {U}}[t,T]}G^{\gamma_t,u}_{t,t+\delta}[V(X^{\gamma_t,u}_{t+\delta})].
$$
                  On the other hand, from (\ref{3.15}), for any
                  $\varepsilon>0$ and $u(\cdot)\in {\mathcal {U}}[t,T]$ there
                  exists an admissible control $\overline{u}(\cdot)\in{\mathcal
                  {U}}[t+\delta,T]$ such that
$$
                     V(X^{\gamma_t,u}_{t+\delta})\leq Y^{X^{
                     \gamma_t,u}_{t+\delta},\overline{u}}(t+\delta)+\varepsilon\ \  P\mbox{-a.s.}
$$
                     For any $u(\cdot)\in {\mathcal {U}}[t,T]$ with $\overline{u}(\cdot)\in{\mathcal {U}}[t+\delta,T]$ from above, we put
$$
                          v(s)= \begin{cases}u(s),\ \ \ \ \ \ s\in[t,t+\delta]; \\
                                   \overline{u}(s),\ \ \ \ \ \
                                   s\in[t+\delta,T];\end{cases}\in {\mathcal
                                   {U}}[t,T],
$$
                   then we have
                   \begin{eqnarray*}
                   &&\left|G^{\gamma_t,u}_{t,t+\delta}\left[Y^{X^{\gamma_t,u}_{t+\delta},\overline{u}}(t+\delta)\right]
                    -G^{\gamma_t,u}_{t,t+\delta}[V(X^{\gamma_t,u}_{t+\delta})]\right|\leq C\varepsilon,
\end{eqnarray*}
where the constant $C$ is independent of
                    admissible control processes.
               Therefore,
\begin{eqnarray*}
                    V(\gamma_t)&\geq&
                    G^{\gamma_t,v}_{t,t+\delta}\left[Y^{X^{\gamma_t,v}_{t+\delta},v}(t+\delta)\right]=G^{\gamma_t,u}_{t,t+\delta}\left[Y^{X^{\gamma_t,u}_{t+\delta},\overline{u}}(t+\delta)\right]\\
                   % &\geq&G^{\gamma_t,u}_{t,t+\delta}[V(t+\delta,X^{\gamma_t,u}_{t+\delta})-\varepsilon]\\
                    &\geq&G^{\gamma_t,u}_{t,t+\delta}[V(X^{\gamma_t,u}_{t+\delta})]-C\varepsilon,\
                    \ \ P\mbox{-a.s.}.
\end{eqnarray*}
%
%                  By Theorem A.3 in \cite{kar}, there exists a sequence $\{u_n\}_{n\geq1}\subset {\cal{U}}[t,T]$ such that $\left\{G^{\gamma_t,u_n}_{t,t+\delta}[V(X^{\gamma_t,u_n}_{t+\delta})]\right\}_{n\geq1}$ is a nondecreasing sequence and
%\begin{eqnarray*}
%\mathop{\esssup}\limits_{u(\cdot)\in{\mathcal
%                              {U}}[t,T]}G^{\gamma_t,u}_{t,t+\delta}[V(X^{\gamma_t,u}_{t+\delta})]=\lim_{n\rightarrow\infty} G^{\gamma_t,u_n}_{t,t+\delta}[V(X^{\gamma_t,u_n}_{t+\delta})], \ \  P\mbox{-a.s.}.
%\end{eqnarray*}
%              Thus,  by the dominated convergence theorem, we get that
%$$
%                      V(\gamma_t)\geq \mathop{\esssup}\limits_{u(\cdot)\in{\mathcal
%                              {U}}[t,T]}G^{\gamma_t,u}_{t,t+\delta}[V(X^{\gamma_t,u}_{t+\delta})]-C_1\varepsilon,\
%                              \ \ \ \mbox{a.s.}
%$$
                  For the arbitrariness of $\varepsilon$, we get (\ref{ddpG}).\ \ $\square$
\par
                  %     In Lemma \ref{lemmavaluev} we have already seen that the
%                       value functional $V(\gamma_t)$ is Lipschitz
%                       continuous in $x$, uniformly in $t$.
                        With the
                       help of Theorem \ref{theoremddp}  now show that the continuity
                       property of $V(\gamma_t)$.
\begin{theorem}\label{theorem3.9}
                          Under  Hypothesis \ref{hypstate} (i), (ii) and Hypothesis \ref{hypcost}, then $V\in C^0(\Lambda)$ and
there is a constant $C>0$ such that for every  $0\leq t\leq s\leq T, \gamma_t,\eta_t\in{\Lambda_t}$,
\begin{eqnarray}\label{hold}
                 |V(\gamma_t)-V(s,\eta_{t,s,A})|&\leq&
                                    C(1+||\gamma_t||_0+||\eta_t||_0)(s-t)^{\frac{1}{2}}+C||\gamma_t-\eta_t||_0.
                                    %\nonumber\\
%                             &\leq&C(1+|x|_C+|y|_C)\delta^{\frac{1}{2}}+C|x-y|_C+C(|x_\delta-x|_C\wedge|y_\delta-y|_C). \ \
\end{eqnarray}
\end{theorem}
{\bf  Proof}. \ \
                 Let $(t,\gamma_t,\eta_t)\in[0,T)\times {\Lambda}\times {\Lambda}$ and $s\in[t,T]$.
                 From Theorem \ref{theoremddp} it follows that for any
                 $\varepsilon>0$ there exists an admissible control
                 $u^\varepsilon\in{\mathcal {U}}[t,s]$ such that
\begin{eqnarray}\label{3.19}
                             \mathbb{E}G^{\gamma_t,u^\varepsilon}_{t,s}[V(X^{\gamma_t,u^\varepsilon}_{s})]
                             +\varepsilon\geq V(\gamma_t)\geq \mathbb{E}G^{\gamma_t,u^\varepsilon}_{t,s}[V(X^{\gamma_t,u^\varepsilon}_{s})].
\end{eqnarray}
Therefore,
\begin{eqnarray}\label{3.20}
                          |V(\gamma_t)-V(\eta_{t,s,A})|\leq
                          |I_1|+|I_2|+\varepsilon,
\end{eqnarray}
                        where
$$
                             I_1=\mathbb{E}G^{\gamma_t,u^\varepsilon}_{t,s}[V(X^{\gamma_t,u^\varepsilon}_{s})]-
                             \mathbb{E}G^{\gamma_t,u^\varepsilon}_{t,s}[V(\eta_{t,s,A})],
$$
$$
                                 I_2=\mathbb{E}G^{\gamma_t,u^\varepsilon}_{t,s}[V(\eta_{t,s,A})]-
                                                V(\eta_{t,s,A}).
$$
                                 By Lemmas \ref{lemma2.5}, \ref{lemmaexist111} and \ref{lemmavaluev} we
                                 have that, for some suitable constant
                                 $C$ independent of the control
                                 $u^\varepsilon$,
\begin{eqnarray*}
                               %(Lemma 2.5 is not enough )
                                |I_1|&\leq& C\mathbb{E}\left|V\left(X^{\gamma_t,u^\varepsilon}_{s}\right)-V(\eta_{t,s,A})\right|\\
                                                &\leq&C\mathbb{E}\left|\left|X^{\gamma_t,u^\varepsilon}_{s}-\eta_{t,s,A}\right|\right|_0\leq C(1+||\gamma_t||_0)(s-t)^{\frac{1}{2}}+C||\gamma_t-\eta_t||_0.
\end{eqnarray*}
%Noting that
%$$
% |x_\delta-y|_C\leq |x_\delta-x|_C+|x-y|_C
%$$
%and
%$$
%|x_\delta-y|_C\leq |x_\delta-y_\delta|_C+|y_\delta-y|_C\leq|x-y|_C+|y_\delta-y|_C,
%$$
%we get that
%$$
% |I^1_{\delta}|\leq C(1+|x|)\delta^{\frac{1}{2}}+C|x-y|_C+C(|x_\delta-x|_C\wedge|y_\delta-y|_C).
%$$
%
%                                   and since $E[\sup_{l\in [0,\delta]}[|X^{\gamma_tu^\varepsilon}(t+l)-x(0)|^2|{\mathcal
%                                                       {F}}_t]\leq
%                                                       C(1+|x|_C^2)\delta$,
%                                   we deduce that
%$$
%                                   |I^1_{\delta}|\leq C(1+|x|)\delta^{\frac{1}{2}}+C\sup_{|\theta-\theta'|\leq \delta}|x(\theta)-x(\theta')|.
%$$
         From the definition of
         $G^{\gamma_t,u^\varepsilon}_{t,s}[\cdot]$ we get that the
         second term $I_2$ can be written as
\begin{eqnarray*}
                                   I_2&=& \mathbb{E}\bigg{[}V(\eta_{t,s,A})+\int^{s}_{t}
                                   q\left(\sigma,X_\sigma^{\gamma_t,u^{\varepsilon}},Y^{\gamma_t,u^{\varepsilon}}(\sigma),Z^{\gamma_t,u^{\varepsilon}}(\sigma),u^{\varepsilon}(\sigma)\right)d\sigma\\
                                   &&-
                                   \int^{s}_{t}Z^{\gamma_t,u^{\varepsilon}}(\sigma)dW(\sigma)\bigg{]}-V(\eta_{t,s,A})\\
                                   &=&\mathbb{E}\int^{s}_{t}
                                   q\left(\sigma,X_\sigma^{\gamma_t,u^{\varepsilon}},Y^{\gamma_t,u^{\varepsilon}}(\sigma),Z^{\gamma_t,u^{\varepsilon}}(\sigma),u^{\varepsilon}(\sigma)\right)d\sigma,
\end{eqnarray*}
                                   where $(Y^{\gamma_t,u^{\varepsilon}}(s),Z^{\gamma_t,u^{\varepsilon}}(s))_{t\leq s\leq
                                   s}$ is the solution of
                                   (\ref{bsdegdpp}) with the terminal
                                   condition $\eta = V(\eta_{t,s,A})$
                                   and the control $u^\varepsilon$.
                                  With the help of the Schwartz
                                   inequality, we then have,  for some suitable constant
                                 $C$ independent of the control
                                 $u^\varepsilon$,
\begin{eqnarray*}
                                    I_2&\leq& (s-t)^{\frac{1}{2}}\bigg{[}\int^{s}_{t}
                                   \mathbb{E}\left|q\left(\sigma,X_\sigma^{\gamma_t,u^{\varepsilon}},Y^{\gamma_t,u^{\varepsilon}}(\sigma),Z^{\gamma_t,u^{\varepsilon}}(\sigma),u^{\varepsilon}(\sigma)\right)\right|^2d\sigma
                                   \bigg{]}^{\frac{1}{2}}\\
                                   &\leq&C(s-t)^{\frac{1}{2}}\bigg{[}\int^{s}_{t}
                                   \mathbb{E}\left(1+\left|\left|X_\sigma^{\gamma_t,u^{\varepsilon}}\right|\right|_0^{2}+\left|Y^{\gamma_t,u^{\varepsilon}}(\sigma)\right|^2+\left|Z^{\gamma_t,u^{\varepsilon}}(\sigma)\right|^2\right)d\sigma
                                   \bigg{]}^{\frac{1}{2}}\\
                                   &\leq&C(1+||\gamma_t||_0+||\eta_t||_0)(s-t)^{\frac{1}{2}}.
\end{eqnarray*}
                                   Hence, from (\ref{3.20}),
$$
                                    |V(\gamma_t)-V(s,\eta_{t,s,A})|\leq
                                    C(1+||\gamma_t||_0+||\eta_t||_0)(s-t)^{\frac{1}{2}}+C||\gamma_t-\eta_t||_0
                                   +\varepsilon,
$$
                                   and letting $\varepsilon\downarrow
                                   0$ we get (\ref{hold}) holds true.
                                   \par
                                   Finally, let $(t,\gamma_t), (s,\gamma'_s)\in[0,T]\times {\Lambda}$ and $s\in[t,T]$, by (\ref{valuelip}) and (\ref{hold})
                                   \begin{eqnarray*}
                                              &&|V(\gamma_t)-V(\gamma'_s)|\leq  |V(\gamma_t)-V(\gamma_{t,s,A})|+ |V(\gamma_{t,s,A})-V(\gamma'_s)|\\
                                              &\leq&  C(1+||\gamma_t||_0)(s-t)^{\frac{1}{2}}+C_p^{\frac{1}{p}}||\gamma_{t,s,A}-\gamma'_s||_0.
                                   \end{eqnarray*}
                                   Thus, $V\in C^0(\Lambda)$.
%$$
%                                   |V(\gamma_t)-V(s,x)|\leq
%                                    C(1+|x|_C)\delta^{\frac{1}{2}}+C\sup_{|\theta-\theta'|\leq \delta}|x(\theta)-x(\theta')|.
%$$
                                  The proof is complete.\ \
                                  $\square$

\section{Viscosity solutions to  PHJB equations: Existence theorem.}%Construction of random invariant manifold

\par
In this section, we consider the  second  order path-dependent  Hamilton-Jacobi-Bellman
                   (PHJB) equation (\ref{hjb1}). As usual, we start with classical solutions.
                                  % In this section, we consider the following second order path-dependent HJB (PHJB) equation:
%\begin{eqnarray}\label{hjb}
%\cases{
%           \partial_tV(\gamma_t)+{\mathbf{H}}(\gamma_t,V(\gamma_t),\partial_xV(\gamma_t),\partial_{xx}V(\gamma_t))= 0,\ \ \  t\in
%                               [0,T),\ \  \gamma_t\in {\Lambda},\cr
%  V(\gamma_T)=\phi(\gamma_T), \ \ \ \gamma_T\in {\Lambda}_T;\cr
%}
%\end{eqnarray}
%      where
%\begin{eqnarray*}
%                                {\mathbf{H}}(\gamma_t,r,p,l)&=&\sup_{u\in{
%                                         {U}}}[
%                        (p,F(\gamma_t,u))_{R^d}+\frac{1}{2}\mbox{tr}[ l G(\gamma_t,u)G^\top(\gamma_t,u)]\\
%                        &&\ \ \ \ \ +q(\gamma_t,r,G^\top(\gamma_t,u)p,u)],  \ \ (\gamma_t,r,p,l)\in {\Lambda}\times R\times R^d\times \Gamma(R^{d}).
%\end{eqnarray*}
%Here we let  $G^\top$ the transpose of the matrix $G$,  $\Gamma(R^d)$  the set of all $(d\times d)$ symmetric matrices and  $(\cdot,\cdot)_{R^d}$ the scalar product of $R^d$.
\par
\begin{definition}\label{definitionccc}     (Classical solution)
              A functional $v\in C_p^{1,2}({\Lambda})$  with  $A^*\partial_xv\in C_p^0({\Lambda})$    is called a classical solution to the PHJB equation (\ref{hjb1}) if it satisfies
              the PHJB¡¡  equation (\ref{hjb1}) point-wisely.
 \end{definition}
% \par
% \begin{definition}\label{definitionusc}     We define
 %\begin{eqnarray*}
%                           &&USC_*(\tilde{\Lambda}^M):=\{u:\tilde{\Lambda}\rightarrow R\  \mbox{such that} \\
%                           &&(i)  \ \mbox{For each fixed} \ \omega_{\hat{t}}\in \tilde{\Lambda}, \bar{u}(t,x):=u((\omega^x_{\hat{t}})_{\hat{t},t})
%                           \ \mbox{is a}\ USC\mbox{-function of } \ (t, x+\omega(\hat{t}))\in [0,T-\hat{t}]\times Q;\\
%                           &&(ii) \ \mbox{For each}\ \omega_t\in \tilde{\Lambda}^M \ \mbox{with}\ t_i\uparrow t,  \limsup_{i\rightarrow\infty}u(\omega_{t_i})
%                           \leq \sup_{|x+\omega(t)|\leq M}u(\omega_{t}^x)\}.
 %\end{eqnarray*}
%              A functional $v\in C^{1,2}(\tilde{\Lambda})$       is called a classical solution to the path-dependent HJB equation (\ref{hjb}) if it satisfies
%              the path-dependent HJB¡¡equation (\ref{hjb}) point-wisely.
% \end{definition}
        \par
        We will prove that the value functional $V$ defined by (\ref{value1}) is a viscosity solution of PHJB equation (\ref{hjb1}).
Define
$$
\Phi=\left\{\varphi\in C_p^{1,2}({\Lambda})| %\phi \ \mbox{is weakly sequentially lower semicontinuous},
A^*\partial_x\varphi\in C_p^0({\Lambda})\right\}.
$$
%We recall, for every $t\in [0,T)$, ${\cal{G}}_t$ is defined in
%\begin{eqnarray*}
%{\mathcal{G}}_t&=&\bigg{\{}g\in C^{1,2}_p(\Lambda^t)| g \ \mbox{satisfies conditon (\ref{4jia081})}\ \mbox{when}\  A \ \mbox{satisfies Hypothesis \ref{hypstate} (i') and} \ X \ \mbox{satisfies}\\
%&& \mbox{(\ref{formular2}) with} \ \vartheta\in L^{p}_{\mathcal{P}}(\Omega\times [0,T];H),
%\varpi\in L^{p}_{\mathcal{P}}(\Omega\times [0,T];L_2(\Xi,H)) \ \mbox{for some}\ p>2\  \mbox{and}\  (t,\gamma_{t})\in [0,T)\times \Lambda,
%\bigg{\}}
%\end{eqnarray*}
%where
%\begin{eqnarray}\label{4jia081}
%                  g(X_s)
%                  &\leq& g(X_t)+\int^{s}_{{t}}(\partial_xg(X_\sigma),\vartheta(\sigma))_H+\frac{1}{2}\mbox{Tr}(\partial_{xx}g(X_\sigma)\varpi(\sigma)
%                  \varpi^*(\sigma))d\sigma\nonumber\\
%                  &&
%                  +\int^{s}_{t}(\partial_xg(X_\sigma),\varpi(\sigma)dW(\sigma))_H,\ s\in [t,T], \ P\mbox{-}a.s.
%\end{eqnarray}
\begin{eqnarray*}
{\mathcal{G}}_t&=&\bigg{\{}g\in C^{0}_p(\Lambda^t)| \exists \  h\in C^1([0,T];\mathbb{R}), \delta_i, \delta'_i\geq0, \gamma^i_{t_i}\in \Lambda, t_i\leq t, i=0,1,2,\ldots, N>0, \ \mbox{with} \  h \geq0,\\
&& \ ||\gamma^i_{t_i}||_0\leq N, \sum_{i=0}^{\infty}(\delta_i+\delta'_i)\leq N \ \mbox{such that}\\
                  && g(\gamma_s)=h(s)\Upsilon^3(\gamma_s)+\sum_{i=0}^{\infty}[\delta_i\overline{\Upsilon}^3(\gamma_s-\gamma^i_{t_i,s,A})
                   +\delta'_i|\gamma_s(s)-e^{(s-t_i)A}\gamma^i_{t_i}(t_i)|^6], \forall\ \gamma_s\in \Lambda^t
\bigg{\}}.
\end{eqnarray*}
For notational simplicity, if $g\in {\mathcal{G}}_t$,
we use $\partial_tg(\gamma_s)$, $\partial_xg(\gamma_s)$ and $\partial_{xx}g(\gamma_s)$ to denote $h_t(s)\Upsilon^3(\gamma_s)+2\sum_{i=0}^{\infty}\delta_i(s-t_i)$,
 $h(s)\partial_x\Upsilon^3(\gamma_s)+\sum_{i=0}^{\infty}\partial_x[\delta_i\overline{\Upsilon}^3(\gamma_s-\gamma^i_{t_i,s,A})
                   +\delta'_i|\gamma_s(s)-e^{(s-t_i)A}\gamma^i_{t_i}(t_i)|^6]$ and \\ $h(s)\partial_{xx}\Upsilon^3(\gamma_s)+\sum_{i=0}^{\infty}\partial_{xx}[\delta_i\overline{\Upsilon}^3(\gamma_s-\gamma^i_{t_i,s,A})
                   +\delta'_i|\gamma_s(s)-e^{(s-t_i)A}\gamma^i_{t_i}(t_i)|^6]$, respectively.
%\begin{eqnarray*}
%{\mathcal{G}}_t&=&\{g\in C^1(\Lambda^t)|  (A_{\mu}\gamma_s(s),\partial_xg(\gamma_s)\leq0\ \mbox{for all}\ \gamma_s\in \Lambda^t\ \mbox{and}\ \mu>0\}.
%\exists \  g_i\in C^1(R), \gamma^i_{t_i}\in \Lambda, t_i\leq t, i=0,1,2,\ldots, M>0, \ \mbox{with}\ g'_i\geq0, \\
%&&\ ||\gamma^i_{t_i}||_0\leq M \ \mbox{such that}\
%                   g(\gamma_s)=h(s)\sum_{i=0}^{\infty}\frac{1}{2^i}g_i(\Upsilon^2(\gamma_s-\gamma^i_{t_i,s,A})), \forall\ \gamma_s\in \Lambda^t
%\}.
%\end{eqnarray*}
\par
                     Now we can give the following definition for  viscosity solutions.
\begin{definition}\label{definition4.1} \ \
 $w\in C^0({\Lambda})$ is called a
                             viscosity subsolution (resp., supersolution)
                             to  (\ref{hjb1}) if the terminal condition,  $w(\gamma_T)\leq \phi(\gamma_T)$(resp., $w(\gamma_T)\geq \phi(\gamma_T)$),
                             $\gamma_T\in {\Lambda}_T$ is satisfied, and for any $\varphi\in \Phi$ and $g\in {\cal{G}}_t$ with $t\in [0,T)$, whenever the function $w-\varphi-g$  (resp.,  $w+\varphi+g$) satisfies
$$
                         0=({w}-\varphi-g)(\gamma_t)=\sup_{(s,\eta_s)\in [t,T]\times\Lambda^t}
                         ({w}- \varphi-g)(\eta_s),
$$
$$
                       \left(\mbox{resp.,}\ \
                         0=({w}+\varphi+g)(\gamma_t)=\inf_{(s,\eta_s)\in [t,T]\times\Lambda^t}
                         ({w}+\varphi+g)(\eta_s),\right)
$$
                   %where $\hat{\gamma}^\mu_{s_\mu}\in {\cal{C}}^\alpha_{\mu,M_0}$, $s_\mu\in [0,T-\Delta]$ and $|\hat{\gamma}^\mu_{s_\mu}(s_\mu)|<
%                   \frac{M_0}{2}$,
              where $\gamma_t\in \Lambda$,      we have
\begin{eqnarray*}
                          &&\partial_t\varphi(\gamma_t)+\partial_tg(\gamma_t)+(A^*\partial_x\varphi(\gamma_t),\gamma_t(t))_H\\
                          &&
                           +{\mathbf{H}}(\gamma_t,\varphi(\gamma_t)+g(\gamma_t),\partial_x\varphi(\gamma_t)+\partial_xg(\gamma_t),\partial_{xx}\varphi(\gamma_t)+\partial_{xx}g(\gamma_t))\geq0,
\end{eqnarray*}
\begin{eqnarray*}
                          (\mbox{resp.,}\ &&-\partial_t\varphi(\gamma_t)-\partial_tg(\gamma_t)-(A^*\partial_x\varphi(\gamma_t),\gamma_t(t))_H\\
                           &&+{\mathbf{H}}(\gamma_t,-\varphi(\gamma_t)-g(\gamma_t),-\partial_x\varphi(\gamma_t)-\partial_xg(\gamma_t),-\partial_{xx}\varphi(\gamma_t)-\partial_{xx}g(\gamma_t))
                          \leq0.)
\end{eqnarray*}
                                $w\in C^0({\Lambda})$ is said to be a
                             viscosity solution to (\ref{hjb1}) if it is
                             both a viscosity subsolution and a viscosity
                             supersolution.
\end{definition}
%\begin{remark}\label{remarkv}
%  %\begin{description}
%%  \rm{
% % \item{(i)}  A viscosity solution $V$ of the PHJB equation (\ref{hjb})  is a
%%                       classical solution (See Definition \ref{definitionccc}) if it further lies in $C^{1,2}_{p,C_0}({\Lambda})$ for some constants $p,C_0>0$.
% %\item{(i)}    In the classical uniqueness proof of viscosity solution to second order HJB equation in infinite
%%                dimensions, the weak compactness of closed balls in separable
%%                Hilbert spaces is used  (See  \cite{swi}). In our case, the  path-dependent HJB
%%                equation is defined on space  ${\Lambda}$, which
%%                does not  have weak compactness.  In order to obtain  the
%%uniqueness of viscosity solutions, the compact subset  sequence $\{{\cal{C}}^\alpha_{\mu,M_0}\}_{\mu>0}$ are used to introduce the new concept of viscosity  solution.
%%\item{(ii)}
%  Assume that  the coefficients $F(\gamma_t,u)=\overline{F}(t,\gamma_t(t),u),
%                     \ G(\gamma_t,u)=\overline{G}(t,\gamma_t(t),u)$,  $ q(\gamma_t,$ $y,z,u)=\overline{q}(t,\gamma_t(t),y,z,u),
%                     \phi(\gamma_T)=\overline{\phi}(\gamma_T(T)), \
%                     (\gamma,y,z,u) \in {\Lambda}_T\times R\times R^d\times U$.
%           Let the  function $V:[0,T]\times R^d\rightarrow R$ be a
%             viscosity solution to (\ref{hjb1}) as a  functional of  $V(\gamma_t):{\Lambda}\rightarrow R$. Then, $V$ is also a classical viscosity
%solution as a function of the  state.
%%}
%%\end{description}
%\end{remark}
\begin{theorem}\label{theoremvexist} \ \
                          Suppose that Hypotheses \ref{hypstate} and  \ref{hypcost}  hold. Then the value
                          functional $V$ defined by (\ref{value1}) is a
                          viscosity solution to (\ref{hjb1}).
\end{theorem}
 {\bf  Proof}. \ \
 First, Let  $\varphi\in \Phi$ and $g\in {\cal{G}}_{\hat{t}}$ with $\hat{t}\in [0,T)$
 such that
\begin{eqnarray}\label{0714}
                         0=(V-\varphi-g)(\hat{\gamma}_{\hat{t}})=\sup_{(s,\eta_s)\in [\hat{t},T]\times\Lambda^{\hat{t}}}
                         (V- \varphi-g)(\eta_s),
\end{eqnarray}
 where $\hat{\gamma}_{\hat{t}}\in \Lambda$.
 Thus,  by the DPP (Theorem \ref{theoremddp}), we obtain that, for all ${\hat{t}}\leq {\hat{t}}+\delta\leq T$,
% We let  $\varphi\in J^+(\hat{\gamma}_{\hat{t}},V)$
%                  with
%                   $(\hat{t},\hat{\gamma}_{\hat{t}})\in [0,T)\times \Lambda$. % and  $\delta_{1}:=(M_0-|\hat{\gamma}_{\hat{t}}(\hat{t})|)>\frac{M_0}{2}$.
 %For any $u\in {\cal{U}}[\hat{t},\hat{t}+\delta]$, we define an ${\cal{F}}_{\hat{t}+\delta}$-measurable set
%\begin{eqnarray*}
%                          A_{\delta,u}:=
%                                              \bigg{\{}||X_{\hat{t}+\delta}^{\hat{\gamma}_{\hat{t}},u}||_0> M_0\bigg{\}}.
%\end{eqnarray*}
% From Lemma \ref{lemmaexist}, it follows that
% \begin{eqnarray}\label{4.3333}
% P(A_{\delta,u})\leq P\bigg{\{}||X_{{\hat{t}}+\delta}^{\hat{\gamma}_{\hat{t}},u}-\hat{\gamma}_{{\hat{t}}}||_0>\frac{\delta_1}{2}\bigg{\}}\leq C_6(1+M_0^6)\frac{\delta^3}{\delta^6_1},
%\end{eqnarray}
%For  $0< \delta\leq T-\hat{t}$, we have $\hat{t}< \hat{t}+\delta \leq T$, then by the DPP (Theorem \ref{theoremddp}), we obtain the following result:
 \begin{eqnarray}\label{4.9}
                           && 0=V(\hat{\gamma}_{\hat{t}})-({{\varphi}}+g) (\hat{\gamma}_{\hat{t}})
                           =\sup_{u(\cdot)\in {\cal{U}}[\hat{t},\hat{t}+\delta]} G^{\hat{\gamma}_{\hat{t}},u}_{{\hat{t}},\hat{t}+\delta}\left[V\left(X^{\hat{\gamma}_{\hat{t}},u}_{{\hat{t}}+\delta}\right)\right]
                           -({{\varphi}}+g) (\hat{\gamma}_{\hat{t}}).
\end{eqnarray}
Then, for any $\varepsilon>0$ and $0<\delta\leq T-\hat{t}$,  we can  find a control  ${u}^{\varepsilon}(\cdot)\equiv u^{{\varepsilon},\delta}(\cdot)\in {\cal{U}}[\hat{t},\hat{t}+\delta]$ such
   that the following result holds:
\begin{eqnarray}\label{4.10}
    -{\varepsilon}\delta%+V(\hat{\gamma}_{\hat{t}})-{{\varphi}} (\hat{\gamma}_{\hat{t}})
    \leq G^{\hat{\gamma}_{\hat{t}},{u}^{\varepsilon}}_{{\hat{t}},\hat{t}+\delta}\left[V\left(X^{\hat{\gamma}_{\hat{t}},{u}^{\varepsilon}}_{{\hat{t}}+\delta}\right)\right]-({{\varphi}}+g) (\hat{\gamma}_{\hat{t}}).
\end{eqnarray}
 We note that
                     $G^{\hat{\gamma}_{\hat{t}},{u}^{\varepsilon}}_{s,\hat{t}+\delta}\left[V\left(X^{\hat{\gamma}_{\hat{t}},{u}^{\varepsilon}}_{{\hat{t}}+\delta}\right)\right]$
                     is defined in terms of the solution of the
                     BSDE:
 \begin{eqnarray}\label{bsde4.10}
\begin{cases}
dY^{\hat{\gamma}_{\hat{t}},{u}^{\varepsilon}}\left(s\right) =
               -q\left(X^{\hat{\gamma}_{\hat{t}},{u}^{\varepsilon}}_s,Y^{\hat{\gamma}_{\hat{t}},{u}^{\varepsilon}}\left(s\right),Z^{\hat{\gamma}_{\hat{t}},
               {u}^{\varepsilon}}\left(s\right),{u}^{\varepsilon}\left(s\right)\right)ds+Z^{\hat{\gamma}_{\hat{t}},{u}^{\varepsilon}}\left(s\right)dW\left(s\right),\   s\in[\hat{t},\hat{t}+\delta], \\
 ~Y^{\hat{\gamma}_{\hat{t}},{u}^{\varepsilon}}\left(\hat{t}+\delta\right)=V\left(X^{\hat{\gamma}_{\hat{t}},
 {u}^{\varepsilon}}_{{\hat{t}}+\delta}\right),
 \end{cases}
\end{eqnarray}
                     by the following formula:
$$
                        G^{\hat{\gamma}_{\hat{t}},{u}^{\varepsilon}}_{s,\hat{t}+\delta}\left[V\left(X^{\hat{\gamma}_{\hat{t}},{u}^{\varepsilon}}_{{\hat{t}}+\delta}\right)\right]
                        =Y^{\hat{\gamma}_{\hat{t}},{u}^{\varepsilon}}(s),\
                        \ \ s\in[\hat{t},\hat{t}+\delta].
$$
%Combining (\ref{4.10}) and (\ref{bsde4.10}), we obtain
%\begin{eqnarray}\label{4.11}
%                           -\varepsilon\delta\leq E\bigg{[}{\varphi}(X^{\hat{\gamma}_{\hat{t}},{u}^\varepsilon}_{{\hat{t}}+\delta})-{{\varphi}} (\hat{\gamma}_{\hat{t}})
%                           +\int^{t+\delta}_{{\hat{t}}}g(s,X^{\gamma_{\hat{t}},{u}^\varepsilon}(s),Y^{\hat{\gamma}_{\hat{t}},{u}^\varepsilon}(s),Z^{\hat{\gamma}_{\hat{t}},
%                           {u}^\varepsilon}(s),{u}^\varepsilon(s))ds\bigg{|}{\cal{F}}_{\hat{t}}\bigg{]}.
%\end{eqnarray}
                         Applying functional It\^{o} formula (\ref{statesop0}) to ${\varphi}\left(X^{\hat{\gamma}_{\hat{t}},{u}^{\varepsilon}}_s\right)$ and  inequality (\ref{jias510815jia}) to
                         $g\left(X^{\hat{\gamma}_{\hat{t}},{u}^{\varepsilon}}_s\right)$,   we get that
%\begin{eqnarray}\label{bsde4.21}
%                             d{\varphi}(X^{\gamma_{\hat{t}},u}_s)
%                             &=&
%                             [{\partial_t{\varphi}}(X^{\gamma_{\hat{t}},u}_s)
%                             +\frac{1}{2}\mbox{tr}[\partial_{xx}{\varphi}(X^{\gamma_{\hat{t}},u}_s)G(X^{\gamma_{\hat{t}},u}_s,u(s))G^{\top}(X^{\gamma_{\hat{t}},u}_s,u(s))]\nonumber\\
%                             &&
%                                         +(\partial_x
%                                         {\varphi}(X^{\gamma_{\hat{t}},u}_s),F(s,X^{\gamma_{\hat{t}},u}_s,u(s)))_{R^d} ] ds
%                                         +(\partial_x{\varphi}(X^{\gamma_{\hat{t}},u}_s))^\top G(X^{\gamma_{\hat{t}},u}_s,u(s))dW(s). \ \ \ \ \ \ \
%\end{eqnarray}
\begin{eqnarray}\label{bsde4.21}
                            && \left({\varphi}+g\right)\left(X^{\hat{\gamma}_{\hat{t}},{u}^{\varepsilon}}_s\right)\leq Y^1\left(s\right)
                             =: \left({\varphi}+g\right)\left(\hat{\gamma}_{\hat{t}}\right)+\int^{s}_{{\hat{t}}} \left({\cal{L}}\left({\varphi}+g\right)\right)\left(X^{\hat{\gamma}_{\hat{t}},{u}^{\varepsilon}}_\sigma,
                                      u^{{\varepsilon}}\left(\sigma\right)\right)d\sigma \nonumber\\
                             &&-\int^{s}_{{\hat{t}}}q\left(X^{\hat{\gamma}_{\hat{t}},{u}^{\varepsilon}}_\sigma,\left({\varphi}+g\right)\left(X^{\hat{\gamma}_{\hat{t}},
                                      {u}^{\varepsilon}}_\sigma\right),\partial_x\left({\varphi}+g\right)\left(X^{\hat{\gamma}_{\hat{t}},{u}^{\varepsilon}}_\sigma\right)
                              G\left(X^{\hat{\gamma}_{\hat{t}},{u}^{\varepsilon}}_\sigma,u^{{\varepsilon}}\left(\sigma\right)\right), u^{\varepsilon}\left(\sigma\right)\right)d\sigma\nonumber\\
                             &&
                             +\int^{s}_{{\hat{t}}}\left(\partial_x\left({\varphi}+g\right)\left(X^{\hat{\gamma}_{\hat{t}},{u}^{\varepsilon}}_\sigma\right),
                             G\left(X^{\hat{\gamma}_{\hat{t}},{u}^{\varepsilon}}_\sigma,u^{\varepsilon}\left(\sigma\right)\right)dW\left(\sigma\right)\right)_H,\ \
\end{eqnarray}
where
\begin{eqnarray*}
                       \left({\cal{L}}\left({\varphi}+g\right)\right)\left(\gamma_t,u\right)
                       &=&\partial_t\left({\varphi}+g\right)\left(\gamma_t\right)+\left(A^*\partial_x {\varphi}\left(\gamma_t\right),\gamma_t\left(t\right)\right)_H
                                         +\left(\partial_x \left({\varphi}+g\right)\left(\gamma_t\right),F\left(\gamma_t,u\right)\right)_{H}\\
                                         &&+\frac{1}{2}\mbox{Tr}[\partial_{xx}\left({\varphi}+g\right)\left(\gamma_t\right)G\left(\gamma_t,u\right)G\left(\gamma_t,u\right)^*]\\
                       &&+
                                         q\left(\gamma_t,\left({\varphi}+g\right)\left(\gamma_t\right),\left({\partial_x\left({\varphi}+g\right)\left(\gamma_t\right)}\right)G\left(\gamma_t,u\right),u\right),
                                           \left(t, \gamma_t,u\right)\in [0,T]\times {\Lambda}\times U.
\end{eqnarray*}
It is clear that
\begin{eqnarray}\label{y_1}
Y_1(\hat{t})=({\varphi}+g)(\hat{\gamma}_{\hat{t}}).
\end{eqnarray}
Set
\begin{eqnarray*}
                          Y^{2,{\hat{\gamma}_{\hat{t}},{u}^{\varepsilon}}}\left(s\right)&:=&
                          Y^1\left(s\right)-Y^{\hat{\gamma}_{\hat{t}},u^{\varepsilon}}\left(s\right)\geq
                          \left({\varphi}+g\right)\left(X_s^{\hat{\gamma}_{\hat{t}},{u}^{\varepsilon}}\right)-Y^{\hat{\gamma}_{\hat{t}},u^{\varepsilon}}\left(s\right),  \ \ s\in[\hat{t},\hat{t}+\delta],\\
                             Y^{3}\left(s\right)&:=&
                          Y^1\left(s\right)-
                          \left({\varphi}+g\right)\left(X_s^{\hat{\gamma}_{\hat{t}},{u}^{\varepsilon}}\right),  \ \ s\in[\hat{t},\hat{t}+\delta],\\
                            Z^{2,{\hat{\gamma}_{\hat{t}},{u}^{\varepsilon}}}\left(s\right)  &:=&\partial_x\left({\varphi}+g\right)\left(X^{\hat{\gamma}_{\hat{t}},{u}^{\varepsilon}}_s\right)
                             G\left(X^{\hat{\gamma}_{\hat{t}},{u}^{\varepsilon}}_s,u^{\varepsilon}\left(s\right)\right)-Z^{\hat{\gamma}_{\hat{t}},u^{\varepsilon}}\left(s\right), \ \ s\in[\hat{t},\hat{t}+\delta].
\end{eqnarray*}
Comparing (\ref{bsde4.10}) and (\ref{bsde4.21}), we have, $P$-a.s.,
\begin{eqnarray*}
                      dY^{2,{\hat{\gamma}_{\hat{t}},{u}^{\varepsilon}}}\left(s\right)
                      &=&\bigg{[}\left({\cal{L}}\left({\varphi}+g\right)\right)\left(X^{\hat{\gamma}_{\hat{t}},{u}^{\varepsilon}}_s,u^{\varepsilon}\left(s\right)\right)
                       -q\bigg{(}X^{\hat{\gamma}_{\hat{t}},{u}^{\varepsilon}}_s,\left({\varphi}+g\right)\left(X^{\hat{\gamma}_{\hat{t}},{u}^{\varepsilon}}_s\right),\\
                       &&
                             \partial_x\left({\varphi}+g\right)\left(X^{\hat{\gamma}_{\hat{t}},{u}^{\varepsilon}}_s\right)
                             G\left(X^{\hat{\gamma}_{\hat{t}},{u}^{\varepsilon}}_s,{u}^{\varepsilon}\left(s\right)\right), {u}^{\varepsilon}\left(s\right)\bigg{)}\\
                             &&+q\left(X^{\hat{\gamma}_{\hat{t}},{u}^{\varepsilon}}_s,Y^{\hat{\gamma}_{\hat{t}},{u}^{\varepsilon}}\left(s\right),
                                  Z^{\hat{\gamma}_{\hat{t}},{u}^{\varepsilon}}\left(s\right),{u}^{\varepsilon}\left(s\right)\right)\bigg{]}ds
                                +Z^{2,{\hat{\gamma}_{\hat{t}},{u}^{\varepsilon}}}\left(s\right)dW\left(s\right)\\
                                &=&\bigg{[}\left({\cal{L}}\left({\varphi}+g\right)\right)\left(X^{\hat{\gamma}_{\hat{t}},{u}^{\varepsilon}}_s,u^{\varepsilon}\left(s\right)\right)+A_1\left(s\right)Y^3\left(s\right)
                                -A_1\left(s\right)Y^{2,{\hat{\gamma}_{\hat{t}},{u}^{\varepsilon}}}\left(s\right)\\
                                &&-\left(A_2\left(s\right),Z^{2,{\hat{\gamma}_{\hat{t}},{u}^{\varepsilon}}}\left(s\right)\right)_{\Xi}\bigg{]}ds+Z^{2,{\hat{\gamma}_{\hat{t}},{u}^{\varepsilon}}}\left(s\right)dW\left(s\right),
\end{eqnarray*}
where $|A_1|\vee|A_2|\leq L$. %($C$ only depends on  Lipschitz constant of $q$).
Applying It\^o  formula  (see also Proposition 2.2 in \cite{el}), we obtain
\begin{eqnarray}\label{4.14}
                      Y^{2,{\hat{\gamma}_{\hat{t}},{u}^{\varepsilon}}}(\hat{t})
                      &=&\mathbb{E}\bigg{[}Y^{2,{\hat{\gamma}_{\hat{t}},{u}^{\varepsilon}}}(t+\delta)\Gamma^{\hat{t}}(\hat{t}+\delta)\nonumber\\
                      &&-
                      \int^{\hat{t}+\delta}_{{\hat{t}}}\Gamma^{\hat{t}}\left(\sigma\right)\left[\left({\cal{L}}\left({\varphi}+g\right)\right)
                       \left(X^{\hat{\gamma}_{\hat{t}},{u}^{\varepsilon}}_\sigma,u^{\varepsilon}\left(\sigma\right)\right)+A_1\left(\sigma\right)Y^3\left(\sigma\right)\right]d\sigma\bigg{|}
                      {\cal{F}}_{\hat{t}}\bigg{]},\ \ \
\end{eqnarray}
where $\Gamma^{\hat{t}}(\cdot)$ solves the linear SDE
$$
               d\Gamma^{\hat{t}}(s)=\Gamma^{\hat{t}}(s)(A_1(s)ds+A_2(s)dW(s)),\ s\in [{\hat{t}},{\hat{t}}+\delta];\ \ \ \Gamma^{\hat{t}}({\hat{t}})=1.
$$
Obviously, $\Gamma^{\hat{t}}\geq 0$. Combining (\ref{4.10}) and (\ref{4.14}), we have
\begin{eqnarray}\label{4.15}
-\varepsilon%+\frac{1}{\delta}[V(\hat{\gamma}_{\hat{t}})-{({\varphi}+g)} (\hat{\gamma}_{\hat{t}})]\nonumber\\
    &\leq& \frac{1}{\delta}\mathbb{E}\bigg{[}-Y^{2,{\hat{\gamma}_{\hat{t}},{u}^{\varepsilon}}}(\hat{t}+\delta)\Gamma^{\hat{t}}(\hat{t}+\delta)+
    \int^{\hat{t}+\delta}_{{\hat{t}}}\Gamma^{\hat{t}}(\sigma)\left[({\cal{L}}({\varphi}+g))\left(X^{\hat{\gamma}_{\hat{t}},{u}^{\varepsilon}}_\sigma,u^{\varepsilon}(\sigma)\right)
    +A_1(\sigma)Y^3(\sigma)\right]d\sigma\bigg{]}\nonumber\\
    &=&-\frac{1}{\delta}\mathbb{E}\bigg{[}Y^{2,{\hat{\gamma}_{\hat{t}},{u}^{\varepsilon}}}(\hat{t}+\delta)\Gamma^{\hat{t}}(\hat{t}+\delta)\bigg{]}+
      \frac{1}{\delta}\mathbb{E}\bigg{[}\int^{\hat{t}+\delta}_{{\hat{t}}}({\cal{L}}({\varphi}+g))(\hat{\gamma}_{{\hat{t}}},u^{\varepsilon}(\sigma))d\sigma\bigg{]}\nonumber\\
                      &&+\frac{1}{\delta}\mathbb{E}\bigg{[}\int^{\hat{t}+\delta}_{{\hat{t}}}\left[({\cal{L}}({\varphi}+g))\left(X^{\hat{\gamma}_{\hat{t}},{u}^{\varepsilon}}_\sigma,u^{\varepsilon}(\sigma)\right)-
                      ({\cal{L}}({\varphi}+g))\left(\hat{\gamma}_{{\hat{t}}},u^{\varepsilon}(\sigma)\right)\right]d\sigma\bigg{]}\nonumber\\
                      &&+\frac{1}{\delta}\mathbb{E}\bigg{[}\int^{\hat{t}+\delta}_{{\hat{t}}}(\Gamma^{\hat{t}}(\sigma)-1)({\cal{L}}({\varphi}+g))\left(X^{\hat{\gamma}_{\hat{t}},{u}^{\varepsilon}}_\sigma,
                                u^{\varepsilon}(\sigma)\right)d\sigma\bigg{]}\nonumber\\
                                &&+\frac{1}{\delta}\mathbb{E}\bigg{[}\int^{\hat{t}+\delta}_{{\hat{t}}}\Gamma^{\hat{t}}(\sigma)A_1(\sigma)Y^3(\sigma)d\sigma\bigg{]}\nonumber\\
    &:=&I+II+III+IV+V.
\end{eqnarray}
Since the coefficients in ${\cal{L}}$  satisfy  linear growth  condition,
 combining the regularity of $\varphi\in \Phi$ and $g\in {\cal{G}}_{\hat{t}}$, there exist a integer
$\bar{p}\geq1$ and a constant $C>0$ independent of $u\in U$ such that, for all $(t,\gamma_t,u)\in [0,T]\times\Lambda\times U$,
\begin{eqnarray}\label{4.4444}|({\varphi}+g)(\gamma_{t})|\vee  |
                      ({\cal{L}}({\varphi}+g))(\gamma_{t},u)|
                      \leq  C(1+||\gamma_{t}||_0)^{\bar{p}}.
\end{eqnarray}
In view of Lemma \ref{lemmaexist111}, we also have
\begin{eqnarray*}
                        %  &&\sup_{u\in{\cal{U}}[\hat{t},\hat{t}+\delta]}E[\sup_{{\hat{t}}\leq s\leq \hat{t}+\delta}|X^{\hat{\gamma}_{\hat{t}},{u}}(s)-\hat{\gamma}_{\hat{t}}(\hat{t})|^8]\leq C\delta^4,\\
                          \sup_{u(\cdot)\in{\cal{U}}[\hat{t},\hat{t}+\delta]}\mathbb{E}|\Gamma^{\hat{t}}(\hat{t}+\delta)-1|^2\leq C\delta.
\end{eqnarray*}
Thus, by (\ref{0714}) and (\ref{bsde4.21}),
%by (\ref{4.3333}) and the definition of $A_{\delta,u}$, we obtain, for some constant $C>0$,
\begin{eqnarray}\label{4.1611}
I&=& -\frac{1}{\delta}\mathbb{E}\left[\left(Y^1(\hat{t}+\delta)-Y^{\hat{\gamma}_{\hat{t}},u^{\varepsilon}}(\hat{t}+\delta)\right)\Gamma^{\hat{t}}(\hat{t}+\delta)\right]\nonumber\\
&\leq&-\frac{1}{\delta}\mathbb{E}\left[\left(({\varphi}+g)\left(X_{\hat{t}+\delta}^{\hat{\gamma}_{\hat{t}},{u}^{\varepsilon}}\right)-Y^{\hat{\gamma}_{\hat{t}},u^{\varepsilon}}(\hat{t}+\delta)\right)
\Gamma^{\hat{t}}(\hat{t}+\delta)\right]\nonumber\\
&=&\frac{1}{\delta}\mathbb{E}\left[\left(V\left(X^{\hat{\gamma}_{\hat{t}},
 {u}^{\varepsilon}}_{{\hat{t}}+\delta}\right)-({\varphi}+g)\left(X_{\hat{t}+\delta}^{\hat{\gamma}_{\hat{t}},{u}^{\varepsilon}}\right)\right)\Gamma^{\hat{t}}(\hat{t}+\delta)\right]\leq 0;
\end{eqnarray}
\begin{eqnarray}\label{4.16}
             II&\leq&\frac{1}{\delta}\bigg{[}\int^{\hat{t}+\delta}_{{\hat{t}}}\sup_{u\in U}({\cal{L}}({\varphi}+g))({\hat{\gamma}_{\hat{t}},{u}})d\sigma\bigg{]}
                      = \partial_t({\varphi}+g)(\hat{\gamma}_{\hat{t}})+(A^*\partial_x{\varphi}(\hat{\gamma}_{\hat{t}}),\hat{\gamma}_{\hat{t}}(\hat{t}))_H\nonumber\\
                           &&+{\mathbf{H}}(\hat{\gamma}_{\hat{t}},({\varphi}+g)(\hat{\gamma}_{\hat{t}}),\partial_x({\varphi}+g)(\hat{\gamma}_{\hat{t}}),
                           \partial_{xx}({\varphi}+g)(\hat{\gamma}_{\hat{t}})).
\end{eqnarray}
  \par
Now we estimate higher order terms  $III$, $IV$ and $V$.
By (\ref{4.4444}) and the dominated convergence theorem, we have
$$
\lim_{\sigma\rightarrow\hat{t}} \mathbb{E}\left|({\cal{L}}({\varphi}+g))\left(X^{{\hat{\gamma}_{\hat{t}},{u}^{\varepsilon}}}_\sigma,{u}^{\varepsilon}(\sigma)\right)-
                      ({\cal{L}}({\varphi}+g))({\hat{\gamma}_{\hat{t}},{u}^{\varepsilon}}(\sigma))\right|=0,
$$
and
\begin{eqnarray*}
&&\lim_{\sigma\rightarrow\hat{t}} \mathbb{E}\left|\Gamma^{\hat{t}}(\sigma)A_1(\sigma)Y^3(\sigma)\right|\\
&\leq& L\lim_{\sigma\rightarrow\hat{t}} \mathbb{E}\left|\Gamma^{\hat{t}}(\sigma)\right|\left[|Y^1(\sigma)-({\varphi}+g)(\hat{\gamma}_{\hat{t}})|
+\left|({\varphi}+g)\left(X^{\hat{\gamma}_{\hat{t}},{u}^{\varepsilon}}_\sigma\right)-({\varphi}+g)(\hat{\gamma}_{\hat{t}})\right|\right] =0,
\end{eqnarray*}
then
\begin{eqnarray}\label{4.18}
\limsup_{\delta\rightarrow0}|III| \leq\frac{1}{\delta}\int^{\hat{t}+\delta}_{{\hat{t}}}\mathbb{E}\left|({\cal{L}}({\varphi}+g))\left(X^{{\hat{\gamma}_{\hat{t}},{u}^{\varepsilon}}}_\sigma,{u}^{\varepsilon}(\sigma)\right)-
                      ({\cal{L}}({\varphi}+g))({\hat{\gamma}_{\hat{t}},{u}^{\varepsilon}}(\sigma))\right|d\sigma=0;
                      %\nonumber\\
%                      &\leq&\frac{C}{\delta}\int^{t+\delta}_{t}E(1+|X^{{{u^{\varepsilon_1}}}}_\sigma(t,x)|_H+|X^{{{u^{\varepsilon_1}}}}(\sigma,t,x)|)^{\bar{p}}
%                      [\sup_{t\leq l\leq t+\sigma}|X^{{{u^{\varepsilon_1}}}}(l,t,x)-x(0)|+|\sigma-t|]d\sigma\nonumber\\
%                      &\leq&{C}(1+|x|_H+|x(0)|)^{\bar{p}}[(1+M_0^2)\delta^{\frac{1}{2}}+\delta+\delta^\alpha];
\end{eqnarray}
\begin{eqnarray}\label{4.180714}
\limsup_{\delta\rightarrow0}|V| \leq\frac{1}{\delta}\int^{\hat{t}+\delta}_{{\hat{t}}}\mathbb{E}\left|\Gamma^{\hat{t}}(\sigma)A_1(\sigma)Y^3(\sigma)\right|d\sigma=0.
\end{eqnarray}
Moreover, for some finite constant $C>0$,
\begin{eqnarray}\label{4.19}
|IV|
&\leq&\frac{1}{\delta}\int^{\hat{t}+\delta}_{{\hat{t}}}\mathbb{E}|\Gamma^{\hat{t}}(\sigma)-1|\left|({\cal{L}}({\varphi}+g))(X^{\hat{\gamma}_{\hat{t}},{u}^{\varepsilon}}_\sigma,
                               u^{\varepsilon}(\sigma)\right|
                       d\sigma\nonumber\\
                       &\leq&\frac{1}{\delta}\int^{\hat{t}+\delta}_{{\hat{t}}}\left(\mathbb{E}\left(\Gamma^{\hat{t}}\left(\sigma\right)-1\right)^2\right)^{\frac{1}{2}}
                       \left(\mathbb{E}\left[\left({\cal{L}}\left({\varphi}+g\right)\right)\left(X^{\hat{\gamma}_{\hat{t}},
                                     {u}^{\varepsilon}}_\sigma,
                               u^{\varepsilon}\left(\sigma\right)\right)\right]^2\right)^{\frac{1}{2}}d\sigma
                    \nonumber\\
                       &\leq& C(1+||\hat{\gamma}_{\hat{t}}||_0)^{\bar{p}}\delta^\frac{1}{2}.
\end{eqnarray}
Substituting  (\ref{4.1611}), (\ref{4.16}) and (\ref{4.19}) into (\ref{4.15}), we have
\begin{eqnarray}\label{4.2000000}
-\varepsilon&\leq& \partial_t({\varphi}+g)(\hat{\gamma}_{\hat{t}}) +(A^*\partial_x{\varphi}(\hat{\gamma}_{\hat{t}}),\hat{\gamma}_{\hat{t}}(\hat{t}))_H
                           +{\mathbf{H}}(\hat{\gamma}_{\hat{t}},({\varphi}+g)(\hat{\gamma}_{\hat{t}}),\partial_x({\varphi}+g)(\hat{\gamma}_{\hat{t}}),
                           \partial_{xx}({\varphi}+g)(\hat{\gamma}_{\hat{t}}))\nonumber\\
    &&+|III|+|V|
    +C(1+||\hat{\gamma}_{\hat{t}}||_0)^{\bar{p}}\delta^\frac{1}{2}.
\end{eqnarray}
Sending $\delta$ to $0$,  by  (\ref{4.18}) and (\ref{4.180714}),  we have
$$
-\varepsilon
    \leq\partial_t\left({\varphi}+g\right)\left(\hat{\gamma}_{\hat{t}}\right)+\left(A^*\partial_x{\varphi}\left(\hat{\gamma}_{\hat{t}}\right),\hat{\gamma}_{\hat{t}}\left(\hat{t}\right)\right)_H
                           +{\mathbf{H}}\left(\hat{\gamma}_{\hat{t}},\left({\varphi}+g\right)\left(\hat{\gamma}_{\hat{t}}\right),\partial_x\left({\varphi}+g\right)\left(\hat{\gamma}_{\hat{t}}\right),
                           \partial_{xx}\left({\varphi}+g\right)\left(\hat{\gamma}_{\hat{t}}\right)\right).
$$
By the arbitrariness of $\varepsilon$, we show
 $V$ is  a viscosity subsolution to (\ref{hjb1}).
 \par
 In a symmetric (even easier) way, we show that  $V$ is also a viscosity supsolution to equation (\ref{hjb1}).
                 This step  completes the proof.\ \ $\Box$
 \par
Now, let us give   the result of  classical  solutions, which show the consistency of viscosity solutions.
 \begin{theorem}\label{theorem3.2}
                      Let $V$ denote the value functional  defined by (\ref{value1}). If  $V\in C^{1,2}_p({\Lambda})$ and $A^*\partial_xV\in C^0_p({\Lambda})$, then
                 $V$ is a classical solution of  equation (\ref{hjb1}).
\end{theorem}
{\bf  Proof }. \ \
                 %We only prove $V(t,x)$ satisfies the  HJB equation (\ref{hjb}). $V^m(t,x)$ satisfies the  HJB equation (\ref{hjbm}) can be proven in the similar
%                 way.
First, using the  definition of $V$ yields  $V(\gamma_T)=\phi(\gamma_T)$ for all $\gamma_T\in \Lambda_T$. Next, for fixed  $(t,\gamma_t,u)\in[0,T)\times {\Lambda}\times U$,
                    from  the DPP (Theorem \ref{theoremddp}), we obtain the following result:
 \begin{eqnarray}\label{4.9000}
                           && 0\geq G^{\gamma_t,u}_{t,t+\delta}[V(X^{{\gamma}_{{t}},u}_{{{t}}+\delta})]
                           -V(\gamma_t),\ \ 0\leq\delta\leq T-t.
\end{eqnarray}
% We note that
%                     $G^{\gamma_t,u}_{t,t+\delta}[V(X^{{\gamma}_{{t}},u}_{{{t}}+\delta})]$
%                     is defined in terms of the solution of the
%                     BSDE:
% \begin{eqnarray}\label{bsde4.10000}
%\cases{
%dY^{{\gamma}_{{t}},{u}}(s) =
%               -q(X^{{\gamma}_{{t}},{u}}_s,Y^{{\gamma}_{{t}},{u}}(s),Z^{{\gamma}_{{t}},
%               {u}}(s),{u})ds+Z^{{\gamma}_{{t}},{u}}(s)dW(s),\ \  s\in[{t},{t}+\delta], \cr
% ~Y^{{\gamma}_{{t}},{u}}({t}+\delta)=V(X^{{\gamma}_{{t}},
% {u}}_{{{t}}+\delta}),\cr
%}
%\end{eqnarray}
%                     by the following formula:
%$$
%                        G^{\gamma_t,u}_{s,t+\delta}[V(X^{{\gamma}_{{t}},u}_{{{t}}+\delta})]
%                        =Y^{{\gamma}_{{t}},{u}}(s),\
%                        \ \ s\in[{t},{t}+\delta].
%$$
Thus
\begin{eqnarray*}
                 V(\gamma_t)\geq Y^{{\gamma}_{{t}},{u}}(t)=\int_{t}^{t+\delta}q(X^{{\gamma}_{{t}},{u}}_s,Y^{{\gamma}_{{t}},{u}}(s),Z^{{\gamma}_{{t}},
               {u}}(s),{u})ds
                                +V(X^{{\gamma}_{{t}},u}_{{{t}}+\delta})-\int_{t}^{t+\delta}Z^{{\gamma}_{{t}},
               {u}}(s)dW(s).
\end{eqnarray*}
              Applying functional It\^{o} formula (\ref{statesop0}) to
                         $V(X^{{\gamma}_{{t}},u}_{{{t}}+\delta})$,   the  inequality above implies that, for all $ 0\leq\delta\leq T-t$,
\begin{eqnarray*}
    0&\geq&\int_{t}^{t+\delta}q(X^{{\gamma}_{{t}},{u}}_s,Y^{{\gamma}_{{t}},{u}}(s),Z^{{\gamma}_{{t}},
               {u}}(s),{u})ds+\int^{t+\delta}_{t} ({\cal{L}}{V})(X_s^{{\gamma}_{{t}},{u}},
                                      u)ds,
                                      \\
                      &&-\int_{t}^{t+\delta}q(X^{{\gamma}_{{t}},{u}}_s,V(X^{{\gamma}_{{t}},{u}}_s),
                      {\partial_xV(X^{{\gamma}_{{t}},{u}}_s)}G(X^{{\gamma}_{{t}},{u}}_s,u),{u})ds, \\
                      &&+\int_{t}^{t+\delta}({\partial_xV(X^{{\gamma}_{{t}},{u}}_s)}G(X^{{\gamma}_{{t}},{u}}_s,u)-Z^{{\gamma}_{{t}},
               {u}}(s))dW(s),\\
               &=&\int^{t+\delta}_{t} [({\cal{L}}{V})(X_s^{{\gamma}_{{t}},{u}},
                                      u)ds -A_1(s)Y^{2}(s)-(A_2(s),Z^{2}(s))_{\Xi}]ds
                                +\int_{t}^{t+\delta}Z^{2}(s)dW(s),
\end{eqnarray*}
where
\begin{eqnarray*}
                          &&Y^{2}(s):=
                          V(X_s^{{\gamma}_{{t}},{u}})-Y^{{\gamma}_{{t}},u}(s),  \ \ s\in[{t},{t}+\delta],\\
                            &&Z^{2}(s)  :=({\partial_xV(X^{{\gamma}_{{t}},{u}}_s)})G(X^{{\gamma}_{{t}},{u}}_s,u)-Z^{{\gamma}_{{t}},
               {u}}(s), \ \ s\in[{t},{t}+\delta],
\end{eqnarray*}
and  $|A_1|\vee|A_2|\leq L$. %($C$ only depends on  Lipschitz constant of $q$).
Applying It\^o formula (see also Proposition 2.2 in \cite{el}),  we obtain
\begin{eqnarray}\label{4.1411}
                      Y^{2}({t})
                      =\mathbb{E}\bigg{[}Y^{2}(t+\delta)\Gamma^{{t}}({t}+\delta)
                      -
                      \int^{{t}+\delta}_{{{t}}}\Gamma^{{t}}(\sigma)({\cal{L}}V)
                       (X^{{\gamma}_{{t}},{u}}_\sigma,u(\sigma))d\sigma\bigg{|}
                      {\cal{F}}_{{t}}\bigg{]},\ \ \
\end{eqnarray}
where $\Gamma^{{t}}(\cdot)$ solves the linear SDE
$$
               d\Gamma^{{t}}(s)=\Gamma^{{t}}(s)(A_1(s)ds+A_2(s)dW(s)),\ s\in [{{t}},{{t}}+\delta];\ \ \ \Gamma^{{t}}({{t}})=1.
$$
Obviously, $\Gamma^{{t}}\geq 0$.
Combining (\ref{4.9000}) and (\ref{4.1411}), we obtain that
\begin{eqnarray*}
0%+\frac{1}{\delta}[V(\hat{\gamma}^{\mu}_{{t}})-{{\varphi}} (\hat{\gamma}^{\mu}_{{t}})]\nonumber\\
    &\geq&
      \frac{1}{\delta}\mathbb{E}\bigg{[}\int^{{t}+\delta}_{{{t}}}({\cal{L}}{V})({\gamma}_{{{t}}},u)d\sigma\bigg{]}+\frac{1}{\delta}\mathbb{E}\bigg{[}\int^{{t}+\delta}_{{{t}}}[({\cal{L}}{V})(X^{{\gamma}_{{t}},{u}}_\sigma,u)-
                      ({\cal{L}}{V})({\gamma}_{{{t}}},u)]d\sigma\bigg{]}\nonumber\\
                      &&+\frac{1}{\delta}\mathbb{E}\bigg{[}\int^{{t}+\delta}_{{{t}}}(\Gamma^{{t}}(\sigma)-1)({\cal{L}}{V})(X^{{\gamma}_{{t}},{u}}_\sigma,
                                u)d\sigma\bigg{]}.
\end{eqnarray*}
By the proving process of the above theorem, letting $\delta\rightarrow0$, we have
\begin{eqnarray*}
                       0\geq({\cal{L}}V)(\gamma_t,u).
\end{eqnarray*}
     Taking the supremum over $u\in U$, %we get that
\begin{eqnarray}\label{idd}
    0\geq  \partial_tV(\gamma_t)+(A^*\partial_{x}V(\gamma_t),\gamma_t(t))_H+{\mathbf{H}}(\gamma_t,V(\gamma_t),\partial_{x}V(\gamma_t),\partial_{xx}V(\gamma_t)).
\end{eqnarray}
          On the other hand,  let $(t,\gamma_t)\in [0,T)\times \Lambda$ be fixed. Then, by (\ref{ddpG}),  there exists an
                         $\tilde{u}(\cdot)\equiv u^{\varepsilon,\delta}(\cdot)\in {\cal{U}}[t,t+\delta]$ for any $\varepsilon>0$ and $0\leq\delta\leq T-t$ such that
\begin{eqnarray}\label{jiajiajia}
    &&-\varepsilon\delta \leq G^{\gamma_t,{\tilde{u}}}_{t,t+\delta}[V(X^{{\gamma}_{{t}},{\tilde{u}}}_{{{t}}+\delta})]
                           -V(\gamma_t).%\\
%    &\geq&E\bigg{[}\int_{t}^{s}q(\sigma,X^{\tilde{u}}(\sigma,t,x),\tilde{u}(\sigma))d\sigma+V(s,X^{\tilde{u}}_s(t,x))-V(t,x)\bigg{]}\\
%                     &=& V_t(t,x)(s-t)+{\mathcal
%                        {S}}(V)(t,x)(s-t)+\int^{s}_{t}
%                        E({\nabla_{x_0}V(t,x)}, \widehat{F}(t,x,{\tilde{u}}(\sigma)))_{R^d} d\sigma\\
%                        &&+\frac{1}{2}\int^{s}_{t}E\mbox{tr}[\nabla_{x_0}^2V(t,x)
%                         \widehat{G}(t,x,{\tilde{u}}(\sigma))\widehat{G}^\top(t,x,{\tilde{u}}(\sigma))]
%                         d\sigma+\int_{t}^{s}Eq(t,x(0),\tilde{u}(\sigma))d\sigma+o(|s-t|)\\
%                     &\geq& V_t(t,x)(s-t)+{\mathcal
%                        {S}}(V)(t,x)(s-t)+{\mathbf{H}}(t,x,\nabla_{x_0}V(t,x),\nabla^2_{x_0}V(t,x))(s-t)+o(|s-t|).
\end{eqnarray}
Thus we have
\begin{eqnarray*}
               V(\gamma_t) -\varepsilon\delta \leq Y^{{\gamma}_{{t}},{\tilde{u}}}(t)= \int_{t}^{t+\delta}q(X^{{\gamma}_{{t}},{\tilde{u}}}_s,Y^{{\gamma}_{{t}},{\tilde{u}}}(s),Z^{{\gamma}_{{t}},
               {\tilde{u}}}(s),{\tilde{u}}(s))ds
                                +V(X^{{\gamma}_{{t}},{\tilde{u}}}_{{{t}}+\delta})-\int_{t}^{t+\delta}Z^{{\gamma}_{{t}},
               {\tilde{u}}}(s)dW(s).
\end{eqnarray*}
              Applying functional It\^{o} formula (\ref{statesop0}) to
                         $V(X^{{\gamma}_{{t}},{\tilde{u}}}_{{{t}}+\delta})$,   the  inequality above implies that, for all $ 0\leq\delta\leq T-t$,
\begin{eqnarray*}
    -\varepsilon\delta&\leq&\int_{t}^{t+\delta}q(X^{{\gamma}_{{t}},{\tilde{u}}}_s,Y^{{\gamma}_{{t}},{\tilde{u}}}(s),Z^{{\gamma}_{{t}},
               {\tilde{u}}}(s),{\tilde{u}}(s))ds+\int^{t+\delta}_{t} ({\cal{L}}{V})(X_s^{{\gamma}_{{t}},{\tilde{u}}},
                                      {\tilde{u}}(s))ds\\
               && -\int_{t}^{t+\delta}q(X^{{\gamma}_{{t}},\tilde{u}}_s,V(X^{{\gamma}_{{t}},\tilde{u}}_s),
               {\partial_xV(X^{{\gamma}_{{t}},\tilde{u}}_s)}G(X^{{\gamma}_{{t}},\tilde{u}}_s,\tilde{u}(s)),\tilde{u}(s))ds\\
               &&+\int_{t}^{t+\delta}({\partial_xV(X^{{\gamma}_{{t}},\tilde{u}}_s)}G(X^{{\gamma}_{{t}},\tilde{u}}_s,\tilde{u}(s))-Z^{{\gamma}_{{t}},
               {\tilde{u}}}(s))dW(s), \\
               &=&\int^{t+\delta}_{t} [({\cal{L}}{V})(X_s^{{\gamma}_{{t}},\tilde{{u}}},
                                      \tilde{u}(s))ds -B_1(s)\tilde{Y}^{2}(s)-(B_2(s),\tilde{Z}^{2}(s))_{\Xi}]ds
                                +\int_{t}^{t+\delta}\tilde{Z}^{2}(s)dW(s),
\end{eqnarray*}
where
\begin{eqnarray*}
                          &&\tilde{Y}^{2}(s):=
                          V(X_s^{{\gamma}_{{t}},\tilde{{u}}})-Y^{{\gamma}_{{t}},\tilde{u}}(s),  \ \ s\in[{t},{t}+\delta],\\
                            &&\tilde{Z}^{2}(s)  :={\partial_xV(X^{{\gamma}_{{t}},{\tilde{u}}}_s)}G(X^{{\gamma}_{{t}},{\tilde{u}}}_s,u)-Z^{{\gamma}_{{t}},
               {\tilde{u}}}(s), \ \ s\in[{t},{t}+\delta],
\end{eqnarray*}
and  $|B_1|\vee|B_2|\leq L$. %($C$ only depends on  Lipschitz constant of $q$).
Applying It\^o  formula  (see also Proposition 2.2 in \cite{el}), we obtain
\begin{eqnarray}\label{4.1422}
                      \tilde{Y}^{2}({t})
                      =\mathbb{E}\bigg{[}\tilde{Y}^{2}(t+\delta)\tilde{\Gamma}^{{t}}({t}+\delta)
                      -
                      \int^{{t}+\delta}_{{{t}}}\tilde{\Gamma}^{{t}}(\sigma)({\cal{L}}V)
                       (X^{{\gamma}_{{t}},{\tilde{u}}}_\sigma,\tilde{u}(\sigma))d\sigma\bigg{|}
                      {\cal{F}}_{{t}}\bigg{]},\ \ \
\end{eqnarray}
where $\tilde{\Gamma}^{{t}}(\cdot)$ solves the linear SDE
$$
               d\tilde{\Gamma}^{{t}}(s)=\tilde{\Gamma}^{{t}}(s)(B_1(s)ds+B_2(s)dW(s)),\ s\in [{{t}},{{t}}+\delta];\ \ \ \tilde{\Gamma}^{{t}}({{t}})=1.
$$
Obviously, $\tilde{\Gamma}^{{t}}\geq 0$.
Combining (\ref{jiajiajia}) and (\ref{4.1422}), we get that
%By the  proving process of the above theorem, we obtain that
\begin{eqnarray*}
-\varepsilon%+\frac{1}{\delta}[V(\hat{\gamma}^{\mu}_{{t}})-{{\varphi}} (\hat{\gamma}^{\mu}_{{t}})]\nonumber\\
    &\leq&
      \partial_tV(\gamma_t)+(A^*\partial_{x}V(\gamma_t),\gamma_t(t))_H+{\mathbf{H}}(\gamma_t,V(\gamma_t),\partial_{x}V(\gamma_t),\partial_{xx}V(\gamma_t))\\
      &&+\frac{1}{\delta}E\bigg{[}\int^{{t}+\delta}_{{{t}}}
      [({\cal{L}}{V})(X^{{\gamma}_{{t}},{\tilde{u}}}_\sigma,{\tilde{u}}(\sigma))-
                      ({\cal{L}}{V})({\gamma}_{{{t}}},{\tilde{u}}(\sigma))]d\sigma\bigg{]}\nonumber\\
                      &&+\frac{1}{\delta}E\bigg{[}\int^{{t}+\delta}_{{{t}}}(\tilde{\Gamma}^{{t}}(\sigma)-1)({\cal{L}}{V})(X^{{\gamma}_{{t}},{\tilde{u}}}_\sigma,
                                {\tilde{u}}(\sigma))d\sigma\bigg{]}.
\end{eqnarray*}
By the proving process of the above theorem, letting $\delta\rightarrow0$, we obtain
\begin{eqnarray*}
                       -\varepsilon\leq \partial_tV(\gamma_t)+(A^*\partial_{x}V(\gamma_t),\gamma_t(t))_H+{\mathbf{H}}(\gamma_t,V(\gamma_t),\partial_{x}V(\gamma_t),\partial_{xx}V(\gamma_t)).
\end{eqnarray*}
               The desired result  is obtained by combining the inequality given above with (\ref{idd}). \ \ $\Box$
                                \par
We conclude this section with   the stability of viscosity solutions.
\begin{theorem}\label{theoremstability}
                      Let $F,G,q,\phi$ satisfy Hypotheses \ref{hypstate} and  \ref{hypcost}, and $v\in C^0(\Lambda)$. Assume
                       \item{(i)}      for any $\varepsilon>0$, there exist $F^\varepsilon, G^\varepsilon, q^\varepsilon, \phi^\varepsilon$ and $v^\varepsilon\in C^0(\Lambda)$ such that  $F^\varepsilon, G^\varepsilon, q^\varepsilon, \phi^\varepsilon$ satisfy  Hypotheses \ref{hypstate} and  \ref{hypcost} and $v^\varepsilon$ is a viscosity subsolution (resp., supsolution) of PHJB equation (\ref{hjb1}) with generators $F^\varepsilon, G^\varepsilon, q^\varepsilon, \phi^\varepsilon$;
                           \item{(ii)} as $\varepsilon\rightarrow0$, $(F^\varepsilon, G^\varepsilon, q^\varepsilon, \phi^\varepsilon,v^\varepsilon)$ converge to
                           $(F, G, q, \phi, v)$  uniformly in the following sense: %for any $(t,\gamma_t,x,y)\in [0,T]\times \Lambda \times R\times R^d$, there exists $\delta$ such that
\begin{eqnarray}\label{sss}
                          \lim_{\varepsilon\rightarrow0}&&\sup_{(t,\gamma_t,x,y,u)\in [0,T]\times \Lambda \times \mathbb{R}\times H\times U}\sup_{\eta_T\in  \Lambda_T}[(|F^\varepsilon-F|+|G^\varepsilon-G|)({\gamma}_t,u)
                         \nonumber\\
                         &&~~~~~~~~~~~~+|q^\varepsilon-q|({\gamma}_t,{x},yG({\gamma}_{{t}},u),u)+|\phi^\varepsilon-\phi|(\eta_T)+|v^\varepsilon-v|({\gamma}_t)]=0.
\end{eqnarray}
                   Then $v$ is a viscosity subsoluiton (resp., supersolution) of PHJB equation (\ref{hjb1}) with generators $F,G,q,\phi$.
\end{theorem}
{\bf  Proof }. \ \ Without loss of generality, we shall only prove the viscosity subsolution property.
First,  from $v^{\varepsilon}$ is a viscosity subsolution of  equation (\ref{hjb1}) with generators $F^{\varepsilon}, G^{\varepsilon},  q^{\varepsilon}, \phi^{\varepsilon}$, it follows that
 $$
                            v^{\varepsilon}(\gamma_T)\leq \phi^{\varepsilon}(\gamma_T),\ \ \gamma_T\in \Lambda_T.
$$
Letting $\varepsilon\rightarrow0$, we  have
 $$
                            v(\gamma_T)\leq \phi(\gamma_T),\ \ \gamma_T\in \Lambda_T.
$$
Next,     Let  $\varphi\in \Phi$ and $g\in {\cal{G}}_{\hat{t}}$ with $\hat{t}\in [0,T)$
 such that
$$
                         0=(V-\varphi-g)(\hat{\gamma}_{\hat{t}})=\sup_{(s,\eta_s)\in[\hat{t},T]\times \Lambda^{\hat{t}}}
                         (V- \varphi-g)(\eta_s),
$$
 where $\hat{\gamma}_{\hat{t}}\in \Lambda$. % we let   $\varphi\in J^+(\hat{\gamma}_{\hat{t}}, v)$ with
%  $(\hat{t},\hat{\gamma}_{\hat{t}})\in [0,T)\times\Lambda$.
 Denote $g_{1}(\gamma_t):=g(\gamma_t)+\overline{\Upsilon}^3(\gamma_t,\hat{\gamma}_{{\hat{t}}})$
 %+|\gamma_t(t)-\hat{\gamma}_{{\hat{t}}}(\hat{t})|^2$
 for all
 $(t,\gamma_t)\in [\hat{t},T]\times\Lambda$. %By Lemma \ref{theoremS},
 Then we have  $g_{1}\in {\cal{G}}_{\hat{t}}$.  Define a sequence of positive numbers $\{\delta_i\}_{i\geq0}$  by %$\delta_0=\beta$ and
        $\delta_i=\frac{1}{2^i}$ for all $i\geq0$. For every $\varepsilon>0$, since $v^{\varepsilon}-\varphi-g_1$ is a  upper semicontinuous functional  and $\overline{\Upsilon}^3(\cdot,\cdot)$ is a gauge-type function, from Lemma \ref{theoremleft} it follows that,
  for every  $(t_0,\gamma^0_{t_0})\in [\hat{t},T]\times \Lambda^{\hat{t}}$ satisfy
 $$
(v^{\varepsilon}-\varphi-g_1)(\gamma^0_{t_0})\geq \sup_{(s,\gamma_s)\in [\hat{t},T]\times \Lambda^{\hat{t}}}(v^{\varepsilon}-\varphi-g_1)(\gamma_s)-\varepsilon,\
\    \mbox{and} \ \ (v^{\varepsilon}-\varphi-g_1)(\gamma^0_{t_0})\geq (v^{\varepsilon}-\varphi-g_1)(\hat{\gamma}_{\hat{t}}),
 $$
  there exist $(t_{\varepsilon},{\gamma}^{\varepsilon}_{t_{\varepsilon}})\in [\hat{t},T]\times \Lambda^{\hat{t}}$ and a sequence $\{(t_i,\gamma^i_{t_i})\}_{i\geq1}\subset [\hat{t},T]\times \Lambda^{\hat{t}}$ such that
  \begin{description}
        \item{(i)} $\overline{\Upsilon}^3(\gamma^0_{t_0},{\gamma}^{\varepsilon}_{t_{\varepsilon}})\leq {\varepsilon}$,  $\overline{\Upsilon}^3(\gamma^i_{t_i},{\gamma}^{\varepsilon}_{t_{\varepsilon}})\leq \frac{\varepsilon}{2^i}$ and $t_i\uparrow t_{\varepsilon}$ as $i\rightarrow\infty$,
        \item{(ii)}  $(v^{\varepsilon}-\varphi-g_1)({\gamma}^{\varepsilon}_{t_{\varepsilon}})-\sum_{i=0}^{\infty}\frac{1}{2^i}\overline{\Upsilon}^3(\gamma^i_{t_i},{\gamma}^{\varepsilon}_{t_{\varepsilon}})\geq (v^{\varepsilon}-\varphi-g_1)(\gamma^0_{t_0})$, and
        \item{(iii)}  $(v^{\varepsilon}-\varphi-g_1)(\gamma_s)-\sum_{i=0}^{\infty}\frac{1}{2^i}\overline{\Upsilon}^3(\gamma^i_{t_i},\gamma_s)
            <(v^{\varepsilon}-\varphi-g_1)({\gamma}^{\varepsilon}_{t_{\varepsilon}})-\sum_{i=0}^{\infty}\frac{1}{2^i}\overline{\Upsilon}^3(\gamma^i_{t_i},{\gamma}^{\varepsilon}_{t_{\varepsilon}})$ for all $(s,\gamma_s)\in [t_{\varepsilon},T]\times \Lambda^{t_{\varepsilon}}\setminus \{(t_{\varepsilon},{\gamma}^{\varepsilon}_{t_{\varepsilon}})\}$.

        \end{description}
  We claim that
\begin{eqnarray}\label{gamma}
d_\infty({\gamma}^{\varepsilon}_{{t}_{\varepsilon}},\hat{\gamma}_{\hat{t}})\rightarrow0  \ \ \mbox{as} \ \ \varepsilon\rightarrow0.
\end{eqnarray}
 Indeed, if not,  by (\ref{s0}), we can assume that there exists an $\nu_0>0$
 such
                    that
$$
                 %   |t_\varepsilon-\hat{t}|^2
%  +
  \overline{\Upsilon}^3({\gamma}^{\varepsilon}_{{t}_{\varepsilon}},\hat{\gamma}_{{\hat{t}}})%+|{\gamma}^{\varepsilon}_{{t}_{\varepsilon}}({t}_{\varepsilon})-\hat{\gamma}_{{\hat{t}}}(\hat{t})|_0^8
  \geq\nu_0.
$$
 Thus,  we obtain that
%Now we have  $({\bar{t}},\bar{\gamma}_{\bar{t}})=(\hat{t},\hat{\gamma}_{\hat{t}})$, for if not,
\begin{eqnarray*}
   &&0=(v- {{\varphi}}-g)(\hat{\gamma}_{\hat{t}})= \lim_{\varepsilon\rightarrow0}(v^\varepsilon-\varphi-g_1)(\hat{\gamma}_{\hat{t}})
   \leq \limsup_{\varepsilon\rightarrow0}\bigg{[}(v^{\varepsilon}-\varphi-g_1)({\gamma}^{\varepsilon}_{t_{\varepsilon}})-\sum_{i=0}^{\infty}\frac{1}{2^i}\overline{\Upsilon}^3(\gamma^i_{t_i},{\gamma}^{\varepsilon}_{t_{\varepsilon}})\bigg{]}\\
   %&=&\overline{\lim_{\varepsilon\rightarrow0}}\bigg{[}(v^\varepsilon-{{\varphi}}-g)({\gamma}^{\varepsilon}_{{t}_{\varepsilon}})%-|t_\varepsilon-\hat{t}|^2
%  -\overline{\Upsilon}^2({\gamma}^{\varepsilon}_{{t}_{\varepsilon}},\hat{\gamma}_{{\hat{t}}})-%|{\gamma}^{\varepsilon}_{{t}_{\varepsilon}}({t}_{\varepsilon})-\hat{\gamma}_{{\hat{t}}}(\hat{t})|_8^8
%  \sum_{i=0}^{\infty}\frac{1}{2^i}\overline{\Upsilon}^2(\gamma^i_{t_i},{\gamma}^{\varepsilon}_{t_{\varepsilon}})\bigg{]}\\
   &\leq&\limsup_{\varepsilon\rightarrow0}{[}(v-{{\varphi}}-g)({\gamma}^{\varepsilon}_{{t}_{\varepsilon}})+(v^\varepsilon-v)({\gamma}^{\varepsilon}_{{t}_{\varepsilon}})
     %-\sum_{i=0}^{\infty}\frac{1}{2^i}\overline{\Upsilon}^2(\gamma^i_{t_i},{\gamma}^{\varepsilon}_{t_{\varepsilon}})
     {]}-\nu_0\leq (v- {{\varphi}}-g)(\hat{\gamma}_{\hat{t}})-\nu_0=-\nu_0,
\end{eqnarray*}
 contradicting $\nu_0>0$.  We notice that, by  the property (i) of $(t_{\varepsilon},{\gamma}^{\varepsilon}_{t_{\varepsilon}})$,
  \begin{eqnarray*}
  2\sum_{i=0}^{\infty}\frac{1}{2^i}({t_{\varepsilon}}-{t}_{i})
  \leq2\sum_{i=0}^{\infty}\frac{1}{2^i}\bigg{(}\frac{\varepsilon}{2^i}\bigg{)}^{\frac{1}{2}}\leq 4\varepsilon^{\frac{1}{2}};
    \end{eqnarray*}
    \begin{eqnarray*}
    &&|\partial_x{\Upsilon}^3({\gamma}^{\varepsilon}_{{t}_{\varepsilon}}-\hat{\gamma}_{{\hat{t}},{t}_{\varepsilon},A})|\leq 18|e^{({t}_{\varepsilon}-\hat{t})A}\hat{\gamma}_{{\hat{t}}}(\hat{t})-{\gamma}^{\varepsilon}_{{t}_{\varepsilon}}({t}_{\varepsilon})|^5;\\
    &&|\partial_{xx}{\Upsilon}^3({\gamma}^{\varepsilon}_{{t}_{\varepsilon}}-\hat{\gamma}_{{\hat{t}},{t}_{\varepsilon},A})|\leq 306|e^{({t}_{\varepsilon}-\hat{t})A}\hat{\gamma}_{{\hat{t}}}(\hat{t})-{\gamma}^{\varepsilon}_{{t}_{\varepsilon}}({t}_{\varepsilon})|^4;
    \end{eqnarray*}
    \begin{eqnarray*}
  \bigg{|}\partial_x\sum_{i=0}^{\infty}\frac{1}{2^i}
                      \Upsilon^3({\gamma}^{\varepsilon}_{t_{\varepsilon}}-\gamma^i_{t_i,t_{\varepsilon},A})
                      \bigg{|}
                      \leq18\sum_{i=0}^{\infty}\frac{1}{2^i}|e^{(t_{\varepsilon}-t_i)A}\gamma^i_{t_i}({t}_{i})-{\gamma}^{\varepsilon}_{t_{\varepsilon}}(t_{\varepsilon})|^5
                     \leq18\sum_{i=0}^{\infty}\frac{1}{2^i}\bigg{(}\frac{\varepsilon}{2^i}\bigg{)}^{\frac{5}{6}}
                      \leq36{\varepsilon}^{\frac{5}{6}}.
                        \end{eqnarray*}
    and
     \begin{eqnarray*}
  \bigg{|}\partial_{xx}\sum_{i=0}^{\infty}\frac{1}{2^i}
                      \Upsilon^3({\gamma}^{\varepsilon}_{t_{\varepsilon}}-\gamma^i_{t_i,t_{\varepsilon},A})
                      \bigg{|}
                      \leq306\sum_{i=0}^{\infty}\frac{1}{2^i}|e^{(t_{\varepsilon}-t_i)A}\gamma^i_{t_i}({t}_{i})-{\gamma}^{\varepsilon}_{t_{\varepsilon}}(t_{\varepsilon})|^4
                     \leq306\sum_{i=0}^{\infty}\frac{1}{2^i}\bigg{(}\frac{\varepsilon}{2^i}\bigg{)}^{\frac{2}{3}}
                      \leq612{\varepsilon}^{\frac{2}{3}}.
                        \end{eqnarray*}
  Then for any $\varrho>0$, by (\ref{sss}) and (\ref{gamma}), there exists $\varepsilon>0$ small enough such that
$$
            \hat{t}\leq {t}_{\varepsilon}< T,  \        \  %|{\gamma}^{\varepsilon}_{{t}_{\varepsilon}}({t}_{\varepsilon})|< \frac{M_0}{2},\ \
             2|{t}_{\varepsilon}-\hat{t}|+2\sum_{i=0}^{\infty}\frac{1}{2^i}({t_{\varepsilon}}-{t}_{i})+|\partial_t{\varphi}({\gamma}^{\varepsilon}_{{t}_{\varepsilon}})-\partial_t{\varphi}(\hat{\gamma}_{\hat{t}})|
             +|\partial_tg({\gamma}^{\varepsilon}_{{t}_{\varepsilon}})-\partial_tg(\hat{\gamma}_{\hat{t}})|\leq \frac{\varrho}{3},$$
$$
|(A^*\partial_x\varphi({\gamma}^{\varepsilon}_{{t}_{\varepsilon}}),{\gamma}^{\varepsilon}_{{t}_{\varepsilon}}({t}_{\varepsilon}))_H-(A^*\partial_x{\varphi}(\hat{\gamma}_{\hat{t}}),\hat{\gamma}_{\hat{t}}(\hat{t}))_H| \leq \frac{\varrho}{3}, \ \mbox{and}\  |I|\leq \frac{\varrho}{3},
$$
where
\begin{eqnarray*}
I&=&{\mathbf{H}}^{\varepsilon}({\gamma}^{\varepsilon}_{{t}_{\varepsilon}},
                           \partial_x\varphi({\gamma}^{\varepsilon}_{{t}_{\varepsilon}})+ \partial_xg_2({\gamma}^{\varepsilon}_{{t}_{\varepsilon}}),\partial_{xx}\varphi({\gamma}^{\varepsilon}_{{t}_{\varepsilon}})+ \partial_{xx}g_2({\gamma}^{\varepsilon}_{{t}_{\varepsilon}}))\\
                       &&-{\mathbf{H}}(\hat{\gamma}_{\hat{t}},\partial_x{\varphi}(\hat{\gamma}_{\hat{t}})+\partial_xg(\hat{\gamma}_{\hat{t}}),\partial_{xx}{\varphi}(\hat{\gamma}_{\hat{t}})+\partial_{xx}g(\hat{\gamma}_{\hat{t}})),\\
g_2({\gamma}^{\varepsilon}_{{t}_{\varepsilon}})&=&g({\gamma}^{\varepsilon}_{{t}_{\varepsilon}})+\overline{\Upsilon}^3({\gamma}^{\varepsilon}_{{t}_{\varepsilon}}-\hat{\gamma}_{{\hat{t}},{t}_{\varepsilon},A})
+\sum_{i=0}^{\infty}\frac{1}{2^i}\overline{\Upsilon}^3({\gamma}^{\varepsilon}_{t_{\varepsilon}}-\gamma^i_{t_i,t_{\varepsilon},A}),
%\\
%\partial_xg_2({\gamma}^{\varepsilon}_{{t}_{\varepsilon}})&=&\partial_xg({\gamma}^{\varepsilon}_{{t}_{\varepsilon}})+\partial_x{\Upsilon}^3({\gamma}^{\varepsilon}_{{t}_{\varepsilon}}-\hat{\gamma}_{{\hat{t}},{t}_{\varepsilon},A})
%+\sum_{i=0}^{\infty}\frac{1}{2^i}\partial_x{\Upsilon}^3({\gamma}^{\varepsilon}_{t_{\varepsilon}}-\gamma^i_{t_i,t_{\varepsilon},A}),\\
%\partial_tg_2({\gamma}^{\varepsilon}_{{t}_{\varepsilon}})&=&\partial_tg({\gamma}^{\varepsilon}_{{t}_{\varepsilon}})+2({t}_{\varepsilon}-\hat{t})
%+2\sum_{i=0}^{\infty}\frac{1}{2^i}(t_{\varepsilon}-t_i),
\end{eqnarray*}
and
\begin{eqnarray*}
                                {\mathbf{H}}^{\varepsilon}(\gamma_t,r,p,l)&=&\sup_{u\in{
                                         {U}}}[
                        (p,F^{\varepsilon}(\gamma_t,u))_{H}+\frac{1}{2}\mbox{Tr}[ l G^{\varepsilon}(\gamma_t,u){G^{\varepsilon}}^*(\gamma_t,u)]\\
                        &&\ \ \ \ \ +q^{\varepsilon}(\gamma_t,r,p{G^{\varepsilon}}(\gamma_t,u),u)],  \ \ (t,\gamma_t,r,p,l)\in [0,T]\times {\Lambda}\times \mathbb{R}\times H\times {\cal{S}}(H).
\end{eqnarray*}
%$$
%                                {\mathbf{H}}^{\varepsilon}(\gamma_t,p)=\inf_{u\in{
%                                         {U}}}[
%                        (p,F^{\varepsilon}(\gamma_t,u))_{R^d}
%                        +q^{\varepsilon}(\gamma_t,u)],  \ \ (\gamma_t,p)\in {\Lambda}\times R^d.
%$$
 Since $v^{\varepsilon}$ is a viscosity subsolution of PHJB equation (\ref{hjb1}) with generators $F^{\varepsilon}, G^{\varepsilon},   q^{\varepsilon}, \phi^{\varepsilon}$, we have
$$
                          \partial_t\varphi({\gamma}^{\varepsilon}_{{t}_{\varepsilon}})+\partial_tg_2({\gamma}^{\varepsilon}_{{t}_{\varepsilon}})
                           +(A^*\partial_x\varphi({\gamma}^{\varepsilon}_{{t}_{\varepsilon}}),{\gamma}^{\varepsilon}_{{t}_{\varepsilon}}({t}_{\varepsilon}))_H+{\mathbf{H}}^{\varepsilon}({\gamma}^{\varepsilon}_{{t}_{\varepsilon}},
                           \partial_x\varphi({\gamma}^{\varepsilon}_{{t}_{\varepsilon}})+\partial_xg_2({\gamma}^{\varepsilon}_{{t}_{\varepsilon}}),
                           \partial_{xx}\varphi({\gamma}^{\varepsilon}_{{t}_{\varepsilon}})+\partial_{xx}g_2({\gamma}^{\varepsilon}_{{t}_{\varepsilon}}))\geq0.
$$
Thus
\begin{eqnarray*}
                       0&\leq&  \partial_t{\varphi}({\gamma}^{\varepsilon}_{{t}_{\varepsilon}})+\partial_tg({\gamma}^{\varepsilon}_{{t}_{\varepsilon}})
                       +2({t}_{\varepsilon}-\hat{t})+2\sum_{i=0}^{\infty}\frac{1}{2^i}({t_{\varepsilon}}-{t}_{i})
                       +(A^*\partial_x\varphi({\gamma}^{\varepsilon}_{{t}_{\varepsilon}}),{\gamma}^{\varepsilon}_{{t}_{\varepsilon}}({t}_{\varepsilon}))_H\\
                       &&+{\mathbf{H}}(\hat{\gamma}_{\hat{t}},\partial_x{\varphi}(\hat{\gamma}_{\hat{t}})+\partial_xg(\hat{\gamma}_{\hat{t}}),\partial_{xx}{\varphi}(\hat{\gamma}_{\hat{t}})+\partial_{xx}g(\hat{\gamma}_{\hat{t}}))+I\\
                       &\leq&\partial_t{\varphi}(\hat{\gamma}_{\hat{t}})+\partial_tg(\hat{\gamma}_{\hat{t}})+(A^*\partial_x{\varphi}(\hat{\gamma}_{\hat{t}}),\hat{\gamma}_{\hat{t}}(\hat{t}))_H
                       +{\mathbf{H}}(\hat{\gamma}_{\hat{t}},\partial_x{\varphi}(\hat{\gamma}_{\hat{t}})+\partial_xg(\hat{\gamma}_{\hat{t}}),\partial_{xx}{\varphi}(\hat{\gamma}_{\hat{t}})+\partial_{xx}g(\hat{\gamma}_{\hat{t}}))+\varrho.
\end{eqnarray*}
Letting $\varrho\downarrow 0$, we show that
$$
\partial_t{\varphi}(\hat{\gamma}_{\hat{t}})+\partial_tg(\hat{\gamma}_{\hat{t}})+(A^*\partial_x{\varphi}(\hat{\gamma}_{\hat{t}}),\hat{\gamma}_{\hat{t}}(\hat{t}))_H
+{\mathbf{H}}(\hat{\gamma}_{\hat{t}},\partial_x{\varphi}(\hat{\gamma}_{\hat{t}})+\partial_xg(\hat{\gamma}_{\hat{t}}),\partial_{xx}{\varphi}(\hat{\gamma}_{\hat{t}})+\partial_{xx}g(\hat{\gamma}_{\hat{t}}))\geq0.
$$
Since ${\varphi}\in \Phi$ and $g\in {\cal{G}}_t$ with $t\in [0,T)$  are arbitrary, we see that $v$ is a viscosity subsolution of PHJB equation (\ref{hjb1}) with generators $F,G,q,\phi$, and thus completes the proof.
\ \ $\Box$

\section{Viscosity solutions to  PHJB equations: Uniqueness theorem.}
\par
             This section is devoted to a  proof of uniqueness of  viscosity
                   solutions to (\ref{hjb1}). This result, together with
                  the results from  the previous section, will be used to characterize
                   the value functional defined by (\ref{value1}).
                   \par
                    Let$\{e_i\}_{i\geq1}$ be an orthonormal basis of $H$ such that $e_i\in {\cal{D}}(A^*)$ for all $i\geq 1$. For every $N\geq 1$, we let $H_N$ denote the vector space generated by $e_1,\ldots,e_N$, $P_N$  denote the orthogonal projection onto $H_N$. Define $Q_N=I-P_N$, we then have
an orthogonal decomposition $H=H_N+H^\bot_N$, where $H^\bot_N=Q_NH$. We will denote by $x_N,y_N, \ldots$ points in $H_N$ and by $x^-_N,y^-_N,\ldots$ points in $H^\bot_N$, and write $x=(x_N,x^-_N)$ for $x\in H$.  Our uniqueness result requires the following assumption of $G$.
                   \begin{hyp}\label{hypstate5666}
                   %Assume $G$ is a measurable functional from $[0,T]\times H\times U$ to $L_2(\Xi,H)$ and there exists a constant $L>0$ such that
%                   $|G(t,x,u)|_{L_2(\Xi,H)}\leq L(1+|x|_{L_2(\Xi,H)})$.
%                   $G:U\rightarrow L_2(\Xi,H)$.
                   For every $(t,\gamma_t)\in [0,T)\times \Lambda$,
                   \begin{eqnarray}\label{g5}
                   \sup_{u\in U}|Q_N{G}(\gamma_t,u)|_{L_2(\Xi;H)}^2\rightarrow0\ \mbox{as}\ N\rightarrow\infty.
                   \end{eqnarray}
\end{hyp}
                   \par
We assume without loss of generality that, there exists a constant $K>0$ such that,
for all $(t,\gamma_t, p,l)\in [0,T]\times{\Lambda}\times H\times {\mathcal{S}}(H)$ and $r_1,r_2\in \mathbb{R}$ such that $r_1<r_2$,
\begin{eqnarray}\label{5.1}
{\mathbf{H}}(\gamma_t,r_1,p,l)-{\mathbf{H}}(\gamma_t,r_2,p,l)\geq K(r_2-r_1).
\end{eqnarray}
We  now state the main result of this section.
\begin{theorem}\label{theoremhjbm}
Suppose Hypotheses \ref{hypstate}, \ref{hypcost} and \ref{hypstate5666}  hold.
                         Let $W_1\in C^0({\Lambda})$ $(\mbox{resp}., W_2\in C^0({\Lambda}))$ be  a viscosity subsolution (resp., supsolution) to equation (\ref{hjb1}) and  let  there exist  constant $L>0$
                          %and $m>0$, %and a  modulus of continuity $\hat{\omega}$,
                        such that, for any  $0\leq t\leq  s\leq T$,
                        $\gamma_t, \eta_s\in{\Lambda}$,
\begin{eqnarray}\label{w}
                                   |W_1(\gamma_t)|\vee |W_2(\gamma_t)|\leq L (1+||\gamma_t||_0);
                                   \end{eqnarray}
\begin{eqnarray}\label{w1}
                              |W_1(\gamma_{t,s,A})-W_1(\eta_t)|\vee|W_2(\gamma_{t,s,A})-W_2(\eta_t)|\leq
                        L(1+||\gamma_t||_0+||\eta_s||_0)|s-t|^{\frac{1}{2}}+L||\gamma_{t}-\eta_s||_0.
\end{eqnarray}
                   Then  $W_1\leq W_2$. %, where $V$ is the value functional defined in (\ref{value1}).
%
%
%
%
% Suppose Hypothesis \ref{hypstate}   holds.
%                         Let $W\in C^0({\Lambda})$ be  a viscosity solution to  (\ref{hjb1}) and  let  there exist a constant $C_1>0$
%                        such that for any  $0\leq t\leq  s\leq T$, $\gamma_t,\gamma'_s\in{\Lambda}$,
%\begin{eqnarray}\label{w}
%                                    &&|W(\gamma_t)|\leq C_1 (1+||\gamma_t||_0), \nonumber\\
%                                    && |W(\gamma_t)-W(\gamma'_s)|
%                  \leq
%                        C_1(1+||\gamma_t||_0+||\gamma'_s||_0)|s-t|^{\frac{1}{2}}+C_1(\mu(|\gamma_{t}-\gamma'_s|)+|\gamma_{t}(t)-\gamma'_s(s)|).
%\end{eqnarray}
%                   Then  $W=V$, where $V$ is the value functional defined in (\ref{value1}).
\end{theorem}
\par
                      Theorems    \ref{theoremvexist} and \ref{theoremhjbm} lead to the result (given below) that the viscosity solution to   PHJB equation given in (\ref{hjb1})
                      corresponds to the value functional  $V$ of our optimal control problem given in (\ref{state1}), (\ref{fbsde1}) and (\ref{value1}).
\begin{theorem}\label{theorem52}\ \
                 Let  Hypotheses \ref{hypstate}, \ref{hypcost} and \ref{hypstate5666}   hold.  Then the value
                          functional $V$ defined by (\ref{value1}) is the unique viscosity
                          solution to (\ref{hjb1}) in the class of functionals satisfying (\ref{w}) and (\ref{w1}).
\end{theorem}
\par
   {\bf  Proof  }. \ \   Theorem \ref{theoremvexist} shows that $V$ is a viscosity solution to equation (\ref{hjb1}).  Thus, our conclusion follows from %Lemma \ref{lemmavaluev},
   Theorems  \ref{theorem3.9} and
    \ref{theoremhjbm}.  \ \ $\Box$
%\par
%{\bf  Proof}. \ \
%               By Theorem 4.4, we know that $V$ is a viscosity solution of (3.8). Thus, our conclusion follows from Theorem 3.2 and Theorem 5.1.\ \ $\Box$

\par
  Next, we prove Theorem \ref{theoremhjbm}.   Let $W_1$ be a viscosity subsolution of PHJB equation (\ref{hjb1}). %We note that it is sufficient to show $W\leq V$ because the inverse of which can be proved in a similar way.
 We  note that for $\delta>0$, the functional
                    defined by $\tilde{W}:=W_1-\frac{\rho}{t}$ is a viscosity subsolution
                   for
 \begin{eqnarray*}
\begin{cases}
{\partial_t} \tilde{W}(\gamma_t)+(A^*\partial_x\tilde{W}(\gamma_t),\gamma_t(t))_H+{\mathbf{H}}(\gamma_t, \tilde{W}(\gamma_t), \partial_x \tilde{W}(\gamma_t),\partial_{xx} \tilde{W}(\gamma_t))
          = \frac{\rho}{t^2}, \ \  (t,\gamma_t)\in [0,T)\times {\Lambda}, \\
\tilde{W}(\gamma_T)=\phi(\gamma_T), \ \  \gamma_T\in \Lambda_T.\\
\end{cases}
\end{eqnarray*}
                As $W_1\leq W_2$ follows from $\tilde{W}\leq W_2$ in
                the limit $\rho\downarrow0$, it suffices to prove
                $W_1\leq W_2$ under the additional assumption given below:
$$
{\partial_t} {W_1}(\gamma_t)+(A^*\partial_xW_1(\gamma_t),\gamma_t(t))_H+{\mathbf{H}}(\gamma_t, {W_1}(\gamma_t), \partial_x {W_1}(\gamma_t),\partial_{xx}{W_1}(\gamma_t))
          \geq c:=\frac{\rho}{T^2}, \  (t,\gamma_t)\in [0,T)\times {\Lambda}.
$$
%$$
%                 \frac{\partial}{\partial t} {W}(t,x)+{\mathcal
%                        {S}}({W})(t,x)+ {\mathbf{H}}(t,x,\nabla_{x_0} {W}(t,x),\nabla^2_{x_0} \tilde{W}(t,x))\geq c,\ \ c:=\frac{\delta}{T^2}, \ \  t\in
%                               [0,T],\ \ x\in {\mathcal{D}}.
%$$
\par
   {\bf  Proof of Theorem \ref{theoremhjbm} } \ \   The proof of this theorem  is rather long. Thus, we split it into several
        steps.
\par
            $Step\  1.$ Definitions of auxiliary functionals.
            \par
 We only need to prove that $W_1(\gamma_t)\leq W_2(\gamma_t)$ for all $(t,\gamma_t)\in
[T-\bar{a},T)\times
       {\Lambda}$.
        Here,
        $$\bar{a}=\frac{1}{8(342L+36)L}\wedge{T}.$$
         Then, we can  repeat the same procedure for the case
        $[T-i\bar{a},T-(i-1)\bar{a})$.  Thus, we assume the converse result that $(\tilde{t},\tilde{\gamma}_{\tilde{t}})\in [T-\bar{a},T)\times
      {\Lambda}$ exists  such that
        $\tilde{m}:=W_1(\tilde{\gamma}_{\tilde{t}})-W_2(\tilde{\gamma}_{\tilde{t}})>0$. % Because $\cup_{\mu>0,M_0>0}{\cal{C}}^\alpha_{\mu,M_0}$  is dense in $(\Lambda,d_\infty)$,
%         by (\ref{w}) there exist $\tilde{t}\in [T-\bar{a},T)$ and $\tilde{\gamma}_{\tilde{t}}\in \cup_{\mu>0,M_0>0}{\cal{C}}^\alpha_{\mu,M_0}$
%         such that
%        $W(\tilde{\gamma}_{\tilde{t}})-V(\tilde{\gamma}_{\tilde{t}})>\tilde{m}$.
\par
         Consider that  $\varepsilon >0$ is  a small number such that
 $$
 W_1(\tilde{\gamma}_{\tilde{t}})-W_2(\tilde{\gamma}_{\tilde{t}})-2\varepsilon \frac{\nu T-\tilde{t}}{\nu
 T}\Upsilon^3(\tilde{\gamma}_{\tilde{t}})%-\frac{\varepsilon}{\tilde{t}-T+\bar{a}}
 >\frac{\tilde{m}}{2},
 $$
      and
\begin{eqnarray}\label{5.3}
                          \frac{\varepsilon}{4\nu T}\leq\frac{c}{2},
\end{eqnarray}
             where
$$
            \nu=1+\frac{1}{8T(342L+36)L}.
$$
 Next,  we define for any  $(t,\gamma_t,\eta_t)\in [T-\bar{a},T]\times{\Lambda}\times{\Lambda}$,
\begin{eqnarray*}
                 \Psi(\gamma_t,\eta_t)=W_1(\gamma_t)-W_2(\eta_t)-{\beta}\Upsilon^3(\gamma_{t},\eta_{t})-\beta^{\frac{1}{3}}|\gamma_{t}(t)-\eta_{t}(t)|^2%-\frac{\varepsilon}{{t}-T+\bar{a}}\\
                 -\varepsilon\frac{\nu T-t}{\nu
                 T}(\Upsilon^3(\gamma_t)+\Upsilon^3(\eta_t)).%-\frac{\varepsilon}{{t}-T+\bar{a}}.
\end{eqnarray*}
%where
%$$
%|\gamma_{t}(t)-\eta_{t}(t)|_N^6=|(\gamma_{t}(t))_N-(\eta_{t}(t))_N|^6+|(\gamma_{t}(t))^-_N-(\eta_{t}(t))^-_N|^6.
%$$
%For every $\varepsilon_1>0$,
 %Noting $\nu$ is independent of  $\beta$, by the definition of  ${\Psi}$,
% there exists an ${M}_1>0$  that is sufficiently  large   that
%           $
%           \sup_{(s,\gamma_s,\eta_s)\in [\tilde{t},T]\times \Lambda^{\tilde{t}}\times \Lambda^{\tilde{t}}}\Psi(\gamma_s,\eta_s)<M_1.
%           $
%           Then $M_1>0$ and  $(\tilde{\gamma}_{\tilde{t}},\tilde{\gamma}_{\tilde{t}})\in {\Lambda}^{\tilde{t}}\times{\Lambda}^{\tilde{t}}$ satisfy
%            $$
%           \Psi(\tilde{\gamma}_{\tilde{t}},\tilde{\gamma}_{\tilde{t}})\geq \sup_{(s,\gamma_s,\eta_s)\in [\tilde{t},T]\times \Lambda^{\tilde{t}}\times \Lambda^{\tilde{t}}}\Psi(\gamma_s,\eta_s)-M_1.
%           $$
        Define a sequence of positive numbers $\{\delta_i\}_{i\geq0}$  by %$\delta_0=\beta$ and
        $\delta_i=\frac{1}{2^i}$ for all $i\geq0$.
           Since $\Psi$ is a  upper semicontinuous function bounded from above and $\overline{\Upsilon}^3$ is a gauge-type function, from Lemma \ref{theoremleft} it follows that,
  for every  $(t_0,\gamma^0_{t_0},\eta^0_{t_0})\in [\tilde{t},T]\times \Lambda^{\tilde{t}}\times \Lambda^{\tilde{t}}$ satisfying
$$
\Psi(\gamma^0_{t_0},\eta^0_{t_0})\geq \sup_{(s,\gamma_s,\eta_s)\in [\tilde{t},T]\times \Lambda^{\tilde{t}}\times \Lambda^{\tilde{t}}}\Psi(\gamma_s,\eta_s)-\frac{1}{\beta},\
\    \mbox{and} \ \ \Psi(\gamma^0_{t_0},\eta^0_{t_0})\geq \Psi(\tilde{\gamma}_{\tilde{t}},\tilde{\gamma}_{\tilde{t}}) >\frac{\tilde{m}}{2},
 $$
  there exist $(\hat{t},\hat{\gamma}_{\hat{t}},\hat{\eta}_{\hat{t}})\in [\tilde{t},T]\times \Lambda^{\tilde{t}}\times \Lambda^{\tilde{t}}$ and a sequence $\{(t_i,\gamma^i_{t_i},\eta^i_{t_i})\}_{i\geq1}\subset
  [\tilde{t},T]\times \Lambda^{\tilde{t}}\times \Lambda^{\tilde{t}}$ such that
  \begin{description}
        \item{(i)} $\Upsilon^3(\gamma^0_{t_0},\hat{\gamma}_{\hat{t}})+\Upsilon^3(\eta^0_{t_0},\hat{\eta}_{\hat{t}})+|\hat{t}-t_0|^2\leq \frac{1}{\beta}$,
         $\Upsilon^3(\gamma^i_{t_i},\hat{\gamma}_{\hat{t}})+\Upsilon^3(\eta^i_{t_i},\hat{\eta}_{\hat{t}})+|\hat{t}-t_i|^2
         \leq \frac{1}{\beta2^i}$ and $t_i\uparrow \hat{t}$ as $i\rightarrow\infty$,
        \item{(ii)}  $\Psi(\hat{\gamma}_{\hat{t}},\hat{\eta}_{\hat{t}})%-\beta[\Upsilon^3(\tilde{\gamma}_{\tilde{t}},\hat{\gamma}_{\hat{t}})+\Upsilon^3(\tilde{\gamma}_{\tilde{t}},\hat{\eta}_{\hat{t}})]
            -\sum_{i=0}^{\infty}\frac{1}{2^i}[\Upsilon^3(\gamma^i_{t_i},\hat{\gamma}_{\hat{t}})
        +\Upsilon^3(\eta^i_{t_i},\hat{\eta}_{\hat{t}})+|\hat{t}-t_i|^2]\geq \Psi(\gamma^0_{t_0},\eta^0_{t_0})$, and
        \item{(iii)}    for all $(s,\gamma_s,\eta_s)\in [\hat{t},T]\times \Lambda^{\hat{t}}\times \Lambda^{\hat{t}}\setminus \{(\hat{t},\hat{\gamma}_{\hat{t}},\hat{\eta}_{\hat{t}})\}$,
        \begin{eqnarray}\label{iii4}
        &&\Psi(\gamma_s,\eta_s)%-\beta[\Upsilon^3(\tilde{\gamma}_{\tilde{t}},\gamma_s)+\Upsilon^3(\tilde{\gamma}_{\tilde{t}},\eta_s)]
        -\sum_{i=0}^{\infty}
        \frac{1}{2^i}[\Upsilon^3(\gamma^i_{t_i},\gamma_s)+\Upsilon^3(\eta^i_{t_i},\eta_s)+|{s}-t_i|^2]\nonumber\\
           &<&\Psi(\hat{\gamma}_{\hat{t}},\hat{\eta}_{\hat{t}})%-\beta[\Upsilon^3(\tilde{\gamma}_{\tilde{t}},\hat{\gamma}_{\hat{t}})
            %+\Upsilon^3(\tilde{\gamma}_{\tilde{t}},\hat{\eta}_{\hat{t}})]
            -\sum_{i=0}^{\infty}\frac{1}{2^i}[\Upsilon^3(\gamma^i_{t_i},\hat{\gamma}_{\hat{t}})
            +\Upsilon^3(\eta^i_{t_i},\hat{\eta}_{\hat{t}})+|\hat{t}-t_i|^2].
        \end{eqnarray}
        \end{description}
             We should note that the point
             $({\hat{t}},\hat{{\gamma}}_{{\hat{t}}},\hat{{\eta}}_{{\hat{t}}})$ depends on  $\beta$ and
              $\varepsilon$.
\par
$Step\ 2.$
There exists ${{M}_0}>0$
    such that
                   \begin{eqnarray}\label{5.10jiajiaaaa}||\hat{\gamma}_{\hat{t}}||_0\vee||\hat{\eta}_{\hat{t}}||_0<M_0,
                   \end{eqnarray} and
 %for every $m\in {\mathbf{N}}$,
 the following result  holds true:
 \begin{eqnarray}\label{5.10}
                      \beta^{\frac{1}{3}} ||\hat{{\gamma}}_{{\hat{t}}}-\hat{{\eta}}_{{\hat{t}}}||_{0}^2
                         +\beta^{\frac{1}{2}}|\hat{{\gamma}}_{{\hat{t}}}(\hat{t})-\hat{{\eta}}_{{\hat{t}}}(\hat{t})|^2
                          \rightarrow0 \ %\ \mbox{uniformly in}\ N \
                                   \mbox{as} \ \beta\rightarrow+\infty.
 \end{eqnarray}
  Let us show the above. First,   noting $\nu$ is independent of  $\beta$, by the definition of  ${\Psi}$,
 there exists an ${M}_0>0$  that is sufficiently  large   that
           $
           \Psi(\gamma_t, \eta_t)<0
           $ for all $t\in [T-\bar{a},T]$ and $||\gamma_t||_0\vee||\eta_t||_0\geq M_0$. Thus, we have $||\hat{\gamma}_{\hat{t}}||_0\vee||\hat{\eta}_{\hat{t}}||_0\vee
           ||{\gamma}^{0}_{t_{0}}||_0\vee||{\eta}^{0}_{t_{0}}||_0<M_0$.
           \par
   Second, by (\ref{iii4}), we have
    %$\Psi(\gamma_s,\eta_s)-\sum_{i=0}^{\infty}[\Upsilon^3(\gamma^i_{t_i},\gamma_s)+\Upsilon^3(\eta^i_{t_i},\eta_s)]
%            <\Psi(\hat{\gamma}_{\hat{t}},\hat{\eta}_{\hat{t}})-\sum_{i=0}^{\infty}[\Upsilon^3(\gamma^i_{t_i},\hat{\gamma}_{\hat{t}})
%            +\Upsilon^3(\eta^i_{t_i},\hat{\eta}_{\hat{t}})]$
   % the definition of $(\hat{t},\hat{\gamma}_{\hat{t}},\hat{\eta}_{\hat{t}})$, we have
 \begin{eqnarray}\label{5.56789}
                        &&2\Psi(\hat{\gamma}_{\hat{t}},\hat{\eta}_{\hat{t}})%-2\beta[\Upsilon^3(\tilde{\gamma}_{\tilde{t}},\hat{\gamma}_{\hat{t}})
            %+\Upsilon^3(\tilde{\gamma}_{\tilde{t}},\hat{\eta}_{\hat{t}})]
            -2\sum_{i=0}^{\infty}\frac{1}{2^i}[\Upsilon^3(\gamma^i_{t_i},\hat{\gamma}_{\hat{t}})
            +\Upsilon^3(\eta^i_{t_i},\hat{\eta}_{\hat{t}})+|\hat{t}-t_i|^2]\nonumber\\
            &\geq&  \Psi(\hat{\gamma}_{\hat{t}},\hat{\gamma}_{\hat{t}})%-\beta[\Upsilon^3(\tilde{\gamma}_{\tilde{t}},\hat{\gamma}_{\hat{t}})
            %+\Upsilon^3(\tilde{\gamma}_{\tilde{t}},\hat{\gamma}_{\hat{t}})]
            -\sum_{i=0}^{\infty}\frac{1}{2^i}[\Upsilon^3(\gamma^i_{t_i},\hat{\gamma}_{\hat{t}})
            +\Upsilon^3(\eta^i_{t_i},\hat{\gamma}_{\hat{t}})+|\hat{t}-t_i|^2]\nonumber\\
            &&
            +\Psi(\hat{\eta}_{\hat{t}},\hat{\eta}_{\hat{t}})%-\beta[\Upsilon^3(\tilde{\eta}_{\tilde{t}},\hat{\gamma}_{\hat{t}})
            %+\Upsilon^3(\tilde{\gamma}_{\tilde{t}},\hat{\eta}_{\hat{t}})]
            -\sum_{i=0}^{\infty}\frac{1}{2^i}[\Upsilon^3(\gamma^i_{t_i},\hat{\eta}_{\hat{t}})
            +\Upsilon^3(\eta^i_{t_i},\hat{\eta}_{\hat{t}})+|\hat{t}-t_i|^2].
 \end{eqnarray}
 This implies that
 \begin{eqnarray}\label{5.6}
                         &&2{\beta}\Upsilon^3(\hat{\gamma}_{\hat{t}},\hat{\eta}_{\hat{t}})+
                         2\beta^{\frac{1}{3}}|\hat{\gamma}_{\hat{t}}(\hat{t})-\hat{\eta}_{\hat{t}}(\hat{t})|^2
                         %+2{\beta}|(\hat{\gamma}_{\hat{t}}(\hat{t})^-_N-(\hat{\eta}_{\hat{t}}(\hat{t})^-_N|^6
                         \nonumber\\
                         &
                         \leq&|W_1(\hat{\gamma}_{\hat{t}})-W_1(\hat{\eta}_{\hat{t}})|
                                   +|W_2(\hat{\gamma}_{\hat{t}})-W_2(\hat{\eta}_{\hat{t}})|+
                                   \sum_{i=0}^{\infty}\frac{1}{2^i}[\Upsilon^3(\eta^i_{t_i},\hat{\gamma}_{\hat{t}})+\Upsilon^3(\gamma^i_{t_i},\hat{\eta}_{\hat{t}})].
                                  % &\leq& 2L(2+||\hat{{\gamma}}^{1}_{{\hat{t}}}||_0+||\hat{{\gamma}}_{{\hat{t}}}^{2}||_0)
%                                   \leq 4L(1+M^m).
 \end{eqnarray}
 On the other hand, by Lemma \ref{theoremS},
 \begin{eqnarray}\label{4.7jiajia130}
 &&\sum_{i=0}^{\infty}\frac{1}{2^i}[\Upsilon^3(\eta^i_{t_i},\hat{\gamma}_{\hat{t}})+\Upsilon^3(\gamma^i_{t_i},\hat{\eta}_{\hat{t}})]\nonumber\\
 %\leq
% 2\sum_{i=0}^{\infty}\frac{1}{2^i}[||\eta^i_{t_i}-\hat{\gamma}_{\hat{t}}||_0^4+||\gamma^i_{t_i}-\hat{\eta}_{\hat{t}}||_0^4]\\
 %&\leq&2^4\sum_{i=0}^{\infty}\frac{1}{2^i}[||\eta^i_{t_i}-\hat{\eta}_{\hat{t}}||_0^4%+||\hat{\eta}_{\hat{t}}-\hat{\gamma}_{\hat{t}}||_0^4
% +||\gamma^i_{t_i}-\hat{\gamma}_{\hat{t}}||_0^4+2||\hat{\gamma}_{\hat{t}}-\hat{\eta}_{\hat{t}}||_0^4]\\
 %&\leq&2^4\sum_{i=0}^{\infty}\frac{1}{2^i}[||\eta^i_{t_i}-\hat{\eta}_{\hat{t}}||_0^4%+||\hat{\eta}_{\hat{t}}-\hat{\gamma}_{\hat{t}}||_0^4
% +||\gamma^i_{t_i}-\hat{\gamma}_{\hat{t}}||_0^4+2^4||{\gamma}^{0}_{t_{0}}-\tilde{\gamma}_{\tilde{t}}||_0^4+2^7||{\eta}^{0}_{t_{0}}-\hat{\eta}_{\hat{t}}||_0^4
% +2^7||{\gamma}^{0}_{t_{0}}-{\eta}^{0}_{t_{0}}||_0^4]\\
 &\leq&2^5\sum_{i=0}^{\infty}\frac{1}{2^i}[\Upsilon^3(\eta^i_{t_i},\hat{\eta}_{\hat{t}})%+||\hat{\eta}_{\hat{t}}-\hat{\gamma}_{\hat{t}}||_0^4
 +\Upsilon^3(\gamma^i_{t_i},\hat{\gamma}_{\hat{t}})+2\Upsilon^3(\hat{\gamma}_{\hat{t}},\hat{\eta}_{\hat{t}})]
 \leq\frac{2^6}{\beta}+{2^7}\Upsilon^3(\hat{\gamma}_{\hat{t}},\hat{\eta}_{\hat{t}}).
 %\frac{54M_1}{\beta}(2^4+\frac{1}{3})=\frac{882M_1}{\beta}.
 \end{eqnarray}
Then we have
 \begin{eqnarray}\label{5.jia6}
                         &&(2{\beta}-2^7)\Upsilon^3(\hat{\gamma}_{\hat{t}},\hat{\eta}_{\hat{t}})+
                          2\beta^{\frac{1}{3}}|\hat{\gamma}_{\hat{t}}(\hat{t})-\hat{\eta}_{\hat{t}}(\hat{t})|^2
                         %+2{\beta}|(\hat{\gamma}_{\hat{t}}(\hat{t})^-_N-(\hat{\eta}_{\hat{t}}(\hat{t})^-_N|^6
                         \nonumber\\
                         &\leq& |W_1(\hat{\gamma}_{\hat{t}})-W_1(\hat{\eta}_{\hat{t}})|
                                   +|W_2(\hat{\gamma}_{\hat{t}})-W_2(\hat{\eta}_{\hat{t}})|+\frac{2^6}{\beta}\nonumber\\
                                   %&\leq&|W_1(\hat{\gamma}_{\hat{t}})-W_1(\hat{\eta}_{\hat{t}})|
%                                   +|W_2(\hat{\gamma}_{\hat{t}})-W_2(\hat{\eta}_{\hat{t}})|
%                                   + 54\times 2^7\sum_{i=0}^{\infty}\frac{1}{2^i}[\frac{1}{2^i\beta}+\frac{1}{\beta}+16M_0^4]\nonumber\\
                                   &\leq& 2L(2+||\hat{\gamma}_{\hat{t}}||_0+||\hat{\eta}_{\hat{t}}||_0)+\frac{2^6}{\beta}
                                   \leq 4L(1+M_0)+\frac{2^6}{\beta}.
 \end{eqnarray}
   Letting $\beta\rightarrow\infty$, we get
                \begin{eqnarray*}
                \Upsilon^3(\hat{\gamma}_{\hat{t}},\hat{\eta}_{\hat{t}})
                \leq \frac{1}{2{\beta}-2^7}\left[4L(1+M_0)+\frac{2^6}{\beta}\right]\rightarrow0\ %\mbox{uniformly in}\ {\hat{\mu}}\
                          \mbox{as} \ \beta\rightarrow+\infty.
                         \end{eqnarray*}
%Noting the norm $ ||\cdot||_{0,8}$ is equivalent to the norm  $||\cdot||_{0}$,
From (\ref{s0}) it follows that
\begin{eqnarray}\label{5.66666123}
||\hat{\gamma}_{\hat{t}}-\hat{\eta}_{\hat{t}}||_0 \rightarrow0\ %\mbox{uniformly in}\ N\
                          \mbox{as} \ \beta\rightarrow+\infty.
 \end{eqnarray}
 %By (\ref{3.6}) and (\ref{w1}), we also have
% \begin{eqnarray}\label{5.65555}
%                         {\beta}d(\hat{{\gamma}}^{1}_{{\hat{t}}},\hat{{\gamma}}_{{\hat{t}}}^{2})
%                         &\leq&|W_1(\hat{{\gamma}}^{1}_{{\hat{t}}})-W_1(\hat{{\gamma}}_{{\hat{t}}}^{2})|
%                                   +|W_2(\hat{{\gamma}}^{1}_{{\hat{t}}})-W_2(\hat{{\gamma}}^{2}_{{\hat{t}}})|\nonumber\\
%                                   &\leq& 2(1+||\hat{{\gamma}}^{1}_{{\hat{t}}}||^m_0\vee||\hat{{\gamma}}^{2}_{{\hat{t}}}||^m_0)\hat{\omega}(
%                                   ||\hat{{\gamma}}^{1}_{{\hat{t}}}-\hat{{\gamma}}^{2}_{{\hat{t}}}||_0) \leq2(1+M^m)\hat{\omega}(
%                                   ||\hat{{\gamma}}^{1}_{{\hat{t}}}-\hat{{\gamma}}^{2}_{{\hat{t}}}||_0).
% \end{eqnarray}
                   % Since $W_1,W_2\in C^0(\Lambda)$,
                   Combining (\ref{s0}),  (\ref{w1}), (\ref{5.6}), (\ref{4.7jiajia130}) and (\ref{5.66666123}), we see
                           that
 \begin{eqnarray}\label{5.10112345}
                        &&{\beta}||\hat{\gamma}_{\hat{t}}-\hat{\eta}_{\hat{t}}||^6_{0}+\beta^{\frac{1}{3}}|\hat{\gamma}_{\hat{t}}(\hat{t})-\hat{\eta}_{\hat{t}}(\hat{t})|^2
                        \leq4{\beta}\Upsilon^3(\hat{\gamma}_{\hat{t}},\hat{\eta}_{\hat{t}})
                        +\beta^{\frac{1}{3}}|\hat{\gamma}_{\hat{t}}(\hat{t})-\hat{\eta}_{\hat{t}}(\hat{t})|^2\nonumber\\
                        &\leq& 4L||\hat{\gamma}_{\hat{t}}-\hat{\eta}_{\hat{t}}||_0
                                   +\frac{2^7}{\beta}+{2^{10}}||\hat{\gamma}_{\hat{t}}-\hat{\eta}_{\hat{t}}||_{0}^6
                                   \rightarrow0 \ %\mbox{uniformly in}\ N \
                                   \mbox{as} \ \beta\rightarrow+\infty.
 \end{eqnarray}
    Then we have
    %\begin{eqnarray}\label{5.1012344556789}
%                      &&\beta^{\frac{1}{3}}
%                      ||\hat{{\gamma}}_{{\hat{t}}}-\hat{{\eta}}_{{\hat{t}}}||^2_{0}
%                                   \rightarrow0 \ %\mbox{uniformly in}\ N \
%                                   \mbox{as} \ \beta\rightarrow+\infty,
% \end{eqnarray}
% and
  \begin{eqnarray}\label{5.10123445567890}
                      &&\beta^{\frac{1}{2}}|\hat{{\gamma}}_{{\hat{t}}}(\hat{t})-\hat{{\eta}}_{{\hat{t}}}(\hat{t})|^2
                     \leq\beta^{\frac{1}{6}}\left(4L||\hat{\gamma}_{\hat{t}}-\hat{\eta}_{\hat{t}}||_0
                                   +\frac{2^7}{\beta}+{2^{10}}||\hat{\gamma}_{\hat{t}}-\hat{\eta}_{\hat{t}}||_{0}^6\right)\nonumber\\
                                   &=&4L\beta^{\frac{1}{6}}||\hat{\gamma}_{\hat{t}}-\hat{\eta}_{\hat{t}}||_0
                                   +\frac{2^7}{\beta^{\frac{5}{6}}}+{2^{10}}\beta^{\frac{1}{6}}||\hat{\gamma}_{\hat{t}}-\hat{\eta}_{\hat{t}}||_{0}^6
                                   \rightarrow0 \ %\mbox{uniformly in}\ N \
                                   \mbox{as} \ \beta\rightarrow+\infty.
 \end{eqnarray}
 Thus, we get (\ref{5.10}) holds true.
 \par
   $Step\ 3.$ There exists % ${{M}_0}>0$  and
   $N_0>0$
    such that
                   $\hat{t}\in [\tilde{t},T)$
                  % , $\hat{{\gamma}}^{1}_{{\hat{t}}}, \hat{{\gamma}}^{2}_{{\hat{t}}}\in
%                \Lambda_{M_0}$ and $|\hat{{\gamma}}^{1}_{{\hat{t}}}({\hat{t}})|\vee|\hat{{\gamma}}^{2}_{{\hat{t}}}({\hat{t}})|<\frac{M_0}{2}$
                for all $\beta\geq N_0$.
 \par
%First,   noting $\nu$ is independent of  $\beta$, by the definition of  ${\Psi}$,
% there exists an ${M}_0>0$  that is sufficiently  large   that
%           $
%           \Psi(\gamma^1_t, \gamma^2_t)<0
%           $ for all $t\in [T-\bar{a},T]$ and $|\gamma^1_t(t)|\vee|\gamma^2_t(t)|\geq\frac{M_0}{2}$. Thus, we have $|\hat{{\gamma}}^{1}_{{\hat{t}}}({\hat{t}})|\vee|\hat{{\gamma}}^{2}_{{\hat{t}}}({\hat{t}})|<\frac{M_0}{2}$.
% \par
%  Next, for the fixed ${M_0}>0$,
By (\ref{5.66666123}), we can let $N_0>0$ be a large number such that
$$
                         L||\hat{\gamma}_{\hat{t}}-\hat{\eta}_{\hat{t}}||_0
                         \leq
                         \frac{\tilde{m}}{4},
$$
               for all $\beta\geq N_0$.
            Then we have $\hat{t}\in [\tilde{t},T)$ for all $\beta\geq N_0$. Indeed, if say $\hat{t}=T$,  we will deduce the following contradiction:
 \begin{eqnarray*}
                         \frac{\tilde{m}}{2}\leq\Psi(\hat{\gamma}_{\hat{t}},\hat{\eta}_{\hat{t}})\leq \phi(\hat{\gamma}_{\hat{t}})-\phi(\hat{\eta}_{\hat{t}})\leq
                        L||\hat{\gamma}_{\hat{t}}-\hat{\eta}_{\hat{t}}||_0
                         \leq
                         \frac{\tilde{m}}{4}.
 \end{eqnarray*}
%
%
%
%
%\par
%Next, noting $\Psi(\gamma^1_T,\gamma^2_T)\leq0<\frac{\tilde{m}}{2}\leq \Psi(\tilde{\gamma}_{\tilde{t}},\tilde{\gamma}_{\tilde{t}})$ for every $\gamma^1_T,\gamma^2_T\in \Lambda_T$,
%   we have $\hat{t}\in (T-\bar{a},T)$.
%\par
%Next, from (\ref{psi0}), it follows that
% \begin{eqnarray*}
% \frac{\tilde{m}}{2}&\leq& W(\hat{{\gamma}}^{1}_{{\hat{t}}})- V(\hat{{\gamma}}^{2}_{{\hat{t}}})\\
% &=&W(\hat{{\gamma}}^{1}_{{\hat{t}}})-W(\hat{{\gamma}}^{1}_{{\hat{t}},T})+W(\hat{{\gamma}}^{1}_{{\hat{t}},T})
%    -V(\hat{{\gamma}}^{2}_{{\hat{t}},T})+V(\hat{{\gamma}}^{2}_{{\hat{t}},T})- V(\hat{{\gamma}}^{2}_{{\hat{t}}})\\
%    &\leq&(C+C_1)(1+2M_0)(T-\hat{t})^{\frac{1}{2}}+L\mu(|\hat{{\gamma}}^{1}_{{\hat{t}}}-\hat{{\gamma}}^{2}_{{\hat{t}}}|).
% \end{eqnarray*}
%On the other hand, from (\ref{5.10}) it follows that there exists a constant $N>0$ such that
% $$
%  ||\hat{{\gamma}}^{1}_{{\hat{t}}}-\hat{{\gamma}}^{2}_{{\hat{t}}}||_0<\frac{\tilde{m}}{4L(1+\mu([0,T]))}, \ \ \mbox{for  all} \ \ \beta>N.
% $$
% Then
% $$
% (C+C_1)(1+2M_0)(T-\hat{t})^{\frac{1}{2}}>\frac{\tilde{m}}{4}.
% $$
% %Letting $\Delta=\frac{\tilde{m}^2}{16(C+C_1)^2(1+2M_0)^2}$,
%  Thus, we have
%  $\hat{t} \in [T-\bar{a},T)$  for all $\beta>N$.
\par
 $Step\ 4.$    Completion of the proof.
\par
          From above all,  for the fixed   $N_0>0$ in step 3, %and $0<\Delta<T$ in step 3,
          we  find
$(\hat{t},\hat{\gamma}_{\hat{t}}), (\hat{t},\hat{\eta}_{\hat{t}})\in [\tilde{t}, T]\times
                 \Lambda^{\tilde{t}}$   satisfying $\hat{t}\in [\tilde{t},T)$  for all $\beta\geq N_0$
           such that
\begin{eqnarray}\label{psi4}
            \Psi_1(\hat{\gamma}_{\hat{t}},\hat{\eta}_{\hat{t}})\geq \Psi(\tilde{\gamma}_{\tilde{t}},\tilde{\gamma}_{\tilde{t}}) \ \ \mbox{and} \ \      \Psi_1(\hat{\gamma}_{\hat{t}},\hat{\eta}_{\hat{t}})\geq
                   {\Psi}_1(\gamma_t,\eta_t),
                   \  (t,\gamma_t,\eta_t)\in [\hat{t},T]\times  \Lambda^{\hat{t}}\times \Lambda^{\hat{t}},
\end{eqnarray}
where we define
$$
     \Psi_1(\gamma_t,\eta_t):=  \Psi(\gamma_t,\eta_t)
        -\sum_{i=0}^{\infty}
        \frac{1}{2^i}\left[{\Upsilon}^3(\gamma^i_{t_i},\gamma_t)+{\Upsilon}^3(\eta^i_{t_i},\eta_t)+|{t}-t_i|^2\right], \ \  \  (t,\gamma_t,\eta_t)\in [\tilde{t},T]\times  \Lambda^{\tilde{t}}\times \Lambda^{\tilde{t}}.
$$
%By (\ref{5.10}), we can let $\beta>N$ be large enough such that
%\begin{eqnarray}\label{betamu}
%2\beta L(1+2L)||\hat{{\gamma}}^{1}_{{\hat{t}}}-\hat{{\gamma}}^{2}_{{\hat{t}}}||_0^2
%                                  +L||\hat{{\gamma}}^{1}_{{\hat{t}}}-\hat{{\gamma}}^{2}_{{\hat{t}}}||_0\leq \frac{c}{8}.
%\end{eqnarray}
We put, for $(t,\gamma_t,\eta_t)\in [\tilde{t},T]\times {\Lambda}\times {\Lambda}$,
\begin{eqnarray*}
                             {W}'_{1}(\gamma_t)=W_1(\gamma_t)-\varepsilon\frac{\nu T-t}{\nu
                 T}\Upsilon^3(\gamma_t)
                % -\varepsilon |t-{\hat{t}}|^2
                 -\varepsilon \overline{\Upsilon}^3(\gamma_t,\hat{\gamma}_{\hat{t}})
                %-\frac{\varepsilon}{t-T+\bar{a}}
                -\sum_{i=0}^{\infty}
        \frac{1}{2^i}\overline{\Upsilon}^3(\gamma^i_{t_i},\gamma_t);
        \end{eqnarray*}
        \begin{eqnarray*}
                             {W}'_{2}(\eta_t)=W_2(\eta_t)+\varepsilon\frac{\nu T-t}{\nu
                 T}\Upsilon^3(\eta_t)
                 %+\varepsilon |t-{\hat{t}}|^2
                 +\varepsilon \overline{\Upsilon}^3(\eta_t,\hat{\eta}_{\hat{t}})+\sum_{i=0}^{\infty}
        \frac{1}{2^i}{\Upsilon}^3(\eta^i_{t_i},\eta_t).
\end{eqnarray*}
                  Now we  define,  for $(t,x_N,y_N,\gamma_t,\eta_t)\in [\hat{t},T]\times H_N\times H_N\times \Lambda^{\hat{t}}\times \Lambda^{\hat{t}}$,
%\begin{eqnarray*}
%                             &&\tilde{W}_{1}(\gamma_t,x_N,x^-_N)=\sup_{z\in H^\bot_N}\sup_{\xi_t\in \Lambda^{\hat{t}},(\xi_t(t))_N=x_N,(\xi_t(t))^-_N=z}
%                             \bigg{[}W'_{1}(\xi_t)
%                             -2^7\beta \Upsilon^3(\gamma_{t},\xi_{t})-2^5\beta|x^-_N-z|^6\bigg{]};\\
%                             %&&\tilde{W}_{1,\mu}(t,\gamma^1,x_0)=\tilde{W}_{1,\mu}(t^1_{x_0},\gamma^1,x_0)-(t^1_{x_0}-t)^{\frac{1}{2}},\ \ t\in [\hat{t}_\mu, t^1_{x_0}];\\
%                             &&\tilde{W}_{2}(\eta_t,y_N,y^-_N)=\inf_{z\in H^\bot_N}\inf_{\xi_t\in\Lambda^{\hat{t}},(\xi_t(t))_N=y_N,(\xi_t(t))^-_N=z}
%                             \bigg{[}W'_{2}(\xi_t)
%                             +2^7\beta \Upsilon^3(\eta_{t},\xi_{t})+2^5\beta|y^-_N-z|^6\bigg{]}.
%\end{eqnarray*}
\begin{eqnarray*}
                             &&\tilde{W}_{1}(\gamma_t,x_N)=\sup_{\xi_t\in \Lambda^{\hat{t}},(\xi_t(t))_N=x_N}
                             \bigg{[}W'_{1}(\xi_t)
                             -2^5\beta \Upsilon^3(\gamma_{t},\xi_{t})-2\beta^{\frac{1}{3}}|(\gamma_t(t))^-_N-(\xi_t(t))^-_N|^2\bigg{]};\\
                             %&&\tilde{W}_{1,\mu}(t,\gamma^1,x_0)=\tilde{W}_{1,\mu}(t^1_{x_0},\gamma^1,x_0)-(t^1_{x_0}-t)^{\frac{1}{2}},\ \ t\in [\hat{t}_\mu, t^1_{x_0}];\\
                             &&\tilde{W}_{2}(\eta_t,y_N)=\inf_{\xi_t\in\Lambda^{\hat{t}},(\xi_t(t))_N=y_N}
                             \bigg{[}W'_{2}(\xi_t)
                             +2^5\beta \Upsilon^3(\eta_{t},\xi_{t})+2\beta^{\frac{1}{3}}|(\eta_t(t))^-_N-(\xi_t(t))^-_N|^2\bigg{]}.
\end{eqnarray*}
Then by Lemma \ref{theoremS}  we obtain that, for all %$t\in [\hat{t}, T]$,
   $(t,x_N,y_N,\xi_t)\in [\hat{t},T]\times H_N\times H_N\times \Lambda$,
\begin{eqnarray}\label{wv}
                         &&\tilde{W}_{1}(\xi_t,x_N)-  \tilde{W}_{2}(\xi_t,y_N)-{\beta^{\frac{1}{3}}|x_N-y_N|^2}\nonumber   \\
                          &=&\sup_{
                          \gamma_t,\eta_t\in\Lambda^{\hat{t}},(\gamma_t(t))_N=x_N,
                          (\eta_t(t))_N=y_N}
                          \bigg{[}W'_{1}(\gamma_t)
                          -2^5\beta \Upsilon^3(\xi_{t},\gamma_{t})-2\beta^{\frac{1}{3}}|(\gamma_t(t))^-_N-(\xi_t(t))^-_N|^2\nonumber\\
                          &&-W'_{2}(\eta_t)-2^5\beta \Upsilon^3(\xi_{t},\eta_{t})-2\beta^{\frac{1}{3}}|(\eta_t(t))^-_N-(\xi_t(t))^-_N|^2
                          -\beta^{\frac{1}{3}}|(\gamma_t(t))_N-(\eta_t(t))_N|^2\bigg{]}\nonumber \\
                          &\leq&\sup_{
                          \gamma_t,\eta_t\in\Lambda^{\hat{t}},(\gamma_t(t))_N=x_N,
                          (\eta_t(t))_N=y_N}\bigg{[}W'_{1}(\gamma_t)
                          -W'_{2}(\eta_t)-\beta \Upsilon^3(\gamma_{t},\eta_{t})-\beta^{\frac{1}{3}}|\gamma_t(t)-\eta_t(t)|^2%-\beta|(\gamma_t(t))^-_N-(\eta_t(t))^-_N|^6
                          \bigg{]}\nonumber \\
                          &\leq&W'_{1}(\hat{{\gamma}}_{{\hat{t}}})-W'_{2}(\hat{{\eta}}_{{\hat{t}}})-\beta \Upsilon^3(\hat{{\gamma}}_{{\hat{t}}},\hat{{\eta}}_{{\hat{t}}})
                          -\beta^{\frac{1}{3}}|\hat{{\gamma}}_{{\hat{t}}}({\hat{t}})-\hat{{\eta}}_{{\hat{t}}}({\hat{t}})|^2,%-\beta|(\hat{{\gamma}}_{{\hat{t}}}({\hat{t}}))^-_N-(\hat{{\eta}}_{{\hat{t}}}({\hat{t}}))^-_N|^6,
\end{eqnarray}
                       where the last inequality becomes equality if and only if $t={\hat{t}}$, $\gamma_t=\hat{{\gamma}}_{{\hat{t}}}, \eta_t=\hat{{\eta}}_{{\hat{t}}}$.
                       %and $x_N=(\hat{{\gamma}}_{{\hat{t}}}({\hat{t}}))_N,y_N=(\hat{{\eta}}_{{\hat{t}}}({\hat{t}}))_N$.
                        The previous inequality becomes equality if  $\xi_{t}=\frac{\gamma_{t}+\eta_{t}}{2}$. % and $z^-_N=\frac{(\gamma_t(t))^-_N+(\eta_t(t))^-_N}{2}$.
                          Then we obtain that, for all %$t\in [\hat{t}, T]$,
  $(t,x_N,y_N,\xi_t)\in [\hat{t},T]\times H_N\times H_N\times \Lambda$,
\begin{eqnarray*}
\tilde{W}_{1}(\xi_t,x_N)-  \tilde{W}_{2}(\xi_t,y_N)-{\beta^{\frac{1}{3}}|x_N-y_N|^2}
       \leq  W'_{1}(\hat{{\gamma}}_{{\hat{t}}})-W'_{2}(\hat{{\eta}}_{{\hat{t}}})
                          -\beta^{\frac{1}{3}}|\hat{{\gamma}}_{{\hat{t}}}({\hat{t}})-\hat{{\eta}}_{{\hat{t}}}({\hat{t}})|^2,
\end{eqnarray*}
    and the equality   holds at ${\hat{t}},(\hat{{\gamma}}_{\hat{t}}({\hat{t}}))_N,(\hat{{\eta}}_{\hat{t}}({\hat{t}}))_N,
    \hat{\xi}_{{\hat{t}}}=\frac{\hat{{\gamma}}_{\hat{t}}+\hat{{\eta}}_{\hat{t}}}{2}$.
 \par
    Define,  for $(t,x_N, y_N)\in [0,T]\times H_N\times H_N$,
               \begin{eqnarray}\label{definewv}
                   \bar{W}_{1}(t,x_N)&=&\begin{cases}\tilde{W}_{1}(\hat{\xi}_{{\hat{t}}},x_N)-(\hat{t}-t)^{\frac{1}{2}}, \ \   t\in[0,\hat{t}),\\
                  \tilde{W}_{1}(\hat{\xi}_{{\hat{t}},t,A}, x_N),  ~~~~~\ \ \ \ \ \ \  t\in[\hat{t},T];            \end{cases}\nonumber\\
                            \bar{W}_{2}(t,y_N)&=&\begin{cases}\tilde{W}_{2}(\hat{\xi}_{{\hat{t}}},y_N)+(\hat{t}-t)^{\frac{1}{2}}, \ \   t\in[0,\hat{t}),\\
                  \tilde{W}_{2}(\hat{\xi}_{{\hat{t}},t,A}, y_N),  ~~~~~\ \ \ \ \ \ \ t\in[\hat{t},T];          \end{cases}
\end{eqnarray}
and
               \begin{eqnarray}\label{definewv111}
                   \breve{{W}}_{1}(t,x_N)=\limsup_{s\rightarrow t}\bar{W}_{1}(s,x_N),\ \ (t,x_N)\in [0,T]\times H_N;
                   %\ \ \mbox{and}\ \
%             \breve{{W}}_{1}(T,x_N)=\bar{W}_1(T,x_N), \ \ x_N\in H_N;
\nonumber\\
                            \breve{W}_{2}(t,y_N)= \liminf_{s\rightarrow t}\bar{W}_{2}(s\mathfrak{\mathfrak{}},y_N),\ \ (t,y_N)\in [0,T]\times H_N.
                           % ,\ \mbox{and}\
%                 \breve{W}_{2}(T,y_N)= \bar{W}_{2}(T,y_N), \  y_N\in  H_N.
\end{eqnarray}
          %     \begin{eqnarray}\label{definewv111}
%                   \breve{{W}}_{1}(t,x_N)=\bar{W}_{1}(t,x_N)\vee\limsup_{s\rightarrow t}\bar{W}_{1}(s,x_N),\ \ (t,x_N)\in [0,T]\times H_N;
%                   %\ \ \mbox{and}\ \
%%             \breve{{W}}_{1}(T,x_N)=\bar{W}_1(T,x_N), \ \ x_N\in H_N;
%\nonumber\\
%                            \breve{W}_{2}(t,y_N)= \bar{W}_{2}(t,y_N)\wedge\liminf_{s\rightarrow t}\bar{W}_{2}(s\mathfrak{\mathfrak{}},y_N),\ \ (t,y_N)\in [0,T]\times H_N.
%                           % ,\ \mbox{and}\
%%                 \breve{W}_{2}(T,y_N)= \bar{W}_{2}(T,y_N), \  y_N\in  H_N.
%\end{eqnarray}
Thus by Lemma \ref{lemma4.3},  $ \breve{W}_{1}(t,x_N)- \breve{W}_{2}(t,y_N)-{\beta^{\frac{1}{3}}|x_N-y_N|^2}$ has a maximum at $({\hat{t}},(\hat{{\gamma}}_{\hat{t}}({\hat{t}}))_N,
(\hat{{\eta}}_{\hat{t}}({\hat{t}}))_N)$
 on $[0,T]\times H_N\times H_N$. %Here $\bar{\mathcal{O}}_{M_0}:=\{x\in R^d:|x|\leq M_0\}$.
Then, by Lemmas \ref{lemma4.3} and \ref{lemma4.344}, the Theorem 8.3 in \cite{cran2} can be used  to obtain sequences  $(x^{k}_N,y^{k}_N)\in H_N\times H_N$, $l_{k},s_{k}\in [0,T]$ such that $(l_{k},x_N^{k})\rightarrow ({\hat{t}},(\hat{{\gamma}}_{{\hat{t}}}({\hat{t}}))_N)$,
                  $(s_{k},y_N^{k})\rightarrow ({\hat{t}},(\hat{\eta}_{{\hat{t}}}({\hat{t}}))_N)$ as $k\rightarrow+\infty$ and the sequences of functions $\varphi_k,\psi_k\in C^{1,2}((T-\bar{a},T)\times H_N)$
                  such that
$$
           \breve{W}_{1}(t, x_N)-\varphi_k(t,x_N)\leq0,\ \ \breve{W}_{2}(t,x_N)+\psi_k(t,x_N)\geq0,
$$
 equalities only  hold true at $(l_{k},x^{k}_N)$,  $(s_{k},y^{k}_N)$, respectively,
\begin{eqnarray*}
                                (\varphi_k)_t(l_{k},x^{k}_N)\rightarrow b_{1}, &&\ \ \ (\psi_k)_t(s_{k},y^{k}_N)\rightarrow b_{2},\\
                                \nabla_x\varphi_k(l_{k},x^{k}_N)\rightarrow 2\beta^{\frac{1}{3}} ((\hat{{\gamma}}_{{\hat{t}}}({{\hat{t}}}))_N-(\hat{{\eta}}_{{\hat{t}}}({{\hat{t}}}))_N), &&\ \ \
                                             -\nabla_x\psi_k(s_{k},y^{k}_N)\rightarrow 2\beta^{\frac{1}{3}} ((\hat{{\gamma}}_{{\hat{t}}}({{\hat{t}}}))_N-(\hat{{\eta}}_{{\hat{t}}}({{\hat{t}}}))_N),\\
                                \nabla_x^2\varphi_k(l_{k},x^{k}_N)\rightarrow X_N, && \ \ \
                                             \nabla_x^2\psi_k(s_{k},y^{k}_N)\rightarrow Y_N,
\end{eqnarray*}
               where $b_{1}+b_{2}=0$ and $X_N,Y_N$ satisfy the following inequality:
\begin{eqnarray}\label{II}
                              {-6\beta^{\frac{1}{3}}}\left(\begin{array}{cc}
                                    I&0\\
                                    0&I
                                    \end{array}\right)\leq \left(\begin{array}{cc}
                                    X_N&0\\
                                    0&Y_N
                                    \end{array}\right)\leq  6\beta^{\frac{1}{3}} \left(\begin{array}{cc}
                                    I&-I\\
                                    -I&I
                                    \end{array}\right).
\end{eqnarray}
%Here
%%$x^5:=\sum_{i=1}^{N}x_i^5e_i$ for $x\in H_N$, and
%$B=P_NP_N$.
\par
We claim that we can assume the  sequences $\{l_{k}\}_{k\geq1}\in [\hat{t},T)$ and $\{s_{k}\}_{k\geq1}\in [\hat{t},T)$. Indeed, if not, for example, there exists a subsequence of $\{l_{k}\}_{k\geq1}$ still denoted by itself such that $l_{k}<\hat{t}$ for all $k\geq0$.
% We claim that there exist sequences $t_{k}\in [\hat{t},T), s_{k}\in [\hat{t},T)$. Indeed, if not, for example, $t_{k}<\hat{t}$.
Since $\breve{W}_{1}(t, x_N)-\varphi_k(t,x_N)$
            has a maximum at $(l_{k},x^{k}_N)$ on $[0,T]\times H_N$, we obtain that
            $$
           (\varphi_k)_t(l_{k},x^{k}_N)={\frac{1}{2}}(\hat{t}-l_{k})^{-\frac{1}{2}}\rightarrow+\infty,\ \mbox{as}\ k\rightarrow+\infty,
           $$
            % By Lemma \ref{lemma4.333333}, we have that
%             $$
%           (\varphi_k)_t(t_{k},x^{k}_0)={\frac{1}{2}}(t^1_{x^{k}_0}\vee t^2_{y^{k}_0}-t_{k})^{-\frac{1}{2}}\rightarrow+\infty,\ \mbox{as}\ k\rightarrow+\infty,
%            $$
            which is contradict to  $(\varphi_k)_t(l_{k},x^{k}_N)\rightarrow b_{1}$, $(\psi_k)_t(s_{k},y^{k}_N)\rightarrow b_{2}$ and $b_{1}+b_{2}=0$.
 \par
  We may without loss of generality assume that $\varphi_k,\psi_k$ grow quadratically at $\infty$. Now we consider the functional,
              for $(t,\gamma_t), (s,\eta_s)\in [T-\bar{a},T]\times{{\Lambda}}$,
\begin{eqnarray}\label{4.1111}
                 \Gamma_{k}(\gamma_t,\eta_s)&=& W'_{1}(\gamma_t)-W'_{2}(\eta_s)-2^5\beta(\Upsilon^3(\gamma_{t},\hat{{\xi}}_{{\hat{t}}})
                +\Upsilon^3(\eta_{s},\hat{{\xi}}_{{\hat{t}}}))-\varphi_k(t,(\gamma_t(t))_N)-\psi_k(s,(\eta_s(s))_N)\nonumber\\
                &&-2\beta^{\frac{1}{3}}|(\gamma_t(t))^-_N-(e^{(t-\hat{t})A}\hat{\xi}_{\hat{t}}(\hat{t}))^-_N|^2-2\beta^{\frac{1}{3}}|(\eta_s(s))^-_N-(e^{(s-\hat{t})A}\hat{\xi}_{\hat{t}}(\hat{t}))^-_N|^2.
\end{eqnarray}
 Define a sequence of positive numbers $\{\delta_i\}_{i\geq0}$  by %$\delta_0=\beta$ and
        $\delta_i=\frac{1}{2^i}$ for all $i\geq0$.  For every $k$ and $\delta>0$,
           % $\Psi$ is a  upper semicontinuous function and $\Upsilon^3$ is a gauge-type function,
           from Lemma \ref{theoremleft} it follows that,
  for every  $(\check{t}_{0},\check{\gamma}^{0}_{\check{t}_{0}}), (\check{s}_{0},\check{\eta}^{0}_{\check{s}_{0}})\in [\hat{t},T]\times  \Lambda^{\hat{t}}$ satisfy
$$
\Gamma_{k}(\check{\gamma}^{0}_{\check{t}_{0}},\check{\eta}^{0}_{\check{s}_{0}})\geq \sup_{(t,\gamma_t),(s,\eta_s)\in [\hat{t},T]\times \Lambda^{\hat{t}}}\Gamma_{k}(\gamma_t,\eta_s)-{\delta},\
%\    \mbox{and} \ \ \Gamma_{k}(\gamma^0_{t_0},\eta^0_{t_0})\geq \Psi(\tilde{\gamma}_{\tilde{t}},\tilde{\gamma}_{\tilde{t}}) >\frac{\tilde{m}}{2},
 $$
  there exist $(\check{t},\check{\gamma}_{\check{t}}), (\check{s},\check{\eta}_{\check{s}})\in [\hat{t},T]\times \Lambda^{\hat{t}}$ and two sequences $\{(\check{t}_{i},\check{\gamma}^{i}_{\check{t}_{i}})\}_{i\geq1}, \{(\check{s}_{i},\check{\eta}^{i}_{\check{s}_{i}})\}_{i\geq1}\subset
  [\hat{t},T]\times \Lambda^{\hat{t}}$ such that
  \begin{description}
        \item{(i)} $\overline{\Upsilon}^3(\check{\gamma}^{0}_{\check{t}_{0}},\check{\gamma}_{\check{t}})+\overline{\Upsilon}^3(\check{\eta}^0_{\check{s}_0},\check{\eta}_{\check{s}})\leq {\delta}$,
         $\overline{\Upsilon}^3(\check{\gamma}^{i}_{\check{t}_{i}},\check{\gamma}_{\check{t}})
         +\overline{\Upsilon}^3(\check{\eta}^{i}_{\check{s}_{i}},\check{\eta}_{\check{s}})\leq \frac{\delta}{2^i}$
          and $\check{t}_{i}\uparrow \check{t}$, $\check{s}_{i}\uparrow \check{s}$ as $i\rightarrow\infty$,
        \item{(ii)}  $\Gamma_k(\check{\gamma}_{\check{t}},\check{\eta}_{\check{s}})%-\beta[\overline{\Upsilon}^3(\tilde{\gamma}_{\tilde{t}},\hat{\gamma}^k_{\hat{t}_k})+\overline{\Upsilon}^3(\tilde{\gamma}_{\tilde{t}},\hat{\eta}^k_{\hat{t}_k})]
            -\sum_{i=0}^{\infty}\frac{1}{2^i}[\overline{\Upsilon}^3(\check{\gamma}^{i}_{\check{t}_{i}},\check{\gamma}_{\check{t}})
        +\overline{\Upsilon}^3(\check{\eta}^{i}_{\check{s}_{i}},\check{\eta}_{\check{s}})]\geq \Gamma_{k}(\check{\gamma}^{0}_{\check{t}_{0}},\check{\eta}^{0}_{\check{s}_{0}})$, and
        \item{(iii)}    for all $(t,\gamma_t,s,\eta_s)\in [\check{t},T]\times \Lambda^{\check{t}}\times[\check{s},T]\times \Lambda^{\check{s}}\setminus \{(\check{t},\check{\gamma}_{\check{t}},\check{s},\check{\eta}_{\check{s}})\}$,
        \begin{eqnarray*}
        \Gamma_{k}(\gamma_t,\eta_s)%-\beta[\overline{\Upsilon}^3(\tilde{\gamma}_{\tilde{t}},\gamma_s)+\overline{\Upsilon}^3(\tilde{\gamma}_{\tilde{t}},\eta_s)]
        -\sum_{i=0}^{\infty}
        \frac{1}{2^i}[\overline{\Upsilon}^3(\check{\gamma}^{i}_{\check{t}_{i}},\gamma_t)+\overline{\Upsilon}^3(\check{\eta}^{i}_{\check{s}_{i}},\eta_s)]
            <\Gamma_{k}(\check{\gamma}_{\check{t}},\check{\eta}_{\check{s}})%-\beta[\overline{\Upsilon}^3(\tilde{\gamma}_{\tilde{t}},\hat{\gamma}^k_{\hat{t}_k})
            %+\overline{\Upsilon}^3(\tilde{\gamma}_{\tilde{t}},\hat{\eta}^k_{\hat{t}_k})]
            -\sum_{i=0}^{\infty}\frac{1}{2^i}[\overline{\Upsilon}^3(\check{\gamma}^{i}_{\check{t}_{i}},\check{\gamma}_{\check{t}})
            +\overline{\Upsilon}^3(\check{\eta}^{i}_{\check{s}_{i}},\check{\eta}_{\check{s}})].
        \end{eqnarray*}
        \end{description}
  By the following Lemma \ref{lemma4.4}, we have
\begin{eqnarray}\label{4.22}
       && \check{t}\rightarrow l_k,\ (\check{\gamma}_{\check{t}}(\check{t}))_N\rightarrow x^k_N,\ \check{s}\rightarrow s_k,\ (\check{\eta}_{\check{s}}(\check{s}))_N\rightarrow y^k_N \ \mbox{as}\ \delta\rightarrow0;
       \end{eqnarray}
and
\begin{eqnarray}\label{4.23}
\lim_{k\rightarrow\infty}\limsup_{\delta\rightarrow0}[%\Upsilon^3(\check{\gamma}_{\check{t}},\hat{\gamma}_{\hat{t}})
%+\Upsilon^3(\check{\eta}_{\check{s}},\hat{\eta}_{\hat{t}})
||\check{\gamma}_{\check{t}}-\hat{\gamma}_{\hat{t},\check{t},A}||_0+ ||\check{\eta}_{\check{s}}-\hat{\eta}_{\hat{t},\check{s},A}||_0]=0.
\end{eqnarray}
From $l_{k},s_{k}\rightarrow {\hat{t}}$ as $k\rightarrow+\infty$ and ${\hat{t}}<T$
%and $|\hat{{\gamma}}^{1}_{{\hat{t}}}({\hat{t}})|\vee|\hat{{\gamma}}^{2}_{{\hat{t}}}({\hat{t}})|<\frac{M_0}{2}$
 for $\beta>N_0$, it follows that, for every fixed $\beta>N_0$,   constant $ K_\beta>0$ exists such that
$$
            |l_{k}|\vee|s_{k}|<T, %\ \    |x^{k}_0|\vee|y^{k}_0|<\frac{M_0}{2},
            \ \ \mbox{for all}    \ \ k\geq K_\beta.
$$
For every fixed $k>K_\beta$, by (\ref{4.22}), there exists constant $\Delta_{k,\beta}>0$ such that
$$
            |\check{t}|\vee|\check{s}|<T, %\ \    |\check{\gamma}^1_{\check{t}}(\check{t})|\vee|\check{\gamma}^2_{\check{s}}(\check{s})|<\frac{M_0}{2},
            \ \ \mbox{for all}    \ \ 0<\delta<\Delta_{k,\beta}.
$$
%Moreover, for every fixed $\mu$, we can choose $k_{\mu}\geq K_{\mu}$ be large enough  such that %$({\check{t}_{k}}_{n,n_k},\check{s}_{n,n_k},\check{x}^1_{n,n_k},\check{x}^2_{n,n_k})$
% \begin{eqnarray}\label{jiajia}
%                 |\Phi_{{k_\mu}}-\Phi_{\mu}|\leq \frac{c}{8},\ \ \
%                 |\Pi_{{k_\mu}}-\Pi_{\mu}|\leq \frac{c}{8}.
% \end{eqnarray}
Now, for every $\beta>N_0$, $k>K_\beta$ and $0<\delta<\Delta_{k,\beta}$, from the definition of viscosity solutions it follows that
\begin{eqnarray}\label{vis1}
                      &&(\varphi_k)_t({\check{t}},(\check{\gamma}_{{\check{t}}}({\check{t}}))_N)
                      -\frac{\varepsilon}{\nu T}\Upsilon^3(\check{\gamma}_{{\check{t}}})%+|\check{\gamma}_{{\check{t}}}({\check{t}})|^2)
                     %-\varepsilon({\check{t}}-T+\bar{a})^{-2}
                     +2\sum_{i=0}^{\infty}\frac{1}{2^i}[(\check{t}-\check{t}_{i})+(\check{t}-t_i)]+(A^*\nabla_{x}(\varphi_k)({\check{t}}, (\check{\gamma}_{{\check{t}}}(\check{{t}}))_N),\check{\gamma}_{{\check{t}}}({\check{t}}))_H
                    \nonumber\\
                     && +2\varepsilon({\check{t}}-{\hat{t}})
                     -4\beta^{\frac{1}{3}}(A^*((\check{\gamma}_{\check{t}}(\check{t}))_N-(e^{(\check{t}-\hat{t})A}\hat{\xi}_{\hat{t}}(\hat{t}))_N),
                     \check{\gamma}_{{\check{t}}}(\check{{t}})%-e^{(\check{t}-\hat{t})A}\hat{\xi}_{\hat{t}}(\hat{t})
                     )_H\nonumber\\
                     &&+{\mathbf{H}}\bigg{(}\check{\gamma}_{{\check{t}}},  W_1(\check{\gamma}_{{\check{t}}}), \nabla_{x}(\varphi_k)({\check{t}}, (\check{\gamma}_{{\check{t}}}(\check{{t}}))_N)
                     +4\beta^{\frac{1}{3}}((\check{\gamma}_{\check{t}}(\check{t}))^-_N-(e^{(\check{t}-\hat{t})A}\hat{\xi}_{\hat{t}}(\hat{t}))^-_N)
                      +2^5\beta\partial_x\Upsilon^3(\check{\gamma}_{{\check{t}}}-\hat{\xi}_{{\hat{t},\check{t},A}})\nonumber\\
                      &&
                      +\varepsilon\partial_x\Upsilon^3(\check{\gamma}_{{\check{t}}}-\hat{\gamma}_{{\hat{t},\check{t},A}})
                      +\varepsilon\frac{\nu T-{\check{t}}}{\nu T}\partial_x\Upsilon^3(\check{\gamma}_{{\check{t}}})  +\partial_x\left[\sum_{i=0}^{\infty}\frac{1}{2^i}
                      \Upsilon^3(\check{\gamma}_{\check{t}}-\check{\gamma}^{i}_{\check{t}_{i},\check{t},A})
                      +\sum_{i=0}^{\infty}\frac{1}{2^i}\Upsilon^3(\check{\gamma}_{\check{t}}-\gamma^i_{t_i,\check{t},A})\right],\nonumber\\
                      &&
                                       \nabla^2_{x}(\varphi_k)({\check{t}}, (\check{\gamma}_{{\check{t}}}(\check{{t}}))_N)
                                       +4\beta^{\frac{1}{3}}Q_NQ_N
                                       +2^5\beta\partial_{xx}\Upsilon^3(\check{\gamma}_{{\check{t}}}-\hat{\xi}_{{\hat{t}},\check{t},A}) +\varepsilon\partial_{xx}\Upsilon^3(\check{\gamma}_{{\check{t}}}-\hat{\gamma}_{{\hat{t}},\check{t},A})\nonumber\\
                                      &&
                                       +\varepsilon\frac{\nu T-{\check{t}}}{\nu T}\partial_{xx}\Upsilon^3(\check{\gamma}_{{\check{t}}})+\partial_{xx}\left[\sum_{i=0}^{\infty}\frac{1}{2^i}
                      \Upsilon^3(\check{\gamma}_{\check{t}}-\check{\gamma}^{i}_{\check{t}_{i},\check{t},A})
                      +\sum_{i=0}^{\infty}\frac{1}{2^i}\Upsilon^3(\check{\gamma}_{\check{t}}-\gamma^i_{t_i,\check{t},A})\right]
                                    \bigg{)}\geq c; \ \
 \end{eqnarray}
 and
 \begin{eqnarray}\label{vis2}
                     &&-(\psi_k)_t({\check{s}},(\check{\eta}_{{\check{s}}}({\check{s}}))_N)
                     +\frac{\varepsilon}{\nu T}\Upsilon^3(\check{\eta}_{{\check{s}}})%+|\check{\eta}_{{\check{s}}}({\check{s}})|^2)
                     %-\varepsilon\frac{\nu T-{\check{s}}}{\nu T}|\check{\gamma}^{2}_{{\check{s}}}({\check{s}})|^2
                      -2\sum_{i=0}^{\infty}\frac{1}{2^i}(\check{s}-\check{s}_{i})-(A^*\nabla_{x}(\psi_k)({\check{s}}, (\check{\eta}_{{\check{s}}}(\check{{s}}))_N),\check{\eta}_{{\check{s}}}({\check{s}}))_H\nonumber\\
                     &&-2\varepsilon({\check{s}}-{\hat{t}})
                     +4\beta^{\frac{1}{3}}(A^*((\check{\eta}_{\check{s}}(\check{s}))_N-(e^{(\check{s}-\hat{t})A}\hat{\xi}_{\hat{t}}(\hat{t}))_N),
                     \check{\eta}_{{\check{s}}}(\check{{s}})%-e^{(\check{s}-\hat{t})A}\hat{\xi}_{\hat{t}}(\hat{t})
                     )_H\nonumber\\
                     &&+{\mathbf{H}}\bigg{(}\check{\eta}_{{\check{s}}}, W_2(\check{\eta}_{{\check{s}}}),
                     -\nabla_{x}(\psi_k)({\check{s}}, (\check{\eta}_{{\check{s}}}(\check{{s}}))_N)-4\beta^{\frac{1}{3}}((\check{\eta}_{\check{s}}(\check{s}))^-_N-(e^{(\check{s}-\hat{t})A}\hat{\xi}_{\hat{t}}(\hat{t}))^-_N)
                      -2^5\beta\partial_x\Upsilon^3(\check{\eta}_{{\check{s}}}-\hat{\xi}_{{\hat{t}},\check{s},A})
                     \nonumber\\
                      && -\varepsilon\partial_x\Upsilon^3(\check{\eta}_{{\check{s}}}-\hat{\eta}_{{\hat{t}},\check{s},A})
                       -\varepsilon\frac{\nu T-{\check{s}}}{\nu T}
                     \partial_x\Upsilon^3(\check{\eta}_{{\check{s}}})-\partial_x\left[\sum_{i=0}^{\infty}\frac{1}{2^i}
                      \Upsilon^3(\check{\eta}_{\check{s}}-\check{\eta}^{i}_{\check{s}_{i},\check{s},A})
                      +\sum_{i=0}^{\infty}\frac{1}{2^i}\Upsilon^3(\check{\eta}_{\check{s}}-\eta^i_{t_i,\check{s},A})\right],
                        \nonumber\\
                         && -\nabla^2_{x}(\psi_k)({\check{s}}, (\check{\eta}_{{\check{s}}}(\check{{s}}))_N)-4\beta^{\frac{1}{3}}Q_NQ_N
                      -2^5\beta\partial_{xx}\Upsilon^3(\check{\eta}_{{\check{s}}}-\hat{\xi}_{{\hat{t}},\check{s},A})
                      -\varepsilon\partial_{xx}\Upsilon^3(\check{\eta}_{{\check{s}}}-\hat{\eta}_{{\hat{t}},\check{s},A})\nonumber\\
                      &&
                      -\varepsilon\frac{\nu T-{\check{s}}}{\nu T}\partial_{xx}\Upsilon^3(\check{\eta}_{{\check{s}}})
                                    -\partial_{xx}\left[\sum_{i=0}^{\infty}\frac{1}{2^i}
                      \Upsilon^3(\check{\eta}_{\check{s}}-\check{\eta}^{i}_{\check{s}_{i},\check{s},A})
                      +\sum_{i=0}^{\infty}\frac{1}{2^i}\Upsilon^3(\check{\eta}_{\check{s}}-\eta^i_{t_i,\check{s},A})\right]
                                     \bigg{)}\leq0. \ \ \
  \end{eqnarray}
 % Here and in the sequel,
% for notational simplicity,
%we use $\partial_x\Upsilon^3(\cdot,\cdot)$ and $\partial_{xx}\Upsilon^3(\cdot,\cdot)$  to denote the first and second spatial derivatives with respect to the first variable, respectively.
  We notice that, by  the property (i) of $(\check{t},\check{\gamma}_{\check{t}},\check{s},\check{\eta}_{\check{s}})$,
  \begin{eqnarray*}
  2\sum_{i=0}^{\infty}\frac{1}{2^i}[(\check{s}-\check{s}_{i})+(\check{t}-\check{t}_{i})]
  \leq4\sum_{i=0}^{\infty}\frac{1}{2^i}\bigg{(}\frac{\delta}{2^i}\bigg{)}^{\frac{1}{2}}\leq 8\delta^{\frac{1}{2}};
    \end{eqnarray*}
    \begin{eqnarray*}
    |\partial_x\Upsilon^3(\check{\gamma}_{{\check{t}}}-\hat{\gamma}_{{\hat{t}},\check{t},A})|
    +|\partial_x\Upsilon^3(\check{\eta}_{{\check{s}}}-\hat{\eta}_{{\hat{t}},\check{s},A})|\leq 18|e^{(\check{t}-\hat{t})A}\hat{\gamma}_{{\hat{t}}}(\hat{t})-\check{\gamma}_{\check{t}}(\check{t})|^5
                      +18|e^{(\check{s}-\hat{t})A}\hat{\eta}_{{\hat{t}}}(\hat{t})-\check{\eta}_{\check{s}}(\check{s})|^5;
    \end{eqnarray*}
    \begin{eqnarray*}
    |\partial_{xx}\Upsilon^3(\check{\gamma}_{{\check{t}}}-\hat{\gamma}_{{\hat{t}},\check{t},A})|
    +|\partial_{xx}\Upsilon^3(\check{\eta}_{{\check{s}}}-\hat{\eta}_{{\hat{t}},\check{s},A})|\leq 306|e^{(\check{t}-\hat{t})A}\hat{\gamma}_{{\hat{t}}}(\hat{t})-\check{\gamma}_{\check{t}}(\check{t})|^4
                      +306|e^{(\check{s}-\hat{t})A}\hat{\eta}_{{\hat{t}}}(\hat{t})-\check{\eta}_{\check{s}}(\check{s})|^4;
    \end{eqnarray*}
     \begin{eqnarray*}
  &&\bigg{|}\partial_x\left[\sum_{i=0}^{\infty}\frac{1}{2^i}
                      \Upsilon^3(\check{\gamma}_{\check{t}}-\check{\gamma}^{i}_{\check{t}_{i},\check{t},A})\right]\bigg{|}
                      +\bigg{|}\partial_x\left[\sum_{i=0}^{\infty}\frac{1}{2^i}
                      \Upsilon^3(\check{\eta}_{\check{s}}-\check{\eta}^{i}_{\check{s}_{i},\check{s},A})\right]\bigg{|}\\
                      &\leq&18\sum_{i=0}^{\infty}\frac{1}{2^i}\left[\left|e^{(\check{t}-\check{t}_{i})A}\check{\gamma}^{i}_{\check{t}_{i}}(\check{t}_{i})-\check{\gamma}_{\check{t}}(\check{t})\right|^5
                      +\left|e^{(\check{s}-\check{s}_{i})A}\check{\eta}^{i}_{\check{s}_{i}}(\check{s}_{i})-\check{\eta}_{\check{s}}(\check{s})\right|^5\right]
                      \leq%24d^{\frac{1}{3}}\sum_{i=0}^{\infty}\frac{1}{2^i}
%                      [|\check{\gamma}^{i}_{\check{t}_{i}}(\check{t}_{i})-\check{\gamma}_{\check{t}}(\check{t})|_{6}^7
%                      +|\check{\eta}^{i}_{\check{s}_{i}}(\check{s}_{i})-\check{\eta}_{\check{s}}(\check{s})|_{6}^7]
                      36\sum_{i=0}^{\infty}\frac{1}{2^i}\bigg{(}\frac{\delta}{2^i}\bigg{)}^{\frac{5}{6}}\leq 72{\delta}^{\frac{5}{6}};
                        \end{eqnarray*}
                        and
                         \begin{eqnarray*}
  &&\bigg{|}\partial_{xx}\left[\sum_{i=0}^{\infty}\frac{1}{2^i}
                      \Upsilon^3(\check{\gamma}_{\check{t}}-\check{\gamma}^{i}_{\check{t}_{i},\check{t},A})\right]\bigg{|}+\bigg{|}\partial_{xx}\left[\sum_{i=0}^{\infty}\frac{1}{2^i}
                      \Upsilon^3(\check{\eta}_{\check{s}}-\check{\eta}^{i}_{\check{s}_{i},\check{s},A})\right]\bigg{|}\\
                      &\leq&306\sum_{i=0}^{\infty}\frac{1}{2^i}\left[\left|e^{(\check{t}-\check{t}_{i})A}\check{\gamma}^{i}_{\check{t}_{i}}(\check{t}_{i})-\check{\gamma}_{\check{t}}(\check{t})\right|^4
                      +\left|e^{(\check{s}-\check{s}_{i})A}\check{\eta}^{i}_{\check{s}_{i}}(\check{s}_{i})-\check{\eta}_{\check{s}}(\check{s})\right|^4\right]
                      \leq612\sum_{i=0}^{\infty}\frac{1}{2^i}\bigg{(}\frac{\delta}{2^i}\bigg{)}^{\frac{2}{3}}\leq 1224{\delta}^{\frac{2}{3}}.
  \end{eqnarray*}
  Combining(\ref{vis1}) and (\ref{vis2}), and letting $\delta\rightarrow0$ and $k\rightarrow\infty$,   we obtain
  \begin{eqnarray}\label{vis112}
                     c+ \frac{\varepsilon}{\nu T}(\Upsilon^3(\hat{{\gamma}}_{{\hat{t}}})+\Upsilon^3(\hat{{\eta}}_{{\hat{t}}})
                     )
                     \leq{\mathbf{H}}_1(W_1(\hat{{\gamma}}_{{\hat{t}}}))-{\mathbf{H}}_2(W_2(\hat{{\eta}}_{{\hat{t}}}))
                     +2\sum_{i=0}^{\infty}\frac{1}{2^i}(\hat{t}-t_i),
\end{eqnarray}
  where, for every $x\in \mathbb{R}$,
\begin{eqnarray*}
                       {\mathbf{H}}_1(x)%(W_1(\hat{{\gamma}}_{{\hat{t}}}))
                       &=&{\mathbf{H}}\bigg{(}\hat{{\gamma}}_{{\hat{t}}},x,
                      %W_1(\hat{{\gamma}}_{{\hat{t}}}),
                                       %6\beta |(\hat{{\gamma}}_{{\hat{t}}}({{\hat{t}}}))_N-(\hat{{\eta}}_{{\hat{t}}}({{\hat{t}}}))_N|^4((\hat{{\gamma}}_{{\hat{t}}}({{\hat{t}}}))_N-(\hat{{\eta}}_{{\hat{t}}}({{\hat{t}}}))_N)
%                                       +6\beta |(\hat{{\gamma}}_{{\hat{t}}}({{\hat{t}}}))^-_N-(\hat{{\eta}}_{{\hat{t}}}({{\hat{t}}}))^-_N|^4((\hat{{\gamma}}_{{\hat{t}}}({{\hat{t}}}))^-_N-(\hat{{\eta}}_{{\hat{t}}}({{\hat{t}}}))^-_N)\\
%                                       &&
                                       2\beta^{\frac{1}{3}}(\hat{{\gamma}}_{{\hat{t}}}({{\hat{t}}})-\hat{{\eta}}_{{\hat{t}}}({{\hat{t}}}))+2^5\beta\partial_x\Upsilon^3(\hat{\gamma}_{{\hat{t}}}-\hat{\xi}_{{\hat{t}}})
                                      +\varepsilon\frac{\nu T-{\hat{t}}}{\nu T}\partial_x\Upsilon^3(\hat{{\gamma}}_{{\hat{t}}})\nonumber\\
                                     &&
                                    +\partial_x\sum_{i=0}^{\infty}\frac{1}{2^i}\Upsilon^3(\hat{\gamma}_{\hat{t}}-\gamma^i_{t_i,\hat{t},A}),
                                     X_N+4\beta^{\frac{1}{3}}Q_NQ_N+2^5\beta\partial_{xx}\Upsilon^3(\hat{\gamma}_{{\hat{t}}}-\hat{\xi}_{{\hat{t}}})\nonumber\\
                                     &&+\varepsilon\frac{\nu T-{\hat{t}}}{\nu T}\partial_{xx}\Upsilon^3(\hat{{\gamma}}_{{\hat{t}}})
                                     +\partial_{xx}\sum_{i=0}^{\infty}\frac{1}{2^i}\Upsilon^3(\hat{\gamma}_{\hat{t}}-\gamma^i_{t_i,\hat{t},A})\bigg{)};\nonumber\\
                      {\mathbf{H}}_2(x)%{(}W_2(\hat{{\eta}}_{{\hat{t}}}))
                      &=&{\mathbf{H}}\bigg{(}\hat{{\eta}}_{{\hat{t}}},x,
                     % W_2(\hat{{\eta}}_{{\hat{t}}}), %6\beta |(\hat{{\gamma}}_{{\hat{t}}}({{\hat{t}}}))_N-(\hat{{\eta}}_{{\hat{t}}}({{\hat{t}}}))_N|^4((\hat{{\gamma}}_{{\hat{t}}}({{\hat{t}}}))_N-(\hat{{\eta}}_{{\hat{t}}}({{\hat{t}}}))_N)
%                                       +6\beta |(\hat{{\gamma}}_{{\hat{t}}}({{\hat{t}}}))^-_N-(\hat{{\eta}}_{{\hat{t}}}({{\hat{t}}}))^-_N|^4((\hat{{\gamma}}_{{\hat{t}}}({{\hat{t}}}))^-_N-(\hat{{\eta}}_{{\hat{t}}}({{\hat{t}}}))^-_N)\nonumber\\
%                     &&
2\beta^{\frac{1}{3}}(\hat{{\gamma}}_{{\hat{t}}}({{\hat{t}}})-\hat{{\eta}}_{{\hat{t}}}({{\hat{t}}}))-2^5\beta\partial_x\Upsilon^3(\hat{\eta}_{{\hat{t}}}-\hat{\xi}_{{\hat{t}}})-\varepsilon\frac{\nu T-{\hat{t}}}{\nu T}\partial_x\Upsilon^3(\hat{{\eta}}_{{\hat{t}}})\nonumber\\
                      &&-\partial_x\sum_{i=0}^{\infty}\frac{1}{2^i}
                                     \Upsilon^3(\hat{\eta}_{\hat{t}}-\eta^i_{t_i,\hat{t},A}),
                      -Y_N-4\beta^{\frac{1}{3}}Q_NQ_N-2^5\beta\partial_{xx}\Upsilon^3(\hat{\eta}_{{\hat{t}}}-\hat{\xi}_{{\hat{t}}})\nonumber\\
                      &&-\varepsilon\frac{\nu T-{\hat{t}}}{\nu T}\partial_{xx}\Upsilon^3(\hat{{\eta}}_{{\hat{t}}})
                      -\partial_{xx}\sum_{i=0}^{\infty}\frac{1}{2^i}
                      \Upsilon^3(\hat{\eta}_{\hat{t}}-\eta^i_{t_i,\hat{t},A})\bigg{)}.
\end{eqnarray*}
 On the other hand, by (\ref{5.1}) and a simple calculation we obtain
\begin{eqnarray}\label{v4}
                {\mathbf{H}}_1( W_1(\hat{{\gamma}}_{{\hat{t}}}))-{\mathbf{H}}_2( W_2(\hat{{\eta}}_{{\hat{t}}}))
                     \leq {\mathbf{H}}_1( W_2(\hat{{\eta}}_{{\hat{t}}}))-{\mathbf{H}}_2(W_2(\hat{{\eta}}_{{\hat{t}}}))
                \leq\sup_{u\in U}(J_{1}+J_{2}+J_{3}),
\end{eqnarray}
 where
\begin{eqnarray}\label{j1}
                               J_{1}&=&\bigg{(} {F}(\hat{{\gamma}}_{{\hat{t}}},u),2\beta^{\frac{1}{3}}(\hat{{\gamma}}_{{\hat{t}}}({{\hat{t}}})-\hat{{\eta}}_{{\hat{t}}}({{\hat{t}}}))
                               +2^5\beta\partial_x\Upsilon^3(\hat{\gamma}_{{\hat{t}}}-\hat{\xi}_{{\hat{t}}})
                               +\varepsilon\frac{\nu T-{\hat{t}}}{\nu T}\partial_x\Upsilon^3(\hat{{\gamma}}_{{\hat{t}}})\nonumber\\
                                     &&
                                     %+2\hat{{\gamma}}_{{\hat{t}}}({\hat{t}}))
                                     +\partial_x\sum_{i=0}^{\infty}\frac{1}{2^i}\Upsilon^3(\hat{\gamma}_{\hat{t}}-\gamma^i_{t_i,\hat{t},A})\bigg{)}_{H}  -\bigg{(} {F}(\hat{{\eta}}_{{\hat{t}}},u),2\beta^{\frac{1}{3}}(\hat{{\gamma}}_{{\hat{t}}}({{\hat{t}}})-\hat{{\eta}}_{{\hat{t}}}({{\hat{t}}}))
                                            -2^5\beta\partial_x\Upsilon^3(\hat{\eta}_{{\hat{t}}}-\hat{\xi}_{{\hat{t}}})\nonumber\\
                                     &&-\varepsilon\frac{\nu T-{\hat{t}}}{\nu T}\partial_x\Upsilon^3(\hat{{\eta}}_{{\hat{t}}})
                                     %+2\hat{{\eta}}_{{\hat{t}}}({\hat{t}}))
                                     -\partial_x\sum_{i=0}^{\infty}\frac{1}{2^i}\Upsilon^3(\hat{\eta}_{\hat{t}}-\eta^i_{t_i,\hat{t},A})\bigg{)}_{H}\nonumber\\
                                  &\leq&2\beta^{\frac{1}{3}}{L}|\hat{{\gamma}}_{{\hat{t}}}({\hat{t}})-\hat{{\eta}}_{{\hat{t}}}({\hat{t}})|
                                 ||\hat{{\gamma}}_{{\hat{t}}}-\hat{{\eta}}_{{\hat{t}}}||_0
                                 +18\beta|\hat{{\gamma}}_{{\hat{t}}}({\hat{t}})-\hat{{\eta}}_{{\hat{t}}}({\hat{t}})|^5L(2+||\hat{{\gamma}}_{{\hat{t}}}||_0
                                 +||\hat{{\eta}}_{{\hat{t}}}||_0)
                                \nonumber\\
                                            &&
                                            +18L\sum_{i=0}^{\infty}\frac{1}{2^i}\left[|e^{(\hat{t}-t_i)A}\gamma^i_{t_i}(t_i)-\hat{\gamma}_{\hat{t}}(\hat{t})|^5
                                            +|e^{(\hat{t}-t_i)A}\eta^i_{t_i}(t_i)-\hat{\eta}_{\hat{t}}(\hat{t})|^5\right]
                                            (1+||\hat{{\gamma}}_{{\hat{t}}}||_0+||\hat{{\eta}}_{{\hat{t}}}||_0)\nonumber\\
                                           && +36\varepsilon \frac{\nu T-{\hat{t}}}{\nu T} L(1+||\hat{{\gamma}}_{{\hat{t}}}||^6_0+||\hat{{\eta}}_{{\hat{t}}}||^6_0);
\end{eqnarray}
\begin{eqnarray}\label{j2}
                               J_{2}&=&\frac{1}{2}\mbox{Tr}\bigg{[} \bigg{(}X_N+4\beta^{\frac{1}{3}}Q_NQ_N+2^5\beta\partial_{xx}\Upsilon^3(\hat{\gamma}_{{\hat{t}}}-\hat{\xi}_{{\hat{t}}})+\varepsilon\frac{\nu T-{\hat{t}}}{\nu T}\partial_{xx}\Upsilon^3(\hat{{\gamma}}_{{\hat{t}}})
                               \nonumber\\
                                            &&+\partial_{xx}\left[\sum_{i=0}^{\infty}\frac{1}{2^i}\Upsilon^3(\hat{\gamma}_{\hat{t}}-\gamma^i_{t_i,\hat{t},A})\right]\bigg{)}{G}(\hat{{\gamma}}_{{\hat{t}}},u)
                                        {G}^*(\hat{{\gamma}}_{{\hat{t}}},u)\bigg{]}-\frac{1}{2}\mbox{Tr}\bigg{[} \bigg{(}-Y_N-4\beta^{\frac{1}{3}}Q_NQ_N\nonumber\\
                                            &&-2^5\beta\partial_{xx}\Upsilon^3(\hat{\eta}_{{\hat{t}}}-\hat{\xi}_{{\hat{t}}})-\varepsilon\frac{\nu T-{\hat{t}}}{\nu T}\partial_{xx}\Upsilon^3(\hat{{\eta}}_{{\hat{t}}})
                                            +\partial_{xx}\left[\sum_{i=0}^{\infty}\frac{1}{2^i}\Upsilon^3(\hat{\eta}_{\hat{t}}-\eta^i_{t_i,\hat{t},A})\right]\bigg{)}
                                        G(\hat{{\eta}}_{{\hat{t}}},u)G^*(\hat{{\eta}}_{{\hat{t}}},u)\bigg{]}\nonumber\\
                               &\leq&3
                               \beta^{\frac{1}{3}}|{G}(\hat{{\gamma}}_{{\hat{t}}},u)-G(\hat{{\eta}}_{{\hat{t}}},u)|_{L_2(\Xi,H)}^2+2\beta^{\frac{1}{3}} (|Q_N{G}(\hat{{\gamma}}_{{\hat{t}}},u)|_{L_2(\Xi,H)}^2+|Q_N{G}(\hat{{\eta}}_{{\hat{t}}},u)|_{L_2(\Xi,H)}^2)\nonumber\\
                               &&
                                      +306\beta|\hat{{\gamma}}_{{\hat{t}}}({\hat{t}})-\hat{{\eta}}_{{\hat{t}}}({\hat{t}})|^4(|{G}(\hat{{\gamma}}_{{\hat{t}}},u)|_{L_2(\Xi,H)}^2+
                                      |G(\hat{{\eta}}_{{\hat{t}}},u)|_{L_2(\Xi,H)}^2)\nonumber\\
                                      &&+153\varepsilon\frac{\nu T-{\hat{t}}}{\nu    T}(|\hat{{\gamma}}_{{\hat{t}}}({\hat{t}})|^4|{G}(\hat{{\gamma}}_{{\hat{t}}},u)|_{L_2(\Xi,H)}^2+|\hat{{\eta}}_{{\hat{t}}}({\hat{t}})|^4|{G}(\hat{{\eta}}_{{\hat{t}}},u)|_{L_2(\Xi,H)}^2)
                                      \nonumber\\
                                      &&
                                      +153\sum_{i=0}^{\infty}\frac{1}{2^i}|e^{(\hat{t}-t_i)A}\gamma^i_{t_i}(t_i)-\hat{\gamma}_{\hat{t}}(\hat{t})|^4
                                        |G(\hat{{\gamma}}_{{\hat{t}}},u)|_{L_2(\Xi,H)}^2\nonumber\\
                                            &&
                                        +153\sum_{i=0}^{\infty}\frac{1}{2^i}|e^{(\hat{t}-t_i)A}\eta^i_{t_i}(t_i)-\hat{\eta}_{\hat{t}}(\hat{t})|^4
                                        |G(\hat{{\eta}}_{{\hat{t}}},u)|_{L_2(\Xi,H)}^2
                                        \nonumber\\
                               &\leq&
                             3
                               \beta^{\frac{1}{3}}{L^2}||\hat{{\gamma}}_{{\hat{t}}}-\hat{{\eta}}_{{\hat{t}}}||_0^2+2\beta^{\frac{1}{3}} (|Q_N{G}(\hat{{\gamma}}_{{\hat{t}}},u)|_{L_2(\Xi,H)}^2+|Q_N{G}(\hat{{\eta}}_{{\hat{t}}},u)|_{L_2(\Xi,H)}^2)\nonumber\\
                                     &&+153\left(
                                     \sum_{i=0}^{\infty}\frac{1}{2^i}\left[\left|e^{(\hat{t}-t_i)A}\gamma^i_{t_i}(t_i)-\hat{\gamma}_{\hat{t}}(\hat{t})\right|^4
                                     +\left|e^{(\hat{t}-t_i)A}\eta^i_{t_i}(t_i)-\hat{\eta}_{\hat{t}}(\hat{t})\right|^4\right]\right)L^2
                                     (1+||\hat{{\gamma}}_{{\hat{t}}}||_0^2
                                             +||\hat{{\eta}}_{{\hat{t}}}||_0^2
                                             )
                                     \nonumber\\
                                             &&+ 306\beta|\hat{{\gamma}}_{{\hat{t}}}({\hat{t}})-\hat{{\eta}}_{{\hat{t}}}({\hat{t}})|^4(2+||\hat{{\gamma}}_{{\hat{t}}}||_0^2
                                             +||\hat{{\eta}}_{{\hat{t}}}||_0^2
                                             )+306\varepsilon \frac{\nu T-{\hat{t}}}{\nu T}L^2(1+||\hat{{\gamma}}_{{\hat{t}}}||^6_0
                                 +||\hat{{\eta}}_{{\hat{t}}}||^6_0) ;
\end{eqnarray}
\begin{eqnarray}\label{j3}
                                 J_{3}&=&q\bigg{(}\hat{{\gamma}}_{{\hat{t}}}, W_2(\hat{{\eta}}_{{\hat{t}}}), \bigg{(}2\beta^{\frac{1}{3}}(\hat{{\gamma}}_{{\hat{t}}}({{\hat{t}}})-\hat{{\eta}}_{{\hat{t}}}({{\hat{t}}}))
                               +2^5\beta\partial_x\Upsilon^3(\hat{\gamma}_{{\hat{t}}}-\hat{\xi}_{{\hat{t}}})
                               +\varepsilon\frac{\nu T-{\hat{t}}}{\nu T}\partial_x\Upsilon^3(\hat{{\gamma}}_{{\hat{t}}}) \nonumber\\
                                     && +\partial_x\left[\sum_{i=0}^{\infty}\frac{1}{2^i}\Upsilon^3(\hat{\gamma}_{\hat{t}}-\gamma^i_{t_i,\hat{t},A})\right]\bigg{)}G(\hat{{\gamma}}_{{\hat{t}}},u),u\bigg{)}-
                                 q\bigg{(}\hat{{\eta}}_{{\hat{t}}}, W_2(\hat{{\eta}}_{{\hat{t}}}), \bigg{(}2\beta^{\frac{1}{3}}(\hat{{\gamma}}_{{\hat{t}}}({{\hat{t}}})-\hat{{\eta}}_{{\hat{t}}}({{\hat{t}}}))
                                           \nonumber\\
                                     && -2^5\beta\partial_x\Upsilon^3(\hat{\eta}_{{\hat{t}}}-\hat{\xi}_{{\hat{t}}})-\varepsilon\frac{\nu T-{\hat{t}}}{\nu T}\partial_x\Upsilon^3(\hat{{\eta}}_{{\hat{t}}})
                                     %+2\hat{{\eta}}_{{\hat{t}}}({\hat{t}}))
                                     -\partial_x\left[\sum_{i=0}^{\infty}\frac{1}{2^i}\Upsilon^3(\hat{\eta}_{\hat{t}}-\eta^i_{t_i,\hat{t},A})\right]\bigg{)}G(\hat{{\eta}}_{{\hat{t}}},u),u\bigg{)}\nonumber\\
                                     &\leq&
                                 L||\hat{{\gamma}}_{{\hat{t}}}-\hat{{\eta}}_{{\hat{t}}}||_0
                               +2\beta^{\frac{1}{3}} L^2|\hat{{\gamma}}_{{\hat{t}}}({\hat{t}})
                                 -\hat{{\eta}}_{{\hat{t}}}({\hat{t}})|\times||\hat{{\gamma}}_{{\hat{t}}}-\hat{{\eta}}_{{\hat{t}}}||_0
                                 +18\beta L^2|\hat{{\gamma}}_{{\hat{t}}}({\hat{t}})-\hat{{\eta}}_{{\hat{t}}}({\hat{t}})|^5
                                 (2+||\hat{{\gamma}}_{{\hat{t}}}||_0
                                             +||\hat{{\eta}}_{{\hat{t}}}||_0
                                             )\nonumber\\
                               &&+18L^2\sum_{i=0}^{\infty}\frac{1}{2^i}\left[\left|\gamma^i_{t_i}(t_i)-e^{(t_i-\hat{t})A}\hat{\gamma}_{\hat{t}}(\hat{t})\right|^5
                                            +\left|\eta^i_{t_i}(t_i)-e^{(t_i-\hat{t})A}\hat{\eta}_{\hat{t}}(\hat{t})\right|^5\right]
                                            (1+||\hat{{\gamma}}_{{\hat{t}}}||_0+||\hat{{\eta}}_{{\hat{t}}}||_0)\nonumber\\
                                            &&+36\varepsilon \frac{\nu T-{\hat{t}}}{\nu T} L^2(1+||\hat{{\gamma}}_{{\hat{t}}}||_0^6
                                             +||\hat{{\eta}}_{{\hat{t}}}||_0^6
                                             );
\end{eqnarray}
 We notice that, by  the property (i) of $(\hat{t},\hat{\gamma}_{\hat{t}},\hat{\eta}_{\hat{t}})$,
  \begin{eqnarray*}
2\sum_{i=0}^{\infty}\frac{1}{2^i}(\hat{t}-t_i)
  \leq2\sum_{i=0}^{\infty}\frac{1}{2^i}\bigg{(}\frac{1}{2^i\beta}\bigg{)}^{\frac{1}{2}}\leq 4{\bigg{(}\frac{1}{{\beta}}\bigg{)}}^{\frac{1}{2}},
    \end{eqnarray*}
     \begin{eqnarray*}
  &&\sum_{i=0}^{\infty}\frac{1}{2^i}[|e^{(\hat{t}-t_i)A}\gamma^i_{t_i}(t_i)-\hat{\gamma}_{\hat{t}}(\hat{t})|^5
                                            +|e^{(\hat{t}-t_i)A}\eta^i_{t_i}(t_i)-\hat{\eta}_{\hat{t}}(\hat{t})|^5]
                      \leq
                      2\sum_{i=0}^{\infty}\frac{1}{2^i}\bigg{(}\frac{1}{2^i\beta}\bigg{)}^{\frac{5}{6}}\leq 4{\bigg{(}\frac{1}{{\beta}}\bigg{)}}^{\frac{5}{6}},
                        \end{eqnarray*}
                        and
                         \begin{eqnarray*}
  \sum_{i=0}^{\infty}\frac{1}{2^i}[|e^{(\hat{t}-t_i)A}\gamma^i_{t_i}(t_i)-\hat{\gamma}_{\hat{t}}(\hat{t})|^4
                                     +|e^{(\hat{t}-t_i)A}\eta^i_{t_i}(t_i)-\hat{\eta}_{\hat{t}}(\hat{t})|^4]
                      \leq
                      2\sum_{i=0}^{\infty}\frac{1}{2^i}\bigg{(}\frac{1}{2^i\beta}\bigg{)}^{\frac{2}{3}}\leq 4{\bigg{(}\frac{1}{{\beta}}\bigg{)}}^{\frac{2}{3}};
  \end{eqnarray*}
  and since $\hat{{\gamma}}_{{\hat{t}}}$ and $\hat{{\eta}}_{{\hat{t}}}$ are independent of $N$, by Hypothesis  \ref{hypstate5666},
  $$
  \sup_{u\in U}[|Q_N{G}(\hat{{\gamma}}_{{\hat{t}}},u)|_{L_2(\Xi,H)}^2+|Q_N{G}(\hat{{\eta}}_{{\hat{t}}},u)|_{L_2(\Xi,H)}^2]\rightarrow0\ \mbox{as}\ N\rightarrow \infty.
  $$
 Combining (\ref{vis112})-(\ref{j3}),  and letting $N\rightarrow\infty$ and then  $\beta\rightarrow\infty$,
  by (\ref{5.10jiajiaaaa}) and (\ref{5.10}),  we obtain
 \begin{eqnarray}\label{vis122}
                     c
                                             &\leq&
                      -\frac{\varepsilon}{\nu T}(\Upsilon^3(\hat{{\gamma}}_{{\hat{t}}})
                     +\Upsilon^3(\hat{{\eta}}_{{\hat{t}}}))+ \varepsilon \frac{\nu T-{\hat{t}}}{\nu T} (342L+36)L(1+||\hat{{\gamma}}_{{\hat{t}}}||_0^6
                                             +||\hat{{\eta}}_{{\hat{t}}}||_0^6
                                             ).
\end{eqnarray}
Recalling $
         \nu=1+\frac{1}{8T(342L+36)L}
$ and $\bar{a}=\frac{1}{8(342L+36)L}\wedge{T}$, by (\ref{s0}) and (\ref{5.3}),the following contradiction is induced:
\begin{eqnarray*}\label{vis122}
                     c\leq
                             \frac{\varepsilon}{4\nu
                              T}\leq \frac{c}{2}.
\end{eqnarray*}
%Letting $N\rightarrow\infty$, by (\ref{g5}), the following contradiction is induced:
%\begin{eqnarray*}\label{vis122}
%                     c\leq
%                              \frac{c}{2}.
%\end{eqnarray*}
%For every fixed $\mu>0$, there exists a constant $\beta>0$ be large enough such that
%$$
%\beta(2L+1)^2(|\hat{{\gamma}}^{1}_{{\hat{t}}}({\hat{t}})-\hat{{\gamma}}^{2}_{{\hat{t}}}({\hat{t}})|^2
%                              +||\hat{{\gamma}}^{1}_{{\hat{t}}}-\hat{{\gamma}}^{2}_{{\hat{t}}}||_0^2)+L||\hat{{\gamma}}^{1}_{{\hat{t}}}-\hat{{\gamma}}^{2}_{{\hat{t}}}||_0\leq \frac{1}{\mu}
%$$
%Then by (\ref{5.3}) and (\ref{betamu}), the following contradiction is induced:
%\begin{eqnarray*}
%                    \frac{c}{2}
%                      \leq\frac{c}{4}.
%\end{eqnarray*}
 The proof is now complete.
 \ \ $\Box$
  \par
 To complete the previous proof, it remains to state and prove the following lemmas.

 \begin{lemma}\label{lemma4.3}\ \
                 The functions $\breve{W}_1$ and $-\breve{W}_2$ defined by (\ref{definewv111}) are upper semicontinuous in $(t,x_N)\in [0,T]\times H_N$, and
                  $ \breve{W}_{1}(t,x_N)- \breve{W}_{2}(t,y_N)-{\beta|x_N-y_N|^6}$ has a maximum at $({\hat{t}},(\hat{{\gamma}}_{\hat{t}}({\hat{t}}))_N,
(\hat{{\eta}}_{\hat{t}}({\hat{t}}))_N)$
 on $[0,T]\times H_N\times H_N$. Moreover,
                 \begin{eqnarray}\label{202002020}
                     \breve{{W}}_{1}(\hat{t},(\hat{\gamma}_{\hat{t}}(\hat{t}))_N)
                     =\bar{W}_{1}(\hat{t},(\hat{\gamma}_{\hat{t}}(\hat{t}))_N),\ \ \
                     \breve{{W}}_{2}(\hat{t},(\hat{\eta}_{\hat{t}}(\hat{t}))_N)
                             =\bar{W}_{2}(\hat{t},(\hat{\eta}_{\hat{t}}(\hat{t}))_N).
\end{eqnarray}
\end{lemma}
\par
   {\bf  Proof  }. \ \
%We shall only prove the result for $\breve{W}_1$, the proof for $-\breve{W}_2$ being similar.
For every $\hat{t}\leq t\leq s\leq T$ and $x_N, y_N\in  H_N$ satisfying $|x_N|\vee|y_N|\leq M_2$ for some constant $M_2>0$, there exists a constant $C>0$ depending only on   $M_2$ such that
\begin{eqnarray*}
                             &&\bar{W}_{1}(t,x_N)-\bar{W}_{1}(s,y_N)=\bar{W}_{1}(t,x_N)-\bar{W}_{1}(s,e^{(s-t)A}x_N)+\bar{W}_{1}(s,e^{(s-t)A}x_N)-\bar{W}_{1}(s,y_N)\\
                             &=&\sup_{\gamma_t\in \Lambda^{\hat{t}},(\gamma_t(t))_N=x_N}
                             \bigg{[}W'_{1}(\gamma_t)
                             -2^5\beta \Upsilon^3(\gamma_{t},\hat{\xi}_{\hat{t}})-2\beta^{\frac{1}{3}}|(\gamma_t(t))^-_N-(e^{(t-\hat{t})A}\hat{\xi}_{\hat{t}}(\hat{t}))^-_N|^2\bigg{]}\\
                             &&-\sup_{\eta_s\in \Lambda^{\hat{t}},(\eta_s(s))_N=e^{(s-t)A}x_N}
                             \bigg{[}W'_{1}(\eta_s)
                             -2^5\beta \Upsilon^3(\eta_{s},\hat{\xi}_{\hat{t}})-2\beta^{\frac{1}{3}}|(\eta_s(s))^-_N-(e^{(s-\hat{t})A}\hat{\xi}_{\hat{t}}(\hat{t}))^-_N|^2\bigg{]}\\
                             &&+\sup_{\eta_s\in \Lambda^{\hat{t}},(\eta_s(s))_N=e^{(s-t)A}x_N}
                             \bigg{[}W'_{1}(\eta_s)
                             -2^5\beta \Upsilon^3(\eta_{s},\hat{\xi}_{\hat{t}})-2\beta^{\frac{1}{3}}|(\eta_s(s))^-_N-(e^{(s-\hat{t})A}\hat{\xi}_{\hat{t}}(\hat{t}))^-_N|^2\bigg{]}\\
                             &&-\sup_{\eta_s\in \Lambda^{\hat{t}},(\eta_s(s))_N=y_N}
                             \bigg{[}W'_{1}(\eta_s)
                             -2^5\beta \Upsilon^3(\eta_{s},\hat{\xi}_{\hat{t}})-2\beta^{\frac{1}{3}}|(\eta_s(s))^-_N-(e^{(s-\hat{t})A}\hat{\xi}_{\hat{t}}(\hat{t}))^-_N|^2\bigg{]}\\
                             &\leq& \sup_{\gamma_t\in \Lambda^{\hat{t}},(\gamma_t(t))_N=x_N}
                             \bigg{[}W'_{1}(\gamma_t)
                             -2^5\beta \Upsilon^3(\gamma_{t},\hat{\xi}_{\hat{t}})-W'_{1}(\gamma_{t,s,A})+2^5\beta \Upsilon^3(\gamma_{t,s,A},\hat{\xi}_{\hat{t}})\nonumber\\
                             &&-2\beta^{\frac{1}{3}}|\gamma_t(t)-e^{(t-\hat{t})A}\hat{\xi}_{\hat{t}}(\hat{t})|^2
                             +2\beta^{\frac{1}{3}}|x_N-(e^{(t-\hat{t})A}\hat{\xi}_{\hat{t}}(\hat{t}))_N|^2
                             +2\beta^{\frac{1}{3}}|e^{(s-\hat{t})A}\gamma_{t}(t)-e^{(s-\hat{t})A}\hat{\xi}_{\hat{t}}(\hat{t})|^2\nonumber\\
                             &&
                             -2\beta^{\frac{1}{3}}|e^{(s-t)A}x_N-(e^{(s-\hat{t})A}\hat{\xi}_{\hat{t}}(\hat{t}))_N|^2\bigg{]}\\
                              && +\sup_{\eta_s\in \Lambda^{\hat{t}},(\eta_s(s))_N=e^{(s-t)A}x_N}
                             \bigg{[}W'_{1}(\eta_s)
                             -2^5\beta \Upsilon^3(\eta_s,\hat{\xi}_{\hat{t}})-W'_{1}(\eta_s+(y_N-e^{(s-t)A}x_N){\mathbf{1}}_{[0,s]})
                             \nonumber\\
                             &&~~~~~~~~~~~~~~~~~+2^5\beta \Upsilon^3(\eta_{s}+(y_N-e^{(s-t)A}x_N){\mathbf{1}}_{[0,s]},\hat{\xi}_{\hat{t}})\bigg{]}\\
                             &\leq& C(s-t)^{\frac{1}{2}}+2\beta^{\frac{1}{3}}|x_N-(e^{(t-\hat{t})A}\hat{\xi}_{\hat{t}}(\hat{t}))_N|^2
                             -2\beta^{\frac{1}{3}}|e^{(s-t)A}x_N-(e^{(s-\hat{t})A}\hat{\xi}_{\hat{t}}(\hat{t}))_N|^2\\
                             &&+C|e^{(s-t)A}x_N-y_N|\\
                              &\leq& C((s-t)^{\frac{1}{2}}+|e^{(s-t)A}x_N-x_N|+|x_N-y_N|+|e^{(s-t)A}\hat{\xi}_{\hat{t}}(\hat{t})-\hat{\xi}_{\hat{t}}(\hat{t})|).
\end{eqnarray*}
%\begin{eqnarray*}
%                             &&\bar{W}_{1}(t,x_N)-\bar{W}_{1}(s,y_N)\\
%                             &=&\sup_{\gamma_t\in \Lambda^{\hat{t}},(\gamma_t(t))_N=x_N}
%                             \bigg{[}W'_{1}(\gamma_t)
%                             -2^5\beta \Upsilon^3(\gamma_{t},\hat{\xi}_{\hat{t}})-2\beta^{\frac{1}{3}}|(\gamma_t(t))^-_N-(e^{(t-\hat{t})A}\hat{\xi}_{\hat{t}}(\hat{t}))^-_N|^2\bigg{]}\\
%                             &&-\sup_{\eta_s\in \Lambda^{\hat{t}},(\eta_s(s))_N=y_N}
%                             \bigg{[}W'_{1}(\eta_s)
%                             -2^5\beta \Upsilon^3(\eta_{s},\hat{\xi}_{\hat{t}})-2\beta^{\frac{1}{3}}|(\eta_s(s))^-_N-(e^{(s-\hat{t})A}\hat{\xi}_{\hat{t}}(\hat{t}))^-_N|^2\bigg{]}\\
%                             &\leq& \sup_{\gamma_t\in \Lambda^{\hat{t}},(\gamma_t(t))_N=x_N}
%                             \bigg{[}W'_{1}(\gamma_t)
%                             -2^5\beta \Upsilon^3(\gamma_{t},\hat{\xi}_{\hat{t}})-W'_{1}(\gamma_{t,s,A}+(y_N-e^{(s-t)A}x_N){\mathbf{1}}_{[0,s]})
%                             \nonumber\\
%                             &&~~~~~~~~~~~~~~~~~+2^5\beta \Upsilon^3(\gamma_{t,s,A}+(y_N-e^{(s-t)A}x_N){\mathbf{1}}_{[0,s]},\hat{\xi}_{\hat{t}})\bigg{]}\\
%                             &\leq& C(|e^{(s-t)A}x_N-y_N|+(s-t)^{\frac{1}{2}}).
%\end{eqnarray*}
Clearly, if $0\leq t\leq s\leq \hat{t}$, we have
\begin{eqnarray*}
                             \bar{W}_{1}(t,x_N)-\bar{W}_{1}(s,y_N)\leq C|x_N-y_N|.
\end{eqnarray*}
and, if $0\leq t\leq  \hat{t}\leq s\leq T$ , we have
\begin{eqnarray*}
                             &&\bar{W}_{1}(t,x_N)-\bar{W}_{1}(s,y_N)\leq\bar{W}_{1}(\hat{t},x_N)-\bar{W}_{1}(s,y_N)\\
                              &\leq& C((s-\hat{t})^{\frac{1}{2}}+|e^{(s-\hat{t})A}x_N-x_N|+|x_N-y_N|+|e^{(s-\hat{t})A}\hat{\xi}_{\hat{t}}(\hat{t})-\hat{\xi}_{\hat{t}}(\hat{t})|).
\end{eqnarray*}
%\begin{eqnarray*}
%                             &&\bar{W}_{1}(t,x_0)-\bar{W}_{1}(t,y_0)\\
%                             &=&\sup_{\gamma_t\in \Lambda^{\hat{t}},\gamma_t(t)=x_0}
%                             \bigg{[}W'_{1}(\gamma_t)
%                             -2^7\beta \Upsilon^3(\gamma_{t},\hat{\xi}_{\hat{t}})\bigg{]}-\sup_{\eta_t\in \Lambda^{\hat{t}},\eta_t(t)=y_0}
%                             \bigg{[}W'_{1}(\eta_t)
%                             -2^7\beta \Upsilon^3(\eta_{t},\hat{\xi}_{\hat{t}})\bigg{]}\\
%                             &\leq& \sup_{\gamma_t\in \Lambda^{\hat{t}},\gamma_t(t)=x_0}
%                             \bigg{[}W'_{1}(\gamma_t)
%                             -2^7\beta \Upsilon^3(\gamma_{t},\hat{\xi}_{\hat{t}})-W'_{1}(\gamma_{t}+(y_0-x_0)1_{[0,t]})
%                             +2^7\beta \Upsilon^3(\gamma_{t}+(y_0-x_0)1_{[0,t]},\hat{\xi}_{\hat{t}})\bigg{]}\\
%                             &\leq& C|x_0-y_0|,
%\end{eqnarray*}
Then we have
\begin{eqnarray}\label{202002024}
                             \bar{W}_{1}(\hat{t},(\hat{\gamma}_{\hat{t}}(\hat{t}))_N)&\leq& \liminf_{s\downarrow \hat{t},y_N\rightarrow(\hat{\gamma}_{\hat{t}}(\hat{t}))_N}[\bar{W}_{1}(s,y_N)
                             +C((s-\hat{t})^{\frac{1}{2}}+|e^{(s-\hat{t})A}(\hat{\gamma}_{\hat{t}}(\hat{t}))_N-(\hat{\gamma}_{\hat{t}}(\hat{t}))_N|\nonumber\\
                             &&~~~~~~~~~~~~~~~~~~+|(\hat{\gamma}_{\hat{t}}(\hat{t}))_N-y_N|+|e^{(s-\hat{t})A}\hat{\xi}_{\hat{t}}(\hat{t})-\hat{\xi}_{\hat{t}}(\hat{t})|)]\nonumber\\
                             &=&\liminf_{s\downarrow \hat{t},y_N\rightarrow(\hat{\gamma}_{\hat{t}}(\hat{t}))_N}\bar{W}_{1}(s,y_N);
\end{eqnarray}
\begin{eqnarray*}
                             |\bar{W}_{1}(t,y_N)-\bar{W}_{1}(t,x_N)|
                             \leq C|x_N-y_N|,\ \ t\in [0,T];
\end{eqnarray*}
and
\begin{eqnarray}\label{2020020240}
                             &&\breve{W}_{1}(t,x_N)=\limsup_{s\rightarrow t}\bar{W}_{1}(s,x_N)\nonumber\\
                             &\geq& \limsup_{s\rightarrow t}[\bar{W}_{1}(t,x_N)-C((s-t)^{\frac{1}{2}}+|e^{(s-t)A}x_N-x_N|+|e^{(s-t)A}\hat{\xi}_{\hat{t}}(\hat{t})-\hat{\xi}_{\hat{t}}(\hat{t})|)]\nonumber\\
                             &=&\bar{W}_{1}(t,x_N).
\end{eqnarray}
Similarly,
\begin{eqnarray}\label{202002024jiaa}
                             -\bar{W}_{2}(\hat{t},\hat{\eta}_{\hat{t}}(\hat{t}))\leq
                            \liminf_{s\downarrow \hat{t},y_N\rightarrow(\hat{\gamma}_{\hat{t}}(\hat{t}))_N}[-\bar{W}_{2}(s,y_N)];
\end{eqnarray}
\begin{eqnarray*}\label{20200205}
                             %\breve{{W}}_{1}(\hat{t},\hat{\gamma}_{\hat{t}}(\hat{t})\geq\bar{W}_{1}(\hat{t},\hat{\gamma}_{\hat{t}}(\hat{t})),\ \ \
                             |\bar{W}_{2}(t,y_N)-\bar{W}_{2}(t,x_N)|
                             \leq C|x_N-y_N|,\  t\in [0,T];
                            % ; \ \ \
%                             -\breve{{W}}_{2}(\hat{t},(\hat{\eta}_{\hat{t}}(\hat{t}))_N)\geq-\bar{W}_{2}(\hat{t},(\hat{\eta}_{\hat{t}}(\hat{t}))_N).
\end{eqnarray*}
%Moreover, by the definitions of $\breve{W}_{1}$ and  $\breve{W}_{2}$,
%\begin{eqnarray}\label{202002024}
%                             \breve{W}_{1}(\hat{t},(\hat{\gamma}_{\hat{t}}(\hat{t}))_N)=\limsup_{s\rightarrow \hat{t}}\bar{W}_{1}(s,(\hat{\gamma}_{\hat{t}}(\hat{t}))_N)
%                             \geq\lim_{s\uparrow \hat{t}}[\bar{W}_{1}(t,(\hat{\gamma}_{\hat{t}}(\hat{t}))_N)-(\hat{t}-s)^{\frac{1}{2}}]=\bar{W}_{1}(t,(\hat{\gamma}_{\hat{t}}(\hat{t}))_N); \ \
%\end{eqnarray}
and
\begin{eqnarray}\label{20200202401}
                            % \breve{W}_{1}(t,x_N)=\limsup_{s\downarrow t}\bar{W}_{1}(s,x_N)
%                             \geq\bar{W}_{1}(t,x_N), \ \
                            -\breve{{W}}_{2}(t,y_N)\geq-\bar{W}_{2}(t,y_N). %,\ \ (t,x_N,y_N)\in [0,T]\times H_N\times H_N.
\end{eqnarray}
In particular,
\begin{eqnarray}\label{20200202482}
                             \breve{W}_{1}(\hat{t},(\hat{\gamma}_{\hat{t}}(\hat{t}))_N)
                             \geq\bar{W}_{1}(t,(\hat{\gamma}_{\hat{t}}(\hat{t}))_N); \ \ \
                            -\breve{{W}}_{2}(\hat{t},(\hat{\eta}_{\hat{t}}(\hat{t}))_N)\geq-\bar{W}_{2}(\hat{t},(\hat{\eta}_{\hat{t}}(\hat{t}))_N).
\end{eqnarray}
Therefore,
\begin{eqnarray}\label{20200202}
                    && \breve{{W}}_{1}(\hat{t},(\hat{\gamma}_{\hat{t}}(\hat{t}))_N)-\breve{{W}}_{2}(\hat{t},(\hat{\eta}_{\hat{t}}(\hat{t}))_N)
                             -\beta^{\frac{1}{3}}|(\hat{\gamma}_{\hat{t}}(\hat{t}))_N-(\hat{\eta}_{\hat{t}}(\hat{t}))_N|^2\nonumber\\
                             &\geq&\bar{W}_{1}(\hat{t},(\hat{\gamma}_{\hat{t}}(\hat{t}))_N)-\bar{W}_{2}(\hat{t},(\hat{\eta}_{\hat{t}}(\hat{t}))_N)
                             -\beta^{\frac{1}{3}}|(\hat{\gamma}_{\hat{t}}(\hat{t}))_N-(\hat{\eta}_{\hat{t}}(\hat{t}))_N|^2. \ \
\end{eqnarray}
On the other hand,  for every $(t,x_N,y_N)\in [0,T]\times H_N\times H_N$,
\begin{eqnarray*}
                            \breve{W}_{1}(t,x_N)- \breve{W}_{2}(t,y_N)-{\beta^{\frac{1}{3}}|x_N-y_N|^2}
                           %=\limsup_{s\rightarrow t}[\bar{W}_{1}(s,x_N)-\bar{W}_{2}(s,y_N)-{\beta|x_N-y_N|^6}]\\
                           \leq\sup_{(l,x,y)\in [0,T]\times H_N\times H_N}[\bar{W}_{1}(l,x)- \bar{W}_{2}(l,y)-{\beta^{\frac{1}{3}}|x-y|^2}].
\end{eqnarray*}
Thus we have
\begin{eqnarray*}
                             &&\bar{W}_{1}(\hat{t},(\hat{\gamma}_{\hat{t}}(\hat{t}))_N)-\bar{W}_{2}(\hat{t},(\hat{\eta}_{\hat{t}}(\hat{t}))_N)
                             -\beta^{\frac{1}{3}}|(\hat{\gamma}_{\hat{t}}(\hat{t}))_N-(\hat{\eta}_{\hat{t}}(\hat{t}))_N|^2\\
                             &=&\sup_{(t,x_N,y_N)\in [0,T]\times H_N\times H_N}[\bar{W}_{1}(t,x_N)- \bar{W}_{2}(t,y_N)-{\beta^{\frac{1}{3}}|x_N-y_N|^2}]\\
                             &\geq&\sup_{(t,x_N,y_N)\in [0,T]\times H_N\times H_N}[\breve{W}_{1}(t,x_N)- \breve{W}_{2}(t,y_N)-{\beta^{\frac{1}{3}}|x_N-y_N|^2}]\\
                             &\geq&
                              \breve{{W}}_{1}(\hat{t},(\hat{\gamma}_{\hat{t}}(\hat{t}))_N)-\breve{{W}}_{2}(\hat{t},(\hat{\eta}_{\hat{t}}(\hat{t}))_N)
                             -\beta^{\frac{1}{3}}|(\hat{\gamma}_{\hat{t}}(\hat{t}))_N-(\hat{\eta}_{\hat{t}}(\hat{t}))_N|^2,
                            % \geq\bar{W}_{1}(\hat{t},\hat{\gamma}_{\hat{t}}(\hat{t}))-\bar{W}_{2}(\hat{t},\hat{\eta}_{\hat{t}}(\hat{t}))
%                             -\beta|\hat{\gamma}_{\hat{t}}(\hat{t})-\hat{\eta}_{\hat{t}}(\hat{t})|_6^6.
\end{eqnarray*}
combining with (\ref{20200202482}), %(\ref{20200205}) and (\ref{20200202}),
we obtain that (\ref{202002020}) holds true, and  $ \breve{W}_{1}(t,x_N)- \breve{W}_{2}(t,y_N)-{\beta^{\frac{1}{3}}|x_N-y_N|^2}$ has a maximum at $({\hat{t}},(\hat{{\gamma}}_{\hat{t}}({\hat{t}}))_N,
(\hat{{\eta}}_{\hat{t}}({\hat{t}}))_N)$
 on $[0,T]\times H_N\times H_N$.
 \par
%Now, let us prove $\breve{{W}}_{1}$ is an upper semicontinuous function.
We shall only prove  that the function $\breve{W}_1$ is an upper semicontinuous function in $(t,x_N)\in [0,T]\times H_N$, the proof for $-\breve{W}_2$ being similar. For every $(t,x_N)\in [0,T]\times H_N$, by the definition of $\breve{{W}}_{1}$,
%\begin{eqnarray}\label{jiahou1}
%                    &&\overline{\lim}_{s\uparrow t}\breve{{W}}_{1}(s,x_N)=\overline{\lim}_{s\uparrow t}\overline{\lim}_{l\downarrow s}\bar{{W}}_{1}(l,x_N)
%                    \leq\overline{\lim}_{s\uparrow t}\overline{\lim}_{l\downarrow s}[\bar{{W}}_{1}(t,x_N)+C(t-l)^{\frac{1}{2}}]\nonumber\\
%                    &=&\bar{{W}}_{1}(t,x_N)\leq
%                    \overline{\lim_{s\downarrow t}}[\bar{{W}}_{1}(s,x_N)+C(s-t)^{\frac{1}{2}}]=\breve{{W}}_{1}(t,x_N).
%\end{eqnarray}
%Moreover,
for every $\varepsilon>0$, there exists a constant $\delta>0$ such that
$$
\breve{{W}}_{1}(t,x_N)\geq \bar{W}_{1}(s,x_N)-\varepsilon \ \ \mbox{for all}\ \ s\in[0\vee(t-\delta),(t+\delta)\wedge T].
$$
Then
\begin{eqnarray*}
                   \breve{{W}}_{1}(t,x_N)\geq\limsup_{s\rightarrow t}\left(%\bar{{W}}_{1}(s,x_N)\vee
                   \limsup_{l\rightarrow s}\bar{{W}}_{1}(l,x_N)\right)-\varepsilon
                   =\limsup_{s\rightarrow t}{\breve{{W}}}_{1}(s,x_N)-\varepsilon.
\end{eqnarray*}
 By the arbitrariness of   $\varepsilon>0$, %combining with (\ref{jiahou1}),
 we obtain that
 $$
 \breve{{W}}_{1}(t,x_N)\geq \limsup_{s\rightarrow t}{\breve{{W}}}_{1}(s,x_N), \ \ \mbox{for all}\ \ (t,x_N)\in[0, T]\times H_N.
 $$
 Therefore,
 \begin{eqnarray*}
                             &&\limsup_{(s,y_N)\rightarrow (t,x_N)}\breve{W}_{1}(s,y_N)
                             =\limsup_{(s,y_N)\rightarrow (t,x_N)}[\breve{W}_{1}(s,y_N)-\breve{W}_{1}(s,x_N)+\breve{W}_{1}(s,x_N)]\\
                             &\leq&\limsup_{(s,y_N)\rightarrow (t,x_N)}C|x_N-y_N|+\breve{W}_{1}(t,x_N)=\breve{W}_{1}(t,x_N), \ \
                             (t,x_N)\in [0,T]\times H_N.
\end{eqnarray*}
The proof is now complete. \ \ $\Box$
 \begin{lemma}\label{lemma4.344}\ \
                 The
                          functions $\breve{W}_{1}$ and $-\breve{W}_{2}$ defined in (\ref{definewv}) satisfy condition (8.5) of Theorem 8.3 in \cite{cran2}. % when $\mu>0$ be large enough.
\end{lemma}
\par
   {\bf  Proof  }. \ \  We only prove  $\breve{W}_{1}$  satisfies condition (8.5) of Theorem 8.3 in \cite{cran2}.
                           The same result for   $-\breve{W}_{2}$ can be obtained by a symmetric  way.
                           \par
 Set $r=\frac{1}{2}|T-{\hat{t}}|%\wedge(\frac{M_0}{2}-|\hat{\gamma}_{{\hat{t}}}({\hat{t}})|))
 $, for a given $L>0$, let $\varphi\in C^{1,2}((T-\bar{a},T)\times H_N)$ be a function such that
            $\breve{W}_{1}(t,x_N)-\varphi(t,x_N)$  has a  maximum at $(\bar{{t}},\bar{x}_{N})\in (T-\bar{a}, T)\times H_N$, moreover, the following inequalities hold true:
\begin{eqnarray*}
                    &&|\bar{{t}}-{\hat{t}}|+|\bar{x}_{N}-\hat{{\gamma}}_{{\hat{t}}}({\hat{t}})|<r,\\
                    &&|\breve{W}_{1}(\bar{{t}},\bar{{x}}_{N})|+|\nabla_x\varphi(\bar{{t}},\bar{{x}}_{N})|
                    +|\nabla^2_x\varphi(\bar{{t}},\bar{{x}}_{N})|\leq L.
\end{eqnarray*}
%Here and in the sequel, $\nabla_x$ and $\nabla^2_{x}$ are the classical partial derivatives in the state variable $x$.
We can modify $\varphi$  such that $\varphi\in C^{1,2}_p((T-\bar{a},T)\times H_N)$, $\breve{W}_1(t,x_N)-\varphi(t,x_N)$  has a strict   maximum at $(\bar{{t}},\bar{x}_N)\in (T-\bar{a}, T)\times H_N$ and the above two inequalities hold true.
If $\bar{{t}}<\hat{t}$, we have $ b=\varphi_{t}(\bar{{t}},\bar{x}_{N})=\frac{1}{2}(\hat{t}-\bar{{t}})^{-\frac{1}{2}} \geq 0$.
%We can modify $\varphi$ such that $\varphi\in C^{1,2}([{\hat{t}_\mu},T]\times H_N)$ and $\bar{W}_{1,\mu}(t,x_0)-\varphi(t,x_0)$  has a strict positive  maximum at $(\bar{{t}},\bar{x}_0)\in ({\hat{t}_\mu}, T)\times H_N$ and the above two inequalities hold true.
If $\bar{{t}}\geq \hat{t}$, we consider the functional
$$
                \Gamma(\gamma_t)= W'_{1}(\gamma_t)-2^5\beta \Upsilon^3(\gamma_{t},\hat{{\xi}}_{{\hat{t}}})
                 -\varphi(t,(\gamma_t(t))_N)-2\beta^{\frac{1}{3}}|(\gamma_t(t))^-_N-(e^{(t-\hat{t})A}\hat{{\xi}}_{{\hat{t}}}(\hat{t}))^-_N|^2,\ (t,\gamma_t)\in [\hat{t},T]\times {{\Lambda}}.
$$
We may assume that $\varphi$ grow quadratically at $\infty$.  By Lemma \ref{theoremleft} we have that,
 for every  $(t_0,\gamma^0_{t_0})\in [\hat{t},T]\times \Lambda^{\hat{t}}$ satisfy
$$
\Gamma(\gamma^0_{t_0})\geq \sup_{(s,\gamma_s)\in [\hat{t},T]\times \Lambda^{\hat{t}}}\Gamma(\gamma_s)-\delta
\geq\breve{W}_{1}(\bar{{t}},\bar{x}_{N})-\varphi(\bar{{t}},\bar{x}_{N})-\delta,
%\    \mbox{and} \ \ \Psi(\gamma^0_{t_0},\eta^0_{t_0})\geq \Psi(\tilde{\gamma}_{\tilde{t}},\tilde{\gamma}_{\tilde{t}}) >\frac{\titlde{m}}{2},
 $$
  there exist $(\breve{t},\breve{\gamma}_{\breve{t}})\in [\hat{t},T]\times \Lambda^{\hat{t}}$ and a sequence $\{(\breve{t}_i,\breve{\gamma}^i_{\breve{t}_i})\}_{i\geq1}\subset
  [\hat{t},T]\times \Lambda^{\hat{t}}$ such that
  \begin{description}
        \item{(i)} $\overline{\Upsilon}^3(\gamma^0_{t_0},\breve{\gamma}_{\breve{t}})\leq \delta$,
         $\overline{\Upsilon}^3(\breve{\gamma}^i_{\breve{t}_i},\breve{\gamma}_{\breve{t}})\leq \frac{\delta}{2^i}$ and $t_i\uparrow \breve{t}$ as $i\rightarrow\infty$,
        \item{(ii)}  $\Gamma(\breve{\gamma}_{\breve{t}})
            -\sum_{i=0}^{\infty}\frac{1}{2^i}\overline{\Upsilon}^3(\breve{\gamma}^i_{\breve{t}_i},\hat{\gamma}_{\hat{t}})
        \geq \Gamma(\gamma^0_{t_0})$, and
        \item{(iii)}    for all $(s,\gamma_s)\in [\breve{t},T]\times \Lambda^{\breve{t}}\setminus \{(\breve{t},\breve{\gamma}_{\breve{t}})\}$,
        \begin{eqnarray*}
        \Gamma(\gamma_s)
        -\sum_{i=0}^{\infty}
        \frac{1}{2^i}\overline{\Upsilon}^3(\breve{\gamma}^i_{\breve{t}_i},\gamma_s)
            <\Gamma(\breve{\gamma}_{\breve{t}})
            -\sum_{i=0}^{\infty}\frac{1}{2^i}\overline{\Upsilon}^3(\breve{\gamma}^i_{\breve{t}_i},\breve{\gamma}_{\breve{t}}).
        \end{eqnarray*}
        \end{description}
%            For every $\delta>0$, by the definition of $\tilde{W}_1$,  there exists $\gamma_{\bar{t}}^{\delta}\in \Lambda_{M_0}(\hat{\gamma}^1_{\hat{t}})$ such that
%\begin{eqnarray}
%            \gamma_{\bar{t}}^{\delta}(\bar{t})=\bar{x}_0,\
%            \ \ W'_1(\bar{t},\gamma_{\bar{t}}^{\delta},\bar{x}_0)-2^7\beta d_1(\gamma_{\bar{t}}^{\delta},\hat{{\xi}}_{{\hat{t}}})\geq \bar{W}_1(\bar{t},\bar{x}_0)-\delta.
%\end{eqnarray}
%  By Lemma \ref{lemmaleft}, there exists
%         $\check{\gamma}_{\check{t}}\in \Lambda_{M_0}(\gamma^{\delta}_{\bar{t}})$
%         such that
%\begin{eqnarray*}
%                 \Upsilon(\check{\gamma}_{\check{t}})\geq\Upsilon(\gamma^{\delta}_{\bar{t}}), \ \ \gamma_t\in \Lambda_{M_0}(\check{\gamma}_{\check{t}}).
%\end{eqnarray*}
 Then
\begin{eqnarray*}
                 &&\breve{W}_1(\breve{t},(\breve{\gamma}_{\breve{t}}(\breve{t}))_N)
                   -\varphi(\breve{t},(\breve{\gamma}_{\breve{t}}(\breve{t}))_N)
                   \geq\bar{W}_1(\breve{t},(\breve{\gamma}_{\breve{t}}(\breve{t}))_N)
                   -\varphi(\breve{t},(\breve{\gamma}_{\breve{t}}(\breve{t}))_N)\\
                   &
                 \geq& W'_1(\breve{\gamma}_{\breve{t}})-2^5\beta \Upsilon^3(\breve{\gamma}_{\breve{t}},\hat{{\xi}}_{{\hat{t}}})
                 -2\beta^{\frac{1}{3}}|(\breve{\gamma}_{\breve{t}}(\breve{t}))^-_N-(e^{(\breve{t}-\hat{t})A}\hat{{\xi}}_{{\hat{t}}}(\hat{t}))^-_N|^2
                -\varphi(\breve{t},(\breve{\gamma}_{\breve{t}}(\breve{t}))_N)\\
                  &\geq& \Gamma(\gamma^0_{t_0})
                 \geq\breve{W}_{1}(\bar{{t}},\bar{x}_{N})-\varphi(\bar{{t}},\bar{x}_{N})-\delta.
               %&\geq&\bar{W}_1(\bar{t},\bar{x}_0) -\varphi(\bar{t},\bar{x}_0)-\delta.
\end{eqnarray*}
Letting $\delta\rightarrow0$, we obtain
\begin{eqnarray}\label{delta0}
       && \breve{t}\rightarrow \bar{t},\ (\breve{\gamma}_{\breve{t}}(\breve{t}))_N\rightarrow \bar{x}_N \ \mbox{as}\ \delta\rightarrow0;
       \end{eqnarray}
 %By Lemma \ref{lemmaleft} and the definition of  $\bar{W}_{1}$, we may assume that
%$$
%                 \Upsilon(\bar{\gamma}_{\bar{t}})\geq\Upsilon(\gamma_t),
%                 \ \  \gamma_t\in  \Lambda_{M_0}(\bar{\gamma}_{\bar{t}}),
%$$
%with equality at $(\bar{t},\bar{\gamma}_{\bar{t}})$ for some $\bar{t}\in  [\hat{t},T]$ and  $\bar{\gamma}_{\bar{t}}\in \Lambda_{M_0}(\hat{{\gamma}}_{\hat{t}})$ satisfying $\bar{\gamma}_{\bar{t}}(\bar{t})=\bar{x}_{0}$.
%\par
Since $\bar{t}\leq {\hat{t}}+\frac{|T-{\hat{t}}|}{2}$, %and $|\bar{x}_{0}|\leq|\hat{{\gamma}}^{1}_{{\hat{t}}}({\hat{t}})|+\frac{\frac{M_0}{2}-|\hat{\gamma}_{{\hat{t}}}({\hat{t}})|}{2}$,
 we get $\bar{t}<T$. %and $|\bar{\gamma}_{\bar{t}}(\bar{t})|=|\bar{x}_{0}|<\frac{M_0}{2}$.
 Then, by (\ref{delta0}),  we have  $\breve{t}<T$ %and $|\breve{\gamma}_{\breve{t}}(\breve{t})|<\frac{M_0}{2}$
  provided by  $\delta>0$ be small enough.
      Thus,  the definition of the viscosity subsolution can be used to obtain the following result:
      \begin{eqnarray}\label{5.15}
                      &&\varphi_t({\breve{t}},(\breve{\gamma}_{{\breve{t}}}({\breve{t}}))_N)
                      -\frac{\varepsilon}{\nu T}\Upsilon^3(\breve{\gamma}_{{\breve{t}}})
                     %-\varepsilon({\breve{t}}-T+\bar{a})^{-2}
                     +2\sum_{i=0}^{\infty}\frac{1}{2^i}[(\breve{t}-{t_i})+(\breve{t}-\breve{t}_i)]+(A^*\nabla_{x}\varphi({\breve{t}}, (\breve{\gamma}_{{\breve{t}}}(\breve{{t}}))_N),\breve{\gamma}_{{\breve{t}}}({\breve{t}}))_H\nonumber\\
                     &&+2\varepsilon({\breve{t}}-{\hat{t}})
                     -4\beta(A^*((\breve{\gamma}_{\breve{t}}(\breve{t}))_N-(e^{(\breve{t}-\hat{t})A}\hat{\xi}_{\hat{t}}(\hat{t}))_N),
                     \breve{\gamma}_{{\breve{t}}}(\breve{{t}}))_H
                     +{\mathbf{H}}\bigg{(}\breve{\gamma}_{{\breve{t}}},  W_1(\breve{\gamma}_{{\breve{t}}}), \nonumber\\
                      && \nabla_{x}\varphi({\breve{t}}, (\breve{\gamma}_{{\breve{t}}}(\breve{{t}}))_N)
                      +2\beta^{\frac{1}{3}}\partial_x|(\breve{\gamma}_{\breve{t}}(\breve{t}))^-_N-(e^{(\breve{t}-\hat{t})A}\hat{\xi}_{\hat{t}}(\hat{t}))^-_N|^2
                     +2^5\beta\partial_x\Upsilon^3(\breve{\gamma}_{{\breve{t}}}-\hat{\xi}_{{\hat{t}},\breve{t},A})
                      +\varepsilon\partial_x\Upsilon^3(\breve{\gamma}_{{\breve{t}}}-\hat{\gamma}_{{\hat{t}},\breve{t},A})\nonumber\\
                      &&
                      +\varepsilon\frac{\nu T-{\breve{t}}}{\nu T}\partial_x\Upsilon^3(\breve{\gamma}_{{\breve{t}}})
                      +\partial_x\left[\sum_{i=0}^{\infty}\frac{1}{2^i}\Upsilon^3(\breve{\gamma}_{\breve{t}}-\breve{\gamma}^i_{\breve{t}_i,\breve{t},A})
                      +\sum_{i=0}^{\infty}\frac{1}{2^i}\Upsilon^3(\breve{\gamma}_{\breve{t}}-{\gamma}^i_{{t}_i,\breve{t},A})\right],
                                  \nabla^2_{x}\varphi({\breve{t}}, (\breve{\gamma}_{{\breve{t}}}(\breve{{t}}))_N)\nonumber\\
                      &&
                                  +2\beta^{\frac{1}{3}}\partial_{xx}|(\breve{\gamma}_{\breve{t}}(\breve{t}))^-_N-(e^{(\breve{t}-\hat{t})A}\hat{\xi}_{\hat{t}}(\hat{t}))^-_N|^2
                                       +2^5\beta\partial_{xx}\Upsilon^3(\breve{\gamma}_{{\breve{t}}}-\hat{\xi}_{{\hat{t}},\breve{t},A})
                                       +\varepsilon\partial_{xx}\Upsilon^3(\breve{\gamma}_{{\breve{t}}}-\hat{\gamma}_{{\hat{t}},\breve{t},A})\nonumber\\
                                     &&
                                       +\varepsilon\frac{\nu T-{\breve{t}}}{\nu T}\partial_{xx}\Upsilon^3(\breve{\gamma}_{{\breve{t}}})
                                     +\partial_{xx}\left[\sum_{i=0}^{\infty}\frac{1}{2^i}\Upsilon^3(\breve{\gamma}_{\breve{t}}-\breve{\gamma}^i_{\breve{t}_i,\breve{t},A})
                                     +\sum_{i=0}^{\infty}\frac{1}{2^i}\Upsilon^3(\breve{\gamma}_{\breve{t}},{\gamma}^i_{{t}_i,\breve{t},A})\right]
                                    \bigg{)}\geq c.
 \end{eqnarray}
Letting $\delta\rightarrow0$, by the definition of ${\mathbf{H}}$, it follows that there exists a constant $\bar{c}$ such that $ b=\varphi_{t}(\bar{t},\bar{x}_{N}) \geq \bar{c}$.
The proof is now complete. \ \ $\Box$
\begin{lemma}\label{lemma4.4}\ \ The maximum points $(\check{t},\check{\gamma}_{\check{t}},\check{s},\check{\eta}_{\check{s}})$
of $\Gamma_{k}(\gamma_s,\eta_s)
        -\sum_{i=0}^{\infty}
        \frac{1}{2^i}[\overline{\Upsilon}^3(\check{\gamma}^{i}_{\check{t}_{i}},\gamma_s)+\overline{\Upsilon}^3(\check{\eta}^{i}_{\check{t}_{i}},\eta_s)]$ defined by (\ref{4.1111}) in $[\check{t},T]\times \Lambda^{\check{t}}\times[\check{s},T]\times \Lambda^{\check{s}}$ satisfy  conditions (\ref{4.22}) and (\ref{4.23}).
%\begin{eqnarray}\label{4.22}
%                         \check{t}\rightarrow t_k,\ \check{\gamma}_{\check{t}}(\check{t})\rightarrow x^k_0,\ \check{s}\rightarrow s_k,\ \check{\eta}_{\check{s}}(\check{s})\rightarrow y^k_0,\ \mbox{as}\ \delta\rightarrow0;
%\end{eqnarray}
%and
%\begin{eqnarray}\label{4.23}
%\lim_{k\rightarrow\infty}\lim_{\delta\rightarrow0}[\overline{\Upsilon}^3(\check{\gamma}_{\check{t}},\hat{\gamma}_{\hat{t}})
%+\overline{\Upsilon}^3(\check{\eta}_{\check{s}},\hat{\eta}_{\hat{t}})]=0.
%\end{eqnarray}
\end{lemma}
\par
   {\bf  Proof  }. \ \
By   (\ref{2020020240}), (\ref{20200202401}) and the definitions of $\bar{W}_{1}$ and $\bar{W}_{2}$, we get that
\begin{eqnarray*}
                 &&\breve{W}_1(\check{t},(\check{\gamma}_{\check{t}}(\check{t}))_N)-  \breve{W}_2(\check{s},(\check{\eta}_{\check{s}}(\check{s}))_N)
                   -\varphi_k(\check{t},(\check{\gamma}_{\check{t}}(\check{t}))_N)
                            -\psi_k(\check{s},(\check{\eta}_{\check{s}}(\check{s}))_N)  \\
                            &\geq&\bar{W}_1(\check{t},(\check{\gamma}_{\check{t}}(\check{t}))_N)-  \bar{W}_2(\check{s},(\check{\eta}_{\check{s}}(\check{s}))_N)
                   -\varphi_k(\check{t},(\check{\gamma}_{\check{t}}(\check{t}))_N)
                            -\psi_k(\check{s},(\check{\eta}_{\check{s}}(\check{s}))_N)  \\
                 &\geq& W'_1(\check{\gamma}_{\check{t}})- {W}'_2(\check{\eta}_{\check{s}})-2^5\beta(\Upsilon^3(\check{\gamma}_{\check{t}},\hat{{\xi}}_{{\hat{t}}})
                +\Upsilon^3(\check{\eta}_{\check{s}},\hat{{\xi}}_{{\hat{t}}}))
                -\varphi_k(\check{t},(\check{\gamma}_{\check{t}}(\check{t}))_N)
                            -\psi_k(\check{s},(\check{\eta}_{\check{s}}(\check{s}))_N)  \\
                            &&-2\beta^{\frac{1}{3}}(|(\check{\gamma}_{\check{t}}(\check{t}))^-_N-(e^{(\check{t}-\hat{t})A}\hat{\xi}_{\hat{t}}(\hat{t}))^-_N|^2
                +|(\check{\eta}_{\check{s}}(\check{s}))^-_N-(e^{(\check{s}-\hat{t})A}\hat{\xi}_{\hat{t}}(\hat{t}))^-_N|^2)\\
                  &\geq&\Gamma_{k}(\check{\gamma}^{0}_{\check{t}_{0}},\check{\eta}^{0}_{\check{s}_{0}})\geq \sup_{(t,\gamma_t),(s,\eta_s)\in [\hat{t},T]\times \Lambda^{\hat{t}}}\Gamma_{k}(\gamma_t,\eta_s)-{\delta}\\
                  &\geq&\breve{W}_1(l_{k},{x}^k_N)- \breve{{W}}_2(s_k,{y}^k_N) -\varphi_k(l_{k},{x}^k_N)-\psi_k(s_k,{y}^k_N)-\delta.
\end{eqnarray*}
Letting $\delta\rightarrow0$, we obtain (\ref{4.22})  and  % by (\ref{wv}) %, (\ref{202002024}), (\ref{202002023}), (\ref{202002024jiaa})  and  (\ref{20200206jinjia})
 %we show that
%\begin{eqnarray}\label{delta}
%       && \check{t}\rightarrow t_k,\ \check{\gamma}_{\check{t}}(\check{t})\rightarrow x^k_0,\ \check{s}\rightarrow s_k,\ \check{\eta}_{\check{s}}(\check{s})\rightarrow y^k_0,\ \mbox{as}\ \delta\rightarrow0;
%       \end{eqnarray}
%{and}
\begin{eqnarray*}
&&\lim_{\delta\rightarrow0}{[}W'_1(\check{\gamma}_{\check{t}})- {W}'_2(\check{\eta}_{\check{s}})-2^5\beta(\Upsilon^3(\check{\gamma}_{\check{t}},\hat{{\xi}}_{{\hat{t}}})
                +\Upsilon^3(\check{\eta}_{\check{s}},\hat{{\xi}}_{{\hat{t}}}))
                -\varphi_{k}(\check{t},(\check{\gamma}_{\check{t}}(\check{t}))_N)
                            \nonumber \\
                            &&~~~~~-\psi_{k}(\check{s},(\check{\eta}_{\check{s}}(\check{s}))_N)
                            -2\beta^{\frac{1}{3}}(|(\check{\gamma}_{\check{t}}(\check{t}))^-_N-(e^{(\check{t}-\hat{t})A}\hat{{\xi}}_{{\hat{t}}}(\hat{t}))^-_N|^2
                +|(\check{\eta}_{\check{s}}(\check{s}))^-_N-(e^{(\check{s}-\hat{t})A}\hat{{\xi}}_{{\hat{t}}}(\hat{t}))^-_N|^2){]}\\
                 &=& \breve{W}_1(l_{k},{x}^k_N)- \breve{{W}}_2(s_k,{y}^k_N) -\varphi_k(l_{k},{x}^k_N)-\psi_k(s_k,{y}^k_N).
\end{eqnarray*}
    Letting $\delta\rightarrow0$ and $k\rightarrow\infty$, by (\ref{wv}) and (\ref{202002024})-(\ref{20200202401}),  we show that
\begin{eqnarray}\label{w1v1}
&&\liminf_{k\rightarrow\infty}\lim_{\delta\rightarrow0}{[}W'_1(\check{\gamma}_{\check{t}})- {W}'_2(\check{\eta}_{\check{s}})-2^5\beta(\Upsilon^3(\check{\gamma}_{\check{t}},\hat{{\xi}}_{{\hat{t}}})
                +\Upsilon^3(\check{\eta}_{\check{s}},\hat{{\xi}}_{{\hat{t}}}))
                            \nonumber \\
                            &&~~~~~
                            -2\beta^{\frac{1}{3}}(|(\check{\gamma}_{\check{t}}(\check{t}))^-_N-(e^{(\check{t}-\hat{t})A}\hat{{\xi}}_{{\hat{t}}}(\hat{t}))^-_N|^2
                +|(\check{\eta}_{\check{s}}(\check{s}))^-_N-(e^{(\check{s}-\hat{t})A}\hat{{\xi}}_{{\hat{t}}}(\hat{t}))^-_N|^2){]}\nonumber\\
                 &=& \liminf_{k\rightarrow\infty}{[}\breve{W}_1(l_{k},{x}^k_N)- \breve{{W}}_2(s_k,{y}^k_N) {]}
                 \geq\bar{W}_1(\hat{t},(\hat{\gamma}_{\hat{t}}(\hat{t}))_N)- \bar{{W}}_2(\hat{t},(\hat{\eta}_{\hat{t}}(\hat{t}))_N)\nonumber\\
                 &=&W'_{1}(\hat{{\gamma}}_{{\hat{t}}})-W'_{2}(\hat{{\eta}}_{{\hat{t}}})-\beta \Upsilon^3(\hat{{\gamma}}_{{\hat{t}}},\hat{{\eta}}_{{\hat{t}}})
                 -\beta^{\frac{1}{3}}|(\hat{{\gamma}}_{{\hat{t}}}(\hat{t}))^-_N-(\hat{{\eta}}_{{\hat{t}}}(\hat{t}))^-_N|^2.
\end{eqnarray}
On the other hand,  without loss of generality, we may assume  $\check{s}\leq \check{t}$,
\begin{eqnarray*}
                      &&W'_1(\check{\gamma}_{\check{t}})- {W}'_2(\check{\eta}_{\check{s}})-2^5\beta(\Upsilon^3(\check{\gamma}_{\check{t}},\hat{{\xi}}_{{\hat{t}}})
                +\Upsilon^3(\check{\eta}_{\check{s}},\hat{{\xi}}_{{\hat{t}}}))\\
                &&-2\beta^{\frac{1}{3}}(|(\check{\gamma}_{\check{t}}(\check{t}))^-_N-(e^{(\check{t}-\hat{t})A}\hat{{\xi}}_{{\hat{t}}}(\hat{t}))^-_N|^2
                +|(\check{\eta}_{\check{s}}(\check{s}))^-_N-(e^{(\check{s}-\hat{t})A}\hat{{\xi}}_{{\hat{t}}}(\hat{t}))^-_N|^2)\\
                      %&=&W(\check{{t}},\check{{x}})-\varepsilon\frac{\nu T-\check{t}}{\nu T}(|\check{x}(0)|^2+|\check{x}|_H^2)-\varepsilon|\check{t}-\hat{t}|^2-\varepsilon|{\check{x}}- \hat{x}_{\check{t}-\hat{t}}|^2_H-\varepsilon|\check{x}(0)-\hat{x}(0)|^4\\
%                      &&-V(\check{{s}},\check{{y}})-\varepsilon\frac{\nu T-\check{s}}{\nu T}(|\check{y}(0)|^2+|\check{y}|_H^2)-\varepsilon|\check{s}-\hat{t}|^2-\varepsilon|{\check{y}}- \hat{y}_{\check{s}-\hat{t}}|^2_H-\varepsilon|\check{y}(0)-\hat{y}(0)|^4\\
%                        && -\alpha|\check{x}-\hat{\xi}_{\check{t}-\hat{t}}|_H^2-\alpha|\check{y}-\hat{\xi}_{\check{s}-\hat{t}}|_H^2-\frac{\varepsilon}{\check{t}-T+\bar{a}}\\
                      &\leq&\Psi_1(\check{\gamma}_{\check{t}},\check{\eta}_{\check{s},\check{t},A})  -W_2(\check{\eta}_{\check{s}})
                      +W_2(\check{\eta}_{\check{s},\check{t},A})
                       -\varepsilon[\overline{\Upsilon}^3(\check{\gamma}_{\check{t}},\hat{\gamma}_{\hat{t}})
+\overline{\Upsilon}^3(\check{\eta}_{\check{s},\check{t},A},\hat{\eta}_{\hat{t}})]\\
&&+2\beta^{\frac{1}{3}}(|(e^{(\check{t}-\check{s})A}\check{\eta}_{\check{s}}(\check{s}))^-_N-(e^{(\check{t}-\hat{t})A}\hat{{\xi}}_{{\hat{t}}}(\hat{t}))^-_N|^2-|(\check{\eta}_{\check{s}}(\check{s}))^-_N-(e^{(\check{s}-\hat{t})A}\hat{{\xi}}_{{\hat{t}}}(\hat{t}))^-_N|^2)\\
&&+\beta^{\frac{1}{3}}|(\check{\gamma}_{\check{t}}(\check{t}))_N-(e^{(\check{t}-\check{s})A}\check{\eta}_{\check{s}}(\check{s}))_N|^2\\
                       &\leq&\Psi_1(\hat{\gamma}_{\hat{t}},\hat{\eta}_{\hat{t}})
                    +2L(1+||\check{\eta}_{\check{s}}||_0)(\check{t}- \check{s})^{\frac{1}{2}}
                      -\varepsilon[\Upsilon^3(\check{\gamma}_{\check{t}},\hat{\gamma}_{\hat{t}})
+\Upsilon^3(\check{\eta}_{\check{s},\check{t},A},\hat{\eta}_{\hat{t}})+2|\check{s}-\hat{t}|^2]\\
&&+2\beta^{\frac{1}{3}}(|(e^{(\check{t}-\check{s})A}\check{\eta}_{\check{s}}(\check{s}))^-_N-(e^{(\check{t}-\hat{t})A}\hat{{\xi}}_{{\hat{t}}}(\hat{t}))^-_N|^2-|(\check{\eta}_{\check{s}}(\check{s}))^-_N-(e^{(\check{s}-\hat{t})A}\hat{{\xi}}_{{\hat{t}}}(\hat{t}))^-_N|^2)\\
&&+\beta^{\frac{1}{3}}|(\check{\gamma}_{\check{t}}(\check{t}))_N-(e^{(\check{t}-\check{s})A}\check{\eta}_{\check{s}}(\check{s}))_N|^2.
\end{eqnarray*}
% if $\check{s}> \check{t}$,
%\begin{eqnarray*}
%                      &&W'_1(\check{\gamma}_{\check{t}})- {W}'_2(\check{\eta}_{\check{s}})-2^7\beta(\Upsilon^3(\check{\gamma}_{\check{t}},\hat{{\xi}}_{{\hat{t}}})
%                +\Upsilon^3(\check{\eta}_{\check{s}},\hat{{\xi}}_{{\hat{t}}}))-\beta |\check{\gamma}_{\check{t}}(\check{t})-\check{\eta}_{\check{s}}(\check{s})|_6^6\\
%                      &\leq&\Psi_1({\check{\gamma}}_{\check{t},\check{s}},{\check{\eta}}_{\check{s}}) +W_1(\check{\gamma}_{\check{t}})
%                      -W_1(\check{\gamma}_{\check{t},\check{s}})
%                       -\varepsilon[\overline{\Upsilon}^3(\check{\gamma}_{\check{t},\check{s}},\hat{\gamma}_{\hat{t}})
%+\overline{\Upsilon}^3(\check{\eta}_{\check{s}},\hat{\eta}_{\hat{t}})]\\
%                       &\leq&\Psi_1(\hat{\gamma}_{\hat{t}},\hat{\eta}_{\hat{t}})
%                      +2L(1+||\check{\gamma}_{\check{t}}||_0)(\check{s}- \check{t})^{\frac{1}{2}}
%                      -\varepsilon[\Upsilon^3(\check{\gamma}_{\check{t}},\hat{\gamma}_{\hat{t}})
%+\Upsilon^3(\check{\eta}_{\check{s}},\hat{\eta}_{\hat{t}})+2|\check{s}-\hat{t}|^2].
%\end{eqnarray*}
 Then letting $\delta\rightarrow0$ and $k\rightarrow\infty$, we obtain
\begin{eqnarray*}
                      &&\limsup_{k\rightarrow\infty}\limsup_{\delta\rightarrow0}\bigg{[}W'_1(\check{\gamma}_{\check{t}})- {W}'_2(\check{\eta}_{\check{s}})-2^5\beta(\Upsilon^3(\check{\gamma}_{\check{t}},\hat{{\xi}}_{{\hat{t}}})
                +\Upsilon^3(\check{\eta}_{\check{s}},\hat{{\xi}}_{{\hat{t}}}))+\varepsilon(\Upsilon^3(\check{\gamma}_{\check{t}},\hat{\gamma}_{\hat{t}})
+\Upsilon^3(\check{\eta}_{\check{s},\check{t},A},\hat{\eta}_{\hat{t}}))
                           \bigg{]}\\
                      &\leq& {W}'_1(\hat{\gamma}_{\hat{t}})-
                    {{W}}'_2(\hat{\eta}_{\hat{t}})-\beta \Upsilon^3(\hat{\gamma}_{\hat{t}},\hat{\eta}_{\hat{t}})-\beta^{\frac{1}{3}}|(\hat{{\gamma}}_{{\hat{t}}}(\hat{t}))^-_N-(\hat{{\eta}}_{{\hat{t}}}(\hat{t}))^-_N|^2.
\end{eqnarray*}
Combining with (\ref{w1v1}), we get
\begin{eqnarray*}
\lim_{k\rightarrow\infty}\limsup_{\delta\rightarrow0}[\Upsilon^3(\check{\gamma}_{\check{t}},\hat{\gamma}_{\hat{t}})
+\Upsilon^3(\check{\eta}_{\check{s},\check{t},A},\hat{\eta}_{\hat{t}})]=0.
\end{eqnarray*}
Noting that, by (\ref{s0}),
$$
||\check{\gamma}_{\check{t}}-\hat{\gamma}_{\hat{t},\check{t},A}||_0+ ||\check{\eta}_{\check{s}}-\hat{\eta}_{\hat{t},\check{s},A}||_0
\leq||\check{\gamma}_{\check{t}}-\hat{\gamma}_{\hat{t},\check{t},A}||_0+ ||\check{\eta}_{\check{s},\check{t},A}-\hat{\eta}_{\hat{t},\check{t},A}||_0\leq4\Upsilon^3(\check{\gamma}_{\check{t}},\hat{\gamma}_{\hat{t}})
+4\Upsilon^3(\check{\eta}_{\check{s},\check{t},A},\hat{\eta}_{\hat{t}})
$$
we obtain (\ref{4.23}) holds true.
  The proof is now complete. \ \ $\Box$

\par

\end{document}